\documentclass[12pt,a4paper]{article}
\usepackage{amssymb,amsmath}
\usepackage{amsfonts}
\usepackage{enumerate}
\usepackage{a4}
\newtheorem{theorem}{Theorem}
\newtheorem{lemma}{Lemma}

\newtheorem{remark}{Remark}

\numberwithin{equation}{section}
\begin{document}
	\baselineskip16pt
	
	\title{\bf Positive values of non-homogeneous quadratic forms of type (1,4): A conjecture of Bambah, Dumir and Hans-Gill }
	\author{Swati Bhardwaj$^1$, Leetika Kathuria$^2$  and Madhu Raka$^1$
		\\ \small{\em $^1$Department of Mathematics, Panjab University, Chandigarh, INDIA}\\
		\small{\em $^2$Mehr Chand Mahajan DAV College for Women, Chandigarh, INDIA}\\
		\footnotesize emails: swatibhardwaj2296@gmail.com,  kathurialeetika@gmail.com, mraka@pu.ac.in\\
		\date{}}
	\maketitle
	\maketitle
	\begin{abstract}
		\noindent Let $Q(x_1, \cdots,x_n)$ be a real indefinite quadratic form of the type $(r,s)$, $n=r+s$,  signature $\sigma=r-s$  and determinant $D\neq 0$. Let $\Gamma_{r,n-r}$ denote the infimum of all numbers $\Gamma$  such that for any  real numbers $c_1, c_2 ,\cdots, c_n$
		there exist integers $x_1, x_ 2,\cdots, x_n$ satisfying
		$$0< Q(x_1+c_1,x_2+c_2,\cdots,x_n+c_n)\leq (\Gamma |D|)^{1/n}.$$
		All the values of $\Gamma_{r,n-r}$ are known except for $\Gamma_{1,4}$. Earlier it was shown that $8\leq \Gamma_{1,4}<12$. It is conjectured that $\Gamma_{1,4}=8$. Here we shall prove that $\Gamma_{1,4}=8$, when (i)  $c_2 \not \equiv 0 \pmod 1$, (ii)  $c_2  \equiv 0 \pmod 1$, $a\geq \frac{1}{2}$, where $a$ is minima of positive definite ternary quadratic forms with determinant $4|D|$, and (iii) in some cases of $c_2  \equiv 0 \pmod 1$, $a< \frac{1}{2}$. We also obtain six critical forms   for which the constant 8 is attained. In the remaining cases we prove that $\Gamma_{1,4}< \frac{32}{3}$.
\vspace{2mm}
\end{abstract}

	\section{Introduction}
	Let $Q(x_1, \cdots,x_n)$ be a real indefinite quadratic form of the type $(r,s)$, $n=r+s$,  signature $\sigma=r-s$  and determinant $D\neq 0$. In 1948, Blaney \cite{Blaney1948} showed that given any real numbers $c_1,c_2,\cdots,c_n$ there exist numbers $\Gamma$ depending only on $r,s$ such that the inequality
	\begin{equation} \label{1.1} 0< Q(x_1+c_1,x_2+c_2,\cdots,x_n+c_n)\leq (\Gamma |D|)^{1/n} \end{equation}
	
	\noindent has a solution in integers $x_1,x_2, \cdots,x_n$.
	Let $\Gamma_{r,s}$ denote the infimum of all such $\Gamma$ for which (\ref{1.1}) has a solution.  The problem is to evaluate $\Gamma_{r,s}$ for different $r,s$ and to determine all those quadratic forms $Q$ and
	$c_1,c_2,\cdots,c_n$ for which equality is needed in (\ref{1.1}) with $\Gamma=\Gamma_{r,s}$. The following table gives history of $\Gamma_{r,s}$ :
	
	$$ \begin{array}{lll} \Gamma_{1,1}&=4 & {\rm ~ Davenport~ and ~Heilbronn}(1947)~ \cite{Davenport1947}\vspace{2mm}\\
	\Gamma_{2,1}&=4 & {\rm ~ Blaney (1950)~\cite{Blaney1950}~ and ~Barnes } (1961)~\cite{Barnes}\vspace{2mm}\\
	\Gamma_{1,2}&=8 & {\rm ~ Dumir} (1967)~\cite{Dumir1967}\vspace{2mm}\\

	\Gamma_{3,1}&=\frac{16}{3} & {\rm ~ Dumir (1968)~ \cite{Dumir1968}}\vspace{2mm}\\
	\Gamma_{2,2}&=16 & {\rm ~ Dumir(1968)~ \cite{Dumir1968}}\vspace{2mm}\\
	\Gamma_{1,3}&=16 & {\rm ~ Dumir~ and ~Hans}\mbox{-}{\rm Gill (1981)~\cite{DumirGill}}\vspace{2mm}\\
	\Gamma_{3,2}&=16 & {\rm ~ Hans}\mbox{-}{\rm Gill~ and ~Raka ~(1980)\cite{GillRaka1980}}\vspace{2mm}\\
	\Gamma_{4,1}&=8 & {\rm ~ Hans}\mbox{-}{\rm Gill~ and ~Raka (1981)~\cite{GillRaka1981}}\vspace{2mm}\\
	\Gamma_{2,3}&=(7/4)^5 & {\rm ~ Bambah, ~ Dumir ~ and ~Hans}\mbox{-}{\rm Gill(1984)~\cite{Bambah1984}}.\end{array}$$
	
	\noindent Bambah, Dumir and Hans-Gill in a series of papers (1981,1983,1984) \cite{Bambah1981, Bambah1983, Bambah1984} proved that $\Gamma_{r,s}=\frac{2^n}{|\sigma|+1}$  for $n\geq 6$, $-1\leq \sigma \leq 3$.\vspace{2mm}
	
	\noindent It was conjectured by Bambah, Dumir and Hans-Gill (1984)~\cite{Bambah1984} that for $n\geq 6$
	\begin{equation}\label{1.2}\Gamma_{r,s}=\left\{\begin{array}{lll} \frac{2^n}{|\sigma|+1} & {\rm if }& |\sigma| \leq 3\vspace{2mm}\\
	\frac{2^n}{\sigma } & {\rm if }& \sigma =4\end{array}\right.\end{equation}
	except for a finite number of exceptions. \vspace{2mm}
	
	\noindent Mary Flahive~(1988)~\cite{Flahive1988} proved the conjecture of Bambah, Dumir and Hans-Gill for $n\geq 21$. Aggarwal and Gupta~(1988, 1991)~ \cite{AggarwalGupta1988}, \cite{AggarwalGupta1991} determined $\Gamma_{r,s}$ for $\sigma=2$ and $n\geq 8$, $\sigma=-3$ and $n\geq 9$ and $\sigma=4$ and $n\geq 6$ and confirming the conjecture of Bambah et al., for these values of $n$ and $\sigma$. Sehmi and Dumir (1994) \cite{DumirSehmi1994} proved that $\Gamma_{2,5}=32$ and $\Gamma_{2,4}= \frac{64}{3}$ was proved by Dumir, Hans-Gill and Sehmi (1995) \cite{DGS1995}.
	Dumir, Hans-Gill and Woods (1994) \cite{DumirWoods} proved that $$ \Gamma_{r,s}=\Gamma_{r',s'} {\rm ~ if ~} r+s=r'+s'=n,~ r-s=\sigma \equiv \sigma'=r'-s' ({\rm mod ~}8 ).$$
	\textbf{Thus $\Gamma_{r,s}$ is determined for all $(r,s)$ except for  $\Gamma_{1,4}$}.  From the work of Jackson (1971) \cite{Jackson} it follows that $\Gamma_{1,4}\leq 32$.
	If  $Q=x_1x_2-\frac{1}{4}(x_2^2+x_3^2+x_4^2+x_5^2)$ and $(c_1,c_2, \cdots,c_5)= (\frac{1}{2}, \cdots,\frac{1}{2})$, the inequality (\ref{1.1}) is not soluble for   $\Gamma_{1,4}< 8$.  Thus  $\Gamma_{1,4}\geq  8$.  Dumir and Sehmi (1994) \cite{DumirSehmi} proved that $8\leq \Gamma_{1,4}< 16$; Raka and Rani (1997) \cite{RakaRani} improved this to $8\leq \Gamma_{1,4}<12$. It is conjectured that $\Gamma_{1,4}=8$.\vspace{2mm}
	
	\noindent The method of proof applied in \cite{DumirSehmi} and \cite{RakaRani} involves very heavy calculations. As now, many software packages like Mathematica are available, we attempted to complete the longstanding conjecture for $\Gamma_{1,4}$, (the only case left in this direction). We have been partially successful in determining the value of $\Gamma_{1,4}$.\vspace{2mm}
	
	\noindent  Marguils (1987) \cite{Marguils} proved the famous Oppenheim Conjecture that if $Q$ is an incommensurable quadratic form in $n\geq 3$ variables, then it  takes arbitrary small non-zero values for integers $x_1,x_2, \cdots,x_n$.   Watson (1960) \cite{Watson} proved that in such case,
	$$ \alpha-\epsilon< Q(x_1+c_1,x_2+c_2,\cdots,x_n+c_n)<\alpha+\epsilon $$ has a solution in integers $x_1,x_2, \cdots,x_n$ for any reals $\alpha, c_1,c_2, \cdots,c_n$ and $\epsilon>0$.
	So we can suppose that $Q$ is a rational form of determinant $D\neq 0$. By Meyer's Theorem, it is a zero form,  being a form in five variables. Using Birch reduction \cite{Birch} and homogeneity we can suppose that $$Q=(x_1 + a_2 x_2 + a_3 x_3 + a_4 x_4 + a_5 x_5) x_2 -\phi(x_3, x_4, x_5),$$
	where $\phi$ is a positive definite ternary quadratic form with determinant $4|D|$. Let $d^5=8|D|.$ Let $a=\min\{\phi(X): 0\neq X \in \mathbb{Z}^3\}.$ Then we get that $0< a\leq d^{5/3}$ from known results. Let $m$
	be an integer satisfying $m<\frac{d}{a}\leq m+1$. It turns out that $m\geq 1$.
	
	\noindent In this paper we prove that $\Gamma_{1,4}=8$ when
	
	\begin{itemize}\item  $c_2 \not \equiv 0 \pmod 1$ (see Theorem \ref{thm1}, Section \ref{sec4})\item $c_2  \equiv 0 \pmod 1$, $a\geq \frac{1}{2}$  (see Theorems \ref{thm2} and \ref{thm3}, Section \ref{sec5})
	\end{itemize}
	We also obtain \textbf{six critical forms} in  Theorems \ref{thm1}-\ref{thm3}  for which the constant 8 is attained in (\ref{1.1}).\vspace{2mm}
	
	\noindent Replacing $\phi(x_3, x_4, x_5)$ by an equivalent form, we can assume that
	$$\phi(x_3,x_4,x_5)=a(x_3+h_4x_4+h_5x_5)^2 + Ax_4^2 + Bx_4x_5 + Cx_5^2,$$  where $0\leq B \leq A \leq C$.
	For $c_2  \equiv 0 \pmod 1$, $a<\frac{1}{2}$, let $K$ be an integer satisfying  $K<\frac{d+(m^2-1)a/4 }{A}\leq K+1$. It turns out that  $K\geq 1$. We further prove that $\Gamma_{1,4}<8$ when
	\begin{itemize}
		\item $c_2  \equiv 0 \pmod 1$, $a<\frac{1}{2}$ and $m\geq 3$ (see Theorem \ref{thm4}, Section \ref{sec6})
		\item $c_2  \equiv 0 \pmod 1$, $a<\frac{1}{2}$, $m=2$ and $K\geq 3$ (see Theorem \ref{thm4'}, Section \ref{sec7}).\end{itemize}
	
	\noindent When  $a<\frac{1}{2}$ and $(m,K)=(2,1)$, we prove that $\Gamma_{1,4}<8.486$ (see Section \ref{sec8}). For $(m,K)=(2,2), (1,2)$ and $(1,1)$, Macbeath's result (see Lemma \ref{lem5}), which is the main tool in our proof,  is not applicable. Here we shall prove that $\Gamma_{1,4}< \frac{32}{3}$ (see  Sections \ref{sec9}-\ref{sec11}). We strongly feel that the conjectured value $\Gamma_{1,4}=8$ is true.

	\section{Some Lemmas}\label{sec2}
	
	We shall call two real quadratic forms $f$ and $g$ in $n$ variables to be equivalent  if there exists a unimodular matrix  $P$ such that $g(X)=f(PX)$ for all $X\in \mathbb{R}^{n}$ and we write $f \sim g.$ In the course of the proof we shall use the following lemmas. Lemma \ref{lem1} and Lemma \ref{lem2} are well known results. See the book by Gruber and Lekkerkerker \cite{Gruber}.
	
	\begin{lemma}\label{lem1}
		Let $\phi(x,y,z)$ be a positive definite quadratic form of determinant $D\neq 0$. Then there exist integers $x,y,z$ such that $$0<\phi(x,y,z)\leq (2|D|)^{1/3},$$ with equality if and only if $\phi \sim \rho(x^2+y^2+z^2+xy+yz+zx),\rho >0$.
	\end{lemma}
	
	\begin{lemma}\label{lem2}
		If $\psi(x,y)$ is a positive definite binary quadratic form of determinant $\delta$, then $$\psi(x,y) \sim Ax^2+Bxy+Cy^2,$$ where $0\leq B \leq A \leq C$ and $AC \leq \frac{4}{3} \delta.$
	\end{lemma}

	\begin{lemma}\label{lem3}
		Let $\alpha,\beta,\gamma$ be real numbers with $\gamma>1$. Let $m$ be the integer defined by $m<\gamma \leq m+1.$ Then for any real $x_0,$ there exists $x \equiv x_0 \pmod 1$ such that
		\begin{equation}\label{eq2.1}
		0<-(x+\alpha)^2+\beta\leq \gamma,
		\end{equation} provided \begin{equation}\label{eq2.2}
		\frac{1}{4} < \beta \leq \gamma+\frac{m^2}{4}.
		\end{equation}
		Further strict inequality in $(\ref{eq2.2})$ implies strict inequality in $(\ref{eq2.1})$.
	\end{lemma}
	This is a result of Dumir(\cite{Dumir1968}).\vspace{2mm}
	
	\noindent We say that the inequality
	\begin{equation}\label{eq1}
	\alpha < Q(x_1,x_2,\cdots,x_n) \leq \beta
	\end{equation}
	
	\noindent is soluble if for any given real numbers $c_1,c_2,\cdots,c_n$, there exist
	$(x_1,x_2,\cdots,x_n) \equiv (c_1,c_2, \cdots,c_n) \pmod 1$ satisfying (\ref{eq1}). If there are integers $u_1,\cdots,u_n$ such that $Q(u_1,\cdots,u_n)=\gamma \neq 0,$ then we say that $Q$ represents $\gamma$.
	
	\begin{lemma}\label{lem4}
		Let $Q(x_1,x_2,\cdots,x_n)$ be a zero form of determinant $D\neq 0$. Let $\alpha, \beta$ be real numbers satisfying \begin{equation}\label{eq2.4}\beta - \alpha \geq  2|D|^{1/n}.\end{equation} Then $(\ref{eq1})$ is soluble.  Further strict inequality in $(\ref{eq2.4})$ implies strict inequality in $(\ref{eq1})$.
	\end{lemma}
	This is a result of Jackson \cite{Jackson}.\vspace{2mm}
	
	\noindent The next  result, which is the main tool in our proof, follows from Macbeath \cite{Macbeath}:
	\begin{lemma}\label{lem5}
		Let $\alpha, \beta, \gamma$ be real numbers with $\alpha > 0, \gamma > 0.$ Let $2h,k$ be integers such that \begin{equation}\label{eq2.5}|h-k^2\alpha|+\frac{1}{2}\leq \gamma.\end{equation} Suppose that either $\alpha \neq h/k^2$ or $\beta \not \equiv h/k \mod{(1/k,2\alpha)} $. Then for any real number $\nu$, there exist integers $x,y$ satisfying \begin{equation}\label{eq2.6}0< \pm x + \beta y\pm \alpha y^2 +\nu \leq \gamma.\end{equation}
		
		\noindent Further strict inequality in $(\ref{eq2.5})$ implies strict inequality in $(\ref{eq2.6})$.
	\end{lemma}

	\begin{lemma}\label{lem7}
		Let $\alpha, \beta, \delta$ and $\nu$ be any real numbers, then \begin{equation}\label{key7}
		0<(x_1+\alpha x_2+\nu)x_2 +\beta\leq \delta
		\end{equation} is soluble for \begin{equation}\label{key71}
		\begin{array}{lll}
		\delta\geq 1/2 & if & c_2 \neq 0\\
		\delta\geq 1 & if & c_2=0.
		\end{array}
		\end{equation}
		\noindent Further strict inequality in $(\ref{key71})$ implies strict inequality in $(\ref{key7})$.
	\end{lemma}
	
	\noindent Proof is trivial and is left to the reader.
	\section{Reduction}\label{sec3}
	As stated in the introduction,  we can suppose that $Q$ is a zero rational form of determinant $D\neq 0$.  Using Birch reduction \cite{Birch} and homogeneity we can suppose that
	$$Q=(x_1 + a_2 x_2 + a_3 x_3 + a_4 x_4 + a_5 x_5) x_2 -\phi(x_3, x_4, x_5),$$
	where $\phi$ is a positive definite ternary quadratic form with determinant $\Delta =4|D|$. Let $d^5=8|D|.$ Let $a= \min\{\phi(X): 0\neq X \in \mathbb{Z}^3\}.$ Then by Lemma \ref{lem1},
	\begin{equation}\label{eq5}
	0< a\leq (2\Delta)^{1/3}=(8|D|)^{1/3}=(d^5)^{1/3}=d^{5/3}
	\end{equation}
	with strict inequality unless $\phi \sim \rho (x^2+y^2+z^2+xy+yz+zx),\rho >0$.\vspace{2mm}
	
	\noindent Without loss of generality we can suppose that the representation of $a$ by $\phi$ is primitive, therefore replacing $\phi$ by an equivalent form we can suppose that $\phi(1,0,0)=a$ and write $$\phi(x_3,x_4,x_5)=a(x_3+h_4x_4+h_5x_5)^2 + \psi(x_4,x_5),$$ where $\psi$ is a positive definite binary form with determinant $\delta = \Delta/a = 4|D|/a$. Therefore, $Q$ becomes $$Q=-a(x_3+h_4x_4+h_5x_5)^2+(x_1+a_2x_2\cdots+a_5 x_5)x_2-\psi(x_4,x_5).$$
	Further, we can suppose that $\frac{-1}{2}<h_i,a_j\leq \frac{1}{2},$ for each $i,j.$
	
	\begin{lemma}\label{lem8}
		If $0<a<d/5$,  the inequality \begin{equation}\label{eq2}
		0<-a(x_3+\cdots)^2+(x_1+a_2x_2\cdots+a_5 x_5)x_2-\psi(x_4,x_5)<d
		\end{equation} is soluble.
	\end{lemma}
	
	\noindent {\bf Proof:} Let $m$ be an integer defined by $m<d/a\leq m+1,$ then by Lemma \ref{lem3}, inequality (\ref{eq2}) is soluble with strict inequality if we can solve
	\begin{equation}\label{eq4}
	a/4 < (x_1+a_2x_2\cdots+a_5 x_5)x_2-\psi(x_4,x_5) < d+m^2a/4.
	\end{equation}
	Since $(x_1+\cdots+a_5 x_5)x_2-\psi(x_4,x_5) $ is a zero form with determinant $-D/a$, (\ref{eq4}) is soluble by Lemma \ref{lem4}, if $$d+\frac{1}{4}(m^2-1)a > 2\left(\frac{|D|}{a}\right)^{1/4} = \left(\frac{2d^5}{a}\right)^{1/4}$$
	which is satisfied for $m \geq 5$. This proves the Lemma. \hfill $\square$

	\noindent Since $\psi$ is a positive definite binary form with determinant $\delta = \Delta/a = 4|D|/a.$ By Lemma \ref{lem2}, we can suppose that \begin{equation}\label{eq6}
	\psi(x_4,x_5) = Ax_4^2 + Bx_4x_5 + Cx_5^2 = A(x_4+\lambda x_5)^2 + tx_5^2,
	\end{equation}
	where
	\begin{equation}\label{eq7}
	0\leq B \leq A \leq C {\rm ~and~} 0<A^2\leq AC\leq \frac{4\delta}{3}=\frac{16|D|}{3a} =\frac{2d^5}{3a},
	\end{equation}
	\begin{equation}\label{eq8}
	0 \leq \lambda \leq \frac{1}{2},~~~~ t=\frac{(\det \psi)}{A}=\frac{\delta}{A}=\frac{4|D|}{aA} =\frac{d^5}{2aA},
	\end{equation}
	and
	\begin{equation}\label{eq09}
	A\leq C=t+A\lambda^2 ~\mbox{which gives}~ 3A/4\leq t ~\mbox{and} ~ \lambda\geq \sqrt{1-t/A}.
	\end{equation}
	\noindent Since $A+h_4^2a$ is a value of $\phi$ and $a$ is the minimum value, we have $$a \leq A+h_4^2a \leq A+a/4,$$ so that \begin{equation}\label{eq9}
	A \geq 3a/4, {\rm ~and~if~} h_4=0 {\rm ~then~} A \geq a.
	\end{equation}
	\noindent We need to show that \begin{equation}\label{eq10}
	0<Q=(x_1+a_2x_2+a_3x_3+...)x_2 - a(x_3+...)^2 -A(x_4+\lambda x_5)^2 - tx_5^2 \leq d
	\end{equation}
	is soluble. Without loss of generality, we can suppose that $\frac{-1}{2}<c_i\leq \frac{1}{2}$ for each $i.$\vspace{2mm}

	\noindent Let $m$ be the integer defined by $m<\frac{d}{a}\leq m+1$.  \vspace{2mm}
	
	\noindent We can rewrite $Q(x_1,x_2,x_3,x_4,x_5)$ as
	\begin{equation*}
	Q=-a(x_3-a_3x_2/2a+h_4x_4+h_5x_5)^2+(x_1+a_2'x_2+a_4'x_4+a_5'x_5)x_2 - A(x_4+\lambda x_5)^2 - tx_5^2,
	\end{equation*}
	where $a_2', a_4', a_5'$ are some real numbers satisfying
	\begin{equation}\begin{array}{ll}\label{eq12}
	a_2'  \equiv  a_2+a_3^2/4a\hspace{-2mm} \pmod 1, &\\
	a_4'  \equiv  a_4-a_3 h_4 \pmod 1, & ~~~~ -1/2 <a_2', a_4', a_5' \leq 1/2. \\
	a_5'  \equiv  a_5-a_3 h_5 \pmod 1, &\\
	\end{array}
	\end{equation}
	
	\noindent Now (\ref{eq10}) is soluble by Lemma \ref{lem3}, if we can solve
	\begin{equation}\label{eq13}
	\begin{array}{lll}
	0<F(x_1,x_2,x_4,x_5)&=&(x_1+a_2'x_2+a_4'x_4+a_5'x_5)x_2 - A(x_4+\lambda x_5)^2 \\
	&&- tx_5^2-a/4 \leq d+(m^2-1)a/4 = \delta_m {\rm ~(say)}.
	\end{array}
	\end{equation}
	
	\noindent Let $K$ be the integer defined by $K< \delta_m/A\leq K+1.$ By definition of $m$, we have
	\begin{equation}\label{eq18'}
	a \geq d/(m+1),~\delta_m = d+(m^2-1)a/4 \geq (m+3)d/4.
	\end{equation}

	\noindent Using (\ref{eq7}), we get
	\begin{equation}\label{eq14}
	\frac{\delta_m}{A} \geq \frac{(m+3)d}{4} \sqrt{\frac{3a}{2d^5}} \geq \frac{m+3}{\sqrt{m+1}}\cdot \frac{\sqrt{3}}{4\sqrt{2}}\cdot \frac{1}{d}.
	\end{equation}

	\noindent We can rewrite $F(x_1,x_2,x_4,x_5)$ as  $$F=-A(x_4-a_4'x_2/2A+\lambda x_5)^2 + (x_1+a_2''x_2+a_5''x_5)x_2 - tx_5^2 -a/4,$$
	where
	\begin{equation}\label{eq16}
	a_2''  \equiv  a_2'+a_4'^2/4A ({\rm mod~ }1),~
	a_5''  \equiv  a_5'-\lambda a_4' ({\rm mod~ }1),  ~ -1/2 <a_2'', a_5'' \leq \frac{1}{2}.
	\end{equation}
	
	\noindent Now (\ref{eq13}) and hence (\ref{eq10}) is soluble, by Lemma \ref{lem3} if we can solve
	\begin{equation}\begin{array}{lll}\label{eq17}
	0<G(x_1,x_2,x_5)& =& (x_1+a_2''x_2+a_5''x_5)x_2 - tx_5^2 - (a+A)/4\\
	&\leq&\delta_m +(K^2-1)A/4 = \delta_{m,K} {\rm ~(say)}.
	\end{array}
	\end{equation}
	
	\noindent By definition of $m$ and $K$, we have
	\begin{equation}\label{eq18}
	A\geq(m+3)d/4(K+1), ~\delta_{m,K} \geq (m+3)(K+3)d/16.
	\end{equation}
	
	\noindent Now since $(x_1+a_2''x_2+a_5''x_5)x_2 - tx_5^2$ is a zero form, (\ref{eq17}) is soluble, by Lemma \ref{lem4} if $$\delta_{m,K} \geq 2(t/4)^{1/3} =(d^5/aA)^{1/3}, $$
	i.e., using (\ref{eq18}), if we have
	\begin{equation}\label{eq19}
	 \left(\frac{m+3}{4}\right)^4\left(\frac{K+3}{4}\right)^3\left(\frac{1}{m+1}\right)\left(\frac{1}{K+1}\right) \geq 1.
	\vspace{2mm}\end{equation}

	\noindent Let $L$ be the integer defined by $L<\delta_{m,K}/t \leq L+1$.
	\noindent We can rewrite $G(x_1,x_2,x_5)$ as $$G= -t(x_5-a_5''x_2/2t)^2+(x_1+a_2'''x_2)x_2 - (a+A)/4,$$ where \begin{equation}\label{eq20}
	a_2'''\equiv a_2''+a_5''/4t \hspace{-2mm} \pmod 1,~~ -1/2<a_2''' \leq 1/2.
	\end{equation}
	Therefore (\ref{eq17}) and hence (\ref{eq10}) is soluble, by Lemma \ref{lem3} if we can solve \begin{equation}\label{eq21}
	0<H(x_1,x_2) = (x_1+a_2'''x_2)x_2 - (a+A+t)/4 \leq \delta_{m,K}+(L^2-1)t/4 =\delta_{m,K,L}~({\rm say}).
	\end{equation}
	By definition of  $m,K$ and $L$, we have
	\begin{equation}\label{eq22}
	t\geq \frac{(m+3)(K+3)d}{16(L+1)},~ \delta_{m,K,L} \geq (m+3)(K+3)(L+3)d/64.
	\end{equation}
	
	\noindent As $(x_1+a_2'''x_2)x_2$ is a zero form, (\ref{eq21}) is soluble, by Lemma \ref{lem4} if
	\begin{equation}\label{eqn23}
	\delta_{m,K,L} \geq 2(1/4)^{1/2}=1.
	\end{equation}
	
	\noindent One observes that strict inequality in (\ref{eq21}) implies strict inequality in (\ref{eq17}) which gives strict inequality in (\ref{eq13}) and hence strict inequality in (\ref{eq10}).

	\section{Proof of $\Gamma_{1,4}= 8$ for $c_2 \not \equiv 0 \pmod 1$}\label{sec4}
	
	\begin{theorem}\label{thm1}
		Let $Q(x_1,x_2, \cdots, x_5)$ be a real indefinite quadratic form  of type $(1,4)$ and of determinant $D\neq 0$. Then  given any real numbers $c_1,c_2, \cdots,c_5$, with $c_2 \not \equiv 0 \pmod1$ there exists $(x_1,x_2,\cdots,x_5) \equiv (c_1,c_2, \cdots,c_5) \pmod 1$ such that \begin{equation}\label{theq1}
		0<Q(x_1+c_1,x_2+c_2, \cdots,x_5+c_5)\leq d=(8 |D|)^{1/5}
		\end{equation}
		with equality if and only if $Q \sim \rho Q_1, \rho >0,$ where $Q_1=(x_1-\frac{1}{4}x_2)x_2-\frac{1}{4}(x_3^2+x_4^2+x_5^2).$ For $Q_1$, equality occurs in $(\ref{theq1})$ if and only if,  $(c_1,c_2,c_3,c_4,c_5)\equiv (\frac{1}{2},\frac{1}{2},\frac{1}{2},\frac{1}{2},\frac{1}{2})\pmod 1.$
	\end{theorem}

	\noindent \textbf{Proof:}
	We use the reduction on $Q$ as done in Section \ref{sec3} and we will show that either (\ref{eq10}) or (\ref{eq13}) or (\ref{eq17}) or (\ref{eq21}) is soluble with strict inequality unless $Q=Q_1$.\\
	\noindent From Lemma {\ref{lem7}}, equation (\ref{eq10}) is soluble  with strict inequality for $d>\frac{1}{2}$. Therefore, let
	\begin{equation}\label{eq*}
	d \leq \frac{1}{2}.\end{equation}
	Since $m$ is the integer defined by $m<\frac{d}{a}\leq m+1$. From (\ref{eq5}) and Lemma \ref{lem8}, we have $5\geq \frac{d}{a}\geq \frac{d}{d^{5/3}}>1$. So $m=1,2,3,{\rm ~or~} 4$ only.\vspace{2mm}
	
	\noindent We will discuss $m=2, 3, 4$ in Lemma \ref{lem7*} and $m=1$ in Lemma \ref{lem8*}.

	\begin{lemma}\label{lem7*} If $c_2 \not \equiv 0 \pmod 1$, and $m=2, 3$ or $4$, then {\rm(\ref{theq1})} is soluble with strict inequality.
	\end{lemma}
	\noindent  {\bf Proof:} By Lemma \ref{lem7}, (\ref{eq13}) is soluble if $\delta_m>1/2.$ Let therefore $\delta_m =d+\frac{(m^2-1)a}{4}\leq 1/2$. Here, as $a\geq \frac{d}{m+1}$, we have $\frac{(m+3)d}{4} \leq \frac{1}{2}$, which gives $d \leq \frac{2}{m+3}$. Also from (\ref{eq7}) we have $$ \frac{\delta_m}{A} \geq \frac{(m+3)d}{4} \sqrt{\frac{3a}{2d^5}} \geq \frac{(m+3)^2}{8} \sqrt{\frac{3}{2}}\frac{1}{\sqrt{m+1}}.
	$$
	So, \begin{equation*}
	\frac{\delta_m}{A} \geq \left\{
	\begin{array}{lll}
	& 3.354 &{\rm~if~} m=4\\
	& 2.755 &{\rm~if~} m=3\\
	& 2.209 &{\rm~if~} m=2.\\
	\end{array}
	\right.
	\end{equation*}
	Since $K< \delta_m/A\leq K+1$, we need to discuss the following cases
	\begin{equation*}\begin{array}{lll}
	{\rm(i)}~~~m=4, K\geq 3  ~~~~&~~~~ {\rm(ii)}~~m=3, K\geq 2 ~~~~&~~~~ {\rm(iii)}~m=2, K\geq 2
	\end{array}
	\end{equation*}
	
	\noindent Now by Jackson's Lemma (Lemma \ref{lem4}), (\ref{eq17}) is soluble with strict inequality if $\delta_{m,K} > 2\left( \frac{t}{4}\right)^{1/3} =\left( \frac{d^5}{aA}\right)^{1/3},$ which is so if
	$\frac{(m+3)(K+3)d}{16} > \left( \frac{d^5}{aA}\right)^{1/3},$ i.e., if
	\begin{equation}\label{eq23}
	\frac{(m+3)^4}{m+1} \frac{(K+3)^4}{K+1} > 4^7.
	\end{equation}
	Also by Lemma \ref{lem7}, (\ref{eq17}) is soluble with strict inequality if $\delta_{m,K}>1/2,$ i.e., if $\frac{(m+3)(K+3)d}{16} > 1/2$.
	From (\ref{eq7}) and (\ref{eq18}), we have
	
	$$\frac{d(m+3)}{4(K+1)}\leq \frac{\delta_m}{K+1} \leq A \leq \sqrt{\frac{2}{3}\frac{d^5}{a}}\leq \sqrt{\frac{2}{3}(m+1)}~ d^{\hspace{1mm}2},$$ which gives
	$d \geq \sqrt{\frac{3}{2(m+1)}}~\frac{m+3}{4}~\frac{1}{(K+1)}.$ Thus $\delta_{m,K} > 1/2$, if
	$$\frac{(m+3)(K+3)d}{16} \geq  \frac{(m+3)(K+3)}{16}\sqrt{\frac{3}{2(m+1)}}~\frac{m+3}{4}~\frac{1}{(K+1)}>1/2$$
	
	\noindent i.e., if \begin{equation}\label{eq24}
	\sqrt{\frac{3}{2(m+1)}}~\left(\frac{m+3}{4}\right)^2~\frac{(K+3)}{4}~\frac{1}{(K+1)} > \frac{1}{2}.
	\end{equation}
	Thus (\ref{eq17}) is soluble with strict inequality if either (\ref{eq23}) or (\ref{eq24}) is satisfied. We find that
	\begin{enumerate}[$\rm(i)$]
		
		\item For $m=4, K\geq 3$, (\ref{eq23}) is satisfied.
		\item For $m=3, K\geq 3$, (\ref{eq23}) is satisfied and for $m=3, K=2$ (\ref{eq24}) is satisfied.
		\item $m=2, K\geq 5$, (\ref{eq23}) is satisfied.

		\item  For $m=2, K=2,3,4;$ (\ref{eq17}) is  soluble when $\delta_{m,K} > 1/2$. So let $\delta_{m,K} \leq 1/2$, so that we have $d \leq \frac{8}{(m+3)(K+3)}.$ Therefore, for $m=2$ and $ K \geq 2,$
		\begin{equation}\label{eq25}
		\frac{\delta_{m,K}}{t} \geq d\frac{(m+3)}{4}\frac{(K+3)}{4}\frac{2aA}{d^5}\geq \frac{(m+3)^4(K+3)^3}{2^{11}(m+1)(K+1)}>4,
		\end{equation}
		which gives $L \geq 4$, (since $L<\delta_{m,K}/t \leq L+1$). \\
		Now for $m=2, K\geq 2,L\geq 4$, we have $\frac{d}{m+1}\leq a \leq d^{\frac{5}{3}}$ which gives $d\geq (m+1)^{-\frac{3}{2}}$, so that
		$$ \delta_{m,K,L} \geq d~ \frac{(m+3)}{4}\frac{(K+3)}{4}\frac{(L+3)}{4}\geq (m+1)^{-\frac{3}{2}} \frac{(m+3)}{4}\frac{(K+3)}{4}\frac{(L+3)}{4}> \frac{1}{2}.$$
		
		\noindent	Therefore, by Lemma \ref{lem7}, (\ref{eq21}) is soluble for $c_2 \not \equiv 0 \pmod 1$.
	\end{enumerate}
	\noindent This completes the proof of Lemma \ref{lem7*}. \hfill $\square$

	\begin{lemma}\label{lem8*} If $c_2 \not \equiv 0 \pmod 1$, and $m=1$ then $(\ref{theq1})$ is soluble with strict inequality unless  $Q \sim \rho Q_1, \rho >0$ and   $(c_1,c_2,c_3,c_4,c_5)\equiv (\frac{1}{2},\frac{1}{2},\frac{1}{2},\frac{1}{2},\frac{1}{2})\pmod 1,$ for which equality is required in $(\ref{theq1})$.
	\end{lemma}
	\noindent  {\bf Proof:}  For $m=1$ we have $1<\frac{d}{a}\leq 2,$ so that from (\ref{eq5}) we have $d/2\leq a \leq d^{5/3}$, i.e., $d\geq 1/\sqrt{8}.$
	Choose $(x_4,x_5) \equiv (c_4,c_5) \pmod 1$ arbitrarily, $ x_2\equiv c_2 \pmod 1$  such that $0 \neq |x_2| \leq \frac{1}{2}$ (since $c_2 \neq 0$), take $x_1=x+c_1$, $x_3=y+c_3$, so that (\ref{eq10}) becomes \begin{equation}\label{eq27}
	0<\pm x +\beta_a y-\alpha y^2+\nu \leq \frac{d}{|x_2|},
	\end{equation}
	where $\alpha= a/|x_2|,~ \beta_a=\pm a_3-2a(c_3+h_4x_4+h_5x_5)/|x_2|$ and $\nu$ is some real number.\vspace{2mm}
	
	\noindent Take $h=1/2, k=1$ in Lemma \ref{lem5} to get
	$$|h-k^2\alpha|+\frac{1}{2}= \Big|\frac{1}{2}-\frac{a}{|x_2|}\Big|+\frac{1}{2}=\left\{ \begin{array}{ll} \frac{a}{|x_2|}& {\rm ~(if~ } \frac{a}{|x_2|}>\frac{1}{2})\vspace{2mm}\\ 1-\frac{a}{|x_2|} & {\rm ~(if~ } \frac{a}{|x_2|}\leq \frac{1}{2})\end{array}\right\} < \frac{d}{|x_2|}$$

	\noindent as $a\leq d^{\frac{5}{3}}<d$ and  $a+d \geq d/2+d = 3d/2>1/2\geq |x_2|$. By Lemma \ref{lem5},  (\ref{eq27}) is soluble with strict inequality  unless $\alpha=1/2$ and $\beta_a \equiv 1/2 \pmod 1$. Taking $x_4=c_4$ and $1+c_4$ simultaneously in $\beta_a$ we get $h_4=0,$ and hence $A\geq a.$ And on taking $x_5=c_5$ and $1+c_5$ simultaneously in $\beta_a$ we get $h_5=0.$\vspace{2mm}
	
	\noindent Also $\beta_a=\pm a_3-c_3 \equiv 1/2 \pmod 1$ gives $(a_3,c_3)=(0,1/2)$ or $(1/2,0).$\vspace{2mm}

	\noindent Now for $a/|x_2|=1/2$, (\ref{eq10}) reduces to
	\begin{equation}\label{eq**}0<\frac{Q}{|x_2|}=\pm(x_1+a_2x_2+\cdots)- \frac{1}{2} x_3^2 -\frac{A}{|x_2|} (x_4+\lambda x_5)^2-\frac{t}{|x_2|}x_5^2\leq \frac{d}{|x_2|}.\end{equation}
	
	\noindent Choose $(x_3,x_5) \equiv (c_3,c_5) \pmod 1$ arbitrarily,  take $x_1=x+c_1$, $x_4=y+c_4$, so that (\ref{eq**}) becomes \begin{equation}\label{eq***}
	0<\pm x +\beta_A y-\alpha_1 y^2+\nu \leq \frac{d}{|x_2|},
	\end{equation}
	where $\alpha_1= A/|x_2|,~ \beta_A=\pm a_4-2A(c_4+\lambda x_5)/|x_2|$ and $\nu$ is some real number.\vspace{2mm}
	
	\noindent Take $h=1/2, k=1$ in Lemma \ref{lem5}. Then since $\frac{A}{|x_2|}\geq \frac{a}{|x_2|}=\frac{1}{2}$ we have
	$$|h-k^2\alpha_1|+\frac{1}{2}= \Big|\frac{1}{2}-\frac{A}{|x_2|}\Big|+\frac{1}{2}= \frac{A}{|x_2|}  < \frac{d}{|x_2|}$$
	
	\noindent  as $A\leq \sqrt{\frac{2d^5}{3a}}<\sqrt{\frac{4}{3}}~d^2<d$ for  $d\leq \frac{1}{2}$. Therefore (\ref{eq***}) is soluble with strict inequality unless $\alpha_1=1/2$ and $\beta_A \equiv 1/2 \pmod 1$. Taking $x_5=c_5$ and $1+c_5$ simultaneously in $\beta_A$ we get $\lambda=0.$ Equation (\ref{eq7}) gives $A\leq C=t.$
	Also $\beta_A=\pm a_4-c_4 \equiv 1/2 \pmod 1$ gives $(a_4,c_4)=(0,1/2)$ or $(1/2,0).$\vspace{2mm}

	\noindent Now for $a/|x_2| = A/|x_2|= 1/2$, choose $(x_3,x_4) \equiv (c_3,c_4) \pmod 1$ arbitrarily,  take $x_1=x+c_1$, $x_5=y+c_5$, so that (\ref{eq10}) becomes \begin{equation}\label{eq29}
	0<\pm x +\beta_t y-\alpha_2 y^2+\nu \leq \frac{d}{|x_2|},
	\end{equation}
	where $\alpha_2= t/|x_2|,~ \beta_t=\pm a_5-2tc_5/|x_2|$ and $\nu$ is some real number.\vspace{2mm}
	
	\noindent  Note that
	\begin{equation}\label{eqa}\frac{t}{|x_2|}=\frac{d^5}{2aA|x_2|}= \frac{d^5}{4a^3}\leq 2d^2\leq \frac{1}{2}. \end{equation}
	Take $h=1/2, k=1$ to find that $$|h-k^2\alpha_2|+\frac{1}{2}=\Big|\frac{1}{2}-\frac{t}{|x_2|}\Big|+\frac{1}{2}=1-\frac{t}{|x_2|}<\frac{d}{|x_2|}$$ as $t+d \geq A+d \geq a+d \geq d/2+d = 3d/2>1/2\geq |x_2|$ and  and hence (\ref{eq29}) is soluble with strict inequality unless $\alpha_2=t/|x_2|=1/2$ and $\beta_t \equiv 1/2 \pmod 1$, which gives $\pm a_5-c_5 \equiv 1/2 \pmod 1$ and hence $(a_5,c_5)=(0,1/2)$ or $(1/2,0).$\vspace{2mm}
	
	\noindent For $t=|x_2|/2,$  we must have equality everywhere in (\ref{eqa}) which gives  $d=1/2$, $a=\frac{d}{2}=\frac{1}{4}$. Also, $\frac{1}{4}=a=\frac{|x_2|}{2}$ gives $|x_2|=\frac{1}{2}$. Thus we have  $a=A=t=1/4$ and $c_2=1/2$. Then, $Q=(x_1+a_2x_2+\cdots+a_5x_5)x_2-\frac{1}{4}(x_3^2+x_4^2+x_5^2)$ and because of symmetry in $x_3,x_4,x_5$ we need to consider the following cases \begin{enumerate}
		\item $(a_3,a_4,a_5)=(0,0,0), (c_3,c_4,c_5)=(\frac{1}{2},\frac{1}{2},\frac{1}{2}) $
		\item $(a_3,a_4,a_5)=(\frac{1}{2},0,0), (c_3,c_4,c_5)=(0,\frac{1}{2},\frac{1}{2}) $
		\item $(a_3,a_4,a_5)=(\frac{1}{2},\frac{1}{2},0), (c_3,c_4,c_5)=(0,0,\frac{1}{2}) $
		\item $(a_3,a_4,a_5)=(\frac{1}{2},\frac{1}{2},\frac{1}{2}), (c_3,c_4,c_5)=(0,0,0) $.
	\end{enumerate}
	\textbf{Case 1}:  $(a_3,a_4,a_5)=(0,0,0), (c_3,c_4,c_5)=(\frac{1}{2},\frac{1}{2},\frac{1}{2}). $\vspace{2mm}
	
	\noindent Here, we need to solve \begin{equation}\label{eq03}
	0<Q=(x_1+a_2x_2)x_2-\frac{1}{4}(x_3^2+x_4^2+x_5^2)\leq 1/2.
	\end{equation}
	As $-1/2<c_1, a_2\leq 1/2$, we have $-3/4<c_1+a_2/2\leq 3/4.$
	\noindent Taking $x_1=x+c_1(x \in \mathbb{Z}),$ $x_2=x_3=x_4=x_5=1/2$ we get $0<Q=(x+c_1)/2+a_2/4-3/16\leq 1/2$ which is soluble if \begin{equation}\label{eq01}
	3/8 < x+c_1+a_2/2 \leq 11/8
	\end{equation} has a solution in integer $x$. Taking $x=0,1,2$ simultaneously, we see that (\ref{eq01}) is soluble when $3/8<c_1+a_2/2\leq 11/8$, or $-5/8<c_1+a_2/2\leq 3/8$, or $-13/8<c_1+a_2/2\leq -5/8$. Thus (\ref{eq01}) is soluble with strict inequality unless $  c_1+a_2/2 \equiv 3/8 \pmod 1.$ \vspace{2mm}
	
	\noindent Similarly taking $ x_1=-x+c_1$,  $x_2=-1/2, x_3=x_4=x_5=1/2$ in (\ref{eq03}) we get that (\ref{eq03}) is soluble if \begin{equation}\label{eq02}
	3/8 < x-c_1+a_2/2 \leq 11/8
	\end{equation} has a solution in integer $x$. Taking $x=0,1,2$ simultaneously, (\ref{eq02}) is soluble when $3/8<-c_1+a_2/2\leq 11/8$, or $-5/8<-c_1+a_2/2\leq 3/8$, or $-13/8<-c_1+a_2/2\leq -5/8$, and hence (\ref{eq02}) is soluble with strict inequality unless $  -c_1+a_2/2 \equiv 3/8 \pmod 1.$\vspace{2mm}
	
	\noindent On solving these we get $c_1=1/2, a_2=-1/4,$ for which $Q$ reduces to $Q_1=(x_1-\frac{1}{4}x_2)x_2 -\frac{1}{4}(x_3^2+x_4^2+x_5^2)$. On taking unimodular transformation $x_1 \rightarrow x_1, x_2 \rightarrow x_2+2x_1, x_3 \rightarrow x_3, x_4 \rightarrow x_4, x_5 \rightarrow x_5,$ we can suppose that  $Q_1=x_1^2-(x_2^2+x_3^2+x_4^2+x_5^2)/4.$ For this quadratic form $Q_1$ and $(c_1,c_2,c_3,c_4,c_5) \equiv (\frac{1}{2},\frac{1}{2},\frac{1}{2},\frac{1}{2},\frac{1}{2})$, equality is required in (\ref{eq03}). This is so because   for $y_i \in \mathbb{Z}$, $16Q_1(y_1+\frac{1}{2},y_2+\frac{1}{2},y_3+\frac{1}{2},y_4+\frac{1}{2},y_5+\frac{1}{2})=4(2y_1+1)^2 - (2y_2+1)^2-(2y_3+1)^2-(2y_4+1)^2-(2y_5+1)^2 \equiv 4-1-1-1-1\equiv 0 \pmod 8.$\vspace{2mm}
	
	\noindent \textbf{Case 2}: $(a_3,a_4,a_5)=(\frac{1}{2},0,0), (c_3,c_4,c_5)=(0,\frac{1}{2},\frac{1}{2}).$\vspace{2mm}
	
	\noindent Here, we need to solve \begin{equation}\label{eq04}
	0<Q=(x_1+a_2x_2+x_3/2)x_2-\frac{1}{4}(x_3^2+x_4^2+x_5^2)\leq 1/2.
	\end{equation}
	Again $-3/4<c_1+a_2/2\leq 3/4.$
	\noindent Taking $x_1=x+c_1(x \in \mathbb{Z}),$ $x_2=x_4=x_5=1/2,$ $ x_3=1$ we get that $0<Q=(x+c_1)/2+a_2/4-1/8\leq 1/2$ which is soluble if \begin{equation}\label{eq05}
	0<x+c_1+a_2/2-1/4\leq 1
	\end{equation} has a solution in integer $x.$ Taking $x=0,1,2$ simultaneously, (\ref{eq05}) is soluble when  $1/4<c_1+a_2/2\leq 5/4$,or $-3/4<c_1+a_2/2\leq 1/4$, or $-7/4<c_1+a_2/2\leq -3/4$. Thus (\ref{eq05}) is soluble with strict inequality unless $  c_1+a_2/2 \equiv 1/4 \pmod 1.$ \vspace{2mm}
	
	\noindent Similarly taking $ x_1=-x+c_1$,  $x_2=-1/2, x_4=x_5=1/2$,$x_3=1$ in (\ref{eq04}) we get that (\ref{eq04}) is soluble if \begin{equation}\label{eq06}
	0<x-c_1+a_2/2-1/4\leq 1
	\end{equation} has a solution in integer $x.$ Taking $x=0,1,2$ simultaneously, (\ref{eq06}) is soluble when $1/4<-c_1+a_2/2<5/4$, or $-3/4<-c_1+a_2/2<1/4$, or $-7/4<-c_1+a_2/2<-3/4$, and hence (\ref{eq04}) is soluble with strict inequality unless $  -c_1+a_2/2 \equiv 1/4 \pmod 1.$\vspace{2mm}
	
	\noindent On solving these we get $c_1=0, a_2=1/2,$ for which $Q$ reduces to $(x_1+\frac{1}{2}x_2+\frac{1}{2}x_3)x_2 -\frac{1}{4}(x_3^2+x_4^2+x_5^2),$ which is equivalent to $Q_1$ by the unimodular transformation, $x_1 \rightarrow x_1-x_2, x_2 \rightarrow x_2, x_3 \rightarrow x_3+x_2, x_4 \rightarrow x_4, x_5 \rightarrow x_5$. \vspace{2mm}

	\noindent \textbf{Case 3}: $(a_3,a_4,a_5)=(\frac{1}{2},\frac{1}{2},0), (c_3,c_4,c_5)=(0,0,\frac{1}{2}). $\vspace{2mm}
	
	\noindent Here, we need to solve $$0<Q=(x_1+a_2x_2+x_3/2+x_4/2)x_2-\frac{1}{4}(x_3^2+x_4^2+x_5^2)\leq 1/2.$$
	\noindent Proceeding as in above cases, we get $c_1=0, a_2=1/4,$ for which $Q$ reduces to $(x_1+\frac{1}{4}x_2+\frac{1}{2}x_3+\frac{1}{2}x_4)x_2 -\frac{1}{4}(x_3^2+x_4^2+x_5^2),$ which is equivalent to $Q_1$ by the unimodular transformation, $x_1 \rightarrow x_1-x_2, x_2 \rightarrow x_2, x_3 \rightarrow x_3+x_2, x_4 \rightarrow x_4+x_2, x_5 \rightarrow x_5 $.\vspace{2mm}
	
	\noindent \textbf{Case 4}: $(a_3,a_4,a_5)=(\frac{1}{2},\frac{1}{2},\frac{1}{2}), (c_3,c_4,c_5)=(0,0,0). $\vspace{2mm}
	
	\noindent Here, we need to solve $$0<Q=(x_1+a_2x_2+x_3/2+x_4/2+x_5/2)x_2-\frac{1}{4}(x_3^2+x_4^2+x_5^2)\leq 1/2.$$
	\noindent Proceeding as in above cases, we get $c_1=0, a_2=0,$ for which $Q$ reduces to $(x_1+\frac{1}{2}x_3+\frac{1}{2}x_4+\frac{1}{2}x_5)x_2 -\frac{1}{4}(x_3^2+x_4^2+x_5^2),$ which is equivalent to $Q_1$ by the unimodular transformation, $x_1 \rightarrow x_1-x_2, x_2\rightarrow x_2, x_3 \rightarrow x_3+x_2, x_4 \rightarrow x_4+x_2, x_5 \rightarrow x_5+x_2$. \vspace{2mm}
	
	\noindent This completes the proof of Lemma \ref{lem8*}.\vspace{1mm} \hfill $\square$
	
	\noindent Now the proof of Theorem \ref{thm1} follows from Lemmas \ref{lem7*} and \ref{lem8*}.

	\section{Proof of $\Gamma_{1,4}= 8$ for $c_2 \equiv 0 \pmod 1$, $a\geq \frac{1}{2}$}\label{sec5}
	\numberwithin{equation}{section}
	\noindent For $c_2 \equiv 0 \pmod 1$, inequality (\ref{eq10}) is soluble by Lemma \ref{lem7} when $d>1.$ Therefore, we can assume that $d \leq 1$. Choose $(x_4,x_5) \equiv (c_4,c_5) \pmod 1$ arbitrarily. Take  $x_2 = \pm 1$, $x_1=x+c_1$ and $x_3=y+c_3$. Then $Q= \pm x + \beta_a y-ay^2+\nu$, where $\beta_a =\pm a_3-2a(c_3+h_4x_4+h_5x_5)$ and $\nu$ is some real number. On taking $h=\frac{1}{2}, k=1$, $\alpha =a$ and $\gamma =d$, (\ref{eq10}) is soluble by Macbeath's Lemma \ref{lem5}, if
	
	\begin{equation}\label{eq5.1}|h-ak^2|+\frac{1}{2}=~\vline~ \frac{1}{2} - a ~\vline +\frac{1}{2} < d.\end{equation}
	
	\noindent For $a\geq \frac{1}{2} $, the above inequality holds as $a\leq d^{\frac{5}{3}}<d$ unless $a=d=1$. Therefore (\ref{eq10}) is soluble by Lemma \ref{lem5} unless $a=d=1$;  or $a=\frac{1}{2}$ and $\beta_a =\pm a_3-2a(c_3+h_4x_4+h_5x_5) \equiv \frac{1}{2} \pmod 1$.  Now, since $-\frac{1}{2} <h_4 \leq \frac{1}{2}$, taking $x_4=c_4$ and $1+c_4$, we get $h_4=0$. Similarly, $h_5=0$. For $h_4=h_5=0$, $\pm a_3-2ac_3 \equiv \frac{1}{2} \pmod 1$ implies that $(a_3,c_3)=(0,\frac{1}{2})$ or $(\frac{1}{2},0)$. Note that for $a <\frac{1}{2}$ it is satisfied unless $a+d \leq 1$.\vspace{2mm}
	
	\noindent  In Section \ref{subsec1}, we discuss the case $a=1, d=1$; and in Section \ref{subsec5.2}, we will discuss the case $a=\frac{1}{2},$ $d\leq 1$, $h_4=h_5=0$ and $(a_3,c_3)=(0,\frac{1}{2})$ or $(\frac{1}{2},0)$. The case $a<\frac{1}{2}$ and $a+d\leq 1$ is too much involved and is discussed in Sections \ref{sec6}-\ref{sec11}.\vspace{3mm}
	
	\noindent The following notations will be used throughout the paper.
	For $n=1,2,3,\cdots$ denote\vspace{1mm}\\
	$~~f_n=n|c_1|+n^2a_2,~~ f_n'=n|c_1|+n^2a_2',~~ f_n''=n|c_1|+n^2a_2'', ~~ f_n'''=n|c_1|+n^2a_2''', $\vspace{1mm}\\
	$~~g_n=-n|c_1|+n^2a_2,~~ g_n'=-n|c_1|+n^2a_2',~~ g_n''=-n|c_1|+n^2a_2'',~~ g_n'''=-n|c_1|+n^2a_2''',$\vspace{1mm}\\
	$~~p_n=nc_1+n^2a_2,~~ p_n'=nc_1+n^2a_2',~~ p_n''=nc_1+n^2a_2'',~~ p_n'''=nc_1+n^2a_2''',$\vspace{1mm}\\
	$~~q_n=-nc_1+n^2a_2,~~ q_n'=-nc_1+n^2a_2', ~~q_n''=-nc_1+n^2a_2'', ~~q_n'''=-nc_1+n^2a_2'''.$\\
	
	\noindent One finds that $-\frac{1}{2}<f_1,f_1',f_1'',f_1'''\leq 1,~ -1<g_1,g_1',g_1'',g_1'''\leq \frac{1}{2}$ and \vspace{1mm}\\$ -1<p_1,p_1',p_1'',p_1''',q_1,q_1',q_1'',q_1'''\leq 1$. Also note that $g_n\leq f_n,~ g_n'\leq f_n',~g_n''\leq f_n''. $
	
	\subsection{$c_2\equiv 0 \pmod 1$, $a=d=1$}\label{subsec1}
	
	\begin{theorem}\label{thm2}
		Let $Q(x_1,x_2, \cdots, x_5)$ be a real indefinite quadratic form of type $(1,4)$ and of determinant $D\neq 0$. Suppose that $c_2 \equiv 0 \pmod1$, and $a=1=d$. Then given any real numbers $c_1,c_2, \cdots,c_5$ there exist $(x_1,x_2,\cdots,x_5) \equiv (c_1,c_2, \cdots,c_5) \pmod 1$ satisfying \begin{equation}\label{eqb} 0<Q(x_1+c_1,x_2+c_2, \cdots,x_5+c_5)\leq d=(8 |D|)^{1/5}.\end{equation}
		Equality occurs in $(\ref{eqb})$ if and only if  $Q \sim \rho Q_2$ or $\rho Q_3, \rho>0$, where $Q_2=x_1x_2-(x_3^2+x_4^2+x_5^2)-(x_3x_4+x_3x_5+x_4x_5)$ and  $ Q_3=(x_1+\frac{1}{2}x_2)x_2-(x_3^2+x_4^2+x_5^2)-(x_3x_4+x_3x_5+x_4x_5)$. For $Q_2$ equality occurs in $(\ref{eqb})$ if and only if $(c_1,c_2,c_3,c_4,c_5)=(0,0,0,0,0)$; and for  $Q_3$ equality occurs in $(\ref{eqb})$ if and only if $(c_1,c_2,c_3,c_4,c_5)=(\frac{1}{2},0,0,0,0)$.
	\end{theorem}
	
	\noindent \textbf{Proof:} Here $ Q=(x_1+a_2x_2+\cdots )x_2-(x_3+h_4x_4+h_5x_5)^2-A(x_4+\lambda x_5)^2-tx_5^2. $ Take  $x_1=x+c_1, x_2=\pm 1, x_3=y+c_3$ and choose $(x_4,x_5)\equiv (c_4,c_5) \pmod 1$ arbitrarily  to get $Q=\pm x+\beta_1 y-y^2+\nu$ where $\beta_1=\pm a_3-2(c_3+h_4c_4+h_5c_5)$ and $\nu$ is some real number. Taking $h=1, k=1$ in Lemma \ref{lem5}, we find that (\ref{eq10}) is soluble with strict inequality unless $\beta_1 \equiv 1 \pmod 1$. Take $x_4=c_4$ and $1+c_4$ simultaneously, we get $2h_4 \equiv 0 \pmod 1$, i.e $h_4 = 0 {\rm ~or~} \frac{1}{2}$. Similarly we get $h_5=0$ or $\frac{1}{2}$.  Note that $h_4 \neq 0$ because if $h_4=0$ then $1=a\leq A+ah_4^2=A<\sqrt{\frac{2d^5}{3a}}=\sqrt{\frac{2}{3}}$, a contradiction. \vspace{2mm}
	
	\noindent Now take  $x_1=x+c_1, x_2=\pm 1, x_4=y+c_4, (x_3, x_5)\equiv (c_3,c_5) \pmod 1$ arbitrarily. Then $Q= \pm x +\beta_2 y-(A+\frac{1}{4})y^2+\nu,$   where $  \beta_2 = \pm a_4-(c_3+c_4/2+h_5c_5)-2A(c_4+\lambda c_5)$ and $\nu$ is some real number. Again, take $h=1, ~k=1$. As $ \frac{3}{4}\leq A<\sqrt{\frac{2}{3}}$, one has $|1-(A+\frac{1}{4})|+\frac{1}{2}<d=1$. Therefore, by Lemma \ref{lem5}, (\ref{eq10}) is soluble with strict inequality  unless $A=\frac{3}{4}$ and $\beta_2 \equiv 1 \pmod 1$.  Further $h_5 \neq 0$ because if $h_5=0$ then $1=a= \min (\phi) \leq A\lambda^2+t+h_5^2 a \leq A/4+d^5/2aA=3/16+2/3<1,$ a contradiction. Therefore, $$h_4=h_5=\frac{1}{2},~a=d=1,~ A=\frac{3}{4},~ t=\frac{d^5}{2aA}=\frac{2}{3}$$
	and $Q=(x_1+\cdots)x_2-(x_3+x_4/2+x_5/2)^2-\frac{3}{4}(x_4+\lambda x_5)^2-\frac{2}{3}x_5^2$. Take $x_1=x+c_1, x_2=\pm 1, x_5=y+c_5, (x_3, x_4)\equiv (c_3,c_4) \pmod 1$ arbitrarily, to get
	$$Q=\pm x+\beta_3 y-Cy^2+\nu, {~\rm ~ with~} C=\frac{3}{4}\lambda^2+\frac{11}{12},$$
	where $\beta_3 =\pm a_5-(c_3+c_4/2+c_5/2)-\frac{3\lambda}{2}(c_4+\lambda c_5)-\frac{4}{3}c_5$ and $\nu$ is some real number. Here $\frac{11}{12} \leq C \leq \frac{3}{16}+\frac{11}{12}.$ Taking $h=1,~ k=1$ we have $|1-C|+\frac{1}{2} < d=1$ and hence from Lemma \ref{lem5}, equation (\ref{eq10}) is soluble with strict inequality  unless $C=1$, i.e., unless $\lambda =\frac{1}{3}$  and $\beta_3 \equiv 1 \pmod 1$.  Therefore equation (\ref{eq10}) is soluble if \begin{equation}\label{eqc}0<Q=(x_1+a_2x_2+\cdots)x_2-(x_3^2+x_4^2+x_5^2)-(x_3x_4+x_3x_5+x_4x_5)\leq 1 \end{equation} is soluble. From $\beta_1, \beta_2, \beta_3 \equiv 0 \pmod 1$, we get $$\pm a_3-2c_3-c_4-c_5 \equiv 0\hspace{-3mm} \pmod 1,$$
	$$\pm a_4-2c_4-c_3-c_5 \equiv 0\hspace{-3mm} \pmod 1,$$
	$$\pm a_5-2c_5-c_3-c_4 \equiv 0\hspace{-3mm} \pmod 1.$$
	Thus, we get $2a_3 \equiv 0 \pmod 1,~ 2a_4 \equiv 0 \pmod 1, {\rm ~and~} 2a_5 \equiv 0 \pmod 1,$ and hence $a_3=0$ or $ \frac{1}{2},~ a_4=0$ or $\frac{1}{2},~ a_5=0$ or $\frac{1}{2}.$ Because of symmetry in $x_3,x_4,x_5$ we have $(a_3,a_4,a_5) = (0,0,0), (\frac{1}{2},0,0), (\frac{1}{2},\frac{1}{2},0)$ or $ (\frac{1}{2},\frac{1}{2},\frac{1}{2}).$
	By using unimodular transformation, $x_i \rightarrow -x_i,$ $1\leq i \leq 5$,  if necessary, we need to discuss the following cases:\vspace{2mm}
	
	\noindent For $(a_3,a_4,a_5)= (0,0,0)$;~ $(c_3,c_4,c_5) \equiv (0,0,0), (\frac{1}{2},\frac{1}{2},\frac{1}{2}), (\frac{1}{4},\frac{1}{4},\frac{1}{4}).$\\
	For $(a_3,a_4,a_5)= (\frac{1}{2},0,0)$;~ $(c_3,c_4,c_5) \equiv (\frac{-3}{8},\frac{1}{8},\frac{1}{8}), (\frac{-1}{8}, \frac{3}{8}, \frac{3}{8}).$\\
	For $(a_3,a_4,a_5)= (\frac{1}{2},\frac{1}{2},0)$;~ $(c_3,c_4,c_5) \equiv (0,0,\frac{1}{2}),(\frac{1}{2},\frac{1}{2},0), (\frac{1}{4},\frac{1}{4},\frac{-1}{4}).$\\
	For $(a_3,a_4,a_5)=(\frac{1}{2},\frac{1}{2},\frac{1}{2})$;~ $(c_3,c_4,c_5) \equiv (\frac{1}{8},\frac{1}{8},\frac{1}{8}), (\frac{3}{8}, \frac{3}{8}, \frac{3}{8}).$\vspace{2mm}
	
	\noindent \textbf{Case 1}: $(a_3,a_4,a_5)= (0,0,0).$\vspace{2mm}
	
	\noindent Here, $(\ref{eqc}) $ reduces to
	\begin{equation}\label{eq005}
	0<Q=(x_1+a_2x_2)x_2 - (x_3^2+x_4^2+x_5^2) -(x_3x_4+x_3x_5+x_4x_5)\leq 1.
	\end{equation}

	\noindent \textbf{Subcase (i)}: $(c_3,c_4,c_5)\equiv (0,0,0).$\vspace{2mm}
	
	\noindent Let $x_1=x+c_1(x \in \mathbb{Z}),$ $x_2=1$, $x_3=x_4=x_5=0$. Then  (\ref{eq005}) reduces to  $0<Q=(x+c_1)+a_2\leq 1.$  Taking $x=0,{\rm~and~}1$  we see that the (\ref{eq005}) is soluble  with strict inequality unless $c_1+a_2 \equiv 0 \pmod 1.$ \vspace{2mm}
	
	\noindent Similarly taking $ x_1=-x+c_1$ (for $x=0,1$),  $x_2=-1, x_3=x_4=x_5=0$ we see that (\ref{eq005}) is soluble with strict inequality unless $  -c_1+a_2 \equiv 0 \pmod 1.$ On solving these we get $(a_2,c_1)=(0,0)$ or  $(a_2,c_1)=(\frac{1}{2},\frac{1}{2})$.\vspace{2mm}
	
	\noindent For $(a_2,c_1)=(0,0)$, the quadratic form $Q$ reduces to $Q_2=x_1x_2 -(x_3^2+x_4^2+x_5^2)-(x_3x_4+x_3x_5+x_4x_5).$ For this quadratic form $Q_2$ and $(c_1,c_2,c_3,c_4,c_5)=(0,0,0,0,0)$ equality is required in (\ref{eq005}) as for $x_i \in \mathbb{Z}$, $Q_2$ takes integral values only. \vspace{2mm}
	
	\noindent For $(a_2,c_1)=(\frac{1}{2},\frac{1}{2})$, the quadratic form $Q$ reduces to $Q_3=(x_1+x_2/2)x_2-(x_3^2+x_4^2+x_5^2)-(x_3x_4+x_3x_5+x_4x_5).$ For this quadratic form $Q_3$ and  $(c_1,c_2,c_3,c_4,c_5)=(\frac{1}{2},0,0,0,0)$ equality is required in (\ref{eq005}). This is so because for $x_i \in \mathbb{Z}$, we see that $8Q_3(x_1+\frac{1}{2},x_2,\cdots,c_5)=4(2x_1+1+x_2)x_2 -8(x_3^2+x_4^2+x_5^2)-8(x_3x_4+x_3x_5+x_4x_5) \equiv 0 \pmod 8.$\vspace{2mm}
	
	\noindent \textbf{Subcase (ii)}: $(c_3,c_4,c_5)\equiv (\frac{1}{2},\frac{1}{2},\frac{1}{2}).$\\
	\noindent Let $x_1=x+c_1(x \in \mathbb{Z}),$ $x_2=1$, $x_3=x_4=x_5=\frac{1}{2}.$ Then (\ref{eq005}) reduces to $0<Q=(x+c_1)+a_2-\frac{3}{2}\leq 1$.  Taking $x=1,2, {\rm~and~}3$ simultaneously, (\ref{eq005}) is soluble  with strict inequality unless $c_1+a_2 \equiv \frac{1}{2} \pmod 1.$ \vspace{2mm}
	
	\noindent Similarly taking $ x_1=-x+c_1$ (for $x=1,2,3$),  $x_2=-1, x_3=x_4=x_5=\frac{1}{2}$, (\ref{eq005}) is soluble  with strict inequality unless $  -c_1+a_2 \equiv \frac{1}{2} \pmod 1.$\vspace{2mm}
	
	\noindent On solving these we get $(a_2,c_1)=(\frac{1}{2},0)$ or  $(a_2,c_1)=(0,\frac{1}{2})$. When $(a_2,c_1)=(\frac{1}{2},0)$, $ Q(0,2,\frac{1}{2},\frac{1}{2},\frac{1}{2})=\frac{1}{2}$ and when $(a_2,c_1)=(0,\frac{1}{2})$, $ Q(\frac{1}{2},4,\frac{1}{2},\frac{1}{2},\frac{1}{2})=\frac{1}{2}$. Therefore (\ref{eq005}) is soluble with strict inequality.\vspace{2mm}

	\noindent \textbf{Subcase (iii)}: $(c_3,c_4,c_5)\equiv (\frac{1}{4},\frac{1}{4},\frac{1}{4}).$\\
	\noindent Let $x_1=x+c_1(x \in \mathbb{Z}),$ $x_2=1$, $x_3=x_4=x_5=\frac{1}{4}.$ Then (\ref{eq005}) reduces to $0<Q=x+c_1+a_2-\frac{3}{8}\leq 1$.  Taking $x=0,1, {\rm~and~}2$ simultaneously, (\ref{eq005}) is soluble  with strict inequality unless $c_1+a_2 \equiv \frac{3}{8} \pmod 1.$ \vspace{2mm}
	
	\noindent Similarly taking $ x_1=-x+c_1$ (for $x=0,1,2$),  $x_2=-1, x_3=x_4=x_5=\frac{1}{4}$, (\ref{eq005}) is soluble  with strict inequality unless $  -c_1+a_2 \equiv \frac{3}{8} \pmod 1.$\vspace{2mm}
	
	\noindent On solving these we get $(a_2,c_1) =(\frac{3}{8},0)$ or  $(a_2,c_1)=(\frac{-1}{8},\frac{1}{2}).$ For $(a_2,c_1) =(\frac{3}{8},0)$, $Q (1,2,\frac{-3}{4},\frac{-3}{4},\frac{-3}{4})=\frac{1}{8}$ and when $(a_2,c_1)=(\frac{-1}{8},\frac{1}{2})$,  $Q(\frac{3}{2},4,\frac{-3}{4},\frac{-3}{4},\frac{-3}{4})=\frac{5}{8}$. Therefore (\ref{eq005}) is soluble with strict inequality.\vspace{2mm}
	
	\noindent \textbf{Case 2}: $(a_3,a_4,a_5)\equiv (\frac{1}{2},0,0).$\vspace{2mm}
	
	\noindent Here, $(\ref{eqc}) $ reduces to
	\begin{equation}\label{eq002}
	0<Q=(x_1+a_2x_2+\frac{1}{2}x_3)x_2 - (x_3^2+x_4^2+x_5^2) -(x_3x_4+x_3x_5+x_4x_5)\leq 1.
	\end{equation}

	\noindent \textbf{Subcase (i)}:  $(c_3,c_4,c_5)\equiv (\frac{-3}{8},\frac{1}{8},\frac{1}{8}).$\\
	Taking $x_1=x+c_1, x_2=1, x_3=\frac{-3}{8}, x_4=x_5=\frac{1}{8}$, and proceeding as in case 1, (\ref{eq002}) is soluble with strict inequality unless $c_1+a_2 \equiv \frac{9}{32} \pmod 1$. And on taking $x_1=-x+c_1, x_2=-1, x_3=\frac{-3}{8}, x_4=x_5=\frac{1}{8}$, (\ref{eq002}) is soluble with strict inequality unless $-c_1+a_2 \equiv \frac{-3}{32} \pmod 1$.
	On solving these we get $(a_2,c_1)=(\frac{3}{32},\frac{3}{16})$, and $(a_2,c_1)= (\frac{19}{32},\frac{-5}{16}) $. When $(a_2,c_1)=(\frac{3}{32},\frac{3}{16})$,  $Q (\frac{3}{16},3,\frac{5}{8},\frac{-7}{8},\frac{-7}{8})=\frac{3}{4}$ and when $(a_2,c_1)= (\frac{19}{32},\frac{-5}{16})$, $Q (\frac{-37}{16},-2,\frac{-3}{8},\frac{-15}{8},\frac{-7}{8})=\frac{9}{32}.$ Therefore (\ref{eq002}) has a solution with strict inequality.\vspace{2mm}

	\noindent \textbf{Subcase (ii)}:  $(c_3,c_4,c_5)\equiv (\frac{-1}{8},\frac{3}{8},\frac{3}{8}).$\\
	Taking $x_1=x+c_1, x_2=1, x_3=\frac{-1}{8}, x_4=x_5=\frac{3}{8}$, (\ref{eq002}) is soluble with strict inequality unless $c_1+a_2 \equiv \frac{13}{32} \pmod 1$. And on taking $x_1=-x+c_1, x_2=-1, x_3=\frac{-1}{8}, x_4=x_5=\frac{3}{8}$, (\ref{eq002}) is soluble with strict inequality unless $-c_1+a_2 \equiv \frac{9}{32} \pmod 1$. On solving these we get $(a_2,c_1)=(\frac{11}{32},\frac{1}{16})$, and $(a_2,c_1)=(\frac{-5}{32},\frac{-7}{16})$. When $(a_2,c_1)=(\frac{11}{32},\frac{1}{16})$,  $Q (\frac{17}{16},3,\frac{-9}{8},\frac{-5}{8},\frac{-5}{8})=\frac{3}{4}$ and when $(a_2,c_1)= (\frac{-5}{32},\frac{-7}{16})$, $Q (\frac{41}{16},3,\frac{-9}{8},\frac{-5}{8},\frac{-5}{8})=\frac{3}{4}.$ Therefore (\ref{eq002}) has a solution with strict inequality.\vspace{2mm}

	\noindent \textbf{Case 3}: $(a_3,a_4,a_5) \equiv (\frac{1}{2},\frac{1}{2},0).$\vspace{2mm}
	
	\noindent Here, $(\ref{eqc}) $ reduces to
	\begin{equation}\label{eq003}
	0<Q= (x_1+a_2x_2+\frac{1}{2}x_3+\frac{1}{2}x_4)x_2 - (x_3^2+x_4^2+x_5^2)-(x_3x_4+x_3x_5+x_4x_5)
	\leq 1.\end{equation}
	
	\noindent \textbf{Subcase (i)}: $(c_3,c_4,c_5) \equiv (0,0,\frac{1}{2}).$\\
	\noindent Taking $x_1=x+c_1, x_2=1, x_3=x_4=0, x_5=\frac{1}{2}$, (\ref{eq003}) is soluble with strict inequality unless $c_1+a_2 \equiv \frac{1}{4} \pmod 1$. And on taking $x_1=-x+c_1, x_2=-1, x_3=x_4=0, x_5=\frac{1}{2}$, (\ref{eq003}) is soluble unless $-c_1+a_2 \equiv \frac{1}{4} \pmod 1$.
	On solving these we get $(a_2,c_1)=(\frac{1}{4},0)$, and $(a_2,c_1)=(\frac{-1}{4},\frac{1}{2})$. When $(a_2,c_1)=(\frac{1}{4},0)$,  $Q(-1,-2,-1,-1,\frac{-1}{2})=\frac{3}{4}$ and when $(a_2,c_1)=(\frac{-1}{4},\frac{1}{2})$, $Q(\frac{-3}{2},-2,1,-1,\frac{-1}{2})=\frac{3}{4}$. Therefore (\ref{eq003}) has a solution with strict inequality.\vspace{2mm}
	
	\noindent \textbf{Subcase (ii)}: $(c_3,c_4,c_5) \equiv (\frac{1}{2},\frac{1}{2},0).$\\
	Taking $x_1=x+c_1, x_2=1, x_3=x_4=\frac{1}{2}, x_5=0$, (\ref{eq003}) is soluble with strict inequality unless $c_1+a_2 \equiv \frac{1}{4} \pmod 1$. And on taking $x_1=-x+c_1, x_2=-1, x_3=x_4=\frac{1}{2}, x_5=0$, (\ref{eq003}) is soluble with strict inequality unless $-c_1+a_2 \equiv \frac{1}{4} \pmod 1$. On solving these we get $(a_2,c_1)=(\frac{1}{4},0)$, and $(a_2,c_1)=(\frac{-1}{4},\frac{1}{2})$. When $(a_2,c_1)=(\frac{1}{4},0)$,  $Q(-2,-2,\frac{-3}{2},\frac{-1}{2},-1)=\frac{3}{4}$ and when $(a_2,c_1)=(\frac{-1}{4},\frac{1}{2})$, $Q(\frac{-3}{2},-2,\frac{-1}{2},\frac{-1}{2},-1)=\frac{1}{4}$. Therefore (\ref{eq003}) has a solution with strict inequality.\vspace{2mm}
	
	\noindent \textbf{Subcase (iii)}: $(c_3,c_4,c_5) \equiv (\frac{1}{4},\frac{1}{4},\frac{-1}{4}).$\\
	Taking $x_1=x+c_1, x_2=1, x_3=x_4=\frac{1}{4}, x_5=\frac{-1}{4}$, (\ref{eq003}) is soluble with strict inequality unless $c_1+a_2 \equiv \frac{-1}{8} \pmod 1$. And on taking $x_1=-x+c_1, x_2=-1, x_3=x_4=\frac{1}{4}, x_5=\frac{-1}{4}$, (\ref{eq003}) is soluble with strict inequality unless $-c_1+a_2 \equiv \frac{3}{8} \pmod 1$. On solving these we get $(a_2,c_1)=(\frac{-3}{8},\frac{1}{4})$, and $(a_2,c_1)=(\frac{1}{8},\frac{-1}{4})$. When $(a_2,c_1)=(\frac{-3}{8},\frac{1}{4})$,  $Q(\frac{13}{4},2,\frac{5}{4},\frac{-3}{4},\frac{-9}{4})=\frac{3}{8}$ and when $(a_2,c_1)=(\frac{1}{8},\frac{-1}{4})$, $Q(\frac{-5}{4},-2,\frac{5}{4},\frac{-3}{4},\frac{-5}{4})=\frac{3}{8}$. Therefore (\ref{eq003}) has a solution with strict inequality.\vspace{2mm}

	\noindent \textbf{Case 4}: $(a_3,a_4,a_5) \equiv (\frac{1}{2},\frac{1}{2},\frac{1}{2}).$\vspace{2mm}
	
	\noindent Here, $(\ref{eqc}) $ reduces to
	\begin{equation}\label{eq004}
	0<Q=(x_1+a_2x_2+\frac{1}{2}x_3+\frac{1}{2}x_4+\frac{1}{2}x_5)x_2 - (x_3^2+x_4^2+x_5^2) -(x_3x_4+x_3x_5+x_4x_5) \leq 1.\end{equation}
	
	\noindent \textbf{Subcase (i)}: $(c_3,c_4,c_5) \equiv (\frac{1}{8},\frac{1}{8},\frac{1}{8}).$\\
	\noindent Taking $x_1=x+c_1, x_2=1, x_3=x_4=x_5=\frac{1}{8}$, (\ref{eq004}) is soluble with strict inequality unless $c_1+a_2 \equiv \frac{-3}{32} \pmod 1$. And on taking $x_1=-x+c_1, x_2=-1, x_3=x_4=x_5=\frac{1}{8}$, (\ref{eq004}) is soluble with strict inequality unless $-c_1+a_2 \equiv \frac{9}{32} \pmod 1$. On solving these we get $(a_2,c_1)=(\frac{3}{32},\frac{-3}{16})$, and $(a_2,c_1)=(\frac{-13}{32},\frac{5}{16})$. When $(a_2,c_1)=(\frac{3}{32},\frac{-3}{16})$,  $Q(\frac{-19}{16},-2,\frac{-7}{8},\frac{-7}{8},\frac{-7}{8})=\frac{25}{32}$ and when $(a_2,c_1)=(\frac{-13}{32},\frac{5}{16})$, $Q(\frac{-27}{16},-2,\frac{9}{8},\frac{-7}{8},\frac{-7}{8})=\frac{25}{32}$. Therefore (\ref{eq004}) has a solution with strict inequality.\vspace{2mm}
	
	\noindent \textbf{Subcase (ii)}: $(c_3,c_4,c_5) \equiv (\frac{3}{8},\frac{3}{8},\frac{3}{8}).$\\
	Taking $x_1=x+c_1, x_2=1, x_3=x_4=x_5=\frac{3}{8}$, (\ref{eq004}) is soluble with strict inequality unless $c_1+a_2 \equiv \frac{9}{32} \pmod 1$. And on taking $x_1=-x+c_1, x_2=-1, x_3=x_4=x_5=\frac{3}{8}$, (\ref{eq004}) is soluble with strict inequality unless $-c_1+a_2 \equiv \frac{-13}{32} \pmod 1$. On solving these we get $(a_2,c_1)=(\frac{-1}{16},\frac{11}{32})$, and $(a_2,c_1)=(\frac{7}{16},\frac{-5}{32})$. When $(a_2,c_1)=(\frac{-1}{16},\frac{11}{32})$,  $Q(\frac{-53}{32},-1,\frac{3}{8},\frac{3}{8},\frac{3}{8})=\frac{3}{16}$ and when $(a_2,c_1)=(\frac{7}{16},\frac{-5}{32})$, $Q(\frac{-37}{32},-1,\frac{3}{8},\frac{3}{8},\frac{3}{8})=\frac{3}{16}$. Therefore (\ref{eq004}) has a solution with strict inequality.\vspace{2mm}
	
	\noindent This completes the proof of Theorem \ref{thm2}.\hfill $\square$
	
	\subsection{$c_2 \equiv 0 \pmod1$ and $a=\frac{1}{2},$ $d\leq 1$, $h_4=h_5=0$, $(a_3,c_3)=(0,\frac{1}{2})$ or $(\frac{1}{2},0)$.}\label{subsec5.2}
	
	\begin{theorem}\label{thm3}
		Let $Q(x_1,x_2, \cdots, x_5)$ be a real indefinite quadratic form of type $(1,4)$ and of determinant $D\neq 0$ as given in $(\ref{eq10})$. Suppose that $c_2 \equiv 0 \pmod1$ and $a=\frac{1}{2},$ $d\leq 1$, $h_4=h_5=0$, $(a_3,c_3)=(0,\frac{1}{2})$ or $(\frac{1}{2},0)$. Then given any real numbers $c_1,c_2, \cdots,c_5$ there exists $(x_1,x_2,\cdots,x_5) \equiv (c_1,c_2, \cdots,c_5) \pmod 1$ such that \begin{equation}\label{eqd}0<Q(x_1+c_1,x_2+c_2, \cdots,x_5+c_5)\leq (8 |D|)^{1/5}.\end{equation} Equality occurs in $(\ref{eqd})$ if and only if  $Q \sim \rho Q_4, \rho Q_5$ or $\rho Q_6$, $\rho >0$, where $Q_4=(x_1+\frac{1}{2}x_3+\frac{1}{2}x_4)x_2-\frac{1}{2}x_3^2-\frac{1}{2}x_4^2-2x_5^2$,   $Q_5=(x_1+\frac{1}{2}x_2+\frac{1}{2}x_3+\frac{1}{2}x_4)x_2-\frac{1}{2}x_3^2-\frac{1}{2}x_4^2-2x_5^2$ and  $Q_6= (x_1+\frac{1}{2}x_2+\frac{1}{2}x_3)x_2-\frac{1}{2}x_3^2-x_4^2-x_5^2$.  Further for $Q_4$, $Q_5$ and $Q_6$ equality holds in $(\ref{eqd})$ if  and only if $(c_1,c_2,c_3,c_4,c_5)=(0,0,0,0,0)$ , $(c_1,c_2,c_3,c_4,c_5)=(\frac{1}{2},0,0,0,0)$ or $(c_1,c_2,c_3,c_4,c_5)=(0,0,0,\frac{1}{2},\frac{1}{2})$ respectively.
	\end{theorem}
	
	\noindent \textbf{Proof:} Here   the inequality   (\ref{eq10}) reduces to
	\begin{equation}\label{eqe}0<Q(x_1, \cdots, x_5)=(x_1+a_2x_2+\cdots+a_5x_5)x_2-\frac{1}{2}x_3^2-A(x_4+\lambda x_5)^2-tx_5^2\leq d.\end{equation}
	
	\noindent To fix the values of $A$ and $t$, we will apply Macbeath's Lemma (Lemma \ref{lem5}) a number of times. \vspace{2mm}
	
	\noindent To fix $A$,  we will choose $(x_3,x_5) \equiv (c_3,c_5) \pmod 1$ arbitrarily , $x_2 = \pm 1$, $x_1=x+c_1$ and $x_4=y+c_4$, $x,y \in \mathbb{Z}$ to get \begin{equation}\label{eqf}Q=\pm x + \beta_{A,Q} y-Ay^2+\nu {\rm ~~ where~~} \beta_{A,Q} = \pm a_4-2A(c_4+\lambda x_5).\end{equation} Then we will choose  integers $2h_A$ and $k_A$ suitably such that \begin{equation}\label{eqg}|h_A-Ak_A^2|+\frac{1}{2}<d.\end{equation} Now by Macbeath's Lemma, (\ref{eqe}) is soluble  unless \begin{equation}\label{eqh}A=\frac{h_A}{k_A^2} {\rm ~~ and ~~}\beta_{A,Q}  \equiv h_A/k_A \mod{(1/k_A,2A)}. \end{equation}
	
	\noindent To fix $C=t+A\lambda^2$, we will choose $(x_3,x_4) \equiv (c_3,c_4) \pmod 1$ arbitrarily , $x_2 = \pm 1$, $x_1=x+c_1$ and $x_5=y+c_5$, $x,y \in \mathbb{Z}$ to get \begin{equation}\label{eqi}Q=\pm x + \beta_{C,Q} y-Cy^2+\nu {\rm ~~ where~~} \beta_{C,Q} = \pm a_5-2A\lambda(x_4+\lambda c_5)-2tc_5.\end{equation} Then we will choose  integers $2h_C$ and $k_C$ suitably such that \begin{equation}\label{eqj}|h_C-Ck_C^2|+\frac{1}{2}<d.\end{equation} Now by Macbeath's Lemma, (\ref{eqe}) is soluble  unless \begin{equation}\label{eqk}C=\frac{h_C}{k_C^2} {\rm ~~ and ~~}\beta_{C,Q} \equiv h_C/k_C \mod{(1/k_C,2C)}. \end{equation}
	
	\noindent When $\lambda=0$, we have $C=t$. Then the above argument will fix $t$. \vspace{2mm}
	
	\noindent Depending upon the ranges of $A$ and $t$, the proof of the theorem is divided into following cases (we call them as Lemmas) :
	\begin{enumerate} \item $A<d$, $t<d$
		\item $A<d$, $d\leq t < d+\frac{1}{2}$
		\item $A<d$, $ t\geq d+\frac{1}{2}$, $d<1$
		\item $A<d$, $ t\geq d+\frac{1}{2}$ and $d=1$
		\item $A\geq d$.
	\end{enumerate}
	\begin{lemma}\label{lem10} If $A<d$ and $t<d$, $(\ref{eqe})$ is soluble with strict inequality.\end{lemma}
	\noindent{\bf Proof:} Since $h_4=0$, (\ref{eq9}) implies $A \geq a=\frac{1}{2}$. Take $h_A=\frac{1}{2}, k_A=1$ so that  (\ref{eqg}) is satisfied. Hence (\ref{eqe}) is soluble  unless $A=\frac{1}{2}$ and $\beta_A=\pm a_4-2A(c_4+\lambda x_5) \equiv \frac{1}{2} \pmod 1$. Since $-\frac{1}{2} < \lambda \leq \frac{1}{2}$, taking $x_5=c_5$ and $1+c_5$, we get $\lambda =0$, and hence $(a_4,c_4)=(0,\frac{1}{2})$ or $(\frac{1}{2},0)$.\vspace{2mm}
	
	\noindent Since $\lambda =0$, (\ref{eq6}) and (\ref{eq7}) implies $C=t+A\lambda^2=t \geq A=\frac{1}{2}$. Take $h_C=\frac{1}{2}, k_C=1$ so that  (\ref{eqj}) is satisfied. Hence (\ref{eqe}) is soluble  unless $C=t=\frac{1}{2}$ and $\beta_{t,Q}=\pm a_5-2tc_5 \equiv \frac{1}{2} \pmod 1$, i.e., unless $(a_5,c_5)=(0,\frac{1}{2})$ or $(\frac{1}{2},0)$.\vspace{2mm}
	
	\noindent Here \vspace{-2mm} $$ Q=(x_1+a_2x_2+\cdots+a_5x_5)x_2-\frac{1}{2}(x_3^2+x_4^2+x_5^2), ~~ d=(1/4)^{\frac{1}{5}}>0.7578.$$
	
	\noindent Because of symmetry in $x_3,x_4,x_5$ we need to consider the following cases
	\begin{enumerate}
		\item $(a_3,a_4,a_5)=(0,0,0), (c_3,c_4,c_5)=(\frac{1}{2},\frac{1}{2},\frac{1}{2})  $
		\item $(a_3,a_4,a_5)=(\frac{1}{2},0,0), (c_3,c_4,c_5)=(0,\frac{1}{2},\frac{1}{2}) $
		\item $(a_3,a_4,a_5)=(\frac{1}{2},\frac{1}{2},0), (c_3,c_4,c_5)=(0,0,\frac{1}{2}) $
		\item $(a_3,a_4,a_5)=(\frac{1}{2},\frac{1}{2},\frac{1}{2}), (c_3,c_4,c_5)=(0,0,0). $
	\end{enumerate}
	
	\noindent We work as in Lemma 11 of Dumir and Sehmi \cite{DumirSehmi} to find that (\ref{eqe}) is soluble with strict inequality. \hfill $\square$
	
	\begin{lemma}\label{lem11} If $A<d$ and $d\leq t<d+\frac{1}{2}$, $(\ref{eqe})$ is soluble with strict inequality.\end{lemma}
	
	\noindent{\bf Proof:} Working as in beginning of Lemma \ref{lem10}, we get   $A=\frac{1}{2}$,  $\lambda =0$, and  $(a_4,c_4)=(0,\frac{1}{2})$ or $(\frac{1}{2},0)$.\vspace{2mm}
	
	\noindent Since $a=\frac{1}{2}=A$, we get $t=\frac{d^5}{2aA}=2d^5$. Therefore $t\geq d$ implies $d^4 \geq \frac{1}{2}$ i.e., $d > 0.8408$. Here $C=t$.
	\noindent  Take $h_t=1, k_t=1$ so that  (\ref{eqj}) is satisfied. Hence (\ref{eqe}) is soluble  unless $t=1$ and $\beta_t=\pm a_5-2tc_5 \equiv 0 \pmod 1$,  i.e., unless $(a_5,c_5)=(0,0), (\frac{1}{2}, \pm \frac{1}{4}),$ or $ (0,\frac{1}{2})$.\vspace{2mm}
	
	\noindent For $t=1,$ we have $d=(1/2)^{1/5}>0.87055$. Thus $ d/a=2d> 1.7411$ which gives $m=1$, thereby $\delta_m/A=d/A=2d>1.7411$, i.e., $K=1$ and hence $\delta_{m,K}=d.$ Here \vspace{-2mm}  $$ Q=(x_1+a_2x_2+\cdots+a_5x_5)x_2-\frac{1}{2}(x_3^2+x_4^2)-x_5^2, ~~ d=(1/2)^{\frac{1}{5}}>0.87055.$$
	\noindent We work as in Lemma 6 of Raka and Rani \cite{RakaRani} to find that (\ref{eqe}) is soluble with strict inequality. \hfill $\square$
	
	\begin{lemma}\label{lem12} If $A<d$, $t\geq d+\frac{1}{2} $ and $d< 1$, $(\ref{eqe})$ is soluble with strict inequality.\end{lemma}
	
	\noindent{\bf Proof:} Working as in beginning of Lemma \ref{lem10}, we get   $a=A=\frac{1}{2}$,  $\lambda =0$, and  $(a_4,c_4)=(0,\frac{1}{2})$ or $(\frac{1}{2},0)$. Already we have $(a_3,c_3)=(0,\frac{1}{2})$ or $(\frac{1}{2},0)$. Because of symmetry in $x_3$ and $x_4$ we need to consider
	\begin{itemize}
		\item $(a_4,c_4)=(0,\frac{1}{2}), (a_3,c_3)$ may be $(0,\frac{1}{2})$ or $(\frac{1}{2},0)$
		\item $(a_4,c_4)=(\frac{1}{2},0), (a_3,c_3)=(\frac{1}{2},0).$
	\end{itemize}
	Since $a=\frac{1}{2}=A$, we get $t=\frac{d^5}{2aA}=2d^5$. Therefore $t\geq d+\frac{1}{2}$ implies $d > 0.9358$ and hence $t > 1.4358 $. Take $h_t=\frac{3}{2}, k_t=1$ so that  (\ref{eqj}) is satisfied. Hence (\ref{eqe}) is soluble  unless $t=\frac{3}{2}$ and $\beta_t=\pm a_5-2tc_5 \equiv \frac{3}{2} \pmod 1$  i.e., unless $(a_5,c_5)= (\frac{1}{2},0),(0,\pm \frac{1}{6})$, $(0,\frac{1}{2})$ or $(\frac{1}{2},\pm \frac{1}{3})$.\vspace{2mm}

	\noindent For $t=\frac{3}{2},$ we have $d=(3/4)^{1/5}>0.94408$. Thus $ d/a=2d> 1.88$ which gives $m=1$, thereby $\delta_m/A=d/A=2d>1.88$ i.e., $K=1$ and hence $\delta_{m,K}=d.$\vspace{2mm}
	
	\noindent \textbf{Case I}: $(a_5,c_5)=(\frac{1}{2},0)$. Taking $x_2=\pm 1$, the following table gives a solution to (\ref{eq17}) with strict inequality, i.e., of \vspace{-2mm} $$0<G=(x_1+a_2''x_2+\frac{1}{2}x_5)x_2 - \frac{3}{2}x_5^2 -\frac{1}{4}< d.\vspace{-2mm}$$
	
	{\scriptsize
		\begin{equation*}\begin{array}{llllll}
		\hline\\
		{\rm ~~~~Range} & x_1 & x_1x_2 & x_5&x_2x_5 & {\rm ~~G}\vspace{2mm}\\
		\hline\\
		
		\frac{1}{4}<f_1''\leq 1&c_1&|c_1|&0&0&f_1''-\frac{1}{4}\vspace{2mm}\\
		
		-\frac{1}{2}<f_1''<d-\frac{3}{4}&\pm 1+c_1&1+|c_1|&0&0&f_1''+\frac{3}{4}\vspace{2mm}\\
		
		{\rm When ~} d-\frac{3}{4}\leq f_1'' \leq \frac{1}{4}: \vspace{2mm}\\
		
		d-\frac{3}{4}\leq g_1''\leq \frac{1}{4} & c_1 & 2|c_1| & 0 & 0 & f_2''-\frac{1}{4} \vspace{2mm}\\
		
		-\frac{3}{4}<g_1''<d-\frac{3}{4}&\pm 1+c_1&1-|c_1|&0&0&g_1''+\frac{3}{4}\vspace{2mm}\\
		
		-1<g_1''\leq -\frac{3}{4}&\pm 1+c_1&2+2|c_1|&\pm 1&1&f_2''+\frac{3}{4}\vspace{2mm}\\
		
		\hline
		\end{array}
		\end{equation*}\vspace{1mm}}
	
	\noindent \textbf{Case II}: $(a_5,c_5)=(0,\pm \frac{1}{6})$. Changing $x_1 \rightarrow -x_1, x_2 \rightarrow -x_2, x_5 \rightarrow -x_5$ if necessary, we can assume that $c_5=\frac{1}{6}$. The following table gives a solution to (\ref{eq17}) with strict inequality, i.e., of
	$$0<G=(x_1+a_2''x_2)x_2-\frac{3}{2}x_5^2-\frac{1}{4}< d.$$
	
	{\scriptsize
		\begin{equation*}\begin{array}{lllll}
		\hline\\
		{\rm ~~~~Range} & x_1 & x_1x_2 & x_5 & {\rm ~~G}\vspace{2mm}\\
		\hline\\
		
		\frac{7}{24}<f_1''\leq 1&c_1&|c_1|&\frac{1}{6}&f_1''-\frac{7}{24}\vspace{2mm}\\
		
		-\frac{1}{2}<f_1''<d-\frac{17}{24}&\pm 1+c_1&1+|c_1|&\frac{1}{6}&f_1''+\frac{17}{24}\vspace{2mm}\\
		
		{\rm When ~} d-\frac{17}{24}\leq f_1'' \leq \frac{7}{4}: \vspace{2mm}\\
		
		d-\frac{17}{24}\leq g_1''\leq \frac{7}{4} &c_1&2|c_1|&\frac{1}{6}&f_2''-\frac{7}{24}\vspace{2mm}\\
		
		-\frac{17}{24}<g_1''<d-\frac{17}{24}&\pm 1+c_1&1-|c_1|&\frac{1}{6}&g_1''+\frac{17}{24}\vspace{2mm}\\
		
		-1<g_1''\leq -\frac{17}{24}&\pm 1+c_1&2+2|c_1|&-\frac{5}{6}&f_2''+\frac{17}{24}\vspace{2mm}\\

		\hline
		\end{array}
		\end{equation*}\vspace{1mm}}
	
	\noindent \textbf{Case III}: $(a_5,c_5)=(\frac{1}{2},\pm \frac{1}{3})$.  Here already we have $(a_4,c_4)=(0,\frac{1}{2})$ or $(\frac{1}{2},0)$. \vspace{2mm}

	\noindent \textbf{Subcase (i)}:  $(a_4,c_4)=(0,\frac{1}{2}), (a_3,c_3)$ may be $(0,\frac{1}{2})$ or $(\frac{1}{2},0)$.\vspace{2mm}
	
	\noindent As $h_4=0, h_5=0$, by (\ref{eq12})  we have  $a'_4=a_4=0, a'_5=a_5=\frac{1}{2}$. Taking $x_5=1/3$, the following table gives a solution to (\ref{eq13}) with strict inequality, i.e., of \vspace{-2mm}
	$$ 0<F=(x_1+a_2'x_2+\frac{1}{2}x_5)x_2-\frac{1}{2}x_4^2-\frac{3}{2}x_5^2-\frac{1}{8}< d.$$

	{\scriptsize
		\begin{equation*}\begin{array}{lllll}
		\hline\\
		{\rm ~~~~Range} & x_1 & x_2 &x_4   & {\rm ~~F}\vspace{2mm}\\
		\hline\\
		
		\frac{1}{4}<p_1'\leq 1&c_1&1&\frac{1}{2}&p_1'-\frac{1}{4}\vspace{2mm}\\
		
		-\frac{3}{4}<p_1'<d-\frac{3}{4}&1+c_1&1&\frac{1}{2}&p_1'+\frac{3}{4}\vspace{2mm}\\
		
		d-\frac{7}{4}\leq p_1'\leq -\frac{3}{4}&-1+c_1&-1&\frac{1}{2}&q_1'+\frac{5}{12}\vspace{2mm}\\
		
		-1<p_1'<d-\frac{7}{4}&2+c_1&1&\frac{1}{2}&p_1'+\frac{7}{4}\vspace{2mm}\\

		{\rm When ~} d-\frac{3}{4}\leq p_1' \leq \frac{1}{4}: \vspace{2mm}\\

		\frac{7}{12}<q_1'\leq 1&c_1&-1&\frac{1}{2}&q_1'-\frac{7}{12}\vspace{2mm}\\
		
		d-\frac{5}{12}\leq q_1'\leq \frac{7}{12} &c_1&2&\frac{3}{2}&p_2'-\frac{13}{12}\vspace{2mm}\\
		
		-\frac{5}{12}<q_1'<d-\frac{5}{12}&-1+c_1&-1&\frac{1}{2}&q_1'+\frac{5}{12}\vspace{2mm}\\
		
		d-\frac{17}{12}\leq q_1'\leq -\frac{5}{12} &c_1&2&\frac{1}{2}&p_2'-\frac{1}{12}\vspace{2mm}\\
		
		-1<q_1'<d-\frac{17}{12}&-2+c_1&-1&\frac{1}{2}&q_1'+\frac{17}{12}\vspace{2mm}\\

		\hline
		\end{array}
		\end{equation*}\vspace{1mm}}

	\noindent \textbf{Subcase (ii)}:  $(a_4,c_4)=(\frac{1}{2},0), (a_3,c_3)=(\frac{1}{2},0).$ \vspace{2mm}
	
	\noindent Here, on taking $x_4=0$, the following table gives a solution to (\ref{eq13}) with strict inequality i.e., of\vspace{-2mm}
	 $$0<F=(x_1+a_2'x_2+\frac{1}{2}x_4+\frac{1}{2}x_5)x_2-\frac{1}{2}x_4^2-\frac{3}{2}x_5^2-\frac{1}{8}<d.$$
	
	{\scriptsize
		\begin{equation*}\begin{array}{lllll}
		\hline\\
		{\rm ~~~~Range} & x_1 & x_2  & x_5 & {\rm ~~F}\vspace{2mm}\\
		\hline\\
		
		\frac{1}{8}<p_1'\leq 1&c_1&1&\frac{1}{3}&p_1'-\frac{1}{8}\vspace{2mm}\\
		
		-\frac{7}{8}<p_1'<d-\frac{7}{8}&1+c_1&1&\frac{1}{3}&p_1'+\frac{7}{8}\vspace{2mm}\\
		
		d-\frac{15}{8}\leq p_1' \leq -\frac{7}{8}&-1+c_1&-1&\frac{1}{3}&q_1'+\frac{13}{24}\vspace{2mm}\\
		
		-1<p_1'<d-\frac{15}{8}&2+c_1&1&\frac{1}{3}&p_1'+\frac{15}{8}\vspace{2mm}\\

		{\rm When ~} d-\frac{7}{8}\leq p_1' \leq \frac{1}{8}: \vspace{2mm}\\
		\frac{11}{24}<q_1'\leq 1&c_1&-1&\frac{1}{3}&q_1'-\frac{11}{24}\vspace{2mm}\\
		d-\frac{13}{24}\leq q_1'\leq \frac{11}{24}&c_1&2&\frac{1}{3}&p_2'+\frac{1}{24}\vspace{2mm}\\
		
		 -\frac{13}{24}<q_1'<d-\frac{13}{24}&-1+c_1&-1&\frac{1}{3}&q_1'+\frac{13}{24}\vspace{2mm}\\
		d-\frac{37}{24}\leq q_1'\leq -\frac{13}{24}&1+c_1&2&-\frac{2}{3}&p_2'+\frac{13}{24}\vspace{2mm}\\
		
		-1<q_1'<d-\frac{37}{24}&-2+c_1&-1&\frac{1}{3}&q_1'+\frac{37}{24}\vspace{2mm}\\		
		\hline
		\end{array}
		\end{equation*}\vspace{1mm}}

	\noindent \textbf{Case IV}: $(a_5,c_5)=(0,\frac{1}{2})$. Here already we have  $(a_4,c_4)=(0,\frac{1}{2})$ or $(\frac{1}{2},0)$.\vspace{2mm}

	\noindent \textbf{Subcase (i)}:  $(a_4,c_4)=(0,\frac{1}{2}), (a_3,c_3)$ may be $(0,\frac{1}{2})$ or $(\frac{1}{2},0)$.\vspace{2mm}
	
	\noindent On taking $x_4=\frac{1}{2}=x_5$, the following table gives a solution to (\ref{eq13}) with strict inequality, i.e., of \vspace{-2mm} $$ 0<F=(x_1+a_2'x_2)x_2-\frac{1}{2}x_4^2-\frac{3}{2}x_5^2-\frac{1}{8}<d.$$
	
	{\scriptsize
		\begin{equation*}\begin{array}{llll}
		\hline\\
		{\rm ~~~~Range} & x_1 & x_1x_2  & {\rm ~~F}\vspace{2mm}\\
		\hline\\
		
		\frac{5}{8}<f_1'\leq 1&c_1&|c_1|&f_1'-\frac{5}{8}\vspace{2mm}\\
		
		-\frac{1}{2}<f_1'<d-\frac{3}{8}&\pm 1+c_1&1+|c_1|&f_1'+\frac{3}{8}\vspace{2mm}\\
		
		{\rm When ~} d-\frac{3}{8}\leq f_1'\leq \frac{5}{8}: \vspace{2mm}\\
		
		-\frac{3}{8}<g_1'\leq\frac{5}{8}&\pm 1+c_1&1-|c_1|&g_1'+\frac{3}{8}\vspace{2mm}\\
		
		-1<g_1'\leq -\frac{3}{8}&c_1&2|c_1|&f_2'-\frac{5}{8}\vspace{2mm}\\
		
		\hline
		\end{array}
		\end{equation*}\vspace{1mm}}
	
	\noindent \textbf{Subcase (ii)}:  $(a_4,c_4)=(\frac{1}{2},0), (a_3,c_3)=(\frac{1}{2},0)$. \vspace{2mm}
	
	\noindent The following table gives a solution to (\ref{eqe}) with strict inequality, i.e., of \vspace{-2mm}
	$$0<Q=(x_1+a_2x_2+\frac{1}{2}x_3+\frac{1}{2}x_4)x_2-\frac{1}{2}(x_3^2+x_4^2)-\frac{3}{2}x_5^2< d.$$

	{\scriptsize
		\begin{equation*}\begin{array}{lllllllll}
		\hline\\
		{\rm ~~~~Range} & x_1 & x_1x_2 &x_3&x_2x_3&x_4 &x_2x_4 & x_5 & {\rm ~~Q}\vspace{2mm}\\
		\hline\\
		
		\frac{3}{8}<f_1\leq 1 &c_1&|c_1|&0&0&0&0&\frac{1}{2}&f_1-\frac{3}{8}\vspace{2mm}\\
		
		-\frac{1}{2}<f_1<d-\frac{5}{8}&\pm 1+c_1&1+|c_1|&0&0&0&0&\frac{1}{2}&f_1+\frac{5}{8}\vspace{2mm}\\

		{\rm When ~} d-\frac{5}{8}\leq f_1\leq \frac{3}{8}: \vspace{2mm}\\

		d-\frac{5}{8}\leq g_1\leq \frac{3}{8}&\pm 2+c_1&4+2|c_1|&\pm 3&6&\pm 2&4&\frac{3}{2}&f_2-\frac{7}{8}\vspace{2mm}\\
		
		-\frac{5}{8}<g_1<d-\frac{5}{8}&\pm 1+c_1&1-|c_1|&0&0&0&0&\frac{1}{2}&g_1+\frac{5}{8}\vspace{2mm}\\
		
		-1<g_1\leq -\frac{5}{8}&\pm 1+c_1&2+2|c_1|&\pm 3&6&\pm 2&4&\frac{1}{2}&f_2+\frac{1}{8}\vspace{2mm}\\

		\hline
		\end{array}
		\end{equation*}}\vspace{2mm}
	
	\noindent This completes the proof of Lemma \ref{lem12}. \hfill $\square$
	
	\begin{lemma}\label{lem13} If $A<d$, $t\geq d+\frac{1}{2}$ and $d= 1$, $(\ref{eqe})$ is soluble with strict inequality unless $Q =\rho Q_4$ or $\rho Q_5, \rho>0$. For $Q_4$ and $Q_5$ equality holds in $(\ref{eqe})$ if  and only if $(c_1,c_2,c_3,c_4,c_5)=(0,0,0,0,0)$ and $(c_1,c_2,c_3,c_4,c_5)=(\frac{1}{2},0,0,0,0)$ respectively.\end{lemma}
	
	\noindent {\bf Proof:} Already we have $a=A=\frac{1}{2},~ \lambda=0=h_4=h_5$. Also $(a_3,c_3)$ and $(a_4,c_4)=$ $(0,\frac{1}{2})$ or $(\frac{1}{2},0)$. Because of symmetry in $x_3$ and $x_4$ we need to consider
	\begin{itemize}
		\item $(a_4,c_4)=(0,\frac{1}{2}), (a_3,c_3)$ may be $(0,\frac{1}{2})$ or $(\frac{1}{2},0)$
		\item $(a_4,c_4)=(\frac{1}{2},0), (a_3,c_3)=(\frac{1}{2},0).$
	\end{itemize} For $d=1$ we have $t=2d^5=2$.  Therefore, (\ref{eqe}) reduces to $$0<Q=(x_1+a_2x_2+a_3x_3+a_4x_4+a_5x_5)x_2-\frac{1}{2}x_3^2-\frac{1}{2}x_4^2-2x_5^2\leq 1.$$
	
	\noindent  Take $h_t=2, k_t=1$ so that  (\ref{eqj}) is satisfied. Hence (\ref{eqe}) is soluble  unless  $\beta_t=\pm a_5-4c_5 \equiv 0 \pmod 1$  i.e., unless $(a_5,c_5)=(0,0), (0,\frac{1}{2}), (0,\pm \frac{1}{4}),$ or $ (\frac{1}{2}, \pm \frac{1}{8})$.\vspace{2mm}
	
	\noindent Here, we have $d/a=2d=2$ which gives $m=1$, thereby $\delta_m/A=d/A=2d=2$ i.e., $K=1$ and hence $\delta_{m,K}=d=1$.\vspace{2mm}
	
	\noindent \textbf{Case I}: $(a_5,c_5)=(0,0)$. \vspace{2mm}

	\noindent \textbf{Subcase (i)}: $(a_4,c_4)=(0,\frac{1}{2}), (a_3,c_3)$ may be $(0,\frac{1}{2})$ or $(\frac{1}{2},0)$.\vspace{2mm}
	
	\noindent Let $x_1=x+c_1(x \in \mathbb{Z}),$ $x_2=1$, $x_4=\frac{1}{2}$, $x_5=0.$ Then (\ref{eq13}) reduces to $0<F=(x+c_1)+a_2'-\frac{1}{4}\leq 1$.  Taking $x=0,1 {\rm~and~}2$ simultaneously, we see that (\ref{eq13}) is soluble  with strict inequality unless $c_1+a_2' \equiv \frac{1}{4} \pmod 1.$ \vspace{2mm}
	
	\noindent Similarly taking $ x_1=-x+c_1$ (for $x=0,1,2$)  $x_2=-1, x_4=\frac{1}{2}, x_5=0$, (\ref{eq13}) is soluble  with strict inequality unless $  -c_1+a_2' \equiv \frac{1}{4} \pmod 1.$ On solving these we get $(a_2',c_1)=(\frac{1}{4},0)$ or  $(-\frac{1}{4},\frac{1}{2})$. When $(a_2',c_1)=(\frac{1}{4},0)$, $ F(1,2,\frac{1}{2},1)=\frac{3}{4}$ and when $(a_2',c_1)=(-\frac{1}{4},\frac{1}{2})$, $ F(\frac{3}{2},2,\frac{3}{2},1)=\frac{7}{8}$. Therefore (\ref{eq13}) is soluble with strict inequality.\vspace{2mm}
	
	\noindent \textbf{Subcase (ii)}: $(a_4,c_4)=(\frac{1}{2},0), (a_3,c_3)=(\frac{1}{2},0)$.\vspace{2mm}
	
	\noindent Here (\ref{eqe}) reduces to $$0<Q=(x_1+a_2x_2+\frac{1}{2}x_3+\frac{1}{2}x_4)x_2-\frac{1}{2}x_3^2-\frac{1}{2}x_4^2-2x_5^2\leq 1.$$
	Take $x_1=c_1$, $x_2=\pm1$, $x_3=x_4=x_5=0$ to get $Q=c_1+a_2$ or $-c_1+a_2$. If $0<c_1+a_2<1$ or $0<-c_1+a_2<1$, we have a solution with strict inequality in (\ref{eqe}). Further take $x_1=\pm1+c_1$, $x_2=\pm1$, $x_3=x_4=x_5=0$ to get $Q=1+c_1+a_2$ or $1-c_1+a_2$. If $-1<c_1+a_2<0$ or
	$-1<-c_1+a_2<0$, we have a solution with strict inequality in (\ref{eqe}). If $c_1+a_2 \equiv 0\pmod 1$ and $-c_1+a_2 \equiv 0\pmod 1$, we get $(a_2,c_1)=(0,0)$ or $(a_2,c_1)=(\frac{1}{2},\frac{1}{2})$.  \vspace{2mm}
	
	\noindent This gives $Q= Q_4= (x_1+\frac{1}{2}x_3+\frac{1}{2}x_4)x_2-\frac{1}{2}x_3^2-\frac{1}{2}x_4^2-2x_5^2$ or $ Q=Q_5=(x_1+\frac{1}{2}x_2+\frac{1}{2}x_3+\frac{1}{2}x_4)x_2-\frac{1}{2}x_3^2-\frac{1}{2}x_4^2-2x_5^2$. For $Q_4$ and $Q_5$ equality holds in (\ref{eqe}) for $(c_1,c_2,c_3,c_4,c_5)=(0,0,0,0,0)$ and $(c_1,c_2,c_3,c_4,c_5)=(\frac{1}{2},0,0,0,0)$ respectively. This is so because for integers $x_i$'s, $8Q_4(x_1,x_2,x_3,x_4,x_5)=8x_1x_2+2x_2^2-(2x_3-x_2)^2-(2x_4-x_2)^2-16x_5^2 \equiv 0 \pmod 8$ and $8Q_5(x_1+\frac{1}{2},x_2,x_3,x_4,x_5)=8x_1x_2+4x_2+6x_2^2-(2x_3-x_2)^2-(2x_4-x_2)^2-16x_5^2 \equiv 0 \pmod 8.$\vspace{2mm}
	
	\noindent \textbf{Case II}: $(a_5,c_5)=(0,\frac{1}{2})$. Here we have $(a_3,c_3)$ and $(a_4,c_4)=(0,\frac{1}{2})$ or $(\frac{1}{2},0)$.\vspace{2mm}
	
	\noindent \textbf{Subcase (i)}:  $(a_4,c_4)=(0,\frac{1}{2})$, $(a_3,c_3)$ may be $(0,\frac{1}{2})$ or $(\frac{1}{2},0)$.\vspace{2mm}
	
	\noindent Let $x_1=x+c_1(x \in \mathbb{Z}),$ $x_2=1$, $x_4=\frac{1}{2}=x_5.$ Then (\ref{eq13}) reduces to $0<F=(x+c_1)+a_2'-\frac{3}{4}\leq 1$.  Taking $x=0,1 {\rm~and~}2$ simultaneously, we see that (\ref{eq13}) is soluble  with strict inequality unless $c_1+a_2' \equiv \frac{3}{4} \pmod 1.$ \vspace{2mm}
	
	\noindent Similarly taking $ x_1=-x+c_1$ (for $x=0,1,2$)  $x_2=-1, x_4=\frac{1}{2}=x_5$, (\ref{eq13}) is soluble  with strict inequality unless $  -c_1+a_2' \equiv \frac{3}{4} \pmod 1.$ On solving these we get $(a_2',c_1)=(\frac{1}{4},\frac{1}{2})$ or  $(-\frac{1}{4},0)$. When $(a_2',c_1)=(\frac{1}{4},\frac{1}{2})$, $ F(\frac{1}{2},2,\frac{3}{2},\frac{1}{2})=\frac{1}{4}$ and when $(a_2',c_1)=(-\frac{1}{4},0)$, $ F(1,2,\frac{1}{2},\frac{1}{2})=\frac{1}{4}$. Therefore (\ref{eq13}) is soluble with strict inequality.\vspace{2mm}

	\noindent \textbf{Subcase (ii)}: $(a_4,c_4)=(\frac{1}{2},0), (a_3,c_3)=(\frac{1}{2},0)$.\vspace{2mm}

	\noindent Let $x_1=x+c_1(x \in \mathbb{Z}),$ $x_2=1$, $x_3=0=x_4,$ $x_5=\frac{1}{2}.$ Then (\ref{eqe}) reduces to $0<Q=(x+c_1)+a_2-\frac{1}{2}\leq 1$.  Taking $x=0 {\rm~and~}1$ simultaneously, we see that (\ref{eqe}) is soluble  with strict inequality unless $c_1+a_2 \equiv \frac{1}{2} \pmod 1.$ \vspace{2mm}
	
	\noindent Similarly taking $ x_1=-x+c_1$ (for $x=0,1,2$)  $x_2=-1, x_3=0=x_4, x_5=\frac{1}{2}$, (\ref{eqe}) is soluble  with strict inequality unless $  -c_1+a_2 \equiv \frac{1}{2} \pmod 1.$ On solving these we get $(a_2,c_1)=(0,\frac{1}{2})$ or  $(\frac{1}{2},0)$. When $(a_2,c_1)=(0,\frac{1}{2})$, $ Q(\frac{1}{2},2,0,0,\frac{1}{2})=\frac{1}{2}$ and when $(a_2,c_1)=(\frac{1}{2},0)$, $ Q(1,2,1,1,\frac{3}{2})=\frac{1}{2}$. Therefore (\ref{eqe}) is soluble with strict inequality.\vspace{2mm}

	\noindent \textbf{Case III}: $(a_5,c_5)=(0, \pm \frac{1}{4})$.\vspace{2mm}

	\noindent Changing $x_1 \rightarrow -x_1, x_2 \rightarrow -x_2, x_5 \rightarrow -x_5$ if necessary, we can assume that $c_5=\frac{1}{4}$.

	\noindent Let $x_1=x+c_1(x \in \mathbb{Z}),$ $x_2=1$, $x_5=\frac{1}{4}.$ Then (\ref{eq17}) reduces to $0<G=(x+c_1)+a_2''-\frac{3}{8}\leq 1$.  Taking $x=0 {\rm~and~}1$ simultaneously, we see that (\ref{eq17}) is soluble  with strict inequality unless $c_1+a_2'' \equiv \frac{3}{8} \pmod 1.$ \vspace{2mm}
	
	\noindent Similarly taking $ x_1=-x+c_1$ (for $x=0,1,2$)  $x_2=-1, x_5=\frac{1}{4}$, (\ref{eq17}) is soluble  with strict inequality unless $  -c_1+a_2'' \equiv \frac{3}{8} \pmod 1.$ On solving these we get $(a_2'',c_1)=(\frac{3}{8},0)$ or  $(-\frac{1}{8},\frac{1}{2})$. When $(a_2'',c_1)=(\frac{3}{8},0)$, $ G(0,2,-\frac{3}{4})=\frac{1}{8}$ and when $(a_2'',c_1)=(-\frac{1}{8},\frac{1}{2})$, $ G(\frac{7}{2},2,-\frac{7}{4})=\frac{1}{8}$. Therefore (\ref{eq17}) is soluble with strict inequality.\vspace{2mm}

	\noindent \textbf{Case IV}: $(a_5,c_5)=(\frac{1}{2}, \pm \frac{1}{8})$.\vspace{2mm}
	
	\noindent Changing $x_1 \rightarrow -x_1, x_2 \rightarrow -x_2, x_5 \rightarrow -x_5$ if necessary, we can assume that $c_5=\frac{1}{8}$.\vspace{2mm}

	\noindent Let $x_1=x+c_1(x \in \mathbb{Z}),$ $x_2=1$, $x_5=\frac{1}{8}.$ Then (\ref{eq17}) reduces to $0<G=(x+c_1)+a_2''-\frac{7}{32}\leq 1$.  Taking $x=0,1 {\rm~and~}2$ simultaneously, we see that (\ref{eq17}) is soluble  with strict inequality unless $c_1+a_2'' \equiv \frac{7}{32} \pmod 1.$ \vspace{2mm}
	
	\noindent Similarly taking $ x_1=-x+c_1$ (for $x=0,1,2$),  $x_2=-1, x_5=\frac{1}{8}$, (\ref{eq17}) is soluble  with strict inequality unless $  -c_1+a_2'' \equiv \frac{11}{32} \pmod 1.$ On solving these we get $(a_2'',c_1)=(-\frac{7}{32},\frac{7}{16})$ or  $(\frac{9}{32},-\frac{1}{16})$. When $(a_2'',c_1)=(-\frac{7}{32},\frac{7}{16})$, $ G(\frac{23}{16},2,\frac{9}{8})=\frac{11}{32}$ and when $(a_2'',c_1)=(\frac{9}{32},-\frac{1}{16})$, $ G(-\frac{1}{16},-2,-\frac{7}{8})=\frac{11}{32}$. Therefore (\ref{eq17}) is soluble with strict inequality.\vspace{2mm}

	\noindent This completes the proof of Lemma \ref{lem13}. \hfill $\square$

	\begin{lemma}\label{lem14} If $A\geq d$, then $(\ref{eqe})$ is soluble with strict inequality unless $Q=\rho Q_6, \rho>0$ where $Q_6= (x_1+\frac{1}{2}x_2+\frac{1}{2}x_3)x_2-\frac{1}{2}x_3^2-x_4^2-x_5^2$. For $Q_6$ equality is required in
		$(\ref{eqe})$ for $(c_1,c_2,c_3,c_4,c_5)=(0,0,0,\frac{1}{2},\frac{1}{2})$.\end{lemma}
	
	\noindent \textbf{Proof:} Since $d\leq A \leq \sqrt{\frac{2d^5}{3a}}= \sqrt{\frac{4d^5}{3}}$, we get $d \geq (\frac{3}{4})^{1/3}>0.9085$. Take $h_A=1, k_A=1$ so that  (\ref{eqg}) is satisfied. Hence (\ref{eqe}) is soluble  unless $A=1$ and $\beta_A=\pm a_4-2A(c_4+\lambda x_5) \equiv 0 \pmod 1$. Since $-\frac{1}{2} < \lambda \leq \frac{1}{2}$, taking $x_5=c_5$ and $1+c_5$, we get $\lambda =0$ or  $\frac{1}{2}$.\vspace{2mm}
	
	\noindent Further  $1=A \leq \sqrt{\frac{4d^5}{3}}$ gives $d \geq (\frac{3}{4})^{1/5}>0.944$.\vspace{4mm}
	
	\noindent \textbf{Case I}: $\lambda=0$.\vspace{2mm}
	
	\noindent Here, $\pm a_4-2A(c_4+\lambda x_5) \equiv 0 \pmod 1$ gives $(a_4,c_4)=(\frac{1}{2},\pm \frac{1}{4}), (0,\frac{1}{2}),$ or $(0,0)$. Changing $x_1 \rightarrow -x_1, x_2 \rightarrow -x_2, x_4 \rightarrow -x_4,$ if necessary, we can assume that $c_4=\frac{1}{4}$.\vspace{2mm}
	
	\noindent Since $\lambda=0,$ using (\ref{eq09}) we get $C=t+A\lambda^2=t\geq A=1$, and from (\ref{eq8}) $t=\frac{d^5}{2aA}=d^5\leq 1$ thereby giving $t=1=d$. Also from (\ref{eq12}) and (\ref{eq16}) we get $a_5''=a_5'=a_5, a_4'=a_4$.\vspace{2mm}
	
	\noindent Take $h_t=1, k_t=1$ so that  (\ref{eqj}) is satisfied. Hence (\ref{eqe}) is soluble  unless  $\beta_t=\pm a_5-2tc_5 \equiv 0 \pmod 1$ i.e., unless $(a_5,c_5)=(\frac{1}{2},\pm \frac{1}{4}), (0,\frac{1}{2})$ or $(0,0)$. Changing $x_1 \rightarrow -x_1, x_2 \rightarrow -x_2, x_5 \rightarrow -x_5,$ if necessary, we can assume that $c_5=\frac{1}{4}$.
	Here, we have $ d/a=2d=2$ which gives $m=1$, thereby $\delta_m=d+(m-1)^2a/4=d=1$.  Here, $Q=(x_1+a_2x_2+\cdots+a_5x_5)x_2-\frac{1}{2}x_3^2-x_4^2-x_5^2$, which is clearly symmetric in $x_4$ and $x_5$, therefore we need to discuss the following cases:
	
	\begin{enumerate}[$\rm(i)$]\item $(a_5,c_5)=(\frac{1}{2}, \frac{1}{4})$, $(a_4,c_4)=(\frac{1}{2},  \frac{1}{4})$ \item $(a_5,c_5)=(\frac{1}{2}, \frac{1}{4})$, $(a_4,c_4)=(0,\frac{1}{2})$
		\item $(a_5,c_5)=(\frac{1}{2}, \frac{1}{4})$, $(a_4,c_4)=(0,0)$ \item $(a_5,c_5)=(0,\frac{1}{2})$, $(a_4,c_4)=(0,\frac{1}{2})$\item $(a_5,c_5)=(0,\frac{1}{2})$, $(a_4,c_4)=(0,0)$\item $(a_5,c_5)=(0,0)$, $(a_4,c_4)=(0,0).$\end{enumerate}
	
	\noindent \textbf{Subcase (i)}: $(a_5,c_5)=(\frac{1}{2}, \frac{1}{4})$, $(a_4,c_4)=(\frac{1}{2},  \frac{1}{4})$.\vspace{2mm}
	
	\noindent  Here, (\ref{eq13}) reduces to
	\begin{equation} \label{alpha}0<F=\big(x_1+a_2'x_2+\frac{1}{2}x_4+\frac{1}{2}x_5\big)x_2-x_4^2-x_5^2-\frac{1}{8}\leq 1.\end{equation}
	
	\noindent Let $x_1=x+c_1(x \in \mathbb{Z}),$ $x_2=1$, $x_4=\frac{1}{4}$, $x_5=\frac{1}{4}.$ Then (\ref{alpha}) reduces to $0<F=x+c_1+a_2'\leq 1$.  Taking $x=0 {\rm~and~}1$ simultaneously, we see that (\ref{alpha}) is soluble  with strict inequality unless $c_1+a_2' \equiv 0 \pmod 1.$ \vspace{2mm}
	
	\noindent Similarly taking $ x_1=-x+c_1$ (for $x=0,1,2$),  $x_2=-1, x_4=\frac{1}{4}$, $x_5=\frac{1}{4}$, (\ref{alpha}) is soluble  with strict inequality unless $  -c_1+a_2' \equiv \frac{1}{2} \pmod 1.$ On solving these we get $(a_2',c_1)=(\frac{1}{4},-\frac{1}{4})$ or  $(-\frac{1}{4},\frac{1}{4})$. When $(a_2',c_1)=(\frac{1}{4},-\frac{1}{4}))$, $ F(-\frac{1}{4},2,\frac{5}{4},\frac{1}{4})=\frac{1}{4}$ and when $(a_2',c_1)=(-\frac{1}{4},\frac{1}{4})$, $ F(\frac{5}{4},2,\frac{5}{4},\frac{5}{4})=\frac{3}{4}$. Therefore (\ref{alpha}) is soluble with strict inequality.\vspace{2mm}
	
	\noindent \textbf{Subcase (ii)}: $(a_5,c_5)=(\frac{1}{2}, \frac{1}{4})$, $(a_4,c_4)=(0,  \frac{1}{2})$.\vspace{2mm}
	
	\noindent  Here, (\ref{eq13}) reduces to
	\begin{equation} \label{beta}0<F=\big(x_1+a_2'x_2+\frac{1}{2}x_5\big)x_2-x_4^2-x_5^2-\frac{1}{8}\leq 1.\end{equation}
	
	\noindent Let $x_1=x+c_1(x \in \mathbb{Z}),$ $x_2=1$, $x_4=\frac{1}{2}$, $x_5=\frac{1}{4}.$ Then (\ref{beta}) reduces to $0<F=x+c_1+a_2'-\frac{5}{16}\leq 1$.  Taking $x=0,1 {\rm~and~}2$ simultaneously, we see that (\ref{beta}) is soluble  with strict inequality unless $c_1+a_2' \equiv \frac{5}{16} \pmod 1.$ \vspace{2mm}
	
	\noindent Similarly taking $ x_1=-x+c_1$ (for $x=0,1,2$),  $x_2=-1, x_4=\frac{1}{2}$, $x_5=\frac{1}{4}$, (\ref{beta}) is soluble  with strict inequality unless $  -c_1+a_2' \equiv \frac{9}{16} \pmod 1.$ On solving these we get $(a_2',c_1)=(\frac{7}{16},-\frac{2}{16})$ or  $(-\frac{1}{16},\frac{6}{16})$. When $(a_2',c_1)=(\frac{7}{16},-\frac{2}{16})$, $ F(-\frac{2}{16},2,\frac{1}{2},\frac{5}{4})=\frac{13}{16}$ and when $(a_2',c_1)=(-\frac{1}{16},\frac{6}{16})$, $ F(\frac{6}{16},2,\frac{1}{2},\frac{1}{4})=\frac{5}{16}$. Therefore (\ref{beta}) is soluble with strict inequality.\vspace{2mm}
	
	\noindent \textbf{Subcase (iii)}: $(a_5,c_5)=(\frac{1}{2}, \frac{1}{4})$, $(a_4,c_4)=(0,  0)$.\vspace{2mm}

	\noindent Let $x_1=x+c_1(x \in \mathbb{Z}),$ $x_2=1$, $x_4=0$, $x_5=\frac{1}{4}.$ Then (\ref{beta}) reduces to $0<F=x+c_1+a_2'-\frac{1}{16}\leq 1$.  Taking $x=0,1 {\rm~and~}2$ simultaneously, we see that (\ref{beta}) is soluble  with strict inequality unless $c_1+a_2' \equiv \frac{1}{16} \pmod 1.$ \vspace{2mm}
	
	\noindent Similarly taking $ x_1=-x+c_1$ (for $x=0,1,2$),  $x_2=-1, x_4=0$, $x_5=\frac{1}{4}$, (\ref{beta}) is soluble  with strict inequality unless $  -c_1+a_2' \equiv \frac{5}{16} \pmod 1.$ On solving these we get $(a_2',c_1)=(\frac{3}{16},-\frac{2}{16})$ or  $(-\frac{5}{16},\frac{6}{16})$. When $(a_2',c_1)=(\frac{3}{16},-\frac{2}{16})$, $ F(-\frac{2}{16},2,0,\frac{5}{4})=\frac{1}{16}$ and when $(a_2',c_1)=(-\frac{5}{16},\frac{6}{16})$, $ F(1+\frac{6}{16},2,1,\frac{1}{4})=\frac{9}{16}$. Therefore (\ref{beta}) is soluble with strict inequality.\vspace{2mm}

	\noindent \textbf{Subcase (iv)}: $(a_5,c_5)=(0,\frac{1}{2}), (a_4,c_4)=(0,\frac{1}{2})$.\vspace{2mm}
	
	\noindent Here, (\ref{eq13}) reduces to
	\begin{equation} \label{eqz}0<F=(x_1+a_2'x_2)x_2-x_4^2-x_5^2-\frac{1}{8}\leq 1.\end{equation}
	
	\noindent Let $x_1=x+c_1(x \in \mathbb{Z}),$ $x_2=1$, $x_4=\frac{1}{2}=x_5.$ Then (\ref{eqz}) reduces to $0<F=(x+c_1)+a_2'-\frac{5}{8}\leq 1$.  Taking $x=0,1 {\rm~and~}2$ simultaneously, we see that (\ref{eqz}) is soluble with strict inequality unless $c_1+a_2' \equiv \frac{5}{8} \pmod 1.$ \vspace{2mm}
	
	\noindent Similarly taking $ x_1=-x+c_1$ (for $x=0,1,2$),  $x_2=-1, x_4=\frac{1}{2}=x_5$, (\ref{eqz}) is soluble  with strict inequality unless $  -c_1+a_2' \equiv \frac{5}{8} \pmod 1.$ On solving these we get $(a_2',c_1)=(\frac{1}{8},\frac{1}{2})$ or  $(-\frac{3}{8},0)$.\vspace{2mm}

	\noindent For $(a'_2,c_1)= (\frac{1}{8},\frac{1}{2})$,  $F(\frac{1}{2},2,\frac{1}{2},\frac{1}{2})=\frac{7}{8}$ gives a solution of (\ref{eqz}). When $(a'_2,c_1)=(-\frac{3}{8},0)$,by (\ref{eq12}) we get  $a_2=-\frac{3}{8}$, if $a_3=0$ and $a_2=\frac{1}{2}$, if $a_3=\frac{1}{2}$. For $a_2=-\frac{3}{8}$, $c_1=0$, $a_3=0$,  from (\ref{eqe}), we get $Q= (x_1-\frac{3}{8}x_2)x_2-\frac{1}{2}x_3^2-x_4^2-x_5^2$. But then $Q(2,2,\frac{3}{2},\frac{1}{2},\frac{1}{2})=\frac{7}{8}$ gives a solution of (\ref{eqe}).\vspace{2mm}
	
	\noindent When $a_2=\frac{1}{2}$, $a_3=\frac{1}{2}$, $c_1=0$, $c_3=0$, $Q= (x_1+\frac{1}{2}x_2+\frac{1}{2}x_3)x_2-\frac{1}{2}x_3^2-x_4^2-x_5^2=Q_6$. For $Q_6$ equality is required in
	(\ref{eqe}) for $(c_1,c_2,c_3,c_4,c_5)=(0,0,0,\frac{1}{2},\frac{1}{2})$. This is so because for integers $x_i$, \vspace{2mm}
	
	\noindent (i) $ 8Q_6(x_1,x_2,x_3,x_4+\frac{1}{2},x_5+\frac{1}{2})= 8x_1x_2+5x_2^2-(2x_3-x_2)^2-2(2x_4+1)^2-2(2x_5+1)^2\equiv 5-1-2-2=0 \pmod 8 $ if $x_2$ is odd.\vspace{2mm}
	
	\noindent (ii) Let $x_2=2x'_2$, then $8Q_6(x_1,x_2,x_3,x_4+\frac{1}{2},x_5+\frac{1}{2})=16(x_1+\frac{1}{2}x_2+\frac{1}{2}x_3)x'_2
	-4x_3^2-2(2x_4+1)^2-2(2x_5+1)^2\equiv -4x_3^2-2-2 \equiv -4-2-2=0 \pmod 8$ if  $x_3$ is odd.
	
	\noindent (iii) Let $x_3=2x'_3$, $x_2=2x'_2$.  Then  $Q_6(x_1,x_2,x_3,x_4+\frac{1}{2},x_5+\frac{1}{2})= 2(x_1+x'_2+x'_3)x'_2-2(x'_3)^2-x_4(x_4+1)-x_5(x_5+1)-\frac{1}{2}=2m-\frac{1}{2}$, for some integer $m$ which can not lie in the interval $(0,1)$ for any $m$. \vspace{2mm}

	\noindent \textbf{Subcase (v)}: $(a_5,c_5)=(0,\frac{1}{2}), (a_4,c_4)=(0,0)$.\vspace{2mm}
	
	\noindent Working as in Subcases (i)-(iv), we find here that $(a_2',c_1)=(\frac{3}{8},0)$ or  $(-\frac{1}{8},\frac{1}{2})$. When $(a_2',c_1)=(\frac{3}{8},0)$, $ F(1,2,1,\frac{3}{2})=\frac{1}{8}$ and when $(a_2',c_1)=(-\frac{1}{8},\frac{1}{2})$, $ F(\frac{5}{2},2,2,\frac{1}{2})=\frac{1}{8}$. Therefore (\ref{eqz}) is soluble with strict inequality.\vspace{2mm}
	
	\noindent \textbf{Subcase (vi)}: $(a_5,c_5)=(0,0)=(a_4,c_4)$.\vspace{2mm}
	
	\noindent Working as in above subcases we find that (\ref{eqz}) is soluble with strict inequality.\vspace{4mm}
	
	\noindent \textbf{Case II:} $\lambda =\frac{1}{2}$.\vspace{2mm}
	
	\noindent Recall here $d >0.94$, $A=1$ and $\beta_A=\pm a_4-2A(c_4+\lambda x_5) \equiv 0 \pmod 1$.  This gives $\pm a_4 -2c_4-x_5 \equiv 0 \pmod 1$. From here we get $2a_4 \equiv 0 \pmod 1$, i.e., $a_4=0$ or $\frac{1}{2}$. For $a_4=0,$ we get $2c_4+c_5\equiv 0 \pmod 1$, and for $a_4=\frac{1}{2}$, we get $2c_4+c_5\equiv \frac{1}{2} \pmod 1$.\vspace{2mm}
	
	\noindent Since $\lambda = \frac{1}{2}$, using (\ref{eq09}) we get $1=A\leq C=t+A\lambda^2=t+\frac{A}{4}=t+\frac{1}{4}$, i.e., $t\geq \frac{3}{4}$ and from (\ref{eq8}) $t=\frac{d^5}{2aA}=d^5\leq d$. Also,  from (\ref{eq12}) and (\ref{eq16}) we get $a_5'=a_5, a_4'=a_4$, and $a_5''\equiv a_5-\frac{1}{2}a_4\pmod 1$, as $h_4=h_5=0$ in this Subsection \ref{subsec5.2}.\vspace{2mm}
	
	\noindent For $C=t+\frac{1}{4}$, take $h_{C}=1, k_{C}=1$ so that  (\ref{eqj}) is satisfied. Hence (\ref{eqe}) is soluble  unless $C=1$ i.e., unless $t=\frac{3}{4}$ and $\beta_{C,Q}=\pm a_5-(x_4+\frac{1}{2}c_5)-2tc_5 \equiv 0 \pmod 1$ which gives $2a_5 \equiv 0 \pmod 1$, i.e., $a_5=0$ or $\frac{1}{2}$. When $a_5=0$, we have $ c_4+2c_5\equiv 0 \pmod 1$ and when $a_5=\frac{1}{2}$, we have $ c_4+2c_5\equiv \frac{1}{2}\pmod 1$.\vspace{2mm}

	\noindent Further, we have $ d/a=2d> 1.88$ which gives $m=1$, thereby $\delta_m/A=d/A=d$ i.e., $K=0$ and hence $\delta_{m,K}=d-\frac{1}{4}.$ \vspace{2mm}
	
	\noindent Here, $Q=(x_1+a_2x_2+\cdots+a_5x_5)x_2-\frac{1}{2}x_3^2-x_4^2-x_5^2-x_4x_5$, which is clearly symmetric in $x_4$ and $x_5$, therefore we need to discuss the  cases when $(a_4,a_5)=(0,0),(\frac{1}{2},0)$ or $(\frac{1}{2},\frac{1}{2})$. When $a_4=a_5=0$, we have $2c_4+c_5\equiv 0 \pmod 1$ and $ c_4+2c_5\equiv 0 \pmod 1$. This gives $c_4\equiv c_5\pmod 1$ and $3c_4\equiv 0 \pmod 1$. So we get either $c_4=c_5=0$ or $c_4=c_5=\pm\frac{1}{3}$. When $a_4=\frac{1}{2}$ and $a_5=0$, we have $2c_4+c_5\equiv \frac{1}{2} \pmod 1$ and $ c_4+2c_5\equiv 0 \pmod 1$. This gives $c_4\equiv c_5+\frac{1}{2}\pmod 1$ and $3c_4\equiv 0 \pmod 1$. So we get either $(c_4,c_5)=(0,\frac{1}{2})$ or $(\pm\frac{1}{3},\mp\frac{1}{6})$. When  $a_4=a_5=\frac{1}{2}$, we get $c_4=c_5= \frac{1}{2}$ or $c_4=c_5=\pm\frac{1}{6}$. By taking the transformation $x_4\rightarrow -x_4, x_5 \rightarrow -x_5$, if necessary, we need consider the following six subcases :\vspace{-2mm}
	
	\begin{enumerate}[$\rm(i)$]\item $(a_5,c_5)=(0,0), (a_4,c_4)=(0,0)$\item $(a_5,c_5)=(0, \frac{1}{3}), (a_4,c_4)=(0,\frac{1}{3})$
		\item $(a_5,c_5)=(0,\frac{1}{2}), (a_4,c_4)=(\frac{1}{2},0)$ \item $(a_5,c_5)=(0,-\frac{1}{6}), (a_4,c_4)=(\frac{1}{2},\frac{1}{3})$\item $(a_5,c_5)=(\frac{1}{2},\frac{1}{2}), (a_4,c_4)=(\frac{1}{2},\frac{1}{2})$\item $(a_5,c_5)=(\frac{1}{2},\frac{1}{6}), (a_4,c_4)=(\frac{1}{2},\frac{1}{6})$.\end{enumerate}
	
	\noindent \textbf{Subcase (i)}: $(a_5,c_5)=(0,0), (a_4,c_4)=(0,0)$.\vspace{2mm}
	
	\noindent In this case we get $a_5''=0$. The following table gives a solution to (\ref{eq17}) with strict inequality, i.e., of\vspace{-3mm}
	 \begin{equation}\label{eqy}0<G=(x_1+a_2''x_2)x_2-\frac{3}{4}x_5^2-\frac{3}{8}<d-\frac{1}{4}.\vspace{-2mm}\end{equation}

	{\scriptsize
		\begin{equation*}\begin{array}{lllll}
		\hline\\
		{\rm ~~~~Range} & x_1 & x_1x_2 & x_5 & {\rm ~~G}\vspace{2mm}\\
		\hline\\
		\frac{3}{8}<f_1''\leq 1 & c_1&|c_1|&0&f_1''-\frac{3}{8}\vspace{2mm}\\
		
		\frac{1}{8}<f_1''\leq \frac{3}{8}& \pm 1+c_1&1+|c_1|&\pm 1&f_1''-\frac{1}{8}\vspace{2mm}\\
		
		-\frac{1}{2}<f_1''<d-\frac{7}{8} &\pm 1+c_1&1+|c_1|&0&f_1''+\frac{5}{8} \vspace{2mm}\\

		{\rm When~}  d-\frac{7}{8}\leq f_1''\leq \frac{1}{8}: \vspace{2mm}\\
		
		\left.\begin{array}{l}d-\frac{7}{8}\leq g_1''\leq \frac{1}{8}\vspace{1mm}\\f_2''>\frac{3}{8}\end{array}\right\}& c_1&2|c_1|&0&f_2''-\frac{3}{8}\vspace{2mm}\\
		
		\left.\begin{array}{l}d-\frac{7}{8}\leq g_1''\leq \frac{1}{8}\vspace{1mm}\\f_2''\leq\frac{3}{8}\end{array}\right\}& c_1&3|c_1|&0&f_3''-\frac{3}{8}\vspace{2mm}\\
		
		-\frac{5}{8}<g_1''<d-\frac{7}{8} &\pm 1+c_1&1-|c_1|&0&g_1''+\frac{5}{8} \vspace{2mm}\\
		
		-\frac{7}{8}<g_1''\leq -\frac{5}{8}&\pm 2+c_1&2-|c_1|&\pm 1&g_1''+\frac{7}{8}\vspace{2mm}\\
		
		-1<g_1''\leq -\frac{7}{8}& \pm 1+c_1&2+2|c_1|&\pm 1& f_2''+\frac{7}{8}\vspace{2mm}\\
		
		\hline
		\end{array}
		\end{equation*}\vspace{1mm}}

	\noindent \textbf{Subcase (ii)}: $(a_5,c_5)=(0, \frac{1}{3}), (a_4,c_4)=(0,\frac{1}{3})$.\vspace{2mm}
	
	\noindent In this case we get $a_5''=0$.  The following table gives a solution to (\ref{eqy}).\vspace{2mm}

	{\scriptsize
		\begin{equation*}\begin{array}{lllll}
		\hline\\
		{\rm ~~~~Range} & x_1 & x_1x_2 & x_5 & {\rm ~~G}\vspace{2mm}\\
		\hline\\
		
		\frac{11}{24}<f_1''\leq 1&c_1&|c_1|&\frac{1}{3}&f_1''-\frac{11}{24}\vspace{2mm}\\
		
		-\frac{7}{24}<f_1''<d-\frac{13}{24}&\pm 1+c_1&1+|c_1|&-\frac{2}{3}&f_1''+\frac{7}{24}\vspace{2mm}\\
		
		-\frac{1}{2}<f_1''\leq -\frac{7}{24}&\pm 1+c_1&1+|c_1|&\frac{1}{3}&f_1''+\frac{13}{24}\vspace{2mm}\\
		
		{\rm When~}  d-\frac{13}{24}\leq f_1''\leq \frac{11}{24}: \vspace{4mm}\\
		
		\left.\begin{array}{l}d-\frac{13}{24}\leq g_1''\leq \frac{11}{24}\vspace{1mm}\\f_2''>\frac{41}{24}\end{array}\right\}& c_1&2|c_1|&\frac{4}{3}&f_2''-\frac{41}{24}\vspace{2mm}\\

		\left.\begin{array}{l}d-\frac{13}{24}\leq g_1''\leq \frac{11}{24}\vspace{1mm}\\f_2''\leq\frac{41}{24}\end{array}\right\}& \pm1+c_1&-3+3|c_1|&\frac{1}{3}&f_3''-\frac{83}{24}\vspace{2mm}\\

		-\frac{7}{24}<g_1''<d-\frac{13}{24}&\pm 1+c_1&1-|c_1|&-\frac{2}{3}&g_1''+\frac{7}{24}\vspace{2mm}\\
		
		-\frac{13}{24}<g_1''\leq -\frac{7}{24}&\pm 1+c_1&1-|c_1|&\frac{1}{3}&g_1''+\frac{13}{24}\vspace{2mm}\\
		
		-\frac{1}{2}<g_1''\leq -\frac{13}{24}&c_1&2|c_1|&\frac{1}{3}&f_2''-\frac{11}{24}\vspace{2mm}\\

		\hline
		\end{array}
		\end{equation*}}\vspace{2mm}

	\noindent \textbf{Subcase (iii)}: $(a_5,c_5)=(0,\frac{1}{2}), (a_4,c_4)=(\frac{1}{2},0)$.\vspace{2mm}
	
	\noindent In this case we get $a_5''=-\frac{1}{4}$. The following table gives a solution to (\ref{eq17}) with strict inequality, i.e., of\vspace{-2mm}
	 \begin{equation}\label{eqx}0<G=(x_1+a_2''x_2-\frac{1}{4}x_5)x_2-\frac{3}{4}x_5^2-\frac{3}{8}<d-\frac{1}{4}.\vspace{-2mm}
	\end{equation}
	
	{\scriptsize
		\begin{equation*}\begin{array}{llllll}
		\hline\\
		{\rm ~~~~Range} & x_1 & x_1x_2 & x_5 &x_2x_5& {\rm ~~G}\vspace{2mm}\\
		\hline\\
		
		\frac{7}{16}<f_1''\leq 1&c_1&|c_1|&\pm \frac{1}{2}&-\frac{1}{2}&f_1''-\frac{7}{16}\vspace{2mm}\\
		
		-\frac{5}{16}<f_1''<d-\frac{9}{16}&\pm 1+c_1&1+|c_1|&\pm \frac{1}{2}&\frac{1}{2}&f_1''+\frac{5}{16}\vspace{2mm}\\
		
		-\frac{1}{2}<f_1''\leq -\frac{5}{16}&\pm 1+c_1&1+|c_1|&\pm \frac{1}{2}&-\frac{1}{2}&f_1''+\frac{9}{16}\vspace{2mm}\\

		{\rm When~}  d-\frac{9}{16}\leq f_1''\leq \frac{7}{16}: \vspace{2mm}\\
		
		d-\frac{9}{16}\leq g_1''\leq \frac{7}{16}&c_1&2|c_1|&\pm \frac{3}{2}&-3&f_2''-\frac{21}{16} \vspace{2mm}\\
		
		-\frac{5}{16}<g_1''<d-\frac{9}{16}&\pm 1+c_1&1-|c_1|&\pm \frac{1}{2}&\frac{1}{2}&g_1''+\frac{5}{16}\vspace{2mm}\\
		
		-\frac{9}{16}<g_1''\leq -\frac{5}{16}&\pm 1+c_1&1-|c_1|&\pm \frac{1}{2}&-\frac{1}{2}&g_1''+\frac{9}{16}\vspace{2mm}\\
		
		-\frac{1}{2}<g_1''\leq -\frac{9}{16} &c_1&2|c_1|&\pm \frac{1}{2}&-1&f_2''-\frac{5}{16} \vspace{2mm}\\

		\hline
		\end{array}
		\end{equation*}}\vspace{2mm}
	
	\noindent \textbf{Subcase (iv)}: $(a_5,c_5)=(0,-\frac{1}{6}), (a_4,c_4)=(\frac{1}{2},\frac{1}{3})$.\vspace{2mm}
	
	\noindent In this case we get $a_5''=-\frac{1}{4}$. The following table gives a solution to (\ref{eqx}).
	
	{\scriptsize
		\begin{equation*}\begin{array}{lllll}
		\hline\\
		{\rm ~~~~Range} & x_1 & x_2 & x_5 & {\rm ~~G}\vspace{2mm}\\
		\hline\\
		
		\frac{17}{48}<p_1''\leq 1&c_1&1&-\frac{1}{6}&p_1''-\frac{17}{48}\vspace{2mm}\\
		
		\frac{5}{48}<p_1''\leq \frac{17}{48}&1+c_1&1&\frac{5}{6}&p_1''-\frac{5}{48}\vspace{2mm}\\
		
		 -\frac{31}{48}<p_1''<d-\frac{43}{48}&1+c_1&1&-\frac{1}{6}&p_1''+\frac{31}{48}\vspace{2mm}\\
		
		-\frac{43}{48}<p_1''\leq -\frac{31}{48}&2+c_1&1&\frac{5}{6}&p_1''+\frac{43}{48}\vspace{2mm}\\
		
		-1<p_1''\leq -\frac{43}{48}&-1+c_1&-1&\frac{5}{6}&q_1''+\frac{5}{16}\vspace{2mm}\\

		{\rm When~}  d-\frac{43}{48}\leq p_1''\leq \frac{5}{48}: \vspace{2mm}\\
		
		\frac{7}{16}<q_1''\leq 1&c_1&-1&-\frac{1}{6}&q_1''-\frac{7}{16}\vspace{2mm}\\
		
		d-\frac{9}{16}\leq q_1''\leq \frac{7}{16}&c_1&2&-\frac{1}{6}&p_2''-\frac{5}{16}\vspace{2mm}\\
		
		-\frac{5}{16}<q_1''<d-\frac{9}{16}&-1+c_1&-1&\frac{5}{6}&q_1''+\frac{5}{16}\vspace{2mm}\\
		
		-\frac{9}{16}<q_1''\leq -\frac{5}{16}&-1+c_1&-1&-\frac{1}{6}&q_1''+\frac{9}{16}\vspace{2mm}\\
		
		d-\frac{25}{16}\leq q_1''\leq -\frac{9}{16}&-1+c_1&2&\frac{5}{6}&p_2''+\frac{11}{16}\vspace{2mm}\\
		
		-1<q_1''\leq d-\frac{25}{16}&-2+c_1&-1&\frac{5}{6}&q_1''+\frac{21}{16}\vspace{2mm}\\

		\hline
		\end{array}
		\end{equation*}\vspace{1mm}}
	
	\noindent \textbf{Subcase (v)}: $(a_5,c_5)=(\frac{1}{2},\frac{1}{2}), (a_4,c_4)=(\frac{1}{2},\frac{1}{2})$.\vspace{2mm}
	
	\noindent In this case we get $a_5''=\frac{1}{4}$. The following table gives a solution to (\ref{eq17}) with strict inequality, i.e., of\vspace{-2mm}
	\begin{equation}\label{eqw}0<G=(x_1+a_2''x_2+\frac{1}{4}x_5)x_2-\frac{3}{4}x_5^2-\frac{3}{8}< d-\frac{1}{4}.\end{equation}

	{\scriptsize
		\begin{equation*}\begin{array}{llllll}
		\hline\\
		{\rm ~~~~Range} & x_1 & x_1x_2 & x_5 &x_2x_5& {\rm ~~G}\vspace{2mm}\\
		\hline\\
		
		\frac{7}{16}<f_1''\leq 1&c_1&|c_1|&\pm \frac{1}{2}&\frac{1}{2}&f_1''-\frac{7}{16}\vspace{2mm}\\
		
		-\frac{5}{16}<f_1''<d-\frac{9}{16}&\pm 1+c_1&1+|c_1|&\pm \frac{1}{2}&-\frac{1}{2}&f_1''+\frac{5}{16}\vspace{2mm}\\
		
		-\frac{1}{2}<f_1''\leq -\frac{5}{16}&\pm 1+c_1&1+|c_1|&\pm \frac{1}{2}&\frac{1}{2}&f_1''+\frac{9}{16}\vspace{2mm}\\
		
		\hline
		\end{array}
		\end{equation*}\vspace{1mm}}

	{\scriptsize
		\begin{equation*}\begin{array}{llllll}
		\hline\\
		{\rm ~~~~Range} & x_1 & x_1x_2 & x_5 &x_2x_5& {\rm ~~G}\vspace{2mm}\\
		\hline\\
		
		{\rm When~}  d-\frac{9}{16}\leq f_1''\leq \frac{7}{16}: \vspace{2mm}\\

		d-\frac{9}{16}\leq g_1''\leq \frac{7}{16}&c_1&2|c_1|&\frac{3}{2}&3&f_2''-\frac{21}{16}\vspace{2mm}\\
		
		-\frac{5}{16}<g_1''<d-\frac{9}{16}&\pm 1+c_1&1-|c_1|&\pm \frac{1}{2}&-\frac{1}{2}&g_1''+\frac{5}{16}\vspace{2mm}\\
		
		-\frac{9}{16}<g_1''\leq -\frac{5}{16}&\pm 1+c_1&1-|c_1|&\pm \frac{1}{2}&\frac{1}{2}&g_1''+\frac{9}{16}\vspace{2mm}\\
		
		-\frac{1}{2}<g_1''\leq -\frac{9}{16}&c_1&2|c_1|&\frac{1}{2}&1&f_2''-\frac{5}{16}\vspace{2mm}\\

		\hline
		\end{array}
		\end{equation*}\vspace{1mm}}
	
	\noindent \textbf{Subcase (vi)}: $(a_5,c_5)=(\frac{1}{2},\frac{1}{6}), (a_4,c_4)=(\frac{1}{2},\frac{1}{6})$.\vspace{2mm}
	
	\noindent
	In this case we get $a_5''=\frac{1}{4}$. The following table gives a solution to (\ref{eqw}).
	
	{\scriptsize
		\begin{equation*}\begin{array}{llllll}
		\hline\\
		{\rm ~~~~Range} & x_1 & x_1x_2 & x_5 &x_2x_5& {\rm ~~G}\vspace{2mm}\\
		\hline\\
		
		\frac{17}{48}<p_1''\leq 1&c_1&c_1&\frac{1}{6}&\frac{1}{6}&p_1''-\frac{17}{48}\vspace{2mm}\\
		
		\frac{5}{48}<p_1''\leq \frac{17}{48}&1+c_1&1+c_1&-\frac{5}{6}&-\frac{5}{6}&p_1''-\frac{5}{48}\vspace{2mm}\\
		
		 -\frac{31}{48}<p_1''<d-\frac{43}{48}&1+c_1&1+c_1&\frac{1}{6}&\frac{1}{6}&p_1''+\frac{31}{48}\vspace{2mm}\\
		
		-\frac{43}{48}<p_1''\leq -\frac{31}{48}&2+c_1&2+c_1&-\frac{5}{6}&-\frac{5}{6}&p_1''+\frac{43}{48}\vspace{2mm}\\
		
		-1<p_1''\leq -\frac{43}{48}&-2+c_1&2-c_1&\frac{7}{6}&-\frac{7}{6}&q_1''+\frac{5}{16}\vspace{2mm}\\

		{\rm When~}  d-\frac{43}{48}\leq p_1''\leq \frac{5}{48}: \vspace{2mm}\\

		\frac{7}{16}<q_1''\leq 1&c_1&-c_1&\frac{1}{6}&-\frac{1}{6}&q_1''-\frac{7}{16}\vspace{2mm}\\
		
		d-\frac{9}{16}\leq q_1''\leq \frac{7}{16}&c_1&2c_1&\frac{1}{6}&\frac{1}{3}&p_2''-\frac{5}{16}\vspace{2mm}\\
		
		 -\frac{5}{16}<q_1''<d-\frac{9}{16}&-1+c_1&1-c_1&-\frac{5}{6}&\frac{5}{6}&q_1''+\frac{5}{16}\vspace{2mm}\\
		
		-\frac{9}{16}<q_1''\leq -\frac{5}{16}&-1+c_1&1-c_1&\frac{1}{6}&-\frac{1}{6}&q_1''+\frac{9}{16}\vspace{2mm}\\
		
		d-\frac{25}{16}\leq q_1''\leq -\frac{9}{16}&1+c_1&2+2c_1&-\frac{5}{6}&-\frac{5}{3}&p_2''+\frac{11}{16}\vspace{2mm}\\
		
		 -1<q_1''<d-\frac{25}{16}&-2+c_1&2-c_1&-\frac{5}{6}&\frac{5}{6}&q_1''+\frac{21}{16}\vspace{2mm}\\
		
		\hline
		\end{array}
		\end{equation*}}\vspace{2mm}
	
	\noindent This completes the proof of Lemma \ref{lem14}.\vspace{2mm} \hfill $\square$
	
	\noindent Now the proof of Theorem \ref{thm3} (where $c_2\equiv 0 \pmod 1$ and $a=\frac{1}{2}$, $d\leq 1$)  follows from Lemmas \ref{lem10} - \ref{lem14}.

	\section{Proof of $\Gamma_{1,4}<8$ for $c_2 \equiv 0 \pmod 1$, $~a< \frac{1}{2}$, $m\geq 3$}\label{sec6}
	
	\numberwithin{equation}{section}
	
	Using  equation (\ref{eq5.1}) and Macbeath's Lemma, we can suppose that  $a+d\leq 1$.
	Working as in Section \ref{sec3}, we will show that either (\ref{eq10}) or (\ref{eq13}) or (\ref{eq17}) or (\ref{eq21}) is soluble with strict inequality. From Lemma \ref{lem7}, equation (\ref{eq10}) is soluble with strict inequality for $d>1$, equation (\ref{eq13}) is soluble with strict inequality for $\delta_m>1$, equation (\ref{eq17}) is soluble with strict inequality for $\delta_{m,K}>1$ and equation (\ref{eq21}) is soluble with strict inequality for $\delta_{m,K,L}>1$.  Therefore, let
	\begin{equation}\label{eq26}
	d \leq 1, ~\delta_m \leq 1, ~ \delta_{m,K}\leq 1, ~\delta_{m,K,L}\leq 1.
	\end{equation}
	\noindent This gives
	
	\begin{equation}\label{eq26'}
	d\leq \frac{4}{m+3},~ d\leq \frac{16}{(m+3)(K+3)},~d\leq \frac{64}{(m+3)(K+3)(L+3)}.
	\end{equation}
	\noindent Since $m$ is an integer defined by $m<\frac{d}{a}\leq m+1$, we have from (\ref{eq5}) and Lemma \ref{lem8}, $5\geq \frac{d}{a} \geq \frac{d}{d^{\frac{5}{3}}}> 1$, unless $d=1, a=d^{\frac{5}{3}}=1$, which can not hold here since $a+d\leq 1$. So we get $m=1,2,3$ or $4$.  Also $a+d \leq 1$ gives $\frac{d}{m+1}+d \leq 1$, i.e., $d\leq \frac{m+1}{m+2}$. For $m=3,4$ we will use the bound  $d\leq \frac{4}{m+3}$ and for $m=1,2$ we will use the bound $d\leq \frac{m+1}{m+2}$. Proceeding as in Section \ref{sec4}, one finds that
	$$ \frac{\delta_m}{A} \geq \frac{(m+3)}{4d} \sqrt{\frac{3}{2}} \frac{1}{\sqrt{m+1}},$$
	which gives
	
	\begin{equation*}
	\frac{\delta_m}{A} \geq \left\{
	\begin{array}{ll}
	1.677 &{\rm~if~} m=4\\
	1.377 &{\rm~if~} m=3\\
	1.178 &{\rm~if~} m=2\\
	1.299 &{\rm ~if~} m=1.
	\end{array}
	\right.
	\end{equation*}
	Since $K< \delta_m/A\leq K+1$, we have $K\geq 1$.\vspace{2mm}

	\noindent Also $\frac{d}{m+1} \leq a \leq d^{\frac{5}{3}}$ gives
	\begin{equation}\label{eq37}
	d \geq \frac{1}{(m+1)^{3/2}}.
	\end{equation}
	
	\noindent Now by Jackson's Lemma (Lemma \ref{lem4}), (\ref{eq17}) is soluble if $\delta_{m,K} > 2\left( \frac{1}{4}t\right)^{1/3} =\left( \frac{d^5}{aA}\right)^{1/3},$ which is so if
	$\frac{(m+3)(K+3)d}{16} > \left( \frac{d^5}{aA}\right)^{1/3},$ i.e., if
	\begin{equation}\label{eq28}
	\frac{(m+3)^4}{m+1} \frac{(K+3)^4}{K+1} > 4^7.
	\end{equation}
	
	\noindent We find that (\ref{eq28}) is satisfied for  $m=4, K \geq 2$; $m=3, K\geq 3$; and  $m=2, K\geq 5$. For $m=1$, we have $\frac{\delta_m}{A}=\frac{d}{A}\leq \frac{4d}{3a} \leq \frac{8}{3}$ (because $\frac{d}{a}\leq 2$ and $A \geq \frac{3a}{4}$), therefore, $K\leq 2$.\vspace{2mm}
	
	\noindent Thus we need to discuss the following cases depending on $m$ and $K$.
	\begin{enumerate}[$\rm(i)$]
		\item $m=4, K=1$
		\item $m=3, K=1,2$
		\item $m=2, K=1,2,3,4$
		\item $m=1, K=1,2$.
	\end{enumerate}
	
	\noindent In this Section we will prove that $\Gamma_{1,4}<8$ for  $m\geq 3$. In the next Section (Section \ref{sec7}), we prove that $\Gamma_{1,4}<8$  for $m=2, K=3,4$.\vspace{2mm}
	
	\noindent Recall the inequalities (\ref{eq10}), (\ref{eq13}) and (\ref{eq17}), with strict inequality  are
	\begin{equation}\label{eq101*}0<Q=(x_1+a_2x_2+...)x_2-a(x_3+h_4x_4+h_5x_5)^2-A(x_4+\lambda x_5)^2 - tx_5^2<d.\end{equation}
	\begin{equation}\label{eq101}
	0<F=(x_1+a_2'x_2+a_4'x_4+a_5'x_5)x_2 - A(x_4+\lambda x_5)^2
	- tx_5^2-\frac{a}{4} < \delta_m.\end{equation}
	\begin{equation}\label{eq101'}
	0<G=(x_1+a''_2x_2+a''_5x_5)x_2
	- tx_5^2-\frac{a+A}{4} < \delta_{m,K}.
	\end{equation}
	\noindent We will show that either (\ref{eq101*}) or (\ref{eq101}) or (\ref{eq101'}) is soluble. To fix the values of $A$ and $t$ we will apply Macbeath's Lemma a number of times. \vspace{4mm}
	
	\noindent To fix $A$  we will choose $x_5 \equiv c_5 \pmod 1$ arbitrarily, $x_2 = \pm 1$, $x_1=x+c_1$ and $x_4=y+c_4$, $x,y \in \mathbb{Z}$ to get \begin{equation}\label{eq102}F=\pm x + \beta_{A,F} y-Ay^2+\nu {\rm ~~ where~~} \beta_{A,F} = \pm a'_4-2A(c_4+\lambda x_5).\end{equation} Then we will choose  integers $2h_A$ and $k_A$ suitably such that \begin{equation}\label{eq103}|h_A-Ak_A^2|+\frac{1}{2}<d+(m^2-1)a/4.\end{equation} Now by Macbeath's Lemma, (\ref{eq101}) is soluble  unless \begin{equation}\label{eq104}A=\frac{h_A}{k_A^2} {\rm ~~ and ~~}\beta_{A,F}  \equiv h_A/k_A \mod{(1/k_A,2A)}. \end{equation}
	
	\noindent This will fix the values of $A$ and $(a_4',c_4)$.  When $\lambda=0$, we fix $t$ and $(a_5',c_5)$ by choosing $x_4 \equiv c_4 \pmod 1$ arbitrarily, $x_2 = \pm 1$, $x_1=x+c_1$ and $x_5=y+c_5$, $x,y \in \mathbb{Z}$ to get \begin{equation}\label{eq105}F=\pm x + \beta_{t,F}~ y-ty^2+\nu {\rm ~~ where~~} \beta_{t,F} = \pm a'_5-2tc_5.\end{equation}
	Similarly to fix $t$ and $(a''_5,c_5)$, we will choose  $x_2 = \pm 1$, $x_1=x+c_1$ and $x_5=y+c_5$, $x,y \in \mathbb{Z}$ to get \begin{equation}\label{eq105'}G=\pm x + \beta_{t,G}~ y-ty^2+\nu {\rm ~~ where~~} \beta_{t,G} = \pm a''_5-2tc_5.\end{equation}
	Then we will choose  integers $2h_t$ and $k_t$ suitably such that \begin{equation}\label{eq106}|h_t-tk_t^2|+\frac{1}{2}<\delta_m {~\rm or ~} |h_t-tk_t^2|+\frac{1}{2}<\delta_{m,K},\end{equation} as the case may be. Now by Macbeath's Lemma, (\ref{eq101}) or (\ref{eq101'}) is soluble  unless \begin{equation}\label{eq107'}t=\frac{h_t}{k_t^2} {\rm ~~ and ~~}\beta_{t,F} {\rm ~~ or ~~} \beta_{t,G} \equiv h_t/k_t \mod{(1/k_t,2t)}. \end{equation}

	\noindent To fix $a$ and $(a_3,c_3)$, we will take
	$$Q(x+c_1,\pm 1,y+c_3,x_4,x_5)=\pm x+\beta_{a,Q}y-ay^2+\nu,$$
	where $\beta_{a,Q}=\pm a_3-2a(c_3+h_4x_4+h_5x_5)$ and find integers $2h_a$ and $k_a$ such that
	\begin{equation}\label{eq107}
	\vline~h_a-ak_a^2~\vline + \frac{1}{2} <d.
	\end{equation}
	Now by Macbeath's Lemma, (\ref{eq101*}) is soluble  unless \begin{equation}\label{eq108}a=\frac{h_a}{k_a^2} {\rm ~~ and ~~}\beta_{a,Q}  \equiv h_a/k_a \mod{(1/k_a,2a)}. \end{equation}

	\noindent The following lemma (Lemma 7 and Remark 1 of Section 6.3 of Raka and Rani \cite{RakaRani}) is helpful in discarding many values of $(b,c)=(a_3,c_3)$ or $(a_4',c_4)$ or $(a_5'',c_5)$.\vspace{2mm}
	
	\begin{lemma}\label{lem15}
		Let $Q$ or $F$ or $G=f(x_1,x_2,y) $+ terms not containing $y$, where $f(x_1,x_2,y)$ $=(x_1+a_2x_2+by)x_2-\alpha y^2, y\equiv c \pmod 1; \alpha =r/2s, r$ and $s$ are positive integers with $(r,s)=1; \pm b-2ac\equiv h/k \pmod{1/k,2a}$ where $2h$ is an integer and $k$ is an integral factor of $s$. Then we just need to consider
		\begin{enumerate}[$\rm(i)$]
			\item $b=0, c= \left\{
			\begin{array}{ll}
			0,1/r,2/r,\cdots & \mbox{if $sh/k$ is an integer }\\
			1/2r,3/2r,\cdots & \mbox{if $sh/k$ is half an integer }
			\end{array}\right.$
			
			\item $b=1/2s, c= \left\{
			\begin{array}{ll}
			0,1/r,2/r,\cdots & \mbox{if $sh/k$ is half an integer }\\
			1/2r,3/2r,\cdots & \mbox{if $sh/k$ is an integer }
			\end{array}\right.$
		\end{enumerate}

		\noindent Further, for $G= \pm x+ \beta_{t,G}~y-ty^2+\nu$ and $t=r/2s$, $(r,s)=1$, $a_5''-2tc_5 \equiv h/k \pmod{1/k,2t}$, if we take $x_2=\pm 1$ only (which we will do except in the case when t=1/6), then for $a_5''=1/2s$, we need to consider only $$c_5= \left\{
		\begin{array}{ll}
		1/r,3/r,5/r,\cdots & \mbox{if $sh/k$ is half an integer }\\
		1/2r,5/2r,9/2r,\cdots & \mbox{if $sh/k$ is an integer }.
		\end{array}\right.$$
	\end{lemma}\vspace{4mm}

	\begin{theorem}\label{thm4}
		Let $Q(x_1,x_2, \cdots, x_5)$ be a real indefinite quadratic form of type $(1,4)$ and of determinant $D\neq 0$. Let $d=(8|D|)^{\frac{1}{5}}$.  Suppose that $c_2 \equiv 0 \pmod1$, $a<\frac{1}{2}$, $d\leq 1$ and $a+d\leq 1$. Let $(m,K)= (4,1),(3,2)$ or $(3,1)$.  Then $(\ref{eq101*})$ and hence  $(\ref{eq10})$ is soluble with strict inequality.
	\end{theorem}

	\noindent \textbf{Proof:} We divide the proof into following subsections:
	
	\subsection{$c_2 \equiv 0 \pmod1$, $a<\frac{1}{2}$, $m=4, K=1$}
	
	\begin{lemma}\label{lem15'}  When $m=4, K=1$, inequality $(\ref{eq101*})$ is soluble.
	\end{lemma}
	
	\noindent {\bf Proof:}
	\noindent Here $d\leq \frac{4}{m+3}$, $4<\frac{d}{a}\leq 5$ and $1<\frac{d+(m^2-1)a/4}{A}\leq 2$. From (\ref{eq7}) and (\ref{eq18}) we have
	\begin{equation}\label{eq109}
	7d/8 \leq A \leq \sqrt{2d^5/3a} \leq \sqrt{10/3}~d^2.
	\end{equation}
	\noindent So we have \begin{equation}\label{eq110}
	0.479< 7\sqrt{30}/80 \leq d \leq 4/7 ~ {\rm ~and ~} 0.4193 < A < 0.597.
	\end{equation}
	
	\noindent Let $h_A=1/2, k_A=1$.  We find that (\ref{eq103}) is satisfied. Therefore, by Macbeath's Lemma, (\ref{eq101}) is soluble unless $A=1/2$ and $\beta_{A,F}\equiv 1/2 \pmod 1$. Taking $x_5=c_5$ and $1+c_5$ in $\beta_{A,F}$ we get $\lambda=0$.  Equation (\ref{eq09}) gives $A\leq t$. Now we have for $d\leq \frac{4}{7}$
	
	\begin{equation}\label{eq111}
	1/2=A \leq t=d^5/a\leq 5d^4< d+\frac{15}{4}a.
	\end{equation}

	\noindent   Taking $h_t=1/2, k_t=1$, and using (\ref{eq111}) we find that (\ref{eq106}) is satisfied. Therefore (\ref{eq101'}) is soluble unless $t=1/2$.\vspace{2mm}
	
	\noindent From (\ref{eq8}), we have
	\begin{equation}\label{eq112}
	a=d^5/2tA=2d^5\leq 2(4/7)^5<0.1219<1/8.
	\end{equation}
	\noindent Also, $d/5\leq a \leq 2d^5$ gives
	\begin{equation}\label{eq113}
	d\geq (1/10)^{1/4}>0.562.
	\end{equation}
	\noindent Taking $h_a=1/2, k_a=2$ and using (\ref{eq112}) and (\ref{eq113}) we find that (\ref{eq107}) is satisfied. Therefore, (\ref{eq101*}) is soluble unless $a=1/8$ and $\beta_{a,Q}=\pm a_3-2a(c_3+h_4x_4+h_5x_5)\equiv 1/4 \pmod {1/2,1/4}$. But we have from (\ref{eq112}), $a\neq 1/8$ which gives that (\ref{eq101*}) is soluble  for all values of $a$.\vspace{2mm}
	
	\noindent This completes the proof of Lemma \ref{lem15'}. \hfill $\square$
	
	\subsection{$c_2 \equiv 0 \pmod1$, $a<\frac{1}{2}$, $m=3, K=2$}\label{subsec6.2}
	
	\begin{lemma}\label{lem16} When $m=3, K=2$, inequality $(\ref{eq101*})$ is soluble.
	\end{lemma}
	
	\noindent {\bf Proof:}
	Here $ 3<\frac{d}{a}\leq 4$ and $2<\frac{d+(m^2-1)a/4}{A}=\frac{d+2a}{A}\leq 3$. From (\ref{eq7}) and (\ref{eq18}) we have
	\begin{equation}\label{eq54}
	d/2 \leq A \leq \sqrt{2d^5/3a} \leq \sqrt{8/3}~d^2.
	\end{equation}
	
	\noindent This together with (\ref{eq26'}) namely  $d\leq \frac{16}{(m+3)(K+3)}$, we get \vspace{-2mm}
	\begin{equation}\label{eq55}
	0.306< \sqrt{6}/8\leq d \leq 8/15. \vspace{-2mm}
	\end{equation}
	
	\noindent Further we have,\vspace{-2mm}
	\begin{equation}\label{eq55'} \begin{array}{ll}
	\frac{\delta_{m,K}}{t}&=\frac{2aA}{d^5} \delta_{m,K}
	\geq \frac{2}{d^5} \frac{d}{m+1}\frac{(m+3)d}{4(K+1)}\frac{(m+3)}{4} \frac{(K+3)}{4}d\vspace{4mm}\\
	&=\frac{1}{32}\frac{(m+3)^2(K+3)}{(m+1)(K+1)}\frac{1}{d^2},
	\end{array}
	\end{equation}
	
	\noindent which is $>1$ for $m=3,K=2$ and $d\leq 8/15.$  From the definition of $L$ we get $L\geq 1.$ For $L\geq 4$ and $d\leq 8/15$, we find that $\delta_{m,K,L}>\frac{(m+3)}{4}\frac{(K+3)}{4}\frac{(L+3)}{4}d>1$. Therefore, we can assume that $L\leq 3$. Further, by definition of $L$, and using (\ref{eq55'})
	we get that
	
	\begin{equation}\label{eq56'}d^2 \geq \frac{1}{32}\frac{(m+3)^2(K+3)}{(m+1)(K+1)(L+1)}.\end{equation}
	
	\noindent   Therefore, using (\ref{eq26'})  namely  $d\leq \frac{64}{(m+3)(K+3)(L+3)}$, we get
	
	\begin{equation}\label{eq56}
	\delta_{3,2,L}\leq 1,{\rm ~and~} \sqrt{15/32(L+1)} \leq d \leq 32/15(L+3).
	\end{equation}
	
	\noindent We distinguish the cases $L=1,2$ and $3$.\vspace{3mm}
	
	\noindent \textbf{Case I}: $\textbf{ m=3, K=2, L=3}$.\vspace{2mm}
	
	\noindent Here $3<\frac{\delta_{m,K}}{t}=\frac{d+2a+3A/4}{t}\leq 4$, $\delta_{m,K}=d+2a+3A/4\geq d+2.\frac{d}{4}+\frac{3}{4}.\frac{d}{2}=\frac{15}{8}d$. From (\ref{eq56}) we have
	\begin{equation}\label{eq66}
	0.3423< \sqrt{15/128} \leq d \leq 32/90<0.3556.
	\end{equation}
	\noindent Further $$0.16< \frac{15}{32}d\leq \frac{d+2a+3A/4}{4}\leq t =\frac{d^5}{2aA}\leq 4d^3< 0.18.$$
	
	\noindent  Taking $h_t=3/2$, $k_t=3$ and using (\ref{eq66})
	we find that $\vline~h_t-tk_t^2~\vline + \frac{1}{2} <\delta_{m,K}$. This is so because $ 9t+\delta_{m,K} \geq \frac{9\delta_{m,K}}{4}+\delta_{m,K}=\frac{13\delta_{m,K}}{4}\geq \frac{195d}{32}>2$ and $ 9t-1< \delta_{m,K}$ as $36d^3-1<\frac{15d}{8}$.  Therefore,  (\ref{eq101'})  is soluble unless $t=1/6$ and $\beta_{t,G}=\pm a_5''-\frac{1}{3}c_5\equiv 1/2 \pmod {\frac{1}{3}}$ i.e., $\pm 3a_5''-c_5\equiv 1/2 \pmod 1 $. By Lemma \ref{lem15} we need to consider $(a_5'',c_5)=(0,1/2)$ or $(1/6,0)$. If $a_5''=1/6$, considering the transformation $x_5 \rightarrow x_5-x_2$, we may assume that $a_5''=1/2.$\vspace{2mm}
	
	\noindent Now, $1/6=t=d^5/2aA\leq 4d^3$ gives $d\geq (1/24)^{1/3}>0.34668$. Also, $1/6=t=d^5/2aA$ gives $A=3d^5/a\leq 12 d^4.$ Hence,
	\begin{equation}\label{eq67}
	\begin{array}{l}
	0.0866< d/4\leq a\leq d/3< 0.11854\vspace{2mm}\\
	0.1733 < d/2 \leq A \leq 12 d^4< 0.19188\vspace{2mm}\\
	0.0649 < \frac{a+A}{4} < 0.07761.
	\end{array}
	\end{equation}
	
	\noindent Here we need to find a solution of\vspace{-2mm}
	\begin{equation}\label{eq68}
	0<G=(x_1+a_2''x_2+a_5''x_5)x_2-x_5^2/6	-(a+A)/4<d+2a+3A/4 \vspace{-2mm}
	\end{equation}
	for $(a_5'',c_5)=(0,1/2)$ and $(1/2,0)$.\vspace{4mm}
	
	\noindent When $(a_5'',c_5)=(0,1/2)$, The following table gives a solution of (\ref{eq68}).\vspace{2mm}

	{\scriptsize
		\begin{equation*}\begin{array}{lllll}
		\hline\\
		{\rm ~~~~Range} & x_1 & x_1x_2 & x_5& {\rm ~~G}\vspace{2mm}\\
		\hline\\
		
		\frac{a+A}{4}+\frac{3}{8}<f_1''\leq 1&c_1&|c_1|&\frac{3}{2}&f_1''-\frac{3}{8}-\frac{a+A}{4}\vspace{2mm}\\
		
		\frac{a+A}{4}+\frac{1}{24}<f_1''\leq \frac{a+A}{4}+\frac{3}{8} & c_1 & |c_1| &\frac{1}{2} & f_1''-\frac{1}{24}-\frac{a+A}{4}\vspace{2mm}\\	
		
		f_1''<d+\frac{9a}{4}+A-\frac{5}{8}&\pm 1+c_1&1+|c_1|& \frac{3}{2}& f_1''+\frac{5}{8}-\frac{a+A}{4}\vspace{2mm}\\

		{\rm When~}  d+\frac{9a}{4}+A-\frac{5}{8}\leq f_1''\leq \frac{a+A}{4}+\frac{1}{24}: \vspace{2mm}\\
		
		d+\frac{9a}{4}+A-\frac{5}{8}\leq g_1'' \leq \frac{a+A}{4}+\frac{1}{24}&c_1&2|c_1|&\frac{1}{2}&f_2''-\frac{1}{24}-\frac{a+A}{4}\vspace{2mm}\\
		
		\frac{a+A}{4}-\frac{5}{8}<g_1''<d+\frac{9a}{4}+A-\frac{5}{8}&\pm 1+c_1&1-|c_1|& \frac{3}{2}& g_1''+\frac{5}{8}-\frac{a+A}{4}\vspace{2mm}\\

		\left.\begin{array}{l}-1<g_1''\leq \frac{a+A}{4}-\frac{5}{8}\vspace{1mm}\\f_2''<d+\frac{9a}{4}+A-\frac{23}{24}\end{array}\right\}&\pm 1+c_1&2+2|c_1|&\frac{5}{2}&f_2''+\frac{23}{24}-\frac{a+A}{4}\vspace{2mm}\\
		
		\left.\begin{array}{l}-1<g_1''\leq \frac{a+A}{4}-\frac{5}{8}\vspace{1mm}\\d+\frac{9a}{4}+A-\frac{23}{24}\leq f_2''\end{array}\right\}& \pm1+c_1&1-|c_1|&\frac{1}{2}&g_1''+\frac{23}{24}-\frac{a+A}{4}\vspace{2mm}\\

		\hline
		\end{array}
		\end{equation*}\vspace{2mm}}
	
	\noindent When $(a_5'',c_5)=(1/2,0)$, the following table gives a solution of (\ref{eq68}).\vspace{2mm}

	{\scriptsize
		\begin{equation*}\begin{array}{llllll}
		\hline\\
		{\rm ~~~~Range} & x_1 & x_1x_2 & x_5&x_2x_5& {\rm ~~G}\vspace{2mm}\\
		\hline\\

		\frac{a+A}{4}<g_1''\leq\frac{1}{2}&c_1&-|c_1|&0&0&g_1''-\frac{a+A}{4}\vspace{2mm}\\	
		
		\frac{a+A}{4}-\frac{1}{3}<g_1''\leq \frac{a+A}{4}&c_1&-|c_1|&\pm 1&1&g_1''+\frac{1}{3}-\frac{a+A}{4}\vspace{2mm}\\
		
		\frac{a+A}{4}-1<g_1''<d+\frac{9a}{4}+A-1&\pm 1+c_1&1-|c_1|&0&0&g_1''+1-\frac{a+A}{4}\vspace{2mm}\\
		
		-1<g_1''\leq \frac{a+A}{4}-1&\pm 1+c_1&1-|c_1|&\pm 1&1&g_1''+\frac{4}{3}-\frac{a+A}{4}\vspace{2mm}\\

		{\rm When~}  d+\frac{9a}{4}+A-1\leq g_1''\leq \frac{a+A}{4}-\frac{1}{3}: \vspace{2mm}\\
		
		d+\frac{9a}{4}+A\leq f_1''\leq 1&\pm 1+c_1&-3+3|c_1|&0&0&f_3''-3-\frac{a+A}{4}\vspace{2mm}\\
		
		\frac{a+A}{4}<f_1''<d+\frac{9a}{4}+A&c_1&|c_1|&0&0&f_1''-\frac{a+A}{4}\vspace{2mm}\\
		
		\frac{a+A}{4}-\frac{1}{3}<f_1''\leq \frac{a+A}{4}&c_1&|c_1|&\pm 1&1&f_1''+\frac{1}{3}-\frac{a+A}{4}\vspace{2mm}\\
		
		-\frac{1}{2}<f_1''\leq \frac{a+A}{4}-\frac{1}{3}&c_1&2|c_1|&\pm 2&4&f_2''+\frac{4}{3}-\frac{a+A}{4}\vspace{2mm}\\		
		
		\hline
		\end{array}
		\end{equation*}\vspace{2mm}}

	\noindent \textbf{Case II}: $\textbf{ m=3, K=2, L=2}$.\vspace{2mm}
	
	\noindent Here $2<\frac{\delta_{m,K}}{t}=\frac{d+2a+3A/4}{t}\leq 3$, $\delta_{m,K}=d+2a+3A/4\geq d+2.\frac{d}{4}+\frac{3}{4}.\frac{d}{2}=\frac{15}{8}d$. From (\ref{eq56}) we have
	\begin{equation}\label{eq63}
	0.39528< \sqrt{5/32} \leq d \leq 32/75<0.42667.
	\end{equation}
	\noindent Further $$0.247< \frac{5}{8}d\leq \frac{d+2a+3A/4}{3}\leq t =\frac{d^5}{2aA}\leq 4d^3< 0.3107.$$
	
	\noindent  Taking $h_t=1$ and $k_t=2$ and using (\ref{eq63})
	we find that $\vline~h_t-tk_t^2~\vline + \frac{1}{2} <\delta_{m,K}$. This is so because $ 4t+\delta_{m,K} \geq \frac{4\delta_{m,K}}{3}+\delta_{m,K}=\frac{7\delta_{m,K}}{3}\geq \frac{35d}{8}>\frac{3}{2}$ and $ 4t-\frac{1}{2}< \delta_{m,K}$ as $16d^3-\frac{1}{2}<\frac{15d}{8}$.  Therefore,  (\ref{eq101'})  is soluble unless $t=1/4$ and $\beta_{t,G}=\pm a_5''-\frac{1}{2}c_5\equiv 1/2 \pmod {\frac{1}{2}}$ i.e., $\pm 2a_5''-c_5\equiv 0 \pmod 1 $. By Lemma \ref{lem15} we need to consider $(a_5'',c_5)=(0,0)$ or $(1/4,1/2)$ only.\vspace{2mm}
	
	\noindent Now, $1/4=t=d^5/2aA\leq 4d^3$ gives $d\geq (1/16)^{1/3}>0.3968$. Also, $1/4=t=d^5/2aA$ gives $A=2d^5/a\leq 8d^4$. Hence,
	\begin{equation}\label{eq64}
	\begin{array}{l}
	0.0992< d/4\leq a\leq d/3< 0.14223\vspace{2mm}\\
	0.1984< d/2 \leq A \leq 8d^4< 0.26513\vspace{2mm}\\
	0.0744< \frac{a+A}{4} < 0.10184.
	\end{array}
	\end{equation}
	
	\noindent Here we need to find a solution of
	\begin{equation}\label{eq65}
	0<G=(x_1+a_2''x_2+a_5''x_5)x_2-x_5^2/4-(a+A)/4<d+2a+3A/4
	\end{equation}
	for $(a_5'',c_5)=(0,0)$ and $(1/4,1/2)$.\vspace{4mm}
	
	\noindent When  $(a_5'',c_5)=(0,0)$, the following table gives a solution of (\ref{eq65}).\vspace{2mm}
	
	{\scriptsize
		\begin{equation*}\begin{array}{lllll}
		\hline\\
		{\rm ~~~~Range} & x_1 & x_1x_2 & x_5 & {\rm ~~G}\vspace{2mm}\\
		\hline\\
		
		\frac{a+A}{4}+\frac{1}{4}<f_1''\leq 1&c_1&|c_1|&\pm 1&f_1''-\frac{1}{4}-\frac{a+A}{4}\vspace{2mm}\\
		
		\frac{a+A}{4}<f_1''\leq \frac{a+A}{4}+\frac{1}{4}&c_1&|c_1|&0&f_1''-\frac{a+A}{4}\vspace{2mm}\\
		
		-\frac{1}{2}<f_1''<d+\frac{9a}{4}+A-\frac{3}{4}&\pm 1+c_1&1+|c_1|&\pm 1&f_1''+\frac{3}{4}-\frac{a+A}{4}\vspace{2mm}\\

		{\rm When~}  d+\frac{9a}{4}+A-\frac{3}{4}\leq f_1''\leq \frac{a+A}{4}: \vspace{2mm}\\

		d+\frac{9a}{4}+A-\frac{3}{4}\leq g_1''\leq \frac{a+A}{4}&c_1 &2|c_1| &0 &f_2''-\frac{a+A}{4} \vspace{2mm}\\
		
		\frac{a+A}{4}-\frac{3}{4}<g_1''<d+\frac{9a}{4}+A-\frac{3}{4}&\pm 1+c_1&1-|c_1|&\pm 1&g_1''+\frac{3}{4}-\frac{a+A}{4}\vspace{2mm}\\
		
		\frac{a+A}{4}-1<g_1''\leq \frac{a+A}{4}-\frac{3}{4} &\pm 1 +c_1&1-|c_1|&0&g_1''+1-\frac{a+A}{4}\vspace{2mm}\\
		
		-1<g_1''\leq \frac{a+A}{4}-1&\pm 1+c_1 &2+2|c_1| &\pm 2 &f_2''+1-\frac{a+A}{4} \vspace{2mm}\\

		\hline
		\end{array}
		\end{equation*}\vspace{1mm}}
	
	\noindent When $(a_5'',c_5)=(1/4,1/2)$, the following table gives a solution of (\ref{eq65}).

	{\scriptsize
		\begin{equation*}\begin{array}{llllll}
		\hline\\
		{\rm ~~~~Range} & x_1 & x_1x_2 & x_5&x_2x_5 & {\rm ~~G}\vspace{2mm}\\
		\hline\\
		
		\frac{a+A}{4}+\frac{3}{16}<f_1''\leq 1&c_1&|c_1|&\frac{1}{2}&-\frac{1}{2}&f_1''-\frac{3}{16}-\frac{a+A}{4}\vspace{2mm}\\
		
		\frac{a+A}{4}-\frac{1}{16}<f_1''\leq \frac{a+A}{4}+\frac{3}{16}&c_1&|c_1|&\frac{1}{2}&\frac{1}{2}&f_1''+\frac{1}{16}-\frac{a+A}{4}\vspace{2mm}\\
		
		-\frac{1}{2}<f_1''<d+\frac{9a}{4}+A-\frac{13}{16}&\pm 1+c_1& 1+|c_1|& \frac{1}{2}&-\frac{1}{2}&f_1''+\frac{13}{16}-\frac{a+A}{4}\vspace{2mm}\\

		{\rm When~}  d+\frac{9a}{4}+A-\frac{13}{16}\leq f_1''\leq \frac{a+A}{4}-\frac{1}{16}: \vspace{2mm}\\

		d+\frac{9a}{4}+A-\frac{13}{16}\leq g_1''\leq \frac{a+A}{4}-\frac{1}{16}&c_1 &2|c_1| &\pm \frac{1}{2} &1 &f_2''+\frac{3}{16}-\frac{a+A}{4}  \vspace{2mm}\\

		\frac{a+A}{4}-\frac{13}{16}<g_1''<d+\frac{9a}{4}+A-\frac{13}{16}&\pm 1+c_1& 1-|c_1|& \frac{1}{2}&-\frac{1}{2}&g_1''+\frac{13}{16}-\frac{a+A}{4}\vspace{2mm}\\
		
		\frac{a+A}{4}-\frac{17}{16}<g_1''\leq \frac{a+A}{4}-\frac{13}{16}&\pm 1+c_1& 1-|c_1|& \frac{1}{2}&\frac{1}{2}&g_1''+\frac{17}{16}-\frac{a+A}{4}\vspace{2mm}\\
		
		-1<g_1''\leq \frac{a+A}{4}-\frac{17}{16}&\pm 2+c_1 &4+2|c_1| &\pm \frac{5}{2} &-5 &f_2''+\frac{19}{16}-\frac{a+A}{4} \vspace{2mm}\\	
		
		\hline
		\end{array}
		\end{equation*}\vspace{1mm}}
	
	\noindent \textbf{Case III}: $\textbf{ m=3, K=2, L=1}$.\vspace{2mm}
	
	\noindent Here $1<\frac{\delta_{m,K}}{t}=\frac{d+2a+3A/4}{t}\leq 2$, $\delta_{m,K}=d+2a+3A/4\geq d+2.\frac{d}{4}+\frac{3}{4}.\frac{d}{2}=\frac{15}{8}d$. From (\ref{eq56}) we have
	\begin{equation}\label{eq57}
	0.4841< \sqrt{15/64} \leq d \leq 8/15<0.53334.
	\end{equation}
	\noindent Further $$0.4538< \frac{15}{16}d\leq \frac{d+2a+3A/4}{2}\leq t =\frac{d^5}{2aA}\leq 4d^3< 0.60684.$$
	
	\noindent Taking $h_t=1/2$ and $k_t=1$ and using (\ref{eq57})
	we find that $\vline~h_t-tk_t^2~\vline + \frac{1}{2} <\delta_{m,K}$, as $ t \leq \delta_{m,K}$ and   $t+\delta_{m,K} \geq \frac{15}{16}d+\frac{15}{8}d=\frac{45}{16}d>1$. Therefore,  (\ref{eq101'})  is soluble unless $t=1/2$ and $\beta_{t,G}=\pm a_5''-c_5\equiv 1/2 \pmod 1$. This gives $(a_5'',c_5)=(0,1/2)$ or $(1/2,0)$.
	
	\noindent Now, $1/2=t=d^5/2aA$ gives $d^5=aA\geq (d/4)\cdot(d/2)$, i.e., $d^3\geq 1/8$,
	and hence we get \vspace{-2mm}
	\begin{equation}\label{eq59}
	d\geq 1/2, a\geq d/4\geq 1/8 ~\mbox{and}~ A\geq d/2 \geq 1/4. \vspace{-1mm}
	\end{equation}
	
	\noindent Also $a<1/2$ and $1/2=t=d^5/2aA$. This gives \vspace{-3mm}
	\begin{equation}\label{eq60}
	A\leq 4d^4\leq 4(8/15)^4<1/3. \vspace{-2mm}
	\end{equation}
	
	\noindent Taking $h_A=1/2, k_A=1$ and using  (\ref{eq60}) we find that $\vline~h_A-Ak_A^2~\vline + \frac{1}{2}=1-A <\delta_m$ is satisfied as $\delta_m+A= d+2a+A\geq d+2d/4+d/2=2d \geq 1$, with strict inequality unless $d=1/2$, $a=\frac{d}{4}=\frac{1}{8}$ and $A=\frac{d}{2}=\frac{1}{4}$. Thus, (\ref{eq101}) and hence (\ref{eq13}) is soluble with strict inequality unless $d=1/2$, $a=\frac{1}{8}$ and $A=\frac{1}{4}$.

	\noindent In this special case (\ref{eq101'}) reduces to\vspace{-2mm}
	\begin{equation}\label{eq62}
	0<G=(x_1+a_2''x_2+a_5''x_5)x_2-\frac{x_5^2}{2}-\frac{3}{32}<\delta_{m,K}=\frac{15}{16} \vspace{-2mm}
	\end{equation}
	for $(a_5'',c_5)=(0,1/2)$ or $(1/2,0)$.\vspace{2mm}
	
	\noindent When $(a_5'',c_5)=(0,1/2)$, the following table gives a solution of (\ref{eq62}).\vspace{2mm}
	
	{\scriptsize
		\begin{equation*}\begin{array}{lllll}
		\hline\\
		{\rm ~~~~Range} & x_1 & x_1x_2 & x_5 & {\rm ~~G}\vspace{2mm}\\
		\hline\\
		
		\frac{7}{32}<f_1''\leq 1&c_1&|c_1|&\frac{1}{2}&f_1''-\frac{7}{32}\vspace{2mm}\\
		
		-\frac{1}{2}<f_1''<\frac{5}{32}&\pm 1+c_1&1+|c_1|&\frac{1}{2}&f_1''+\frac{25}{32}\vspace{2mm}\\

		{\rm When~}  \frac{5}{32}\leq f_1''\leq \frac{7}{32}: \vspace{2mm}\\

		\frac{5}{32}\leq g_1''\leq \frac{7}{32}&c_1&2|c_1|&\frac{1}{2}&f_2''-\frac{7}{32}\vspace{2mm}\\

		-\frac{25}{32}<g_1''<\frac{5}{32}&\pm 1+c_1&1-|c_1|&\frac{1}{2}&g_1''+\frac{25}{32}\vspace{2mm}\\
		
		-1<g_1''\leq -\frac{25}{32}&\pm 1+c_1&2+2|c_1|&\frac{3}{2}&f_2''+\frac{25}{32}\vspace{2mm}\\

		\hline
		\end{array}
		\end{equation*}\vspace{1mm}}

	\noindent When $(a_5'',c_5)=(1/2,0)$, taking $x_5=0$,
	the following table gives a solution of (\ref{eq62}).\vspace{2mm}

	{\scriptsize
		\begin{equation*}\begin{array}{llll}
		\hline\\
		{\rm ~~~~Range} & x_1 & x_1x_2  & {\rm ~~G}\vspace{2mm}\\
		\hline\\
		
		\frac{3}{32}<f_1''\leq 1&c_1&|c_1|&f_1''-\frac{3}{32}\vspace{2mm}\\
		
		-\frac{1}{2}<f_1''<\frac{1}{32}&\pm 1+c_1&1+|c_1|&f_1''+\frac{29}{32}\vspace{2mm}\\

		{\rm When~}  \frac{1}{32}\leq f_1''\leq \frac{3}{32}: \vspace{2mm}\\

		\frac{1}{32}\leq g_1''\leq \frac{3}{32} &c_1&3|c_1|&f_3''-\frac{3}{32}\vspace{2mm}\\
		
		-\frac{29}{32}<g_1''<\frac{1}{32}&\pm 1+c_1&1-|c_1|&g_1''+\frac{29}{32}\vspace{2mm}\\
		
		-1<g_1''\leq -\frac{29}{32} &\pm 1+c_1&3+3|c_1|&f_3''+\frac{93}{32}\vspace{2mm}\\

		\hline
		\end{array}
		\end{equation*}\vspace{2mm}}
	
	\noindent This completes the proof of Lemma \ref{lem16}. \hfill $\square$
	
	\subsection{$c_2 \equiv 0 \pmod1$, $a<\frac{1}{2}$, $m=3, K=1$}

	\begin{lemma}\label{lem17} When $m=3, K=1$, inequality $(\ref{eq101*})$ is soluble.
	\end{lemma}
	
	\noindent {\bf Proof:} Here $d\leq \frac{4}{m+3}=\frac{2}{3},~ 3<\frac{d}{a}\leq 4$ and $1<\frac{\delta_m}{A}=\frac{d+(m^2-1)a/4}{A}=\frac{d+2a}{A}\leq 2$. From (\ref{eq7}) and (\ref{eq18}) we have\vspace{-2mm}
	\begin{equation}\label{eq69}
	3d/4 \leq A \leq \sqrt{2d^5/3a} \leq \sqrt{8/3}~d^2.
	\end{equation}
	
	\noindent This  gives\vspace{-2mm}
	\begin{equation}\label{eq70}
	0.459<3\sqrt{6}/16\leq d \leq 2/3.
	\end{equation}
	\noindent and hence,
	\begin{equation}\label{eq71}
	0.34425 <3d/4\leq A\leq \sqrt{8/3}d^2< 0.7258.
	\end{equation}
	
	\noindent Let $h_A=1/2, k_A=1$. Since $A+d+2a\geq 3d/4+d+d/2=9d/4>1$ and $A<\delta_m$, we find that (\ref{eq103}) is satisfied. Therefore, by Macbeath's Lemma, (\ref{eq101}) is soluble unless $A=1/2$ and $\beta_{A,F}\equiv 1/2 \pmod 1$. Taking $x_5=c_5$ and $1+c_5$ in $\beta_{A,F}$ we get $\lambda=0$, and so $\pm a_4'-c_4\equiv 1/2\pmod1$ which gives $(a_4',c_4)=(0,1/2) {\rm ~or~} (1/2,0)$.  Equation (\ref{eq09}) gives $A\leq t$. Now we have for $d\leq \frac{2}{3}$,\vspace{-2mm}
	
	\begin{equation}\label{eq0111}
	1/2=A \leq t=d^5/a\leq 4d^4< d+2a.
	\end{equation}

	\noindent   Taking $h_t=1/2, k_t=1$, and using (\ref{eq0111}) we find that (\ref{eq106}) is satisfied. Therefore (\ref{eq101'}) is soluble unless $t=1/2$ and $\beta_{t,G}\equiv 1/2 \pmod 1$ i.e., $\pm a_5'-c_5\equiv 1/2\pmod 1$ which gives $(a_5',c_5)=(0,1/2) {\rm ~or~}(1/2,0)$.\vspace{2mm}
	
	\noindent Because of symmetry in $x_4 {\rm ~and~} x_5$, we need to consider \begin{enumerate}[$\rm (i)$]
		\item $(a_4',a_5')=(0,0), (c_4,c_5)=(1/2,1/2)$
		\item $(a_4',a_5')=(0,1/2), (c_4,c_5)=(1/2,0)$
		\item $(a_4',a_5')=(1/2,1/2), (c_4,c_5)=(0,0)$.
	\end{enumerate}

	\noindent Here, $\frac{1}{2}=t=\frac{d^5}{a} {\rm ~gives~} a=2d^5$ and $t\leq 4d^4$ gives $\frac{1}{2}=\frac{d^5}{a}\leq 4d^4$ and hence $d^4\geq \frac{1}{8}$. And  $1<\frac{d+2a}{A}\leq 2$ gives $d+4d^5\leq 1$. Thus,\vspace{-2mm}
	\begin{equation}\label{EQ2}
	0.5946<\left(\frac{1}{8}\right)^{1/4}\leq d <0.624 \vspace{-2mm}
	\end{equation}
	
	\noindent so that,\vspace{-2mm}
	\begin{equation}\label{EQ1}
	d+2a\geq d+d/2=3d/2>0.8919{\rm ~and~}
	a\geq d/4>0.1486.
	\end{equation}

	\noindent \textbf{Case (I):}  $(a_4',a_5')=(0,0), (c_4,c_5)=(1/2,1/2)$. The following table gives a solution of (\ref{eq101}) i.e., to\vspace{-2mm}
	$$0<F=(x_1+a_2'x_2)x_2-\frac{1}{2}x_4^2-\frac{1}{2}x_5^2-\frac{a}{4}<d+2a. \vspace{-2mm}$$

	{\scriptsize
		\begin{equation*}\begin{array}{llllll}
		\hline\\
		{\rm ~~~~Range} & x_1 & x_1x_2 & x_4&x_5& {\rm ~~F}\vspace{2mm}\\
		\hline\\

		\frac{1}{4}+\frac{a}{4}<f_1'\leq 1&c_1&|c_1|&\frac{1}{2}&\frac{1}{2}&f_1'-\frac{1}{4}-\frac{a}{4}\vspace{2mm}\\
		
		-\frac{1}{2}<f_1'<-\frac{3}{4}+d+\frac{9a}{4}&\pm 1+c_1&1+|c_1|&\frac{1}{2}&\frac{1}{2}&f_1'+\frac{3}{4}-\frac{a}{4}\vspace{2mm}\\

		{\rm When~}  -\frac{3}{4}+d+\frac{9a}{4}\leq f_1'\leq \frac{1}{4}+\frac{a}{4}: \vspace{2mm}\\

		-\frac{3}{4}+d+\frac{9a}{4}\leq g_1'\leq \frac{1}{4}+\frac{a}{4}&c_1 &2|c_1| &\frac{1}{2}&\frac{1}{2}&f_2'-\frac{1}{4}-\frac{a}{4} \vspace{2mm}\\
		
		-\frac{3}{4}+\frac{a}{4}<g_1'<-\frac{3}{4}+d+\frac{9a}{4}&\pm 1+c_1&1-|c_1|&\frac{1}{2}&\frac{1}{2}&g_1'+\frac{3}{4}-\frac{a}{4}\vspace{2mm}\\

		\left.\begin{array}{l}-1<g_1'\leq -\frac{3}{4}+\frac{a}{4}\vspace{1mm}\\-\frac{3}{4}+\frac{a}{4}<f_3'\end{array}\right\}&\pm 1+c_1&3+3|c_1|&\frac{3}{2}&\frac{3}{2}&f_3'+\frac{3}{4}-\frac{a}{4}\vspace{2mm}\\
		
		\left.\begin{array}{l}-1<g_1'\leq -\frac{3}{4}+\frac{a}{4}\vspace{1mm}\\f_3'\leq -\frac{3}{4}+\frac{a}{4}\end{array}\right\}& \pm1+c_1&2+2|c_1|&\frac{3}{2}&\frac{1}{2}&f_2'+\frac{3}{4}-\frac{a}{4}\vspace{2mm}\\

		\hline
		\end{array}
		\end{equation*}\vspace{1mm}}

	\noindent \textbf{Case (II):} $(a_4',a_5')=(0,1/2), (c_4,c_5)=(1/2,0)$. Here, (\ref{eq101}) reduces to\vspace{-2mm}
	$$0<F=(x_1+a_2'x_2+\frac{1}{2}x_5)x_2-\frac{1}{2}x_4^2-\frac{1}{2}x_5^2-\frac{a}{4}<d+2a. \vspace{-2mm}$$
	
	\noindent Taking $x_5=0$, the following table gives its solution.
	
	{\scriptsize
		\begin{equation*}\begin{array}{lllll}
		\hline\\
		{\rm ~~~~Range} & x_1 & x_1x_2 & x_4& {\rm ~~F}\vspace{2mm}\\
		\hline\\
		
		\frac{1}{8}+\frac{a}{4}<f_1'\leq 1&c_1&|c_1|&\frac{1}{2}&f_1'-\frac{1}{8}-\frac{a}{4}\vspace{2mm}\\
		
		-\frac{1}{2}<f_1'<-\frac{7}{8}+d+\frac{9a}{4}&\pm 1+c_1&1+|c_1|&\frac{1}{2}&f_1'+\frac{7}{8}-\frac{a}{4}\vspace{2mm}\\

		{\rm When~}  -\frac{7}{8}+d+\frac{9a}{4}\leq f_1'\leq \frac{1}{8}+\frac{a}{4}: \vspace{2mm}\\

		-\frac{7}{8}+d+\frac{9a}{4}\leq g_1'\leq \frac{1}{8}+\frac{a}{4}&c_1 &2|c_1| &\frac{1}{2}  &f_2'-\frac{1}{8}-\frac{a}{4} \vspace{2mm}\\
		
		-\frac{7}{8}+\frac{a}{4}<g_1'<-\frac{7}{8}+d+\frac{9a}{4}&\pm 1+c_1&1-|c_1|&\frac{1}{2}&g_1'+\frac{7}{8}-\frac{a}{4}\vspace{2mm}\\
		
		-1<g_1'\leq -\frac{7}{8}+\frac{a}{4}&\pm 1+c_1 &2+2|c_1| &\frac{3}{2} &f_2'+\frac{7}{8}-\frac{a}{4} \vspace{2mm}\\

		\hline
		\end{array}
		\end{equation*}\vspace{1mm}}
	
	\noindent \textbf{Case (III):} $(a_4',a_5')=(1/2,1/2), (c_4,c_5)=(0,0)$. Here, (\ref{eq101}) reduces to\vspace{-2mm}
	 \begin{equation}\label{EQN1}0<F=(x_1+a_2'x_2+\frac{1}{2}x_4+\frac{1}{2}x_5)x_2-\frac{1}{2}x_4^2-\frac{1}{2}x_5^2
	-\frac{a}{4}<d+2a.\end{equation}

	\noindent We note that (\ref{EQN1}) has no solution for $c_1=1/2$ and $a_2'=1/2$ and $d+\frac{9a}{4}\leq 1$, because then $F+\frac{a}{4}$ takes integral values only. Therefore we have to go to $Q$.\vspace{2mm}
	
	\noindent Choose $x_1=x+c_1, x_2=\pm 1, x_4=y+c_4$ and $(x_3,x_5)\equiv (c_3,c_5)\pmod 1$ arbitrarily so that $Q=\pm x+\beta_{Q}' y-(\frac{1}{2}+ah_4^2)y^2+\nu$ where, $\beta_{Q}'=\pm a_4-2a(x_3+h_4c_4+h_5x_5)$ and $\nu$ is some real number. Since, $ah_4^2+\frac{1}{2}\leq \frac{a}{4}+\frac{1}{2}=\frac{2d^5}{4}+\frac{1}{2}<d$, therefore on taking $h=1/2, k=1$ we find that $\vline~\frac{1}{2}-(\frac{1}{2}+ah_4^2)~\vline+\frac{1}{2}<d$  and hence (\ref{eq101*}) is soluble unless $\frac{1}{2}+ah_4^2=\frac{1}{2}$ i.e., unless $h_4=0$ which gives  $a_4=a'_4=1/2$. Similarly, we get $h_5=0$ and $a_5=1/2$. \vspace{2mm}
	
	\noindent Now, to fix values of $a$, choose $x_1=x+c_1, x_2=\pm 1, x_3=y+c_3$ and $(x_4,x_5)\equiv (c_4,c_5)\pmod 1$ arbitrarily so that $Q=\pm x+\beta_{a,Q} y-ay^2+\nu$ where, $\beta_{a,Q}=\pm a_3-2ac_3$ and $\nu$ is some real number.\vspace{2mm}
	
	\noindent Using (\ref{EQ2}), we get $0.1486< d/4\leq a=2d^5<0.18922=\lambda_0(say)$. We divide the range of $a$ into $6$ subintervals $[\lambda_n,\lambda_{n-1}]$, $n=1,2,\cdots 6$ and in each subinterval choose suitable integers $2h_n$ and $k_n$ such that
	\begin{equation}\label{EQ3}
	\vline~h_n-ak_n^2~\vline +\frac{1}{2}<d.
	\end{equation}
	\noindent The choice of $(h_n,k_n)$ is done in the following way.  For a fixed value of the pair $(h_n,k_n)$ $2h_n,k_n \in \mathbb{Z}$, (\ref{EQ3}) is satisfied if $a$ lies in the interval $I_n=\big(a_n(d),b_n(d)\big)$,  where $a_n(d)= (h_n+1/2-d)/k^2_n$, $ b_n(d)=(h_n-1/2+d)/k^2_n$ and $h_n/k^2_n $ is the center of the interval $I_n$, $I_n\supset [\lambda_n,\lambda_{n-1}]$. The intervals $I_n$ are overlapping. Therefore by Macbeath Lemma, (\ref{eq101*}) is soluble for $ a\in I_n$, in particular for $a\in [\lambda_n,\lambda_{n-1}]$ (unless $a=\frac{h_n}{k_n^2}$).
	\begin{enumerate}
		\item Take $(h_1,k_1)=(3,4)$. The center of $I_1$ namely $h_1/k^2_1 =\frac{3}{16}<\lambda_0$. Note that $a=2d^5<\lambda_0<b_1(d)=\frac{1}{16}(2.5+d) $ for $0.5946<d<0.624$. Also $a_1(d)=\frac{1}{16}(3.5-d)< 2d^5$ for $d\geq \gamma_1$ where $\gamma_1$ is a positive real root of $32d^5+d-3.5=0$ satisfying $0.6178<\gamma_1<0.6179$. So $a_1(d)=\frac{1}{16}(3.5-d)< 0.1802=\lambda_1$ (say) for $d\geq \gamma_1>0.6178$.
		
		\item Take $(h_2,k_2)=(4.5,5)$.
		Note that $a=2d^5<\lambda_1<b_2(d)=\frac{1}{25}(4+d) $ for $0.5946<d<0.6179$. Also $a_2(d)=\frac{1}{25}(5-d)< 2d^5$ for $d\geq \gamma_2$ where $\gamma_2$ is a positive real root of $50d^5+d-5=0$ satisfying $0.6145<\gamma_2<0.61465$. So $a_2(d)=\frac{1}{25}(5-d)< 0.1755=\lambda_2$ (say) for $d\geq \gamma_2>0.6145$.
		
		\item Take $(h_3,k_3)=(1.5,3)$.
		Note that $a=2d^5<\lambda_2<b_3(d)=\frac{1}{9}(1+d) $ for $0.5946<d<0.61465$.
		Also $a_3(d)=\frac{1}{9}(2-d)<2d^5$ for $d\leq \gamma_3$, where $\gamma_3$ is a positive real root of $18d^5+d-2=0$ satisfying $0.6<\gamma_3<0.6001$. So $a_3(d)=\frac{1}{9}(2-d)< 0.156=\lambda_3$ (say) for $d\geq \gamma_3>0.6$.
		
		\item  Take $(h_4,k_4)=(2.5,4)$.   Note that $a=2d^5<\lambda_3<b_4(d)=\frac{1}{16}(2+d) $ for $0.5946<d<0.6001$. Also $a_4(d)=\frac{1}{16}(3-d)< 2d^5$ for $d\leq \gamma_4$, where $\gamma_4$ is a positive real root of $32d^5+d-3=0$ satisfying $0.5958<\gamma_4<0.5959$. So $a_4(d)=\frac{1}{16}(3-d)< 0.1503=\lambda_4$ (say) for $d\geq \gamma_4>0.5958$.
		
		\item  Take $(h_5,k_5)=(15,10)$.   Note that $a=2d^5<\lambda_4<b_5(d)=\frac{1}{100}(14.5+d) $ for $0.5946<d<0.5959$. Also $a_5(d)=\frac{1}{100}(15.5-d)< 2d^5$ for $d\leq \gamma_5$, where $\gamma_5$ is a positive real root of $200d^5+d-15.5=0$ satisfying $0.5949<\gamma_5<0.59494$. So $a_5(d)=\frac{1}{100}(15.5-d)< 0.1491=\lambda_5$ (say) for $d\geq \gamma_5>0.5949$.
		
		\item  Take $(h_6,k_6)=(9.5,8)$.   Note that $a=2d^5<\lambda_5<b_6(d)=\frac{1}{64}(9+d) $ for $0.5946<d<0.59494$. Also $a_6(d)=\frac{1}{64}(10-d)< 0.1486=\lambda_6$ for $d>0.5946$.
		
	\end{enumerate}
	
	\noindent Therefore we find that  (\ref{eq101*}) is soluble unless $a=h_n/k_n^2$, $n=1,2,\cdots,6$.  For $n=2$, i.e., for $a=9/50$, we find that (\ref{EQ3}) is satisfied by another pair $(h,k)=(11.5,8)$ as well. For $n=4$ i.e., for $a=2.5/16$, $\lambda_3=0.156<2.5/16$, so this value need not be considered.  For $n=6$, $h_n/k_n^2=\frac{9.5}{64}<\lambda_6$. Therefore, by Macbeath's Lemma, (\ref{eq101*}) is soluble unless $a=h_n/k_n^2,$ for $n=1,3,5$ i.e., unless $a=3/16, 1/6 {\rm ~or~} 3/20$ and $\beta_Q\equiv h_n/k_n\pmod{1/k_n,2a}$. In each of these cases, we will give solutions of (\ref{eq101*}) depending upon values of $a_2$ and $c_1$.\vspace{4mm}

	\noindent\textbf{Subcase (i): $a=\frac{3}{16}$. } \vspace{2mm}
	
	\noindent Here  $a=2d^5=\frac{3}{16}$  gives
	$d=\left(\frac{3}{32}\right)^{1/5}<0.62287$.
	\noindent Using Lemma \ref{lem15} we need to consider $(a_3,c_3)=(0,0), (0,1/3), (1/16,1/6) {\rm ~and~} (1/16,1/2).$\vspace{2mm}
	
	\noindent When $(a_3,c_3)=(0,0)$, on taking $x_4=0=x_5$, the following table gives a solution to (\ref{eq101*}) i.e., of
	 \begin{equation}\label{EQN2}0<Q=(x_1+a_2x_2+\frac{1}{2}x_4+\frac{1}{2}x_5)x_2-\frac{3}{16}x_3^2-\frac{1}{2}
	x_4^2-\frac{1}{2}x_5^2<d.\end{equation}
	
	{\scriptsize
		\begin{equation*}\begin{array}{lllll}
		\hline\\
		{\rm ~~~~Range} & x_1 & x_1x_2 &x_3& {\rm ~~Q}\vspace{2mm}\\
		\hline\\
		
		\frac{3}{4}<f_1\leq 1&c_1&|c_1|&\pm 2&f_1-\frac{3}{4}\vspace{2mm}\\
		
		\frac{3}{16}<f_1\leq \frac{3}{4}&c_1&|c_1|&\pm 1&f_1-\frac{3}{16}\vspace{2mm}\\
		
		-\frac{1}{4}<f_1\leq \frac{3}{16}&\pm 1+c_1&1+|c_1|&\pm 2&f_1+\frac{1}{4}\vspace{2mm}\\
		
		-\frac{1}{2}<f_1\leq -\frac{1}{4}&\pm 1+c_1&1+|c_1|&\pm 1&f_1+\frac{13}{16}\vspace{2mm}\\

		\hline
		\end{array}
		\end{equation*}\vspace{1mm}}
	
	\noindent When $(a_3,c_3)=(0,1/3)$, on taking $x_4=0=x_5$, the following table gives a solution to (\ref{EQN2}).
	
	{\scriptsize
		\begin{equation*}\begin{array}{lllll}
		\hline\\
		{\rm ~~~~Range} & x_1 & x_1x_2 &x_3& {\rm ~~Q}\vspace{2mm}\\
		\hline\\
		
		\frac{25}{48}<f_1\leq 1&c_1&|c_1|&-\frac{5}{3}&f_1-\frac{25}{48}\vspace{2mm}\\
		
		\frac{1}{48}<f_1\leq \frac{25}{48}&c_1&|c_1|&\frac{1}{3}&f_1-\frac{1}{48}\vspace{2mm}\\
		
		-\frac{23}{48}<f_1\leq \frac{1}{48}&\pm 1+c_1&1+|c_1|&-\frac{5}{3}&f_1+\frac{23}{48}\vspace{2mm}\\
		
		-\frac{1}{2}<f_1\leq -\frac{23}{48}&\pm 1+c_1&1+|c_1|&\frac{1}{3}&f_1+\frac{47}{48}\vspace{2mm}\\
		
		\hline
		\end{array}
		\end{equation*}\vspace{1mm}}
	
	\noindent When $(a_3,c_3)=(1/16,1/6)$, on taking $x_2=1$ and $x_4=0=x_5$, the following table gives a solution to (\ref{eq101*}) i.e., of\vspace{-3mm}
	 \begin{equation}\label{EQN3}0<Q=(x_1+a_2x_2+\frac{1}{16}x_3+\frac{1}{2}x_4+\frac{1}{2}x_5)x_2-\frac{3}{16}x_3^2-\frac{1}{2}
	x_4^2-\frac{1}{2}x_5^2<d.\end{equation}
	
	{\scriptsize
		\begin{equation*}\begin{array}{llll}
		\hline\\
		{\rm ~~~~Range} & x_1  &x_3& {\rm ~~Q}\vspace{2mm}\\
		\hline\\
		
		\frac{143}{192}<p_1\leq 1&c_1&\frac{13}{6}&p_1-\frac{143}{192}\vspace{2mm}\\
		
		\frac{35}{192}<p_1\leq \frac{143}{192}&c_1&\frac{7}{6}&p_1-\frac{35}{192}\vspace{2mm}\\	
		
		-\frac{49}{192}<p_1\leq \frac{35}{192}&1+c_1&\frac{13}{6}&p_1+\frac{49}{192}\vspace{2mm}\\
		
		-\frac{157}{192}<p_1\leq -\frac{49}{192}&1+c_1&\frac{7}{6}&p_1+\frac{157}{192}\vspace{2mm}\\
		
		-1<p_1\leq -\frac{157}{192}&1+c_1&\frac{1}{6}&p_1+\frac{193}{192}\vspace{2mm}\\
		
		\hline
		\end{array}
		\end{equation*}\vspace{1mm}}
	
	\noindent When $(a_3,c_3)=(1/16,1/2)$, on taking $x_4=0=x_5$, the following table gives a solution to (\ref{EQN3}).
	
	{\scriptsize
		\begin{equation*}\begin{array}{llllll}
		\hline\\
		{\rm ~~~~Range} & x_1 & x_1x_2 &x_3&x_2x_3& {\rm ~~Q}\vspace{2mm}\\
		\hline\\
		
		\frac{33}{64}<f_1\leq 1&c_1&|c_1|&\frac{3}{2}&-\frac{3}{2}&f_1-\frac{33}{64}\vspace{2mm}\\
		
		\frac{5}{64}<f_1\leq \frac{33}{64}&c_1&|c_1|&\frac{1}{2}&-\frac{1}{2}&f_1-\frac{5}{64}\vspace{2mm}\\
		
		-\frac{31}{64}<f_1\leq \frac{5}{64}&\pm 1+c_1&1+|c_1|&\frac{3}{2}&-\frac{3}{2}&f_1+\frac{31}{64}\vspace{2mm}\\
		
		-\frac{1}{2}<f_1\leq -\frac{31}{64}&\pm 1+c_1&1+|c_1|&\frac{1}{2}&-\frac{1}{2}&f_1+\frac{59}{64}\vspace{2mm}\\

		\hline
		\end{array}
		\end{equation*}\vspace{1mm}}
	
	\noindent\textbf{Subcase (ii): $a=\frac{1}{6}$.}  \vspace{2mm}
	
	\noindent Here  $a=2d^5=\frac{1}{6}$  gives
	$d=\left(\frac{1}{12}\right)^{1/5}<0.60837$.
	Using Lemma \ref{lem15} we need to consider $(a_3,c_3)=(0,1/2) {\rm ~and~} (1/6,0).$

	\noindent When $(a_3,c_3)=(0,1/2)$, on taking $x_5=0$, the following table gives a solution to (\ref{eq101*}) i.e., of
	 $$0<Q=(x_1+a_2x_2+\frac{1}{2}x_4+\frac{1}{2}x_5)x_2-\frac{1}{6}x_3^2-\frac{1}{2}x_4^2-\frac{1}{2}
	x_5^2<d.$$
	
	{\scriptsize
		\begin{equation*}\begin{array}{lllllll}
		\hline\\
		{\rm ~~~~Range} & x_1 & x_1x_2 &x_3& x_4&x_2x_4& {\rm ~~Q}\vspace{2mm}\\
		\hline\\
		
		d+\frac{3}{8}\leq f_1\leq 1&c_1&2|c_1|&\frac{5}{2}&\pm1&-2&f_2-\frac{61}{24}\vspace{2mm}\\
		
		\frac{3}{8}<f_1<d+\frac{3}{8}&c_1&|c_1|&\frac{3}{2}&0&0&f_1-\frac{3}{8}\vspace{2mm}\\
		
		\frac{1}{24}<f_1\leq \frac{3}{8}&c_1&|c_1|&\frac{1}{2}&0&0&f_1-\frac{1}{24}\vspace{2mm}\\
		
		-\frac{1}{2}<f_1< d-\frac{5}{8}&\pm 1+c_1&1+|c_1|&\frac{3}{2}&0&0&f_1+\frac{5}{8}\vspace{2mm}\\
		
		{\rm When~}  d-\frac{5}{8}\leq f_1\leq \frac{1}{24}: \vspace{2mm}\\
		
		d-\frac{5}{8}\leq g_1\leq \frac{1}{24}& c_1&2|c_1|&\frac{3}{2}&\pm1&2&f_2+\frac{1}{8}\vspace{2mm}\\
		
		-\frac{5}{8}<g_1< d-\frac{5}{8}&\pm 1+c_1&1-|c_1|&\frac{3}{2}&0&0&g_1+\frac{5}{8}\vspace{2mm}\\
		
		-\frac{23}{24}<g_1\leq -\frac{5}{8}&\pm 1+c_1&1-|c_1|&\frac{1}{2}&0&0&g_1+\frac{23}{24}\vspace{2mm}\\
		
		-1<g_1\leq -\frac{23}{24}&\pm 2+c_1&4+2|c_1|&\frac{9}{2}&\pm1&2&f_2+\frac{9}{8}\vspace{2mm}\\
		
		\hline
		\end{array}
		\end{equation*}\vspace{1mm}}	
	
	\noindent When $(a_3,c_3)=(1/6,0)$, on taking $x_5=0$, the following table gives a solution to (\ref{eq101*}) i.e., of
	$$0<Q=(x_1+a_2x_2+\frac{1}{6}x_3+\frac{1}{2}x_4+\frac{1}{2}x_5)x_2-\frac{1}{6}x_3^2-
	\frac{1}{2}x_4^2-\frac{1}{2}x_5^2<d.$$
	
	{\scriptsize
		\begin{equation*}\begin{array}{llllllll}
		\hline\\
		{\rm ~~~~Range} & x_1 & x_1x_2 &x_3&x_2x_3& x_4&x_2x_4& {\rm ~~Q}\vspace{2mm}\\
		\hline\\
		
		d+\frac{1}{3}\leq f_1\leq 1&\pm 1+c_1&-2+2|c_1|&\pm 1&-2&0&0&f_2-\frac{5}{2}\vspace{2mm}\\
		
		\frac{1}{3}<f_1<d+\frac{1}{3}&c_1&|c_1|&\pm 1&-1&0&0&f_1-\frac{1}{3}\vspace{2mm}\\
		
		0<f_1\leq \frac{1}{3}&c_1&|c_1|& 0&0&0&0&f_1\vspace{2mm}\\
		
		-\frac{1}{2}\leq f_1<d-\frac{2}{3}&\pm 1+c_1&1+|c_1|&\pm 1&-1&0&0&f_1+\frac{2}{3}\vspace{2mm}\\

		{\rm When~}  d-\frac{2}{3}\leq f_1\leq 0: \vspace{2mm}\\

		d-\frac{2}{3}\leq g_1\leq 0&c_1&2|c_1|&0&0&\pm 1&2&f_2+\frac{1}{2}\vspace{2mm}\\
		
		-\frac{2}{3}<g_1<d-\frac{2}{3}&\pm 1+c_1&1-|c_1|&\pm 1&-1&0&0&g_1+\frac{2}{3}\vspace{2mm}\\
		
		-1<g_1\leq -\frac{2}{3}&\pm 1+c_1&1-|c_1|& 0&0&0&0&g_1+1\vspace{2mm}\\

		\hline
		\end{array}
		\end{equation*}\vspace{1mm}}	
	
	\noindent\textbf{Subcase (iii): $a=\frac{3}{20}$.}  \vspace{2mm}
	
	\noindent Here  $a=2d^5=\frac{3}{20}$  gives
	$d=\left(\frac{3}{40}\right)^{1/5}<0.59568$.
	Using Lemma \ref{lem15} we need to consider $(a_3,c_3)=(0,0), (0,1/3), (1/20,1/6) {\rm ~and~} (1/20,1/2).$\vspace{2mm}
	
	\noindent For each of these values of $(a_3,c_3)$ we take $x_2=\pm 1, x_4=0=x_5$ and give suitable values to  $x_3$, which covers the range of $f_1''$ or, $p_1''$ and hence (\ref{eq101'}) is soluble.\vspace{2mm}
	
	\noindent This completes the proof of Lemma \ref{lem17}.\vspace{2mm} \hfill $\Box$
	
	\noindent Now the proof of Theorem \ref{thm4} follows from Lemmas \ref{lem15'} - \ref{lem17}.

	\section{Proof of $\Gamma_{1,4}<8$ for $c_2 \equiv 0 \pmod 1$, $~a< \frac{1}{2}$,  $m=2, K=3,4$} \label{sec7}
	\numberwithin{equation}{section}
	\begin{theorem}\label{thm4'}
		Let $Q(x_1,x_2, \cdots, x_5)$ be a real indefinite quadratic form of type $(1,4)$ and of determinant $D\neq 0$. Let $d=(8|D|)^{\frac{1}{5}}$.  Suppose that $c_2 \equiv 0 \pmod1$, $a<\frac{1}{2}$, $d\leq 1$ and $a+d\leq 1$. Let $(m,K)=(2,4),(2,3)$.  Then $(\ref{eq101*})$ and hence $(\ref{eq10})$ is soluble with strict inequality.
	\end{theorem}
	
	\noindent \textbf{Proof:} We divide the proof into following subsections:

	\subsection{$c_2 \equiv 0 \pmod1$, $a<\frac{1}{2}$, $m=2, K=4$}
	
	\begin{lemma}\label{lem18}  When $m=2,K=4$, inequality $(\ref{eq101*})$ is soluble.
	\end{lemma}
	
	\noindent {\bf Proof:} Here $d\leq \frac{4}{m+3}, 2<\frac{d}{a}\leq 3$ and $4<\frac{d+(m^2-1)a/4}{A}\leq 5$. From (\ref{eq7}) and (\ref{eq18}) we have
	\begin{equation}\label{eq72}
	d/4 \leq A \leq \sqrt{2d^5/3a} \leq \sqrt{2}~d^2.
	\end{equation}
	
	\noindent Also we have $d/3 \leq a\leq d^{5/3}$. This together with $d\leq \frac{16}{(m+3)(K+3)}$, gives
	\begin{equation}\label{eq73}
	0.1924< 1/\sqrt{27}\leq d \leq 16/35.
	\end{equation}
	\noindent Further we have,  \begin{equation}\label{eq55''} \begin{array}{ll}
	\frac{\delta_{m,K}}{t}&=\frac{2aA}{d^5} \delta_{m,K}
	\geq \frac{2}{d^5} \frac{d}{m+1}\frac{(m+3)d}{4(K+1)}\frac{(m+3)}{4} \frac{(K+3)}{4}d\vspace{4mm}\\
	&=\frac{1}{32}\frac{(m+3)^2(K+3)}{(m+1)(K+1)}\frac{1}{d^2},
	\end{array}
	\end{equation}
	
	\noindent which is $>1$ for $m=2,K=4$ and $d\leq 16/35.$  From the definition of $L$ we get $L\geq 1.$ For $L\geq 7$ and $d\geq 1/\sqrt{27}$, we find that $\delta_{m,K,L}>\frac{(m+3)}{4}\frac{(K+3)}{4}\frac{(L+3)}{4}d>1$. Therefore, we can assume that $L\leq 6$. Further, from (\ref{eq55''}) and (\ref{eq26'}) we have
	
	\begin{equation*}d^2 \geq \frac{1}{32}\frac{(m+3)^2(K+3)}{(m+1)(K+1)(L+1)} {\rm ~ and, ~} d\leq \frac{64}{(m+3)(K+3)(L+3)}\end{equation*}
	\noindent and hence,
	\begin{equation}\label{eq74}
	\delta_{2,4,L}\leq 1,{\rm ~and~} \sqrt{35/96(L+1)} \leq d \leq 64/35(L+3).
	\end{equation}
	\noindent The cases $L=4,5$ and $6$ do not arise because then we have from (\ref{eq74})
	\begin{equation*}
	L=4 : ~~ 0.27003< \sqrt{7/96} \leq d \leq 64/245<0.26123
	\end{equation*}
	\begin{equation*}
	L=5 :~~ 0.2465< \sqrt{35/576} \leq d \leq 8/35<0.22858
	\end{equation*}
	\begin{equation*}
	L= 6 : ~~ 0.2282< \sqrt{5/96} \leq d \leq 64/315<0.20318
	\end{equation*}
	
	\noindent which is absurd. We distinguish the cases $L=1,2$ and $3$.\vspace{2mm}

	\noindent \textbf{Case I:} $\textbf{m=2, K=4, L=3}$.\vspace{2mm}
	
	\noindent Here $\delta_{m,K}=d+\frac{3a}{4}+\frac{15A}{4}\geq d+\frac{d}{4}+\frac{15d}{16}=\frac{35d}{16}$. From (\ref{eq74}) we have
	\begin{equation}\label{eq79}
	0.3019< \sqrt{35/384} \leq d \leq 32/105<0.30477.
	\end{equation}
	\noindent Further $$0.1651< \frac{35d}{64}\leq \frac{\delta_{m,K}}{4}=\frac{d+3a/4+15A/4}{4}\leq t =\frac{d^5}{2aA}\leq 6d^3< 0.16986.$$
	
	\noindent Taking $h_t=3/2$, $k_t=3$ and using (\ref{eq79})
	we find that $\vline~h_t-tk_t^2~\vline + \frac{1}{2} <\delta_{m,K}$, as $9t<\delta_{m,K}+1$ and $\delta_{m,K}+9t \geq \frac{35d}{16}+\frac{315d}{64}>2$. Therefore,  (\ref{eq101'}) and hence (\ref{eq17}) is soluble unless $t=1/6$ and $\beta_{t,G}\equiv \pm a_5''-\frac{1}{3}c_5= \frac{1}{2}\pmod {\frac{1}{3}}$. This gives, by Lemma \ref{lem15}, $(a_5'',c_5)=(0,1/2)$ or $(1/6,0)$. If $a_5''=1/6,$ considering the transformation $x_5\rightarrow x_5-x_2$, we may assume that $a_5''=1/2$.\vspace{2mm}

	\noindent Now, $1/6=t=d^5/2aA\leq 6d^3$ gives $d\geq (1/36)^{1/3}>0.3028$, $ \delta_{2,4}\geq \frac{35d}{16}> 0.66249.$  Also, $1/6=t=d^5/2aA$ gives $A=3d^5/a\leq 9d^4$. Hence,
	\begin{equation}\label{eq80}
	\begin{array}{l}
	0.1009< d/3\leq a\leq d/2< 0.152385\vspace{2mm}\\
	0.0757< d/4 \leq A \leq 9d^4< 0.07765\vspace{2mm}\\
	0.04415< \frac{a+A}{4} < 0.05751.
	\end{array}
	\end{equation}
	
	\noindent Thus we need to find the solutions of
	\begin{equation}\label{eq81}
	0<G=(x_1+a_2''x_2+a_5''x_5)x_2-x_5^2/6-(a+A)/4<\delta_{2,4}=d+\frac{3a}{4}+\frac{15A}{4}
	\end{equation}
	for $(a_5'',c_5)=(0,1/2)$ and $(1/2,0)$. We work as in the case $m=3,K=2,L=3$ (Section \ref{subsec6.2}) and find that (\ref{eq81}) is soluble.\\

	\noindent \textbf{Case II:} $\textbf{m=2, K=4, L=2}$.\vspace{2mm}
	
	\noindent Here $\delta_{m,K}=d+\frac{3a}{4}+\frac{15A}{4}\geq d+\frac{d}{4}+\frac{15d}{16}=\frac{35d}{16}$. From (\ref{eq74}) we have
	\begin{equation}\label{eq78}
	0.3486< \sqrt{35/288} \leq d \leq 64/175<0.36572.
	\end{equation}
	\noindent Further $$0.25418< \frac{35d}{48}\leq \frac{\delta_{m,K}}{3}=\frac{d+3a/4+15A/4}{3}\leq t =\frac{d^5}{2aA}\leq 6d^3< 0.2935.$$
	
	\noindent Taking $h_t=1$, $k_t=2$ and using (\ref{eq78})
	we find that $\vline~h_t-tk_t^2~\vline + \frac{1}{2} <\delta_{m,K}$, as $\delta_{m,K}+\frac{1}{2}-4t \geq \frac{35}{16}d+\frac{1}{2}-24d^3>0$. Note that here $t>\frac{1}{4}$. Therefore,  (\ref{eq101'}) and hence (\ref{eq17}) is soluble with strict inequality.\vspace{2mm}

	\noindent \textbf{Case III}: $\textbf{m=2, K=4, L=1}$.\vspace{2mm}
	
	\noindent Here $\delta_{m,K}=d+\frac{3a}{4}+\frac{15A}{4}\geq d+\frac{d}{4}+\frac{15d}{16}=\frac{35d}{16}$. From (\ref{eq74}) we have
	\begin{equation}\label{eq75}
	0.4269< \sqrt{35/192} \leq d \leq 16/35=0.45715.
	\end{equation}
	\noindent Further $$0.4669< \frac{35d}{32}\leq \frac{\delta_{m,K}}{2}=\frac{d+3a/4+15A/4}{2}\leq t =\frac{d^5}{2aA}\leq 6d^3< 0.57323.$$
	
	\noindent Taking $h_t=1/2$, $k_t=1$ and using (\ref{eq75})
	we find that $\vline~h_t-tk_t^2~\vline + \frac{1}{2} <\delta_{m,K}$, as $ t < \delta_{m,K}$ and   $t+\delta_{m,K} \geq \frac{35}{32}d+\frac{35}{16}d=\frac{105}{32}d>1$. Therefore,  (\ref{eq101'})  is soluble unless $t=1/2$ and $\beta_{t,G}=\pm a_5''-c_5\equiv 1/2 \pmod 1$. This gives $(a_5'',c_5)=(0,1/2)$ or $(1/2,0)$.\vspace{2mm}

	\noindent Now, $1/2=t=d^5/2aA\leq 6d^3$ gives $d\geq (1/12)^{1/3}>0.4367$. Also $\delta_{2,4}\geq \frac{35d}{16}>0.9554$, $1/2=t=d^5/2aA$ gives $A=d^5/a\leq 3d^4$. Hence,
	\begin{equation}\label{eq76}
	\begin{array}{l}
	0.1455< d/3\leq a\leq d/2< 0.228575\vspace{2mm}\\
	0.1091< d/4 \leq A \leq 3d^4< 0.13103\vspace{2mm}\\
	0.06365< \frac{a+A}{4} < 0.08991.
	\end{array}
	\end{equation}
	
	\noindent Thus we need to find a solution of
	\begin{equation}\label{eq77}
	0<G=(x_1+a_2''x_2+a_5''x_5)x_2-x_5^2/2-(a+A)/4<d+3a/4+15A/4
	\end{equation}
	for $(a_5'',c_5)=(0,1/2)$ and $(1/2,0)$.\vspace{2mm}
	
	\noindent When  $(a_5'',c_5)=(0,1/2)$, the following table gives a solution of (\ref{eq77}).\vspace{2mm}
	
	{\scriptsize
		\begin{equation*}\begin{array}{lllll}
		\hline\\
		{\rm ~~~~Range} & x_1 & x_1x_2 & x_5 & {\rm ~~G}\vspace{2mm}\\
		\hline\\
		
		\frac{1}{8}+\frac{a+A}{4}<f_1''\leq 1&c_1&|c_1|&\frac{1}{2}&f_1''-\frac{1}{8}-\frac{a+A}{4}\vspace{2mm}\\
		
		-\frac{1}{2}<f_1''<-\frac{7}{8}+d+a+4A&\pm 1+c_1&1+|c_1|&\frac{1}{2}&f_1''+\frac{7}{8}-\frac{a+A}{4}\vspace{2mm}\\

		{\rm When~}  -\frac{7}{8}+d+a+4A\leq f_1''\leq \frac{1}{8}+\frac{a+A}{4}: \vspace{2mm}\\

		-\frac{7}{8}+d+a+4A\leq g_1''\leq \frac{1}{8}+\frac{a+A}{4}&c_1&2|c_1|&\frac{1}{2}&f_2-\frac{1}{8}-\frac{a+A}{4}\vspace{2mm}\\
		
		-\frac{7}{8}+\frac{a+A}{4}<g_1''<-\frac{7}{8}+d+a+4A&\pm 1+c_1&1-|c_1|&\frac{1}{2}&g_1''+\frac{7}{8}-\frac{a+A}{4}\vspace{2mm}\\
		
		-1<g_1''\leq -\frac{7}{8}+\frac{a+A}{4}&\pm 1+c_1&2+2|c_1|&\frac{3}{2}&f_2''+\frac{7}{8}-\frac{a+A}{4}\vspace{2mm}\\

		\hline
		\end{array}
		\end{equation*}\vspace{2mm}}
	
	\noindent When  $(a_5'',c_5)=(1/2,0)$, on taking $x_5=0$, the following table gives a solution of (\ref{eq77}).\vspace{2mm}
	
	{\scriptsize
		\begin{equation*}\begin{array}{llll}
		\hline\\
		{\rm ~~~~Range} & x_1 & x_1x_2  & {\rm ~~G}\vspace{2mm}\\
		\hline\\
		
		\frac{a+A}{4}<f_1''\leq 1&c_1&|c_1|&f_1''-\frac{a+A}{4}\vspace{2mm}\\
		
		-\frac{1}{2}<f_1''<-1+d+a+4A&\pm 1+c_1&1+|c_1|&f_1''+1-\frac{a+A}{4}\vspace{2mm}\\

		{\rm When~}  -1+d+a+4A\leq f_1''\leq \frac{a+A}{4}: \vspace{2mm}\\

		-1+d+a+4A\leq g_1''\leq \frac{a+A}{4}&c_1&2|c_1|&f_2''-\frac{a+A}{4}\vspace{2mm}\\
		
		-1+\frac{a+A}{4}<g_1''<-1+d+a+4A&\pm 1+c_1&1-|c_1|&g_1''+1-\frac{a+A}{4}\vspace{2mm}\\
		
		-1<g_1''\leq -1+\frac{a+A}{4}&\pm 1+c_1&3+3|c_1|&f_3''+3-\frac{a+A}{4}\vspace{2mm}\\

		\hline
		\end{array}
		\end{equation*}\vspace{2mm}}
	
	\noindent This completes the proof of Lemma \ref{lem18}.\hfill $\square$

	\subsection{$c_2 \equiv 0 \pmod1$, $a<\frac{1}{2}$, $m=2, K=3$}\label{subsec7.2}

	\begin{lemma}\label{lem19}  When $m=2,K=3$, inequality (\ref{eq101*}) is soluble.
	\end{lemma}
	
	\noindent {\bf Proof:} Here $d\leq \frac{4}{m+3}, 2<\frac{d}{a}\leq 3$ and $3<\frac{d+(m^2-1)a/4}{A}\leq 4$. From (\ref{eq7}) and (\ref{eq18}) we have
	\begin{equation}\label{eq82}
	5d/16 \leq A \leq \sqrt{2d^5/3a} \leq \sqrt{2}~d^2.
	\end{equation}

	\noindent This together with $d\leq \frac{16}{(m+3)(K+3)}$, gives
	\begin{equation}\label{eq83}
	0.2209< 5\sqrt{2}/32\leq d \leq 8/15< 0.53334.
	\end{equation}
	
	\noindent Further we have, from (\ref{eq55''})
	\begin{equation*}
	\frac{\delta_{m,K}}{t}\geq \frac{1}{32}\frac{(m+3)^2(K+3)}{(m+1)(K+1)}\frac{1}{d^2}
	\end{equation*}
	
	\noindent which is $>1$ for $m=2,K=3$ and $d\leq 8/15.$  From the definition of $L$ we get $L\geq 1.$ For $L\geq 7$ and $d\geq 5\sqrt{2}/32$, we find that $\delta_{m,K,L}>\frac{(m+3)}{4}\frac{(K+3)}{4}\frac{(L+3)}{4}d>1$. \vspace{2mm}
	
	\noindent Therefore, we can assume that $L\leq 6$. Further, from (\ref{eq55''}) and (\ref{eq26'}) we have
	
	\begin{equation*}d^2 \geq \frac{1}{32}\frac{(m+3)^2(K+3)}{(m+1)(K+1)(L+1)} {\rm ~ and, ~} d\leq \frac{64}{(m+3)(K+3)(L+3)}\end{equation*}
	
	\noindent and hence,
	
	\begin{equation}\label{eq84}
	\delta_{2,3,L}\leq 1,{\rm ~and~} \sqrt{25/64(L+1)} \leq d \leq 32/15(L+3).
	\end{equation}
	
	\noindent We distinguish the cases $L=1,2,3,4,5$ and $6$.\vspace{4mm}

	\noindent \textbf{Case I:} $\textbf{m=2, K=3, L=6}.$\vspace{2mm}
	
	\noindent Here $6<\frac{\delta_{2,3}}{t}=\frac{d+3a/4+2A}{t}\leq 7$, $\delta_{2,3}\geq \frac{15}{8}d$, $\delta_{2,3,6}=\delta_{2,3}+\frac{L^2-1}{4}t =\delta_{2,3}+\frac{35}{4}t$. From (\ref{eq84}) we have
	\begin{equation}\label{eq125}
	0.23622< \sqrt{25/448} \leq d \leq 32/135<0.23704.
	\end{equation}
	
	\noindent Further, since $a\geq d/3$ and $A\geq 5d/16$, we get
	
	\begin{equation}\label{eq126}
	\begin{array}{l}
	0.06327< \frac{15d}{56}\leq \frac{\delta_{2,3}}{7}\leq t =\frac{d^5}{2aA}\leq \frac{24d^3}{5}< 0.06394\vspace{2mm}\\
	\frac{a+A+t}{4}\geq \frac{307}{1344}d\geq 0.05395\vspace{2mm}\\
	\delta_{2,3,6}=\delta_{2,3}+\frac{35}{4}t\geq \delta_{2,3}+\frac{75}{32}d\geq \frac{135}{32}d> 0.9965.
	
	\end{array}
	\end{equation}
	
	\noindent Here, we will give a solution of (\ref{eq21}) with strict inequality, i.e., of
	\begin{equation}\label{eq127}
	0<H=(x_1+a_2'''x_2)x_2-(a+A+t)/4<d+3a/4+2A+\frac{35}{4}t= \delta_{2,3,6}.
	\end{equation}

	\noindent Also $\delta_{2,3,6}=d+3a/4+2A+\frac{35}{4}t \leq 1$ gives $$\begin{array}{l}\frac{(a+A+t)}{4}=\frac{1}{3}(\frac{3a}{4}+2A+\frac{35t}{4})-\frac{5A}{12}-\frac{32t}{12} \leq \frac{1-d}{3} - \frac{25d}{192}-\frac{5d}{7}=\frac{1}{3}-\frac{1583d}{1344} < 0.05513.\end{array}$$

	\noindent The following table gives a
	solution of (\ref{eq127}).\vspace{2mm}
	
	{\scriptsize
		\begin{equation*}\begin{array}{llll}
		\hline\\
		{\rm ~~~~Range} & x_1 & x_1x_2  & {\rm ~~H}\vspace{2mm}\\
		\hline\\
		
		\frac{a+A+t}{4}<f_1'''\leq 1&c_1&|c_1|&f_1'''-\frac{a+A+t}{4}\vspace{2mm}\\
		
		-\frac{1}{2}<f_1'''<-1+d+a+\frac{9A}{4}+\frac{25t}{4}&\pm 1+c_1&1+|c_1|&f_1'''+1-\frac{a+A+t}{4}\vspace{2mm}\\

		{\rm When~}  -1+d+a+\frac{9A}{4}+\frac{25t}{4}\leq f_1'''\leq \frac{a+A+t}{4}: \vspace{2mm}\\

		-1+d+a+\frac{9A}{4}+\frac{25t}{4}\leq g_1'''\leq \frac{a+A+t}{4}&c_1&2|c_1|&f_2'''-\frac{a+A+t}{4}\vspace{2mm}\\
		
		-1+\frac{a+A+t}{4}<g_1'''<-1+d+a+\frac{9A}{4}+\frac{25t}{4}&\pm 1+c_1&1-|c_1|&g_1'''+1-\frac{a+A+t}{4}\vspace{2mm}\\
		
		-1<g_1'''\leq -1+\frac{a+A+t}{4}&\pm 2+c_1&6-3|c_1|&g_3'''+6-\frac{a+A+t}{4}\vspace{2mm}\\
		
		\hline
		\end{array}
		\end{equation*}\vspace{1mm}}
	
	\noindent \textbf{Case II:} $\textbf{m=2, K=3, L=5}.$\vspace{2mm}
	
	\noindent Here $5<\frac{\delta_{2,3}}{t}=\frac{d+3a/4+2A}{t}\leq 6$, $\delta_{2,3}\geq \frac{15}{8}d$, $\delta_{2,3,5}=\delta_{2,3}+\frac{L^2-1}{4}t =\delta_{2,3}+6t$. From (\ref{eq84}) we have
	\begin{equation}\label{eq122}
	0.25515< \sqrt{25/384} \leq d \leq 4/15<0.26667.
	\end{equation}
	
	\noindent Further, since $a\geq d/3$ and $A\geq 5d/16$, we get \vspace{-2mm}
	
	\begin{equation}\label{eq123}
	\begin{array}{l}
	0.079734< \frac{5d}{16}\leq \frac{\delta_{2,3}}{6}\leq t =\frac{d^5}{2aA}\leq \frac{24d^3}{5}< 0.091026\vspace{2mm}\\
	\frac{a+A+t}{4}\geq \frac{23d}{96}> 0.06112\vspace{2mm}\\
	\delta_{2,3,5}=\delta_{2,3}+6t\geq \delta_{2,3}+\delta_{2,3}\geq \frac{15}{4}d> 0.9568.
	
	\end{array}
	\end{equation}
	
	\noindent Here, we will give a solution of (\ref{eq21})with strict inequality, i.e., of \vspace{-2mm}
	\begin{equation}\label{eq124}
	0<H=(x_1+a_2'''x_2)x_2-(a+A+t)/4<d+3a/4+2A+6t= \delta_{2,3,5}.\vspace{-2mm}
	\end{equation}
	
	\noindent Also $\delta_{2,3,5}=d+3a/4+2A+6t \leq 1$ gives \vspace{-2mm}
	$$\begin{array}{l}\frac{(a+A+t)}{4}=\frac{1}{3}(\frac{3a}{4}+2A+6t)-\frac{5A}{12}-\frac{7t}{4} \leq \frac{1-d}{3} - \frac{25d}{192}-\frac{35d}{64}=\frac{1}{3}-\frac{194d}{192} < 0.0756.\end{array}\vspace{-2mm}$$

	\noindent The following table gives a solution of (\ref{eq124}).\vspace{1mm}
	
	{\scriptsize
		\begin{equation*}\begin{array}{llll}
		\hline\\
		{\rm ~~~~Range} & x_1 & x_1x_2  & {\rm ~~H}\vspace{2mm}\\
		\hline\\
		
		\frac{a+A+t}{4}<f_1'''\leq 1&c_1&|c_1|&f_1'''-\frac{a+A+t}{4}\vspace{2mm}\\
		
		-\frac{1}{2}<f_1'''<-1+d+a+\frac{9A}{4}+\frac{25t}{4}&\pm 1+c_1&1+|c_1|&f_1'''+1-\frac{a+A+t}{4}\vspace{2mm}\\

		{\rm When~}  -1+d+a+\frac{9A}{4}+\frac{25t}{4}\leq f_1'''\leq \frac{a+A+t}{4}: \vspace{2mm}\\

		-1+d+a+\frac{9A}{4}+\frac{25t}{4}\leq g_1'''\leq \frac{a+A+t}{4}&c_1&2|c_1|&f_2'''-\frac{a+A+t}{4}\vspace{2mm}\\
		
		-1+\frac{a+A+t}{4}<g_1'''<-1+d+a+\frac{9A}{4}+\frac{25t}{4}&\pm 1+c_1&1-|c_1|&g_1'''+1-\frac{a+A+t}{4}\vspace{2mm}\\
		
		-1<g_1'''\leq -1+\frac{a+A+t}{4}&\pm 1+c_1&3+3|c_1|&f_3'''+3-\frac{a+A+t}{4}\vspace{2mm}\\
		
		\hline
		\end{array}
		\end{equation*}\vspace{1mm}}
	
	\noindent \textbf{Case III:} $\textbf{m=2, K=3, L=4}.$\vspace{2mm}
	
	\noindent Here $4<\frac{\delta_{m,K}}{t}=\frac{d+3a/4+2A}{t}\leq 5$, $\delta_{m,K}=d+\frac{3a}{4}+2A\geq d+\frac{d}{4}+\frac{5d}{8}=\frac{15d}{8}$. From (\ref{eq84}) we have
	\begin{equation}\label{eq128}
	0.279508< \sqrt{5/64} \leq d \leq 32/105<0.30477.
	\end{equation}
	
	\noindent Further, since $a\geq d/3$ and $A\geq 5d/16$ we get, $ t =\frac{d^5}{2aA}\leq \frac{24}{5}d^3$.\vspace{2mm}

	\noindent Also since $\delta_{2,3,4}=\delta_{2,3}+\frac{L^2-1}{4}t\leq 1$, we get $t\leq \frac{4}{15}(1-\delta_{2,3}) \leq \frac{4}{15}(1-\frac{15d}{8})$. Therefore,
	\begin{equation}\label{eq891}
	\begin{array}{ll}
	0.1048< \frac{15d}{40}\leq \frac{\delta_{2,3}}{5}\leq t\leq &\min\{\frac{24}{5}d^3,\frac{4}{15}(1-\frac{15d}{8})\}
	=f_0(d) {\rm ~(say)~}\vspace{2mm}\\
	&=\left\{
	\begin{array}{ll}
	\frac{24}{5}d^3 & {\rm ~if~} d\leq \alpha\\
	\frac{4}{15}(1-\frac{15d}{8}) &{\rm ~if~} d\geq \alpha
	\end{array}
	\right.\vspace{2mm}\\
	&<0.121=\lambda_0 ~(<\frac{1}{8}){\rm ~(say)~}
	\end{array}
	\end{equation}
	\noindent  where $\alpha $ is a root of $\frac{24}{5}d^3+\frac{1}{2}d-\frac{4}{15}=0 $ satisfying $0.2926<\alpha<0.2927$.\vspace{2mm}

	\noindent We divide the range of $t$ into $53$ subintervals $[\lambda_n,\lambda_{n-1}], n=1,2,\cdots 53$ and in each subinterval choose suitable integers $2h_n$ and $k_n$ such that
	
	\begin{equation}\label{eq129}
	\vline~ h_n-tk_n^2~\vline +\frac{1}{2} <\frac{15}{8}d\leq d+\frac{3a}{4}+2A
	\end{equation}
	is satisfied. The choice of $(h_n,k_n)$ is done in a manner similar to that in Section 6.3 of \cite{RakaRani}. For a fixed value of the pair $(h_n,k_n)$, $2h_n,k_n \in \mathbb{Z}$, (\ref{eq129}) is satisfied if $t$ lies in the interval $I_n=(a_n,b_n)\supseteq[\lambda_n,\lambda_{n-1}]$, where $a_n=a_n(d)= (h_n+1/2-\frac{15}{8}d)/k^2_n$, $ b_n=b_n(d)=(h_n-1/2+\frac{15}{8}d)/k^2_n>\lambda_{n-1}$ and $h_n/k^2_n $ is the center of the interval $I_n$. For $(h_1,k_1)=(0.5,2)$ where, center of $I_1$ namely $h_1/k^2_1=1/8$ lies to the left of $\lambda_0$ (or is very close to it) and also $f_0(d)<b_1(d)$. So for $t\in I_1=(a_1,b_1)$, (\ref{eq129}) is true.\\
	\noindent Let now $t\leq a_1(d)$, also $t\leq \frac{24}{5}d^3$, so that $t\leq \min\{\frac{24}{5}d^3, a_1(d)\}=f_1(d),$ say. If $\gamma_1$ is the positive real root of the cubic equation $\frac{24}{5}d^3-a_1(d)=0$; we compute some good enough bounds of $\gamma_1$ namely, $0.2882<\gamma_1<0.28823$ and so $t\leq f_1(d)\leq f_1(\gamma_1)<0.115=\lambda_1$ (say). This means that (\ref{eq129}) is satisfied for the region $t\in [\lambda_1, \lambda_{0}]$ by the pair $(h_1,k_1)$.\\
	\noindent We take $(h_2,k_2)=(1,3)$ so that $h_2/k_2^2<\lambda_1$ and $f_1(d)<b_2(d).$ Repeating this process 53 times we find that  (\ref{eq101'}) is soluble unless $t=h_n/k_n^2$, $1\leq n\leq 53$.  The values of $(h_n, k_n), \lambda_n$ are given in the following table. In the column headed by `Remarks', ``tbd" stands for `to be discussed'; ``na" stands for `not applicable' as $h_n/k_n^2$ is out of the interval $[\lambda_n,\lambda_{n-1}]$. If for some $n$, $t=h_n/k_n^2$, (\ref{eq129}) is satisfied by  another pair $(h_n',k_n')$, it is listed in this column.\vspace{2mm}
	
	\noindent We have taken help of software package Mathematica to find the values of pairs $(h_n, k_n), (h_n', k_n')$ and $\lambda_n$, here and in later sections also. We would also like to point out that the choice of $(h_n, k_n), (h_n', k_n')$ and $\lambda_n$, taken by us is not unique.\vspace{2mm}
	
	{\scriptsize
		\begin{equation*}\begin{array}{lllllllll}
		\hline\\
		n & (h_n,k_n) & \lambda_n& \mbox{Remarks}&~~ &n & (h_n,k_n) & \lambda_n& \mbox{Remarks}\vspace{2mm}\\
		\hline\\
		1&(0.5,2)&0.115& \mbox{na} &&2&(1,3)&0.1079&\mbox{tbd}\vspace{2mm}\\
		3&(110.5,32)&0.107882& \mbox{na}&&4&(117.5,33)&0.107871&\mbox{na}\vspace{2mm}\\
		5&(259,49)&0.10786& \mbox{na}&&6&(248.5,48)&0.107844&(314.5,54)\vspace{2mm}\\	
		7&(280.5,51)&0.107832& \mbox{tbd}&&8&(164,39)&0.107805&(606.5,75)\vspace{2mm}\\
		9&(291.5,52)&0.107793& (498.5,68)&&10&(15.5,12)&0.10767&\mbox{tbd}\vspace{2mm}\\
		11&(13,11)&0.107207& (27.5,16)&&12&(21,14)&0.107002& \mbox{tbd}\vspace{2mm}\\
		\hline
		\end{array}
		\end{equation*}}
	{\scriptsize
		\begin{equation*}\begin{array}{lllllllll}
		\hline\\
		n & (h_n,k_n) & \lambda_n& \mbox{Remarks}&~~ &n & (h_n,k_n) & \lambda_n& \mbox{Remarks}\vspace{2mm}\\
		\hline\\
		13&(78,27)&0.1069& (154.5,38)&&14&(109.5,32)&0.106907& (131,35)\vspace{2mm}\\
		15&(216.5,45)&0.10689999& (509,69)&&16&(1069,100)&0.1068973&(1090.5,101)\vspace{2mm}\\
		17&(2247.5,145)&0.10689525& \mbox{tbd}&&18&(944.5,94)&0.1068892& (924.5,93)\vspace{2mm}\\
		19&(278,51)&0.1068715& (790.5,86)&&20&(138.5,36)&0.106847& (188.5,42)\vspace{2mm}\\
		21&(56.5,23)&0.106754& (61.5,24)&&22&(24,15)&0.106547& (38.5,19)\vspace{2mm}\\
		23&(18,13)&0.106352&(83.5,28) &&24&(333.5,56)&0.106338&(345.5,57) \vspace{2mm}\\
		25&(77.5,27)&0.1062737& (153.5,38)&&26&(42.5,20)&0.1061842& \mbox{tbd}\vspace{2mm}\\
		27&(161.5,39)&0.106163&(178.5,41) &&28&(102,31)&0.1061122&(130,35)\vspace{2mm}\\
		29&(95.5,30)&0.1060821&(130,35) &&30&(115.5,33)&0.1060367&\mbox{tbd} \vspace{2mm}\\
		31&(344.5,57)&0.1060246& (332.5,56)&&32&(332.5,56)&0.1060185&(344.5,57) \vspace{2mm}\\
		33&(187,42)&0.105995&(196,43) &&34&(196,43)&0.10599&(187,42) \vspace{2mm}\\
		35&(122.5,34)&0.105982&(108.5,32) &&36&(108.5,32)&0.105947&(122.5,34) \vspace{2mm}\\
		37&(61,24)&0.105909&(145,37) &&38&(56,23)&0.105812&(83,28) \vspace{2mm}\\
		39&(89,29)&0.105807&\mbox{na} &&40&(71.5,26)&0.105732&\mbox{tbd} \vspace{2mm}\\
		41&(129.5,35)&0.105694&(137,36) &&42&(243.5,48)&0.105691&\mbox{na} \vspace{2mm}\\
		43&(214,45)&0.105667&(243.5,48) &&44&(626.5,77)&0.105663&\mbox{na} \vspace{2mm}\\
		45&(875,91)&0.105661&\mbox{na} &&46&(1521.5,120)&0.105659&(1278.5,110) \vspace{2mm}\\
		47&(77,27)&0.10564&\mbox{na} &&48&(66,25)&0.10562&\mbox{na} \vspace{2mm}\\
		49&(30.5,17)&0.10545&(27,16) &&50&(27,16)&0.105372&(30.5,17) \vspace{2mm}\\
		51&(51,22)&0.1053205&(253,49) &&52&(38,19)&0.105195&\mbox{tbd}\vspace{2mm}\\
		53&(8.5,9)&0.1047&(10.5,10) \vspace{2mm}\\
		\hline
		\end{array}
		\end{equation*}\vspace{1mm}}
	
	\noindent Thus by Macbeath Lemma, (\ref{eq101'}) is soluble unless $t=1/9,$ $11/102,$ $31/288,$ $3/28,$ $31/290,$ $17/160,$ $7/66,$ $11/104,$  $ 2/19$ and $\beta_{t,G} = \pm a_5''-2tc_5\equiv h_n/k_n \pmod {1/k_n, 2t}$. For each of these values to $t$ we will give solutions of (\ref{eq101'}) depending upon values of $a_2''$ and $c_1$. We discuss here $t=1/9,$ the remaining values of $t$ are settled similarly while taking $x_2=\pm 1$.\vspace{2mm}
	
	\noindent When $t=\frac{1}{9}$, since $\frac{d}{a}\leq 3$, and $\frac{1}{9}=t=\frac{d^5}{2aA}\leq \frac{24d^3}{5}$, we get $d^3\geq \frac{5}{216}$ i.e., $d\geq (\frac{5}{216})^{\frac{1}{3}}>0.28499$. Also, $\frac{1}{9}=t=\frac{d^5}{2aA}$ gives $A=\frac{9d^5}{2a}\leq \frac{27d^4}{2}$. And from (\ref{eq128}), we have $d\leq \frac{32}{105}$. Therefore,
	\begin{equation}\label{eq130}
	\begin{array}{l}
	0.09499< d/3\leq a\leq d/2< 0.152385\vspace{2mm}\\
	0.089059< 5d/16 \leq A \leq 27d^4/2< 0.116473\vspace{2mm}\\
	0.046012< \frac{a+A}{4} < 0.0672145.
	\end{array}
	\end{equation}

	\noindent Using Lemma \ref{lem15} we need to consider $(a_5'',c_5)=(0,0), (0,1/2), {\rm ~and~} (1/18,1/4).$\vspace{2mm}
	
	\noindent When $(a_5'',c_5)=(0,0)$, the following table gives a solution to (\ref{eq101'}) i.e., of\vspace{-2mm}
	\begin{equation}\label{EQN49}0<G=(x_1+a_2''x_2)x_2-\frac{1}{9}x_5^2-\frac{a+A}{4}
	<d+\frac{3a}{4}+2A.\end{equation}
	
	{\scriptsize
		\begin{equation*}\begin{array}{lllll}
		\hline\\
		{\rm ~~~~Range} & x_1 & x_1x_2 & x_5& {\rm ~~G}\vspace{2mm}\\
		\hline\\
		
		\frac{4}{9}+\frac{a+A}{4}<f_1''\leq 1&c_1&|c_1|& \pm 2&f_1''-\frac{4}{9}-\frac{a+A}{4}\vspace{2mm}\\
		
		\frac{a+A}{4}<f_1''\leq \frac{4}{9}+\frac{a+A}{4}&c_1&|c_1|& 0&f_1''-\frac{a+A}{4}\vspace{2mm}\\	
		
		-\frac{5}{9}+\frac{a+A}{4}<f_1''<-\frac{5}{9}+d+a+\frac{9A}{4}&\pm 1+c_1&1+|c_1|&\pm 2&f_1''+\frac{5}{9}-\frac{a+A}{4}\vspace{2mm}\\
		
		-\frac{1}{2}<f_1''\leq -\frac{5}{9}+\frac{a+A}{4}&\pm 1+c_1&1+|c_1|&0&f_1''+1-\frac{a+A}{4}\vspace{2mm}\\

		{\rm When~}  -\frac{5}{9}+d+a+\frac{9A}{4}\leq f_1''\leq \frac{a+A}{4}: \vspace{2mm}\\

		-\frac{5}{9}+d+a+\frac{9A}{4}\leq g_1''\leq \frac{a+A}{4}&c_1&2|c_1|&0&f_2''-\frac{a+A}{4}\vspace{2mm}\\
		
		-\frac{5}{9}+\frac{a+A}{4}<g_1''<-\frac{5}{9}+d+a+\frac{9A}{4}&\pm 1+c_1&1-|c_1|&\pm 2&g_1''+\frac{5}{9}-\frac{a+A}{4}\vspace{2mm}\\
		
		-1+\frac{a+A}{4}<g_1''\leq -\frac{5}{9}+\frac{a+A}{4}&\pm 1+c_1&1-|c_1|&0&g_1''+1-\frac{a+A}{4}\vspace{2mm}\\
		
		-1<g_1''\leq -1+\frac{a+A}{4}&\pm 1+c_1&2+2|c_1|&\pm 3&f_2''+1-\frac{a+A}{4}\vspace{2mm}\\	 
		\hline
		\end{array}
		\end{equation*}\vspace{1mm}}
	
	\noindent When $(a_5'',c_5)=(0,1/2)$, the following table gives a solution to (\ref{EQN49}).\vspace{1mm}
	
	{\scriptsize
		\begin{equation*}\begin{array}{lllll}
		\hline\\
		{\rm ~~~~Range} & x_1 & x_1x_2 & x_5& {\rm ~~G}\vspace{2mm}\\
		\hline\\
		
		\frac{25}{36}+\frac{a+A}{4}<f_1''\leq 1&c_1&|c_1|&\frac{5}{2} &f_1''-\frac{25}{36}-\frac{a+A}{4}\vspace{2mm}\\
		
		\frac{1}{4}+\frac{a+A}{4}<f_1''\leq \frac{25}{36}+\frac{a+A}{4}&c_1&|c_1|& \frac{3}{2}&f_1''-\frac{1}{4}-\frac{a+A}{4}\vspace{2mm}\\	
		
		\frac{1}{36}+\frac{a+A}{4}<f_1''\leq \frac{1}{4}+\frac{a+A}{4}&c_1&|c_1|&\frac{1}{2}&f_1''-\frac{1}{36}-\frac{a+A}{4}\vspace{2mm}\\
		
		-\frac{11}{36}+\frac{a+A}{4}<f_1''\leq \frac{1}{36}+\frac{a+A}{4}&\pm 1+c_1&1+|c_1|&\frac{5}{2}&f_1''+\frac{11}{36}-\frac{a+A}{4}\vspace{2mm}\\
		
		-\frac{1}{2}<f_1''\leq -\frac{11}{36}+\frac{a+A}{4}&\pm 1+c_1&1+|c_1|&\frac{3}{2}&f_1''+\frac{3}{4}-\frac{a+A}{4}\vspace{2mm}\\
		
		\hline
		\end{array}
		\end{equation*}\vspace{1mm}}
	
	\noindent When $(a_5'',c_5)=(1/18,1/4)$, the following table gives a solution to (\ref{eq101'}) i.e., of
	\begin{equation}\label{EQN51}0<G=(x_1+a_2''x_2+\frac{1}{18}x_5)x_2-\frac{1}{9}x_5^2
	-\frac{a+A}{4}<d+\frac{3a}{4}+2A.\end{equation}
	
	{\scriptsize
		\begin{equation*}\begin{array}{lllll}
		\hline\\
		{\rm ~~~~Range} & x_1 & x_2 & x_5& {\rm ~~G}\vspace{2mm}\\
		\hline\\
		
		\frac{7}{16}+\frac{a+A}{4}<p_1''\leq 1&c_1&1&\frac{9}{4} &p_1''-\frac{7}{16}-\frac{a+A}{4}\vspace{2mm}\\
		
		-\frac{1}{144}+\frac{a+A}{4}<p_1''\leq \frac{7}{16}+\frac{a+A}{4}&c_1&1& \frac{1}{4}&p_1''+\frac{1}{144}-\frac{a+A}{4}\vspace{2mm}\\	
		
		-\frac{11}{48}+\frac{a+A}{4}<p_1''\leq -\frac{1}{144}+\frac{a+A}{4}&2+c_1&1&\frac{17}{4}&p_1''+\frac{11}{48}-\frac{a+A}{4}\vspace{2mm}\\
		
		-\frac{9}{16}+\frac{a+A}{4}<p_1''\leq -\frac{11}{48}+\frac{a+A}{4}& 1+c_1&1&\frac{9}{4}&p_1''+\frac{9}{16}-\frac{a+A}{4}\vspace{2mm}\\
		
		-\frac{145}{144}+\frac{a+A}{4}<p_1''\leq -\frac{9}{16}+\frac{a+A}{4}& 1+c_1&1&\frac{1}{4}&p_1''+\frac{145}{144}-\frac{a+A}{4}\vspace{2mm}\\
		
		-1<p_1''\leq -\frac{145}{144}+\frac{a+A}{4}& 3+c_1&1&\frac{17}{4}&p_1''+\frac{59}{48}-\frac{a+A}{4}\vspace{2mm}\\
		
		\hline
		\end{array}
		\end{equation*}\vspace{1mm}}
	
	\noindent For each of the remaining values of $t$ and corresponding $(a_5'',c_5)$ on taking $x_2=\pm 1$ and giving suitable values to  $x_5$, we find that the range of $f_1''$ or, $p_1''$ is covered and hence (\ref{eq101'}) is soluble.\vspace{4mm}
	
	\noindent \textbf{Case IV:} $\textbf{m=2, K=3, L=3}.$\vspace{2mm}
	
	\noindent Here $3<\frac{\delta_{m,K}}{t}=\frac{d+3a/4+2A}{t}\leq 4$, $\delta_{2,3}=d+\frac{3a}{4}+2A\geq d+\frac{d}{4}+\frac{5d}{8}=\frac{15d}{8}$. From (\ref{eq84}) we have
	
	\begin{equation}\label{eq98}
	0.3125 = \sqrt{25/256} \leq d \leq 16/45=0.35556.
	\end{equation}
	
	\noindent Further, since $a\geq d/3$ and $A\geq 5d/16$, we get $ t =\frac{d^5}{2aA}\leq \frac{24d^3}{5}$. Also since $\delta_{2,3,3}=\delta_{2,3}+\frac{L^2-1}{4}t\leq 1$, we get $t\leq \frac{1}{2}(1-\delta_{2,3}) \leq \frac{1}{2}(1-\frac{15d}{8})$. Therefore,
	
	\begin{equation*}\label{eq891'}
	\begin{array}{ll}
	0.1464< \frac{15d}{32}\leq \frac{\delta_{2,3}}{4}\leq t\leq &\min\{\frac{24d^3}{5},\frac{1}{2}(1-\frac{15d}{8})\}
	=f_0(d) {\rm ~(say)~}\vspace{2mm}\\
	&=\left\{
	\begin{array}{ll}
	\frac{24d^3}{5} & {\rm ~if~} d\leq \alpha\\
	\frac{1}{2}(1-\frac{15d}{8}) &{\rm ~if~} d\geq \alpha
	\end{array}
	\right. ~\leq 0.1845=\lambda_0 {\rm (say)},\vspace{2mm}\\
	\end{array}\vspace{-2mm}
	\end{equation*}
	
	\noindent  where $\alpha $ is a positive real root of $\frac{24}{5}d^3+\frac{15}{16}d-\frac{1}{2}=0 $ satisfying $0.337<\alpha<0.3375$.\vspace{2mm}

	\noindent We divide the range of $t$ into $7$ subintervals $[\lambda_n,\lambda_{n-1}], n=1,2,\cdots 7$ and in each subinterval choose suitable integers $2h_n$ and $k_n$ such that
	$\vline~ h_n-tk_n^2~\vline +\frac{1}{2} <\frac{15}{8}d\leq d+\frac{3a}{4}+2A.$\vspace{1mm}

	{\footnotesize
		\begin{equation*}\begin{array}{lllllllll}
		\hline\\
		n & (h_n,k_n) & \lambda_n& \mbox{Remarks}&~~ &n & (h_n,k_n) & \lambda_n& \mbox{Remarks}\vspace{1mm}\\
		\hline\\
		
		1&(4.5,5)&0.1757&(6.5,6) &&2&(1.5,3)&0.1558& \mbox{tbd} \vspace{2mm}\\
		
		3&(2.5,4)&0.1506&\mbox{na} &&4&(15,10)&0.1498&\mbox{tbd} \vspace{2mm}\\
		
		5&(9.5,8)&0.1471&(12,9) &&6&(42.5,17)&0.147&\mbox{tbd} \vspace{2mm}\\	
		
		7&(33,15)&0.1464&(37.5,16) &&&&&\vspace{2mm}\\
		
		\hline
		\end{array}
		\end{equation*}\vspace{1mm}}
	
	\noindent The choice of $(h_n,k_n)$ is done in a manner similar to that in Lemma \ref{lem18} ( $m=2,K=3, L=4$ ).\vspace{2mm}
	
	\noindent We find that  (\ref{eq101'}) is soluble unless  $t=1/6,$ $3/20, 5/34$ and $\beta_{t,G} = \pm a_5''-2tc_5\equiv h_n/k_n \pmod {1/k_n, 2t}$. For $t=1/6$, we  give here solutions of (\ref{eq101'}) depending upon values of $a_2''$ and $c_1$. The value $t=3/20$ and $5/34$ can be similarly disposed off by taking $x_2=\pm 1$.\vspace{2mm}

	\noindent When $t=\frac{1}{6}$, we have $\frac{1}{6}=t\leq \frac{24d^3}{5}$. This gives $d^3\geq \frac{5}{144}$, i.e., $d\geq (\frac{5}{144})^{1/3}>0.32623$. Also, $\frac{1}{6}=t=\frac{d^5}{2aA}$ gives $A=\frac{6d^5}{2a}\leq 9d^4$. And from (\ref{eq98}), we have $d\leq \frac{16}{45}$. Therefore,
	\begin{equation}\label{eq100}
	\begin{array}{l}
	0.10874< d/3\leq a\leq d/2< 0.17778\vspace{2mm}\\
	0.10194< 5d/16 \leq A \leq 9d^4< 0.14385\vspace{2mm}\\
	0.05267 < \frac{a+A}{4} < 0.08041.
	\end{array}
	\end{equation}

	\noindent Using Lemma \ref{lem15} we need to consider $(a_5'',c_5)=(0,1/2) {\rm ~and~} (1/6,0).$ If $a_5''=1/6,$ considering the transformation $x_5\rightarrow x_5-x_2$, we may assume that $a_5''=1/2$.\vspace{2mm}
	
	\noindent When $(a_5'',c_5)=(0,\frac{1}{2})$, the following table gives a solution to (\ref{eq101'}), i.e., of
	 \begin{equation}\label{EQN40}0<G=(x_1+a_2''x_2)x_2-\frac{1}{6}x_5^2-\frac{a+A}{4}<d+\frac{3a}{4}+2A.\end{equation}
	
	{\scriptsize
		\begin{equation*}\begin{array}{lllll}
		\hline\\
		{\rm ~~~~Range} & x_1 & x_1x_2 & x_5& {\rm ~~G}\vspace{2mm}\\
		\hline\\
		
		\frac{3}{8}+\frac{a+A}{4}<f_1''\leq 1&c_1&|c_1|& \frac{3}{2}&f_1''-\frac{3}{8}-\frac{a+A}{4}\vspace{2mm}\\
		
		\frac{1}{24}+\frac{a+A}{4}<f_1''\leq \frac{3}{8}+\frac{a+A}{4}&c_1&|c_1|& \frac{1}{2}&f_1''-\frac{1}{24}-\frac{a+A}{4}\vspace{2mm}\\	
		
		-\frac{1}{2}\leq f_1''<-\frac{5}{8}+d+a+\frac{9A}{4}&\pm 1+c_1&1+|c_1|&\frac{3}{2}&f_1''+\frac{5}{8}-\frac{a+A}{4}\vspace{2mm}\\

		{\rm When~}  -\frac{5}{8}+d+a+\frac{9A}{4}\leq f_1''\leq \frac{1}{24}+\frac{a+A}{4}: \vspace{2mm}\\

		-\frac{5}{8}+d+a+\frac{9A}{4}\leq g_1''\leq \frac{1}{24}+\frac{a+A}{4}&c_1&2|c_1|&\frac{1}{2}&f_2''-\frac{1}{24}-\frac{a+A}{4}\vspace{2mm}\\
		
		-\frac{5}{8}+\frac{a+A}{4}<g_1''<-\frac{5}{8}+d+a+\frac{9A}{4}&\pm 1+c_1&1-|c_1|&\frac{3}{2}&g_1''+\frac{5}{8}-\frac{a+A}{4}\vspace{2mm}\\
		
		-\frac{23}{24}+\frac{a+A}{4}<g_1''\leq -\frac{5}{8}+\frac{a+A}{4}&\pm 1+c_1&1-|c_1|&\frac{1}{2}&g_1''+\frac{23}{24}-\frac{a+A}{4}\vspace{2mm}\\
		
		-1<g_1''\leq -\frac{23}{24}+\frac{a+A}{4}&\pm 1+c_1&2+2|c_1|&\frac{5}{2}&f_2''+\frac{23}{24}-\frac{a+A}{4}\vspace{2mm}\\

		\hline
		\end{array}
		\end{equation*}\vspace{1mm}}
	
	\noindent When $(a_5'',c_5)=(\frac{1}{2},0)$, the following table gives a solution to (\ref{eq101'}) i.e., of
	\begin{equation}\label{EQN41}0<G=(x_1+a_2''x_2+\frac{1}{2}x_5)x_2-\frac{1}{6}x_5^2-
	\frac{a+A}{4}<d+\frac{3a}{4}+2A.\end{equation}
	
	{\scriptsize
		\begin{equation*}\begin{array}{llllll}
		\hline\\
		{\rm ~~~~Range} & x_1 & x_1x_2 & x_5&x_2x_5& {\rm ~~G}\vspace{2mm}\\
		\hline\\
		
		\frac{a+A}{4}<g_1''\leq \frac{1}{2}&c_1&-|c_1|&0&0&g_1''-\frac{a+A}{4}\vspace{2mm}\\	
		
		-\frac{1}{3}+\frac{a+A}{4}<g_1''\leq \frac{a+A}{4}&\pm 1+c_1&1-|c_1|&\pm 1&-1&g_1''+\frac{1}{3}-\frac{a+A}{4}\vspace{2mm}\\	
		
		-1+\frac{a+A}{4}<g_1''<-1+d+a+\frac{9A}{4}&\pm 1+c_1&1-|c_1|&0&0&g_1''+1-\frac{a+A}{4}\vspace{2mm}\\
		
		-1<g_1''\leq -1+\frac{a+A}{4}&\pm 2+c_1&2-|c_1|&\pm 1&-1&g_1''+\frac{4}{3}-\frac{a+A}{4}\vspace{2mm}\\

		{\rm When~}  -1+d+a+\frac{9A}{4}\leq g_1''\leq -\frac{1}{3}+\frac{a+A}{4}: \vspace{2mm}\\

		\left.\begin{array}{l}d+a+\frac{9A}{4}\leq f_1''\leq 1\vspace{1mm}\\f_2''< \frac{7}{6}+d+a+\frac{9A}{4}\end{array}\right\}&c_1&2|c_1|&\pm 1&-2&f_2''-\frac{7}{6}-\frac{a+A}{4}\vspace{2mm}\\
		
		\left.\begin{array}{l}d+a+\frac{9A}{4}\leq f_1''\leq 1\vspace{1mm}\\f_2''\geq \frac{7}{6}+d+a+\frac{9A}{4}\end{array}\right\}& \pm 1+c_1&3+3|c_1|&\pm 3&-9&f_3''-3-\frac{a+A}{4}\vspace{2mm}\\

		\frac{a+A}{4}<f_1''<d+a+\frac{9A}{4}&c_1&|c_1|&0&0&f_1''-\frac{a+A}{4}\vspace{2mm}\\
		
		-\frac{1}{3}+\frac{a+A}{4}<f_1''<-\frac{1}{3}+d+a+\frac{9A}{4}&\pm 1+c_1&1+|c_1|&\pm 1&-1&f_1''+\frac{1}{3}-\frac{a+A}{4}\vspace{2mm}\\
		
		-\frac{1}{2}<f_1''\leq -\frac{1}{3}+\frac{a+A}{4}&\pm 3+c_1&6+2|c_1|&\pm 3&-6&f_2''+\frac{3}{2}-\frac{a+A}{4}\vspace{2mm}\\

		\hline
		\end{array}
		\end{equation*}\vspace{1mm}}
	
	\noindent \textbf{Case V:} $\textbf{m=2, K=3, L=2}.$\vspace{2mm}
	
	\noindent Here $2<\frac{\delta_{m,K}}{t}=\frac{d+3a/4+2A}{t}\leq 3$, $\delta_{2,3}=d+\frac{3a}{4}+2A\geq d+\frac{d}{4}+\frac{5d}{8}=\frac{15d}{8}$. From (\ref{eq84}) we have
	
	\begin{equation}\label{eq94}
	0.3608< \sqrt{25/192} \leq d \leq 32/75<0.42667.
	\end{equation}
	
	\noindent Further, since $a\geq d/3$ and $A\geq 5d/16$ we get, $ t =\frac{d^5}{2aA}\leq \frac{24d^3}{5}.$\vspace{2mm}
	
	\noindent Also since $\delta_{2,3,2}=\delta_{2,3}+\frac{L^2-1}{4}t\leq 1$, we get $t\leq \frac{4(1-\delta_{2,3})}{3} \leq \frac{4(1-\frac{15d}{8})}{3}$. Therefore,
	
	\begin{equation*}\label{eq890}
	\begin{array}{ll}
	0.2255< \frac{5d}{8}\leq \frac{\delta_{2,3}}{3}\leq t\leq &\min\{\frac{24d^3}{5},\frac{4(1-\frac{15d}{8})}{3}\}
	=f_0(d) {\rm ~(say)~}\vspace{2mm}\\
	&=\left\{
	\begin{array}{ll}
	\frac{24d^3}{5} & {\rm ~if~} d\leq \alpha\\
	\frac{4(1-\frac{15d}{8})}{3} &{\rm ~if~} d\geq \alpha
	\end{array}
	\right.\vspace{2mm}\\
	&\leq 0.3212
	\end{array}
	\end{equation*}
	
	\noindent  where $\alpha $ is a root of $\frac{24}{5}d^3+\frac{5}{2}d-\frac{4}{3}=0 $ satisfying $0.405<\alpha<0.406$. Take $h_t=8$, $k_t=5$ if $\frac{1}{25}(8.5-\frac{15}{8}d)< t \leq f_0(d)$ and $h_t=1$,$k_t=2$ when $
	\frac{\delta_{2,3}}{3}\leq t\leq \min\{\frac{24d^3}{5},\frac{1}{25}(8.5-\frac{15}{8}d)\}$ to see that $\vline~ h_t-tk_t^2~\vline +\frac{1}{2} <d+\frac{3a}{4}+2A=\delta_{2,3}$.
	
	When $t=\frac{h}{k^2}=\frac{8}{25}$, we find that $|3-\frac{8}{25} 9|=0.12<\delta_{2,3}-\frac{1}{2}.$\vspace{2mm}

	\noindent Therefore,  (\ref{eq101'})  is soluble unless $t=1/4$ and $\beta_{t,G}=\pm a_5''-c_5\equiv 1/2 \pmod {1/2}$. This gives $(a_5'',c_5)=(0,0)$ or $(1/4,1/2)$.\vspace{2mm}
	
	\noindent Now, $1/4=t=d^5/2aA\leq 24d^3/5$ gives $d\geq (5/96)^{1/3}>0.3734$. Also, $1/4=t=d^5/2aA$ gives $A=2d^5/a\leq 6d^4$. Hence,
	
	\begin{equation}\label{eq96}
	\begin{array}{l}
	0.12446< d/3\leq a\leq d/2< 0.21334\vspace{2mm}\\
	0.1166< 5d/16 \leq A \leq 6d^4< 0.19885\vspace{2mm}\\
	0.06026 < \frac{a+A}{4} < 0.10305.
	\end{array}
	\end{equation}
	
	\noindent Thus we need to find a solution of
	\begin{equation}\label{eq97}
	0<G=(x_1+a_2''x_2+a_5''x_5)x_2-x_5^2/4-(a+A)/4<d+3a/4+2A.
	\end{equation}
	for $(a_5'',c_5)=(0,0)$ and $(1/4,1/2)$.\vspace{4mm}
	
	\noindent When  $(a_5'',c_5)=(0,0)$, the following table gives a solution of (\ref{eq97}).\vspace{2mm}
	
	{\scriptsize
		\begin{equation*}\begin{array}{lllll}
		\hline\\
		{\rm ~~~~Range} & x_1 & x_1x_2 & x_5 & {\rm ~~G}\vspace{2mm}\\
		\hline\\
		
		\frac{1}{4}+\frac{a+A}{4}<f_1''\leq 1&c_1&|c_1|&\pm 1&f_1''-\frac{1}{4}-\frac{a+A}{4}\vspace{2mm}\\
		
		\frac{a+A}{4}<f_1''\leq \frac{1}{4}+\frac{a+A}{4}&c_1&|c_1|&0&f_1''-\frac{a+A}{4}\vspace{2mm}\\
		
		-\frac{1}{2}<f_1''<-\frac{3}{4}+d+a+\frac{9A}{4}&\pm 1+c_1&1+|c_1|&\pm 1&f_1''+\frac{3}{4}-\frac{a+A}{4}\vspace{2mm}\\

		{\rm When~}  -\frac{3}{4}+d+a+\frac{9A}{4}\leq f_1'\leq \frac{a+A}{4}: \vspace{2mm}\\
		
		-\frac{3}{4}+d+a+\frac{9A}{4}\leq g_1''\leq \frac{a+A}{4}&c_1&3|c_1|&0&f_3''-\frac{a+A}{4}\vspace{2mm}\\
		
		-\frac{3}{4}+\frac{a+A}{4}<g_1''<-\frac{3}{4}+d+a+\frac{9A}{4}&\pm 1+c_1&1-|c_1|&\pm 1&g_1''+\frac{3}{4}-\frac{a+A}{4}\vspace{2mm}\\
		
		-1+\frac{a+A}{4}<g_1''\leq -\frac{3}{4}+\frac{a+A}{4}&\pm 1+c_1&1-|c_1|&0&g_1''+1-\frac{a+A}{4}\vspace{2mm}\\
		
		-1<g_1'' \leq -1+\frac{a+A}{4}&\pm 1+c_1&3+3|c_1|&0&f_3''+3-\frac{a+A}{4}\vspace{2mm}\\
		
		\hline
		\end{array}
		\end{equation*}\vspace{2mm}}
	
	\noindent When  $(a_5'',c_5)=(1/4,1/2)$, the following table gives a solution of (\ref{eq97}).
	
	{\scriptsize
		\begin{equation*}\begin{array}{llllll}
		\hline\\
		{\rm ~~~~Range} & x_1 & x_1x_2 & x_5&x_2x_5 & {\rm ~~G}\vspace{2mm}\\
		\hline\\
		
		-\frac{1}{16}+\frac{a+A}{4}<g_1''\leq \frac{1}{2}&c_1&-|c_1|&\pm \frac{1}{2}&\frac{1}{2}&g_1''+\frac{1}{16}-\frac{a+A}{4}\vspace{2mm}\\
		
		-\frac{13}{16}+\frac{a+A}{4}<g_1''<-\frac{13}{16}+d+a+\frac{9A}{4}&\pm 1+c_1&1-|c_1|&\pm \frac{1}{2}&-\frac{1}{2}&g_1''+\frac{13}{16}-\frac{a+A}{4}\vspace{2mm}\\
		
		-1<g_1''\leq -\frac{13}{16}+\frac{a+A}{4}&\pm 1+c_1&1-|c_1|&\pm \frac{1}{2}&\frac{1}{2}&g_1''+\frac{17}{16}-\frac{a+A}{4}\vspace{2mm}\\
		
		{\rm When~}  -\frac{13}{16}+d+a+\frac{9A}{4}\leq g_1'\leq -\frac{1}{16}+\frac{a+A}{4}: \vspace{2mm}\\

		\frac{3}{16}+d+a+\frac{9A}{4}\leq f_1''\leq 1&\pm 1+c_1&-2+2|c_1|&\pm \frac{1}{2}&-1&f_2''-\frac{37}{16}-\frac{a+A}{4}\vspace{2mm}\\

		\frac{3}{16}+\frac{a+A}{4}<f_1''<\frac{3}{16}+d+a+\frac{9A}{4}&c_1&|c_1|&\pm \frac{1}{2}&-\frac{1}{2}&f_1''-\frac{3}{16}-\frac{a+A}{4}\vspace{2mm}\\
		
		-\frac{1}{16}+\frac{a+A}{4}<f_1''\leq \frac{3}{16}+\frac{a+A}{4}&c_1&|c_1|&\pm \frac{1}{2}&\frac{1}{2}&f_1''+\frac{1}{16}-\frac{a+A}{4}\vspace{2mm}\\
		
		-\frac{1}{2}<f_1'' \leq -\frac{1}{16}+\frac{a+A}{4}&\pm 1+c_1&2+2|c_1|&\pm \frac{3}{2}&-3&f_2''+\frac{11}{16}-\frac{a+A}{4}\vspace{2mm}\\

		\hline
		\end{array}
		\end{equation*}\vspace{1mm}}
	
	\noindent \textbf{Case VI:} $\textbf{m=2, K=3, L=1}.$\vspace{2mm}
	
	\noindent Here $1<\frac{\delta_{m,K}}{t}=\frac{d+3a/4+2A}{t}\leq 2$, $\delta_{m,K}=d+\frac{3a}{4}+2A\geq d+\frac{d}{4}+\frac{5d}{8}=\frac{15d}{8}$. From (\ref{eq84}) we have
	
	\begin{equation}\label{eq85}
	0.4419 \sqrt{25/128} \leq d \leq 8/15<0.53334.
	\end{equation}
	
	\noindent Further, since $a\geq d/3$ and $A\geq 5d/16$ we get, \vspace{-2mm}
	$$0.4142< \frac{15d}{16}\leq \frac{d+3a/4+2A}{2}\leq t =\frac{d^5}{2aA}\leq \frac{24d^3}{5}< 0.7282.$$
	
	\noindent Taking $h_t=1/2$ and $k_t=1$ and using (\ref{eq85})
	we find that $\vline~h_t-tk_t^2~\vline + \frac{1}{2} <\delta_{m,K}$, as $ t < \delta_{m,K}$ and   $t+\delta_{m,K} \geq \frac{15}{16}d+\frac{15}{8}d=\frac{45}{16}d>1$. Therefore,  (\ref{eq101'}) and hence (\ref{eq17}) is soluble unless $t=1/2$ and $\beta_{t,G}=\pm a_5''-c_5\equiv 1/2 \pmod 1$. This gives $(a_5'',c_5)=(0,1/2)$ or $(1/2,0)$.\vspace{2mm}

	\noindent Now, $1/2=t=d^5/2aA\leq 24d^3/5$ gives $d\geq (5/48)^{1/3}>0.4705$. Also, $1/2=t=d^5/2aA$ gives $A=d^5/a\leq 3d^4$. Hence,
	\begin{equation}\label{eq86}
	\begin{array}{l}
	0.1568< d/3\leq a\leq d/2< 0.26667\vspace{2mm}\\
	0.147< 5d/16 \leq A \leq 3d^4< 0.24274\vspace{2mm}\\
	0.07595< \frac{a+A}{4} < 0.12736.
	\end{array}
	\end{equation}
	
	\noindent Thus we need to find a solution of
	\begin{equation}\label{eq87}
	0<G=(x_1+a_2''x_2+a_5''x_5)x_2-x_5^2/2-(a+A)/4<d+3a/4+2A \vspace{-2mm}
	\end{equation}
	for $(a_5'',c_5)=(0,1/2)$ and $(1/2,0)$.\vspace{4mm}
	
	\noindent When  $(a_5'',c_5)=(0,1/2)$, the following table gives a solution of (\ref{eq87}).\vspace{2mm}
	
	{\scriptsize
		\begin{equation*}\begin{array}{lllll}
		\hline\\
		{\rm ~~~~Range} & x_1 & x_1x_2 & x_5 & {\rm ~~G}\vspace{2mm}\\
		\hline\\
		
		\frac{1}{8}+\frac{a+A}{4}<f_1''\leq 1&c_1&|c_1|&\frac{1}{2}&f_1''-\frac{1}{8}-\frac{a+A}{4}\vspace{2mm}\\
		
		-\frac{1}{2}<f_1''<-\frac{7}{8}+d+a+\frac{9A}{4}&\pm 1+c_1&1+|c_1|&\frac{1}{2}&f_1''+\frac{7}{8}-\frac{a+A}{4}\vspace{2mm}\\

		{\rm When~}  -\frac{7}{8}+d+a+\frac{9A}{4}\leq f_1''\leq \frac{1}{8}+\frac{a+A}{4}: \vspace{2mm}\\

		-\frac{7}{8}+d+a+\frac{9A}{4}\leq g_1''\leq \frac{1}{8}+\frac{a+A}{4}&c_1&2|c_1|&\frac{1}{2}&f_2''-\frac{1}{8}-\frac{a+A}{4}\vspace{2mm}\\
		
		-\frac{7}{8}+\frac{a+A}{4}<g_1''<-\frac{7}{8}+d+a+\frac{9A}{4}&\pm 1+c_1&1-|c_1|&\frac{1}{2}&g_1''+\frac{7}{8}-\frac{a+A}{4}\vspace{2mm}\\
		
		-1<g_1'' \leq -\frac{7}{8}+\frac{a+A}{4}&\pm 1+c_1&2+2|c_1|&\frac{3}{2}&f_2''+\frac{7}{8}-\frac{a+A}{4}\vspace{2mm}\\

		\hline
		\end{array}
		\end{equation*}\vspace{1mm}}

	\noindent When  $(a_5'',c_5)=(1/2,0)$, equation (\ref{eq87}) has no solution for $c_1=\frac{1}{2}$ and $a_2''=\frac{1}{2}$. This is so because for $x_1=x+\frac{1}{2}$,
	\begin{equation*} \frac{(a+A)}{2}< (2x+1+x_2)x_2 +x_5(x_2-x_5)<2\big(\delta_{2,3}+\frac{(a+A)}{4}\big)\end{equation*}
	does not have a solution for some values of $a,~A$ and $d$. For $a=\frac{d}{3}$, $A=\frac{5d}{16}$ and $d<0.49$, we note that $2\big(\delta_{2,3}+\frac{(a+A)}{4}\big)<2$ and $(2x+1+x_2)x_2 +x_5(x_2-x_5)$ takes integral values only and it can not be equal to $1$ for all integers $x,x_2,x_5$. \vspace{2mm}
	
	\noindent Therefore, we need to go to $F$. Here since $a\geq d/3, A\geq 5d/16$ and $t=1/2$ we have $5d/16\leq A=d^5/a\leq 3d^4$ which gives \vspace{-2mm}
	\begin{equation}\label{eq88}
	d\geq (5/48)^{1/3}>0.4705. \vspace{-2mm}
	\end{equation}
	
	\noindent Also since $\delta_{2,3}\leq 1$, we get $A\leq \left(1-(d+3a/4)\right)/2 \leq (1-5d/4)/2$. Therefore, \vspace{-2mm}
	\begin{equation}\label{eq89}
	\begin{array}{ll}
	0.147< 5d/16 \leq A\leq &min\{3d^4,(1-5d/4)/2\}=f_0(d) {\rm ~(say)~}\vspace{2mm}\\
	&=\left\{
	\begin{array}{ll}
	3d^4 & {\rm ~if~} d\leq 1/2\\
	(1-5d/4)/2 &{\rm ~if~} d\geq 1/2
	\end{array}
	\right.\vspace{2mm}\\
	&\leq \frac{3}{16}=\lambda_{0}(\rm say).
	\end{array}\vspace{-2mm}
	\end{equation}
	
	\noindent We divide the range of $A$ into $7$ subintervals $[\lambda_n,\lambda_{n-1}], n=1,2,\cdots 7$ and in each subinterval choose suitable integers $2h_n$ and $k_n$ such that
	$\vline~ h_n-Ak_n^2~\vline +\frac{1}{2} <\frac{5}{4}d \leq d+\frac{3a}{4}$.\vspace{1mm}
	
	{\scriptsize
		\begin{equation*}\begin{array}{lllllllll}
		\hline\\
		n & (h_n,k_n) & \lambda_n& \mbox{Remarks}&~~ &n & (h_n,k_n) & \lambda_n&\mbox{Remarks}\vspace{2mm}\\
		\hline\\
		
		1&(3,4)&0.1801& \mbox{tbd} &&2&(4.5,5)&0.1757&(6.5,6) \vspace{2mm}\\
		
		3&(1.5,3)&0.1566& \mbox{tbd} &&4&(2.5,4)&0.15079&(1.5,3) \vspace{2mm}\\
		
		5&(15,10)&0.1491& \mbox{tbd} &&6&(12,9)&0.14705& (9.5,8)\vspace{2mm}\\	
		
		7&(42.5,17)&0.1467& \mbox{tbd} &&&&&\vspace{2mm}\\
		
		\hline
		\end{array}
		\end{equation*}\vspace{1mm}}
	
	\noindent The choice of $(h_n,k_n)$ is done in a manner similar to that in the case  $m=2,K=3, L=4$  (now for $A$ instead of $t$). We find that  (\ref{eq101}) is soluble unless  $A=3/16, 1/6,$ $3/20, 5/34$ and $\beta_{A,F} = \pm a'_4-2A(c_4+\lambda x_5)\equiv h_n/k_n \pmod {1/k_n, 2A}$. For each of these values to $A$ we get solutions of (\ref{eq101}) depending upon values of $a_2'$ and $c_1$:  \vspace{2mm}
	
	\noindent Note that we have $A\leq \frac{3}{16}$. Also $F=(x_1+\cdots)x_2-A(x_4+\lambda x_5)^2-\frac{1}{2}x_5^2-\frac{a}{4}=\pm x+\beta y-(A\lambda^2+\frac{1}{2})y^2+\nu,$ and $ d+\frac{3}{4}a \geq d+\frac{3}{4}.\frac{d}{3}=\frac{5d}{4}$. Therefore, on taking $h=\frac{1}{2}, k=1$ we get $|\frac{1}{2}-(A\lambda^2+\frac{1}{2})|+\frac{1}{2}=A\lambda^2+\frac{1}{2}\leq \frac{1}{2}+\frac{3}{64}<\frac{5d}{4}\leq d+\frac{3a}{4}=\delta_{2}.$ Thus (\ref{eq101}) is soluble unless $A\lambda^2 +\frac{1}{2}=\frac{1}{2} $ which gives $\lambda=0$. Then from (\ref{eq16}), we get $a_5''\equiv a_5' \pmod 1$ which gives $a_5'=\frac{1}{2}$.\vspace{2mm}
	
	\noindent We discuss here  $A=1/6$ only.  The other three $A=3/16,3/20, 5/34$ are settled in the similar way by just taking $x_2=\pm 1$.\vspace{2mm}

	\noindent Let $A=\frac{1}{6}$. Here, since $\frac{d}{a}\leq 3$, and $\frac{1}{2}=t=\frac{d^5}{2aA}=\frac{6d^5}{2a}$, we get $d^4\geq \frac{1}{18}$ i.e., $d\geq (\frac{1}{18})^{\frac{1}{4}}>0.48549$. Also from (\ref{eq85}), we have $d\leq \frac{8}{15}$. This, together with $\frac{d}{3} \leq a \leq \frac{4A}{3}$ gives
	\begin{equation}\label{EQN23}
	0.16183< a < 0.22223.
	\end{equation}
	
	\noindent Using Lemma \ref{lem15} we need to consider $(a_4',c_4)=(0,1/2) {\rm ~and~} (1/6,0).$ If $a_4'=1/6,$ considering the transformation $x_4\rightarrow x_4-x_2$, we may assume that $a_4'=1/2$.\vspace{4mm}
	
	\noindent When $(a_4',c_4)=(0,\frac{1}{2})$, taking $x_5=0$, the following table gives a solution to (\ref{eq101}) i.e. of
	\begin{equation}\label{EQN26}0<F=(x_1+a_2'x_2+\frac{1}{2}x_5)x_2-\frac{1}{6}x_4^2-\frac{1}{2}
	x_5^2-\frac{a}{4}<d+\frac{3a}{4}.\end{equation}
	
	{\scriptsize
		\begin{equation*}\begin{array}{lllll}
		\hline\\
		{\rm ~~~~Range} & x_1 & x_1x_2 & x_4& {\rm ~~F}\vspace{2mm}\\
		\hline\\
		
		\frac{3}{8}+\frac{a}{4}<f_1'\leq 1&c_1&|c_1|&\pm \frac{3}{2}&f_1'-\frac{3}{8}-\frac{a}{4}\vspace{2mm}\\
		
		\frac{1}{24}+\frac{a}{4}<f_1'\leq \frac{3}{8}+\frac{a}{4}&c_1&|c_1|&\pm \frac{1}{2}&f_1'-\frac{1}{24}-\frac{a}{4}\vspace{2mm}\\	
		
		-\frac{5}{8}+\frac{a}{4}<f_1'<-\frac{5}{8}+d+a&\pm 1+c_1&1+|c_1|&\pm \frac{3}{2}&f_1'+\frac{5}{8}-\frac{a}{4}\vspace{2mm}\\
		
		-\frac{1}{2}<f_1'\leq -\frac{5}{8}+\frac{a}{4}&\pm 1+c_1&1+|c_1|&\pm \frac{1}{2}&f_1'+\frac{23}{24}-\frac{a}{4}\vspace{2mm}\\

		{\rm When~}  -\frac{5}{8}+d+a\leq f_1'\leq \frac{1}{24}+\frac{a}{4}: \vspace{2mm}\\

		-\frac{5}{8}+d+a\leq g_1'\leq \frac{1}{24}+\frac{a}{4}&c_1&2|c_1|&\pm \frac{1}{2}&f_2'-\frac{1}{24}-\frac{a}{4}\vspace{2mm}\\
		
		-\frac{5}{8}+\frac{a}{4}<g_1'<-\frac{5}{8}+d+a&\pm 1+c_1&1-|c_1|&\pm \frac{3}{2}&g_1'+\frac{5}{8}-\frac{a}{4}\vspace{2mm}\\
		
		-\frac{23}{24}+\frac{a}{4}<g_1'\leq -\frac{5}{8}+\frac{a}{4}&\pm 1+c_1&1-|c_1|&\pm \frac{1}{2}&g_1'+\frac{23}{24}-\frac{a}{4}\vspace{2mm}\\
		
		-1<g_1'\leq -\frac{23}{24}+\frac{a}{4}&\pm 1+c_1&2+2|c_1|&\pm \frac{5}{2}&f_2'+\frac{23}{24}-\frac{a}{4}\vspace{2mm}\\	
		\hline
		\end{array}
		\end{equation*}\vspace{2mm}}
	
	\noindent When $(a_4',c_4)=(\frac{1}{2},0)$, the following table gives a solution to (\ref{eq101}) i.e., of
	\begin{equation}\label{EQN27}0<F=(x_1+a_2'x_2+\frac{1}{2}x_4+\frac{1}{2}x_5)x_2-\frac{1}{6}
	x_4^2-\frac{1}{2}x_5^2-\frac{a}{4}<d+\frac{3a}{4}.\end{equation}
	
	{\scriptsize
		\begin{equation*}\begin{array}{llllllll}
		\hline\\
		{\rm ~~~~Range} & x_1 & x_1x_2 & x_4&x_2x_4&x_5&x_2x_5& {\rm ~~F}\vspace{2mm}\\
		\hline\\

		\frac{a}{4}<g_1'\leq \frac{1}{2}&c_1&-|c_1|&0&0&0&0&g_1'-\frac{a}{4}\vspace{2mm}\\	
		
		-\frac{1}{3}+\frac{a}{4}<g_1'\leq \frac{a}{4}&\pm 1+c_1&1-|c_1|&\pm 1&-1&0&0&g_1'+\frac{1}{3}-\frac{a}{4}\vspace{2mm}\\
		
		-1+\frac{a}{4}<g_1'<-1+d+a&\pm 1+c_1&1-|c_1|&0&0&0&0&g_1'+1-\frac{a}{4}\vspace{2mm}\\
		
		-1<g_1'\leq -1+\frac{a}{4}&\pm 2+c_1&2-|c_1|&\pm 1&-1&0&0&g_1'+\frac{4}{3}-\frac{a}{4}\vspace{2mm}\\
		
		\hline
		\end{array}
		\end{equation*}}

	{\scriptsize
		\begin{equation*}\begin{array}{llllllll}
		\hline\\
		{\rm ~~~~Range} & x_1 & x_1x_2 & x_4&x_2x_4&x_5&x_2x_5& {\rm ~~F}\vspace{2mm}\\
		\hline\\

		{\rm When~}  -1+d+a\leq g_1'\leq -\frac{1}{3}+\frac{a}{4}: \vspace{2mm}\\

		d+a\leq f_1'\leq 1&\pm 1+c_1&-2+2|c_1|&0&0&\pm 1&2&f_2'-\frac{3}{2}-\frac{a}{4}\vspace{2mm}\\
		
		\frac{a}{4}<f_1'<d+a&c_1&|c_1|&0&0&0&0&f_1'-\frac{a}{4}\vspace{2mm}\\	
		
		-\frac{1}{3}+\frac{a}{4}<f_1'\leq \frac{a}{4}&\pm 1+c_1&1+|c_1|&\pm 1&-1&0&0&f_1'+\frac{1}{3}-\frac{a}{4}\vspace{2mm}\\
		
		-\frac{1}{2}<f_1'\leq -\frac{1}{3}+\frac{a}{4}&c_1&2|c_1|&\pm 3&6&0&0&f_2'+\frac{3}{2}-\frac{a}{4}\vspace{2mm}\\

		\hline
		\end{array}
		\end{equation*}}
	
	\noindent This completes the proof of Lemma \ref{lem19}.\vspace{2mm}\hfill $\square$
	
	\noindent Now Theorem \ref{thm4'} follows from Lemmas \ref{lem18} and \ref{lem19}.
	
	\section{Evaluation of $\Gamma_{1,4}$ when $(m,K)=(2,1)$}\label{sec8}
	\numberwithin{equation}{section}
	
	\noindent When $(A,t,\lambda,a_5'',c_5) \neq (1/3,1/4,1/2,0,0)$ or $d+6d^5> \frac{2}{3}$, we prove that (\ref{eq101*}) and hence (\ref{eq10}) is soluble with strict inequality with $d=(8|D|)^{1/5}$ (see Theorem \ref{thm6}). Note that the real root $\gamma$ of the increasing function  $d+6d^5- \frac{2}{3}=0$ satisfies $ 0.49261<\gamma <0.49262$. \vspace{2mm}
	
	\noindent When $(A,t,\lambda,a_5'',c_5) = (1/3,1/4,1/2,0,0)$ and $d+6d^5\leq  \frac{2}{3}$, i.e., $d\leq \gamma$, we find that $0<F< \delta_2$ and $0<G< \delta_{2,1}$ do not have a solution (see Lemma \ref{prob}). (Here we have $d>0.48549$ using the given values of $a, A$ and $t$.) So we have to go to  $Q$. But now, since $d$ is not always bigger than $1/2$, Macbeath's Lemma (Lemma \ref{lem5}) is not applicable. Therefore, we can not proceed further in this case. We find that $d+6d^5$ is always $> \frac{2}{3}$  if  we work with a bigger constant for $\Gamma_{1,4}$ namely with $8.486$ instead of $8$. Thus in this case we could prove that $\Gamma_{1,4}<8.486$ only.
	\begin{theorem}\label{thm6}
		Let $Q(x_1,x_2, \cdots, x_5)$ be a real indefinite quadratic form of type $(1,4)$ and of determinant $D\neq 0$. Let $d=(8|D|)^{\frac{1}{5}}$.  Suppose that $c_2 \equiv 0 \pmod1$, $a<\frac{1}{2}$, $d\leq 1$ and $a+d\leq 1$. Let $(m,K)= (2,1)$.  Then $(\ref{eq101*})$ and hence $(\ref{eq10})$ is soluble with strict inequality, unless when $A=1/3, t=1/4, \lambda=1/2, (a_5'',c_5)=(0,0)$ and $d+6d^5< \frac{2}{3}$.
	\end{theorem}
	\noindent {\bf Proof:}
	\noindent Here, $2 < \frac{d}{a}\leq 3$ and $1<\frac{\delta_{2}}{A}=\frac{d+3a/4}{A} \leq 2$. Since, $\frac{d}{3}+d \leq a+d \leq 1$ we get $d \leq \frac{3}{4}$. Also $d/3 \leq a\leq d^{5/3}$. Now,\vspace{-2mm}
	\begin{equation}\label{Eq1}
	\frac{5d}{8}=\frac{(m+3)d}{4(K+1)} \leq A \leq \sqrt{\frac{2d^5}{3a}} \leq \sqrt{\frac{2d^4.3}{3}} = \sqrt{2}~d^2.
	\end{equation}
	\noindent  This gives $d \geq \frac{5\sqrt{2}}{16} \sim 0.4419$ and $0.27618<A<0.7955$. We divide the proof into two cases:
	
	$$\mbox{(i)~} A> 1-(d+\frac{3a}{4})  \mbox{~~~and}~~~~ \mbox{(ii)~}A\leq 1-(d+3a/4).$$
	
	\subsection{\normalfont $A> 1-(d+\frac{3a}{4})$}

	\noindent Considering $h_A=1/2$ and $k_A=1$, we find that $\vline~ h_A-Ak_A^2~\vline + \frac{1}{2} \leq d+\frac{3a}{4}$ is satisfied as $1-A<d+3a/4$. Therefore, by Macbeath's Lemma (Lemma \ref{lem5}), (\ref{eq101}) is soluble unless $A=\frac{1}{2}$ and $\beta_{A,F}=\pm a_4'-2A(c_4+\lambda x_5)\equiv 1/2 \pmod {1}$. Taking $x_5=c_5$ and $1+c_5$ simultaneously, we get $\lambda=0$. Therefore, by Lemma \ref{lem15} we need to consider $(a_4',c_4)=(0,1/2)$ or $(1/2,0)$.\vspace{3mm}
	
	\noindent Since, $A\leq C=t+A\lambda^2$ and $\lambda=0$, we have $$1/2=A\leq t =d^5/2aA=d^5/a\leq 3d^4$$ so that $$3/4 \geq d \geq (1/6)^{1/4} > 0.63894.$$
	
	\noindent For, $1/2 \leq t \leq 3d^4 \leq 0.94922$ taking $h_1=1, k_1=1$ we find that $|h_1-tk_1^2|+1/2<\delta_{2,1} =d+3a/4$ is satisfied when $t \in [0.65, 0.95]$ and the centre $h_1/k_1^2=1$ lies outside the interval, and on taking $h_2=1/2, k_2=1$ we find that $|h_2-tk_2^2|+1/2<d+3a/4$ is satisfied when $t \in [0.5, 0.65]$. Thus, (\ref{eq101'}) is soluble unless $t=\frac{1}{2}$ and $\beta_{t,G}=\pm a_5''-2tc_5 \equiv 1/2 \pmod 1$,  i.e., unless $(a_5'',c_5)=(0,1/2)$ or $(1/2,0)$. \vspace{2mm}
	
	\noindent Also, since $\lambda=0$ we get from (\ref{eq16}) that $a_5''=a_5'$. Thus, because of symmetry, we need to consider
	\begin{equation*}\begin{array}{l}
	\mbox{(i)}~~a_4'=a_5'=0, c_4=c_5=1/2~~~~~~~~~~~~~~~~~~~~~~~~~~~~~~~~~~~~~~~~~~~~~~~~~\vspace{2mm}\\
	\mbox{(ii)}~~a_4'=0, a_5'=1/2, c_4=1/2, c_5=0\vspace{2mm}\\
	\mbox{(iii)}~~a_4'=a_5'=1/2, c_4=c_5=0.\vspace{2mm}\\
	\end{array}\vspace{2mm}
	\end{equation*}	
	\noindent {\bf Case (i)}: $a_4'=0=a_5', c_4=1/2=c_5$. The following table gives a solution to (\ref{eq101}), i.e., of
	$$0<F=(x_1+a_2'x_2)x_2-\frac{1}{2}x_4^2-\frac{1}{2}x_5^2-\frac{a}{4}< d+\frac{3a}{4}.$$
	{\scriptsize
		\begin{equation*}\begin{array}{llllll}
		\hline\\
		{\rm ~~~~Range} & x_1 & x_1x_2 & x_4& x_5 & {\rm ~~F}\vspace{2mm}\\
		\hline\\
		
		\frac{1}{4}+\frac{a}{4}<f_1'\leq 1&c_1&|c_1|&\frac{1}{2}&\frac{1}{2}&f_1'-\frac{1}{4}-\frac{a}{4}\vspace{2mm}\\
		
		-\frac{1}{2}\leq f_1'<-\frac{3}{4}+d+a&\pm 1+c_1&1+|c_1|&\frac{1}{2}&\frac{1}{2}&f_1'+\frac{3}{4}-\frac{a}{4}\vspace{2mm}\\
		{\rm When ~} -\frac{3}{4}+d+a\leq f_1' \leq \frac{1}{4}+\frac{a}{4}: \vspace{2mm}\\
		
		\left.\begin{array}{l}-\frac{3}{4}+d+a\leq g_1' \leq \frac{1}{4}+\frac{a}{4}\vspace{1mm}\\f_2'<\frac{1}{4}+d+a\end{array}\right\}
		&c_1&2|c_1|&\frac{1}{2}&\frac{1}{2}&f_2'-\frac{1}{4}-\frac{a}{4}\vspace{2mm}\\
		
		\left.\begin{array}{l}-\frac{3}{4}+d+a\leq g_1' \leq \frac{1}{4}+\frac{a}{4}\vspace{1mm}\\f_2'\geq \frac{1}{4}+d+a\end{array}\right\}
		&c_1&3|c_1|&\frac{3}{2}&\frac{3}{2}&f_3'-\frac{9}{4}-\frac{a}{4}\vspace{2mm}\\

		-\frac{3}{4}+\frac{a}{4}<g_1'<-\frac{3}{4}+d+a&\pm 1+c_1&1-|c_1|&\frac{1}{2}&\frac{1}{2}&g_1'+\frac{3}{4}-\frac{a}{4}\vspace{2mm}\\

		\left.\begin{array}{l}-1< g_1' \leq -\frac{3}{4}+\frac{a}{4}\vspace{1mm}\\f_2'<-\frac{3}{4}+d+a\end{array}\right\}
		&\pm 1+c_1&2+2|c_1|&\frac{3}{2}&\frac{1}{2}&f_2'+\frac{3}{4}-\frac{a}{4}\vspace{2mm}\\
		
		\left.\begin{array}{l}-1< g_1' \leq -\frac{3}{4}+\frac{a}{4}\vspace{1mm}\\f_2'\geq -\frac{3}{4}+d+a\end{array}\right\}
		&\pm 1+c_1&3+3|c_1|&\frac{3}{2}&\frac{3}{2}&f_3'+\frac{3}{4}-\frac{a}{4}\vspace{2mm}\\
		
		\hline
		\end{array}
		\end{equation*}\vspace{1mm}}
	
	\noindent {\bf Case (ii)}: $a_4'=0, a_5'=1/2, c_4=1/2, c_5=0$. The following table gives a solution to (\ref{eq101}), i.e., of\vspace{-2mm}
	$$0<F=(x_1+a_2'x_2+\frac{1}{2}x_5)x_2-\frac{1}{2}x_4^2-\frac{1}{2}x_5^2-\frac{a}{4}< d+\frac{3a}{4}. \vspace{-2mm}$$

	{\scriptsize
		\begin{equation*}\begin{array}{lllllll}
		\hline\\
		{\rm ~~~~Range} & x_1 & x_1x_2 & x_4& x_5&x_2x_5 & {\rm ~~F}\vspace{2mm}\\
		\hline\\
		\frac{1}{8}+d+a\leq f_1'\leq 1&\pm 1+c_1&-2+2|c_1|&\frac{3}{2}&\pm 1&2&f_2'-\frac{21}{8}-\frac{a}{4}\vspace{2mm}\\		
		 \frac{1}{8}+\frac{a}{4}<f_1'<\frac{1}{8}+d+a&c_1&|c_1|&\frac{1}{2}&0&0&f_1'-\frac{1}{8}-\frac{a}{4}\vspace{2mm}\\
		
		-\frac{1}{2}\leq f_1'<-\frac{7}{8}+d+a&\pm 1+c_1&1+|c_1|&\frac{1}{2}&0&0&f_1'+\frac{7}{8}-\frac{a}{4}\vspace{2mm}\\

		{\rm When ~}-\frac{7}{8}+d+a\leq f_1' \leq \frac{1}{8}+\frac{a}{4} : \vspace{2mm}\\
		
		\left.\begin{array}{l}-\frac{7}{8}+d+a\leq g_1' \leq \frac{1}{8}+\frac{a}{4}\vspace{1mm}\\f_2'<-\frac{3}{8}+d+a\end{array}\right\}
		&c_1&2|c_1|&\pm \frac{1}{2}&\pm 1&2&f_2'+\frac{3}{8}-\frac{a}{4}\vspace{2mm}\\
		
		\left.\begin{array}{l}-\frac{7}{8}+d+a\leq g_1' \leq \frac{1}{8}+\frac{a}{4}\vspace{1mm}\\f_2'\geq -\frac{3}{8}+d+a\end{array}\right\}
		&c_1&2|c_1|&\frac{1}{2}&0&0&f_2'-\frac{1}{8}-\frac{a}{4}\vspace{2mm}\\

		-\frac{7}{8}+\frac{a}{4}<g_1'<-\frac{7}{8}+d+a&\pm 1+c_1&1-|c_1|&\frac{1}{2}&0&0&g_1'+\frac{7}{8}-\frac{a}{4}\vspace{2mm}\\	
		
		\left.\begin{array}{l}-1\leq g_1' \leq -\frac{7}{8}+\frac{a}{4}\vspace{1mm}\\f_2'\leq-\frac{7}{8}+\frac{a}{4}\end{array}\right\}
		&\pm 1+c_1&2+2|c_1|&\frac{3}{2}&\pm 1&2&f_2'+\frac{11}{8}-\frac{a}{4}\vspace{2mm}\\
		
		\left.\begin{array}{l}-1\leq g_1' \leq -\frac{7}{8}+\frac{a}{4}\vspace{1mm}\\f_2'> -\frac{7}{8}+\frac{a}{4}\end{array}\right\}
		&\pm 1+c_1&2+2|c_1|&\frac{3}{2}&0&0&f_2'+\frac{7}{8}-\frac{a}{4}\vspace{2mm}\\
		
		\hline
		\end{array}
		\end{equation*}\vspace{1mm}}
	
	\noindent {\bf Case (iii)}: $a_4'=1/2=a_5', c_4=0=c_5$. Here $0<F=(x_1+a_2'x_2+\frac{1}{2}x_4+\frac{1}{2}x_5)x_2-\frac{1}{2}x_4^2-\frac{1}{2}x_5^2-\frac{a}{4}< d+\frac{3a}{4}$ has no solution when $c_1=\frac{1}{2}$, $a_2'=\frac{1}{2}$ and $a+d=2d^5+d<1$. This is so because $ (2x+1+x_2+x_4+x_5)x_2-x_4^2-x_5^2$ never takes the value $1$ for integers $x,x_2,x_4,x_5$. So we have to go to $Q$. We divide the range of $a$ into five subintervals and in each subinterval choose suitable integers $2h_n$ and $k_n$ such that $\vline~ h_n-ak_n^2~\vline + \frac{1}{2} <d.$\vspace{1mm}
	
	{\footnotesize
		\begin{equation*}\begin{array}{lllllllll}
		\hline\\
		n & (h_n, k_n) & \lambda_n&\mbox{Remarks}&~& n & (h_n, k_n) & \lambda_n&\mbox{Remarks}\vspace{2mm}\\
		\hline\\
		
		1&(0.5,1)&0.311&\mbox{na}&&2&(5,4)&0.301&\mbox{na}\vspace{2mm}\\
		
		3&(7.5,5)&0.295&\mbox{tbd}&&4&(1,2)&0.215&\mbox{tbd}\vspace{2mm}\\
		
		5&(2,3)&0.2129&\mbox{na}&&&&&\\

		\hline
		\end{array}
		\end{equation*}\vspace{1mm}}
	
	\noindent The choice of $(h_n,k_n)$ is done in a manner similar to that in the Case  $m=2,K=3, L=4$ (now for $a$ instead of $t$). Thus, (\ref{eq101*}) is soluble unless $a=\frac{3}{10}$ or $\frac{1}{4}$ and $\beta_{a,Q}=\pm a_3-2a(c_3+h_4x_4+h_5x_5) \equiv h/k \pmod {1/k, 2a}$.\vspace{4mm}
	
	\noindent \textbf{{Subcase (i): $a=\frac{3}{10}$.}}\vspace{2mm}
	
	\noindent Here, $t=\frac{d^5}{2aA}, t=\frac{1}{2}, A=\frac{1}{2}$ and $a=\frac{3}{10}=2d^5$ gives $d=\left(  \frac{3}{20} \right)^{1/5} > 0.68425$.
	
	\noindent $\beta_{a,Q}=\pm a_3-2a(c_3+h_4x_4+h_5x_5) \equiv h/k \pmod {1/k, 2a}$ gives $\pm a_3-\frac{3}{5}(c_3+h_4x_4+h_5x_5) \equiv 7.5/5 \pmod {1/5,3/5}$. Taking $x_4=c_4$ and $1+c_4$ simultaneously, we get $3h_4\equiv 0\pmod 1$, i.e., $h_4=0$ or $1/3$. Similarly taking $x_5=c_5$ and $1+c_5$ simultaneously, we get $3h_5\equiv 0\pmod 1$ i.e., $h_5=0$ or $1/3$. If $h_4=1/3$, choose $(x_3,x_5) \equiv (c_3,c_5) \pmod 1$ arbitrarily and take $x_2=\pm 1$, $x_1=x+c_1, x_4=y+c_4$, so that (\ref{eq101*}) becomes\vspace{-2mm}
	\begin{equation}\label{key1}
	0<\pm x+\beta_{\alpha} y - \alpha y^2+\nu <d,
	\end{equation}
	
	\noindent where $\alpha=1/2+\frac{1}{30}=\frac{8}{15}$, $\beta_{\alpha}=\pm a_4-2ah_4(x_3+h_4c_4+h_5x_5)-2Ac_4$ and $\nu$ is some real number. Taking $h=\frac{1}{2}, k=1$, we find that $|h-\alpha k^2|+1/2<d$, so that by Macbeath Lemma, (\ref{eq101*}) is soluble.  Similarly, for $h_5=1/3$, (\ref{eq101*}) is soluble. Therefore, we must have $ h_4=h_5=0$. Now as $a_4'=1/2=a_5'$, using (\ref{eq12}), we get $a_4=1/2=a_5$.\vspace{2mm}

	\noindent Now $\pm a_3-\frac{3}{5}c_3 \equiv  7.5/5 \pmod {1/5,3/5}$ gives $\pm 5a_3-3c_3\equiv \frac{1}{2} \pmod 1$ which gives $(a_3,c_3)=(0,1/6)$, $(0,1/2)$, $(1/10,0)$ and $(1/10,1/3)$. In all of these, on taking $x_2=\pm 1$ and choosing $x_3\equiv c_3 \pmod 1$ suitably, we find that (\ref{eq101*}) is soluble.\vspace{4mm}

	\noindent \textbf{{Subcase (ii): $a=\frac{1}{4}$.}}\vspace{2mm}
	
	\noindent Here, $t=\frac{d^5}{2aA}, t=\frac{1}{2}, A=\frac{1}{2}$ and $a=\frac{1}{4}=2d^5$ gives $d=\left(  \frac{1}{8} \right)^{1/5} > 0.6597$. Also $h=1, k=2$.
	$\beta_{a,Q}=\pm a_3-2a(c_3+h_4x_4+h_5x_5) \equiv h/k \pmod {1/k, 2a}$ gives $\pm a_3-(c_3+h_4x_4+h_5x_5)/2 \equiv 1/2 \pmod {1/2}$. Taking $x_4=c_4$ and $1+c_4$ simultaneously, we get $h_4=0$ and taking $x_5=c_5$ and $1+c_5$ simultaneously, we get $h_5=0$. Therefore, we get $\pm 2a_3-c_3 \equiv 0 \pmod {1}$. Applying Lemma \ref{lem15}, we get $(a_3,c_3)=(0,0)$ or $(1/4,1/2)$. Also  as $a_4'=1/2=a_5'$, and $h_4=0=h_5$, we get using (\ref{eq12}) that $(a_4,c_4)=(a_5,c_5)=(1/2,0)$.\vspace{4mm}

	\noindent When $(a_3,c_3)=(0,0)$, on taking $x_5=0$, the following table gives a solution to (\ref{eq101*}) i.e., of\vspace{-2mm}
	\begin{equation}
	 0<Q=(x_1+a_2x_2+\frac{1}{2}x_4+\frac{1}{2}x_5)x_2-\frac{1}{4}x_3^2-\frac{1}{2}x_4^2-\frac{1}{2}x_5^2<d.
	\end{equation}
	
	{\scriptsize
		\begin{equation*}\begin{array}{lllllll}
		\hline\\
		{\rm ~~~~Range} & x_1 & x_1x_2 &x_3& x_4&x_2x_4& {\rm ~~Q}\vspace{2mm}\\
		\hline\\

		d+\frac{1}{4}\leq f_1\leq 1&\pm 1+c_1&-2+2|c_1|&\pm 2&\pm 1&2&f_2'-\frac{5}{2}\vspace{2mm}\\	
		
		\frac{1}{4}<f_1< d+\frac{1}{4}&c_1&|c_1|&\pm 1&0&0&f_1-\frac{1}{4}\vspace{2mm}\\
		
		0<f_1\leq\frac{1}{4}&c_1&|c_1|&0&0&0&f_1\vspace{2mm}\\	
		
		-\frac{1}{2}\leq f_1<d-\frac{3}{4}&\pm 1+c_1&1+|c_1|&\pm 1&0&0&f_1+\frac{3}{4}\vspace{2mm}\\
		
		{\rm When ~} d-\frac{3}{4}\leq f_1 \leq 0 :\vspace{2mm}\\
		
		d-\frac{3}{4}\leq g_1\leq 0&c_1&2|c_1|&0&\pm 1&2&f_2'+\frac{1}{2}\vspace{2mm}\\
		
		-\frac{3}{4}<g_1<d-\frac{3}{4}&\pm 1+c_1&1-|c_1|&\pm 1&0&0&g_1+\frac{3}{4}\vspace{2mm}\\
		
		-1<g_1\leq-\frac{3}{4}&\pm 1+c_1&1-|c_1|&0&0&0&g_1+1\vspace{2mm}\\
		
		\hline
		\end{array}
		\end{equation*}\vspace{4mm}}
	
	\noindent When $(a_3,c_3)=(1/4,1/2)$, on taking $x_4=0=x_5$, the following table gives a solution to (\ref{eq101*}) i.e., of \vspace{-2mm}
	\begin{equation}
	 0<Q=(x_1+a_2x_2+\frac{1}{4}x_3+\frac{1}{2}x_4+\frac{1}{2}x_5)x_2-\frac{1}{4}x_3^2-\frac{1}{2}x_4^2-\frac{1}{2}x_5^2<d.
	\end{equation}

	{\scriptsize
		\begin{equation*}\begin{array}{llllll}
		\hline\\
		{\rm ~~~~Range} & x_1 & x_1x_2 &x_3&x_2x_3& {\rm ~~Q}\vspace{2mm}\\
		\hline\\
		-\frac{1}{16}<g_1\leq \frac{1}{2}&\pm 1+c_1&1-|c_1|&\pm \frac{3}{2}&-\frac{3}{2}&g_1+\frac{1}{16}\vspace{2mm}\\
		
		-\frac{13}{16}<g_1< d-\frac{13}{16}&\pm 1+c_1&1-|c_1|&\pm \frac{1}{2}&-\frac{1}{2}&g_1+\frac{13}{16}\vspace{2mm}\\
		
		-1<g_1\leq -\frac{13}{16}&\pm 2+c_1&2-|c_1|&\pm \frac{3}{2}&-\frac{3}{2}&g_1+\frac{17}{16}\vspace{2mm}\\

		{\rm When ~}-\frac{13}{16}+d\leq g_1 \leq -\frac{1}{16}	: \vspace{2mm}\\

		\frac{3}{16}+d\leq f_1\leq 1  &\pm 1+c_1&-2+2|c_1|&\pm \frac{1}{2}&-1&f_2-\frac{37}{16}\vspace{2mm}\\
		
		\frac{3}{16}<f_1< d+\frac{3}{16}&c_1&|c_1|&\pm \frac{1}{2}&-\frac{1}{2}&f_1-\frac{3}{16}\vspace{2mm}\\
		
		-\frac{1}{16}<f_1\leq\frac{3}{16}&\pm 1+c_1&1+|c_1|&\pm \frac{3}{2}&-\frac{3}{2}&f_1+\frac{1}{16}\vspace{2mm}\\
		
		-\frac{1}{2}<f_1 \leq -\frac{1}{16}&\pm 1+c_1&2+2|c_1|&\pm \frac{7}{2}&7&f_2+\frac{11}{16}\vspace{2mm}\\
		
		\hline
		\end{array}
		\end{equation*}}	
	
	\subsection{\normalfont $A\leq 1-(d+\frac{3a}{4})$}
	
	\noindent Here, we have $\frac{5d}{8}\leq A\leq 1-\left( d+\frac{3a}{4}  \right)$ i.e., $d\leq \frac{8}{15}$. Already, $d\geq 0.4419$, so that
	
	$$0.4419\leq d\leq \frac{8}{15}, \mbox{~~~and}$$
	\begin{equation}
	0.2071<\frac{15d}{32}=\frac{3}{4}\cdot \frac{5d}{8}\leq \frac{3A}{4}\leq t=\frac{d^5}{2aA}\leq \frac{12d^3}{5}<0.3643=\lambda_0 \mbox{(say)}.
	\end{equation}

	\noindent We divide the range of $t$ into $22$ subintervals $[\lambda_n,\lambda_{n-1}], n=1,2,\cdots 22$ and in each subinterval choose suitable integers $2h_n$ and $k_n$ such that
	$\vline~ h_n-tk_n^2~\vline +\frac{1}{2} <d+\frac{3a}{4}$.
	The choice of $(h_n,k_n)$ is done in a manner already described.\vspace{1mm}

	{\scriptsize
		\begin{equation*}\begin{array}{lllllllllll}
		\hline\\
		n & h_n & k_n & \lambda_n&\mbox{Remarks}&~~~~ &n & h_n & k_n & \lambda_n&\mbox{Remarks}\vspace{2mm}\\
		\hline\\
		
		1&0.5&1&0.34531&\mbox{na}&&2&3&3&0.3182&\mbox{tbd}\vspace{2mm}\\
		
		3&5&4&0.3051&\mbox{tbd}&&4&7.5&5&0.2952&\mbox{tbd}\vspace{2mm}\\
		
		5&14.5&7&0.294&\mbox{na}&&6&10.5&6&0.2885&\mbox{tbd}\vspace{2mm}\\
		
		7&2.5&3&0.2771&(4.5,4)&&8&1&2&0.2317&\mbox{tbd}\vspace{2mm}\\
		
		9&28&11&0.23081&(52,15)&&10&39&13&0.2307&\mbox{tbd}\vspace{2mm}\\
		
		11&23&10&0.2302&\mbox{na}&&12&2&3&0.2156&(3.5,4)\vspace{2mm}\\
		
		13&17.5&9&0.2155&\mbox{na}&&14&10.5&7&0.21311&\mbox{tbd}\vspace{2mm}\\
		
		15&36&13&0.21268&(69,18)&&16&61.5&17&0.21264&\mbox{na}\vspace{2mm}\\
		
		17&85&20&0.21236&\mbox{tbd}&&18&143.5&26&0.2122&(178.5,29)\vspace{2mm}\\
		
		19&30.5&12&0.21181&\mbox{na}&&20&13.5&8&0.2101&(93,21)\vspace{2mm}\\
		
		21&17&9&0.2098&(111,23)&&22&7.5&6&0.2071&\mbox{tbd}\vspace{2mm}\\

		\hline
		\end{array}
		\end{equation*}\vspace{1mm}}
	
	\noindent  Thus (\ref{eq101*}) is soluble unless $t=1/3,$ $5/16,$ $3/10,$ $7/24,$ $1/4,$ $3/13,$ $3/14,$ $17/80,$ $5/24,$ and $\beta_{t,G} = \pm a_5''-2tc_5\equiv h_n/k_n \pmod {1/k_n, 2t}$.
	Here, we have\vspace{-2mm}
	\begin{equation*}\begin{array}{l}
	(5t/12)^{1/3}\leq d \leq \min{\{8/15, 32t/15\}}\vspace{2mm}\\
	5d/8\leq A\leq \min \{1-5d/4, 4t/3\}\vspace{2mm}\\
	d/3\leq a\leq \min \{d/2, 4A/3\}.\vspace{2mm}\\
	\end{array}\vspace{-2mm}
	\end{equation*}

	\noindent  For $t=5/16,$ $3/10,$ $7/24,$ $3/13,$ $3/14,$ $17/80,$ and $5/24,$ we compute the values of $(a_5'',c_5)$, using Lemma \ref{lem15}, and find that (\ref{eq101'}) is soluble by taking $x_2=\pm 1$ and suitably choosing $x_5 \equiv c_5 \pmod 1$.\vspace{2mm}
	
	\noindent We discuss $t=1/3$ and $t=1/4$ separately in the following lemmas:\vspace{2mm}
	
	\begin{lemma}\label{lem21}
		When $t=1/3,$ $(\ref{eq101'})$ and hence $(\ref{eq17})$ is soluble with strict inequality.
	\end{lemma}\vspace{-3mm}
	
	\noindent {\bf Proof:} Here, we have
	\begin{equation}\label{Eqn1''}
	\begin{array}{l}
	0.5178<(5t/12)^{1/3}\leq d\leq 8/15<0.5334\vspace{2mm}\\
	0.3236< 5d/8 \leq A \leq 1-5d/4 < 0.35275\vspace{2mm}\\
	0.1726< d/3\leq a\leq d/2< 0.2667\vspace{2mm}\\
	0.12405 < \frac{a+A}{4} < 0.15487.
	\end{array}
	\end{equation}
	
	\noindent We have $\beta_{t,G}=\pm a_5'-\frac{2}{3}c_5\equiv 1 ({\rm mod~} \frac{1}{3},\frac{2}{3})$, i.e., $\pm 3a_5'-2c_5 \equiv 0 \pmod 1$.
	
	\noindent Using Lemma \ref{lem15} we need to consider $(a_5'',c_5)=(0,0),$ $(0,1/2)$ and $(1/6,1/4).$ In each of these, we give a solution to (\ref{eq101'}), i.e., of\vspace{-2mm}
	\begin{equation}\label{EQ58'}0<G=(x_1+a_2''x_2+a_5''x_5)x_2-\frac{1}{3}x_5^2-\frac{a+A}{4}
	<d+\frac{3a}{4}.\end{equation}
	\noindent \textbf{Case (i): $(a_5'',c_5)=(0,0)$}
	
	{\scriptsize
		\begin{equation*}\begin{array}{lllll}
		\hline\\
		{\rm ~~~~Range} & x_1 & x_1x_2 & x_5& {\rm ~~G}\vspace{2mm}\\
		\hline\\
		
		\frac{1}{3}+\frac{a+A}{4}<f_1''\leq 1&c_1&|c_1|&1&f_1''-\frac{1}{3}-\frac{a+A}{4}\vspace{2mm}\\
		
		\frac{a+A}{4}<f_1''\leq \frac{1}{3}+\frac{a+A}{4}&c_1&|c_1|&0&f_1''-\frac{a+A}{4}\vspace{2mm}\\
		
		-\frac{2}{3}+\frac{a+A}{4}<f_1''<-\frac{2}{3}+d+\frac{A}{4}+a&\pm 1+c_1&1+|c_1|&1&f_1''+\frac{2}{3}-\frac{a+A}{4}\vspace{2mm}\\
		
		-\frac{1}{2}<f_1''\leq -\frac{2}{3}+\frac{a+A}{4}&\pm 1+c_1&1+|c_1|&0&f_1''+1-\frac{a+A}{4}\vspace{2mm}\\

		{\rm When~}  -\frac{2}{3}+d+\frac{A}{4}+a\leq f_1''\leq \frac{a+A}{4}: \vspace{2mm}\\

		-\frac{2}{3}+d+\frac{A}{4}+a\leq g_1''\leq \frac{a+A}{4}&c_1&2|c_1|&0&f_2''-\frac{a+A}{4}\vspace{2mm}\\	
		
		-\frac{2}{3}+\frac{a+A}{4}<g_1''<-\frac{2}{3}+d+\frac{A}{4}+a&\pm 1+c_1&1-|c_1|&1&g_1''+\frac{2}{3}-\frac{a+A}{4}\vspace{2mm}\\
		
		-1+\frac{a+A}{4}<g_1''\leq -\frac{2}{3}+\frac{a+A}{4}&\pm 1+c_1&1-|c_1|&0&g_1''+1-\frac{a+A}{4}\vspace{2mm}\\
		
		-1<g_1''\leq -1+\frac{a+A}{4}&\pm 2+c_1&4+2|c_1|&3&f_2''+1-\frac{a+A}{4}\vspace{2mm}\\

		\hline
		\end{array}
		\end{equation*}\vspace{1mm}}

	\noindent \textbf{Case (ii): $(a_5'',c_5)=(0,1/2)$}
	
	{\scriptsize
		\begin{equation*}\begin{array}{lllll}
		\hline\\
		{\rm ~~~~Range} & x_1 & x_1x_2 & x_5& {\rm ~~G}\vspace{2mm}\\
		\hline\\
		
		\frac{1}{12}+\frac{a+A}{4}<g_1''\leq \frac{1}{2}&c_1&-|c_1|&\frac{1}{2}&g_1''-\frac{1}{12}-\frac{a+A}{4}\vspace{2mm}\\
		
		-\frac{1}{4}+\frac{a+A}{4}<g_1''\leq \frac{1}{12}+\frac{a+A}{4}&\pm 1+c_1&1-|c_1|&\frac{3}{2}&g_1''+\frac{1}{4}-\frac{a+A}{4}\vspace{2mm}\\
		
		-\frac{11}{12}+\frac{a+A}{4}<g_1''<-\frac{11}{12}+d+\frac{A}{4}+a&\pm 1+c_1&1-|c_1|&\frac{1}{2}&g_1''+\frac{11}{12}-\frac{a+A}{4}\vspace{2mm}\\
		
		-1<g_1''\leq-\frac{11}{12}+\frac{a+A}{4}&\pm 2+c_1&2-|c_1|&\frac{3}{2}&g_1''+\frac{5}{4}-\frac{a+A}{4}\vspace{2mm}\\	
		\hline
		\end{array}
		\end{equation*}}

	{\scriptsize
		\begin{equation*}\begin{array}{lllll}
		\hline\\
		{\rm ~~~~Range} & x_1 & x_1x_2 & x_5& {\rm ~~G}\vspace{2mm}\\
		\hline\\
		
		{\rm When ~} -\frac{11}{12}+d+\frac{A}{4}+a\leq g_1''\leq -\frac{1}{4}+\frac{a+A}{4}:\vspace{2mm}\\

		\frac{1}{12}+d+\frac{A}{4}+a\leq f_1''\leq 1&\pm c_1&-2|c_1|&\frac{1}{2}&g_2''-\frac{1}{12}-\frac{a+A}{4}\vspace{2mm}\\	
		
		 \frac{1}{12}+\frac{a+A}{4}<f_1''<\frac{1}{12}+d+\frac{A}{4}+a&c_1&|c_1|&\frac{1}{2}&f_1''-\frac{1}{12}-\frac{a+A}{4}\vspace{2mm}\\
		
		-\frac{1}{4}+\frac{a+A}{4}<f_1''\leq \frac{1}{12}+\frac{a+A}{4}&\pm 1+c_1&1+|c_1|&\frac{3}{2}&f_1''+\frac{1}{4}-\frac{a+A}{4}\vspace{2mm}\\
		
		-\frac{1}{2}<f_1''\leq -\frac{1}{4}+\frac{a+A}{4}&\pm 1+c_1&2+2|c_1|&\frac{3}{2}&f_2''+\frac{5}{4}-\frac{a+A}{4}\vspace{2mm}\\

		\hline
		\end{array}
		\end{equation*}\vspace{1mm}}
	
	\noindent \textbf{Case (iii): $(a_5'',c_5)=(1/6,1/4)$}

	{\scriptsize
		\begin{equation*}\begin{array}{lllll}
		\hline\\
		{\rm ~~~~Range} & x_1 & x_2 & x_5& {\rm ~~G}\vspace{2mm}\\
		\hline\\
		
		\frac{5}{16}+\frac{a+A}{4}<p_1''\leq 1&c_1&1&\frac{5}{4}&p_1''-\frac{5}{16}-\frac{a+A}{4}\vspace{2mm}\\
		
		-\frac{1}{48}+\frac{a+A}{4}<p_1''\leq \frac{5}{16}+\frac{a+A}{4}&c_1&1&\frac{1}{4}&p_1''+\frac{1}{48}-\frac{a+A}{4}\vspace{2mm}\\
		
		 -\frac{11}{16}+\frac{a+A}{4}<p_1''<-\frac{11}{16}+d+\frac{A}{4}+a&1+c_1&1&\frac{5}{4}&p_1''+\frac{11}{16}-\frac{a+A}{4}\vspace{2mm}\\
		
		-\frac{49}{48}+\frac{a+A}{4}< p_1''\leq -\frac{11}{16}+\frac{a+A}{4}&1+c_1&1&\frac{1}{4}&p_1''+\frac{49}{48}-\frac{a+A}{4}\vspace{2mm}\\
		
		-\frac{27}{16}+d+\frac{A}{4}+a\leq p_1''\leq -\frac{49}{48}+\frac{a+A}{4}&-1+c_1&-1&\frac{5}{4}&q_1''+\frac{13}{48}-\frac{a+A}{4}\vspace{2mm}\\
		
		 -1<p_1''<-\frac{27}{16}+d+\frac{A}{4}+a&2+c_1&1&\frac{5}{4}&p_1''+\frac{27}{16}-\frac{a+A}{4}\vspace{2mm}\\

		{\rm When~}  -\frac{11}{16}+d+\frac{A}{4}+a\leq p_1''\leq -\frac{1}{48}+\frac{a+A}{4}: \vspace{2mm}\\
		
		\frac{1}{16}+d+\frac{A}{4}+a\leq q_1''\leq 1&c_1&2&\frac{9}{4}&p_2''-\frac{15}{16}-\frac{a+A}{4}\vspace{2mm}\\
		
		 \frac{1}{16}+\frac{a+A}{4}<q_1''<\frac{1}{16}+d+\frac{A}{4}+a&c_1&-1&\frac{1}{4}&q_1''-\frac{1}{16}-\frac{a+A}{4}\vspace{2mm}\\
		
		 -\frac{13}{48}+\frac{a+A}{4}<q_1''\leq\frac{1}{16}+\frac{a+A}{4}&-1+c_1&-1&\frac{5}{4}&q_1''+\frac{13}{48}-\frac{a+A}{4}\vspace{2mm}\\
		
		-\frac{15}{16}+d+\frac{A}{4}+a\leq q_1''\leq -\frac{13}{48}+\frac{a+A}{4}&c_1&2&\frac{1}{4}&p_2''+\frac{1}{16}-\frac{a+A}{4}\vspace{2mm}\\
		
		 -\frac{15}{16}+\frac{a+A}{4}<q_1''<-\frac{15}{16}+d+\frac{A}{4}+a&-1+c_1&-1&\frac{1}{4}&q_1''+\frac{15}{16}-\frac{a+A}{4}\vspace{2mm}\\
		
		 -1<q_1''\leq-\frac{15}{16}+\frac{a+A}{4}&-2+c_1&-1&\frac{5}{4}&q_1''+\frac{61}{48}-\frac{a+A}{4}\vspace{2mm}\\
		
		\hline
		\end{array}
		\end{equation*}\vspace{1mm}}
	
	\noindent This proves Lemma \ref{lem21}. \hfill $\square$
	
	\begin{remark}\label{rem1}\normalfont Recall the inequality (\ref{eq101}) is\vspace{-2mm}
		$$ 0<F=(x_1+a_2'x_2+a_4'x_4+a_5'x_5)x_2-A(x_4+\lambda x_5)^2-tx_5^2-\frac{a}{4}<\delta_m=d+\frac{(m^2-1)a}{4}.\vspace{-2mm}$$
		
		\noindent In addition to (\ref{eq102}) and (\ref{eq103}), sometimes we choose $x_4\equiv c_4\pmod 1$ arbitrarily, $x_2=\pm1$, $x_1=x+c_1$ and $x_5=y+c_5$ to get\vspace{-2mm}
		$$ F(x+c_1, \pm1,x_4,y+c_5)=\pm x+\beta'_{F,C}y-Cy^2+\nu, \vspace{-2mm}$$ where $C=A\lambda^2+t$, $\beta'_{F,C}=\pm a'_5 -2A\lambda x_4-2Cc_5$. Then for a solution of (\ref{eq101}), we find integers $2h_C$ and $k_C$ such that \vspace{-2mm}
		
		\begin{equation}\label{eq45} |h_C-Ck_C^2|+\frac{1}{2}<\delta_m\end{equation}
		\noindent and apply Macbeath's Lemma. We can also rewrite $F$ as
		$$F=-C(x_5+A\lambda x_4/C-a_5'x_2/2C)^2+(x_1+a_2^*x_2+a_4^*x_4)x_2-Atx_4^2/C-a/4,$$ where
		\begin{equation}\label{G*}
		\begin{array}{l}
		a_2^*\equiv a_2'+a_5'^2/4C ~({\rm mod~} 1),\vspace{2mm}\\
		a_4^*\equiv a_4'-\lambda Aa_5'/C ~({\rm mod~} 1), -1/2 <a_2^*,a_4^*\leq 1/2.
		\end{array}
		\end{equation}
		
		\noindent If $M$ is an integer $\geq 1$ such that $M<\delta_{m}/C \leq M+1,$ then
		(\ref{eq101}) is soluble, by Lemma \ref{lem3} if
		\begin{equation}\begin{array}{lll}\label{e17}
		0<G^*(x_1,x_2,x_4)& =& (x_1+a_2^*x_2+a_4^*x_4)x_2 - Atx_4^2/C - a/4-C/4\vspace{2mm}\\
		&<&\delta_m +(M^2-1)C/4 = \delta_{m,M}^* {\rm ~(say)}
		\end{array}
		\end{equation}
		Take $x_2=\pm1$, $x_1=x+c_1$ and $x_4=y+c_4$ to get
		$$ G^*(x+c_1, \pm1,y+c_4)=\pm x+\beta^*_{G}y-\frac{At}{C}y^2+\nu,$$ where $\beta^*_{G}=\pm a^*_4 -2\frac{At}{C} c_4$. Then for a solution of (\ref{e17}) and hence of (\ref{eq101}), we find integers $2h^*$ and $k^*$ such that
		\begin{equation}\label{eq45'} |h^*-\frac{At}{C}k^{{*}^2}|+\frac{1}{2}<\delta_{m,M}^*\end{equation}
		and apply Macbeath's Lemma. We will use these in Lemma \ref{lem22} and also in  later Sections \ref{sec9}-\ref{sec11}.\end{remark}
	
	\begin{lemma}\label{lem22}
		When $t=1/4,$ $(\ref{eq101})$ and hence $(\ref{eq13})$ is soluble with strict inequality except when $(a_5'',c_5)=(0,0)$, $A=1/3$, $\lambda=1/2$.
	\end{lemma}

	\noindent {\bf Proof:} Here, we have
	\begin{equation}\label{Eqn1'}
	\begin{array}{l}
	0.4705<(5t/12)^{1/3}\leq d\leq 8/15<0.5334\vspace{2mm}\\
	0.1568< d/3\leq a\leq d/2< 0.2667\vspace{2mm}\\
	0.294< 5d/8 \leq A \leq 4t/3 < 0.3334\vspace{2mm}\\
	0.1127 \leq \frac{a+A}{4} \leq 0.1501.
	\end{array}
	\end{equation}
	\noindent We have $\beta_{t,G}=\pm a_5'-\frac{2}{4}c_5\equiv \frac{1}{2} ({\rm mod~} \frac{1}{2},\frac{1}{2})$, i.e., $\pm 2a_5'-c_5 \equiv 0\pmod 1$.
	\noindent Using Lemma \ref{lem15}, we need to consider $(a_5'',c_5)=(0,0)$ and $(1/4,1/2).$\vspace{2mm}

	\noindent When $(a_5'',c_5)=(1/4,1/2)$, the following table gives a solution to (\ref{eq101'}), i.e., of
	 \begin{equation}\label{EQ58}0<G=(x_1+a_2''x_2+\frac{1}{4}x_5)x_2-\frac{1}{4}x_5^2-\frac{a+A}{4}<d+\frac{3a}{4}.\end{equation}
	
	{\scriptsize
		\begin{equation*}\begin{array}{llllll}
		\hline\\
		{\rm ~~~~Range} & x_1 & x_1x_2 & x_5&x_2x_5& {\rm ~~G}\vspace{2mm}\\
		\hline\\
		
		\frac{3}{16}+d+\frac{A}{4}+a\leq f_1''\leq 1&\pm1+c_1&-2+2|c_1|&\pm \frac{1}{2}&-1&f_2''-\frac{37}{16}-\frac{a+A}{4}\vspace{2mm}\\
		
		\frac{3}{16}+\frac{a+A}{4}<f_1''<\frac{3}{16}+d+\frac{A}{4}+a&c_1&|c_1|&\pm \frac{1}{2}&-\frac{1}{2}&f_1''-\frac{3}{16}-\frac{a+A}{4}\vspace{2mm}\\
		
		-\frac{1}{16}+\frac{a+A}{4}<f_1''\leq\frac{3}{16}+\frac{a+A}{4}&c_1&|c_1|&\pm \frac{1}{2}&\frac{1}{2}&f_1''+\frac{1}{16}-\frac{a+A}{4}\vspace{2mm}\\	
		
		-\frac{1}{2}<f_1''<-\frac{13}{16}+d+\frac{A}{4}+a&\pm 1+c_1&1+|c_1|&\pm \frac{1}{2}&-\frac{1}{2}&f_1''+\frac{13}{16}-\frac{a+A}{4}\vspace{2mm}\\
		
		\hline
		\end{array}
		\end{equation*}\vspace{1mm}}

	{\scriptsize
		\begin{equation*}\begin{array}{llllll}
		\hline\\
		{\rm ~~~~Range} & x_1 & x_1x_2 & x_5&x_2x_5& {\rm ~~G}\vspace{2mm}\\
		\hline\\
		
		{\rm When ~} -\frac{13}{16}+d+\frac{A}{4}+a\leq f_1'' \leq -\frac{1}{16}+\frac{a+A}{4}:\vspace{2mm}\\
		
		\left.\begin{array}{l}-\frac{13}{16}+d+\frac{A}{4}+a\leq g_1'' \leq -\frac{1}{16}+\frac{a+A}{4}\vspace{1mm}\\f_2''> -\frac{3}{16}+\frac{a+A}{4}\end{array}\right\}	 &c_1&2|c_1|&\pm \frac{1}{2}&1&f_2''+\frac{3}{16}-\frac{a+A}{4}\vspace{2mm}\\
		
		\left.\begin{array}{l}-\frac{13}{16}+d+\frac{A}{4}+a\leq g_1'' \leq -\frac{1}{16}+\frac{a+A}{4}\vspace{1mm}\\f_2''\leq -\frac{3}{16}+\frac{a+A}{4}\end{array}\right\}	 &\pm 1+c_1&2+2|c_1|&\pm \frac{7}{2}&7&f_2''+\frac{11}{16}-\frac{a+A}{4}\vspace{2mm}\\
		
		-\frac{13}{16}+\frac{a+A}{4}<g_1''<-\frac{13}{16}+d+\frac{A}{4}+a&\pm 1+c_1&1-|c_1|&\pm \frac{1}{2}&-\frac{1}{2}&g_1''+\frac{13}{16}-\frac{a+A}{4}\vspace{2mm}\\
		
		-\frac{17}{16}+\frac{a+A}{4}<g_1''\leq-\frac{13}{16}+\frac{a+A}{4}&\pm 1+c_1&1-|c_1|&\pm \frac{1}{2}&\frac{1}{2}&g_1''+\frac{17}{16}-\frac{a+A}{4}\vspace{2mm}\\
		
		\left.\begin{array}{l}-1\leq g_1'' \leq -\frac{17}{16}+\frac{a+A}{4}\vspace{1mm}\\f_2''<-\frac{27}{16}+d+\frac{A}{4}+a\end{array}\right\}	 &\pm 1+c_1&2+2|c_1|&\pm \frac{5}{2}&5&f_2''+\frac{27}{16}-\frac{a+A}{4}\vspace{2mm}\\
		
		\left.\begin{array}{l}-1\leq g_1'' \leq -\frac{17}{16}+\frac{a+A}{4}\vspace{1mm}\\f_2''\geq -\frac{27}{16}+d+\frac{A}{4}+a\end{array}\right\}	 &\pm 2+c_1&4+2|c_1|&\pm \frac{9}{2}&9&f_2''+\frac{19}{16}-\frac{a+A}{4}\vspace{2mm}\\
		
		\hline
		\end{array}
		\end{equation*}\vspace{1mm}}

	\noindent When $(a_5'',c_5)=(0,0)$, (\ref{eq17}) has no solution for $(a_2'',c_1)=(0,0)$ and for $d+a+A/4< \frac{3}{4}$. This is so because $0<G=x_1x_2-\frac{1}{4}x_5^2-\frac{a+A}{4}\leq d+\frac{3a}{4}$ with $x_1,x_2,x_5 \in \mathbb{Z}$ is not soluble, as $0<a+A\leq  4x_1x_2-x_5^2\leq  4(d+a+A/4)<3$ has no solution. We note that $4x_1x_2-x_5^2$ never takes the value $1$ or $2$ for integers $x_1,x_2,x_5$, as $-x_5^2 \not \equiv 1 {\rm ~ or ~} 2 \pmod 4$. So we have to go to $F$. We will apply Macbeath Lemma on (\ref{eq101}) and find integers $2h_n$ and $k_n$ such that $\vline~ h_n-Ak_n^2~\vline +\frac{1}{2} <d+\frac{3a}{4}$, in a manner described earlier.\vspace{2mm}

	{\scriptsize
		\begin{equation*}\begin{array}{lllllllllll}
		\hline\\
		n & h_n & k_n & \lambda_n&\mbox{Remarks}&~~~~ &n & h_n & k_n & \lambda_n&\mbox{Remarks}\vspace{1mm}\\
		\hline\\
		
		1&3&3&0.3226&\mbox{tbd}&&2&8&5&0.3183&(11.5,6)\vspace{1mm}\\
		
		3&5&4&0.3074&\mbox{tbd}&&4&11&6&0.3034&(15,7)\vspace{1mm}\\
		
		5&24.5&9&0.3031&\mbox{na}&&6&7.5&5&0.2971&\mbox{tbd}\vspace{1mm}\\
		
		7&14.5&7&0.2945&(24,9)&&8&29.5&10&0.2943&\mbox{na}\vspace{1mm}\\
		
		9&85&17&0.2939&\mbox{tbd}&&&&&&\vspace{1mm}\\

		\hline
		\end{array}
		\end{equation*}\vspace{1mm}}
	
	\noindent  Thus (\ref{eq101}) is soluble unless $A=1/3,$ $5/16,$ $3/10,$ $5/17,$  and $\beta_{A,F} = \pm a_4'-2A(c_4+\lambda x_5)\equiv h_n/k_n \pmod {1/k_n, 2A}$. This gives $A=1/3, \lambda=0, 1/2$; $A=5/16, \lambda=0, 1/5, 2/5$; $A=3/10, \lambda=0,1/3$ and $A=5/17, \lambda=0, 1/10, 2/10, 3/10, 4/10, 5/10$. But from (\ref{eq8}) and (\ref{eq09}), we have \vspace{-2mm}
	$$1/2\geq \lambda \geq \sqrt{1-t/A},\vspace{-2mm}$$ which discards all the values of $A$ and $\lambda$ except $A=5/17, \lambda=2/5, 1/2$ and $A=1/3, \lambda=1/2$.\vspace{2mm}
	
	\noindent For $A=5/17, \lambda=2/5$, $t=\frac{1}{4}$, we have $C=A\lambda^2+t=\frac{101}{340}$. We find that (\ref{eq45}) is satisfied on taking $(h_C,k_C)=(24,9)$. Therefore by Macbeath's Lemma, (\ref{eq101}) is soluble.\vspace{2mm}
	
	\noindent For $A=5/17, \lambda=1/2$, $t=\frac{1}{4},$ we have $C=\frac{11}{34}$ and $\frac{At}{C}= \frac{5}{22}$. Here $1< \frac{\delta_m}{C}<2$, as $\delta_m=d+\frac{3}{4}a\leq 2A$. Therefore $M=1$, see Remark \ref{rem1}. On taking $(h^*,k^*)=(2,3)$, we find that (\ref{eq45'}) is satisfied. Then by Macbeath's Lemma, (\ref{e17}) and hence (\ref{eq101}) is soluble.\hfill $\square$
	
	\begin{lemma}{\label{prob}}
		When $t=1/4,$ $A=1/3$, $(a_5'',c_5)=(0,0)$ and $\lambda=1/2$, $(\ref{eq101})$ and hence $(\ref{eq10})$ is soluble with strict inequality only when $d+a=d+6d^5>\frac{2}{3}$.
	\end{lemma}
	
	\noindent {\bf Proof:} For $t=1/4$ and $A=1/3$ we get  $F=(x_1+a_2'x_2+a_4'x_4+a_5'x_5)x_2-\frac{1}{3}(x_4^2+x_4x_5+x_5^2)-\frac{a}{4}$, which is symmetric in $x_4$ and $x_5$. As $c_5=0$, by symmetry, we can take  $c_4=0$. Here $h_A=k_A=3$. Also, $\beta_{A,F}=\pm a_4'-2A(c_4+\frac{1}{2} x_5)\equiv h_A/k_A \pmod{1/k_A,2A}$. Taking $c_4=0=x_5$, we get $3a_4'\equiv 0\pmod 1$ which gives $a_4'=0$ or $1/3$. By symmetry, $a_5'=0$ or $1/3$. From (\ref{eq16}), we have $0=a_5''\equiv a_5'-a_4'/2\pmod 1$. Also $|a_5'-a_4'/2|\leq 3/4$ so that we must have $a_4'=0=a_5'$. Note that $t=d^5/2aA\leq 9d^4/2$ gives $d\geq (2t/9)^{1/4}=(1/18)^{1/4}>0.48549.$ Inequality (\ref{eq101}) reduces to
	\begin{equation}\label{probeq}
	0<F=(x_1+a_2'x_2)x_2-\frac{1}{3}(x_4^2+x_4x_5+x_5^2)-\frac{a}{4}<d+\frac{3a}{4}.
	\end{equation}
	
	\noindent But (\ref{probeq}) has no solution for $(a_2',c_1)=(0,0)$ when $d+a=d+6d^5\leq\frac{2}{3}$. This is so because $3x_1x_2-(x_4^2+x_4x_5+x_5^2)$ never takes the value $1$ for integers $x_1,x_2,x_4,x_5$, and so $0< \frac{a}{4}< x_1x_2-\frac{1}{3}(x_4^2+x_4x_5+x_5^2)< d+a\leq \frac{2}{3}$ has no solution for integers $x_1,x_2,x_4,x_5$. Therefore we assume that $d+a=d+6d^5>\frac{2}{3}$, which gives $d \geq 0.4927$.\vspace{2mm}
	
	\noindent The following table gives a solution to (\ref{probeq}),\vspace{2mm}
	
	{\scriptsize
		\begin{equation*}\begin{array}{lllllll}
		\hline\\
		{\rm ~~~~Range} & x_1 & x_1x_2&x_4 & x_5&x_4x_5 & {\rm ~~F}\vspace{2mm}\\
		\hline\\

		\frac{1}{3}+d+a\leq f_1'\leq 1&\pm 1+c_1&-2+2|c_1|&1&0&0&f_2'-\frac{7}{3}-\frac{a}{4}\vspace{2mm}\\
		
		d+a\leq f_1'<\frac{1}{3}+d+a&c_1&|c_1|&1&0&0&f_1'-\frac{1}{3}-\frac{a}{4}\vspace{2mm}\\
		
		\frac{a}{4}<f_1'<d+a&c_1&|c_1|&0&0&0&f_1'-\frac{a}{4}\vspace{2mm}\\
		
		-\frac{1}{2}<f_1'<-\frac{2}{3}+d+a&\pm 1+c_1&1+|c_1|&1&0&0&f_1'+\frac{2}{3}-\frac{a}{4}\vspace{2mm}\\
		{\rm When ~}-\frac{2}{3}+d+a\leq f_1' \leq \frac{a}{4} : \vspace{2mm}\\
		
		-\frac{2}{3}+\frac{a}{4}<g_1'<-\frac{2}{3}+d+a&\pm 1+c_1&1-|c_1|&1&0&0&g_1'+\frac{2}{3}-\frac{a}{4}\vspace{2mm}\\
		
		-1+\frac{a}{4}<g_1'\leq-\frac{2}{3}+\frac{a}{4}&\pm 1+c_1&1-|c_1|&0&0&0&g_1'+1-\frac{a}{4}\vspace{2mm}\\
		
		\hline
		\end{array}
		\end{equation*}\vspace{1mm}}

	\noindent except in the following two cases:
	\begin{enumerate}
		\item $-\frac{2}{3}+d+a\leq f_1'\leq \frac{a}{4}$,~ $-\frac{2}{3}+d+a\leq g_1'\leq \frac{a}{4}$
		\item $-\frac{2}{3}+d+a\leq f_1'\leq \frac{a}{4}, ~-1< g_1'\leq -1+\frac{a}{4}$.	
	\end{enumerate}
	
	\noindent \textbf{Case 1:} $-\frac{2}{3}+d+a\leq f_1',g_1'\leq \frac{a}{4}$.\vspace{2mm}
	
	\noindent Here as $g_1'\geq -\frac{2}{3}+d+a > 0$, we have
	$ g_1'=-|c_1|+a_2' > 0$. Which  gives $a_2' > 0$.\vspace{2mm}

	\noindent   Find an integer $n\geq 1$ such that \vspace{-2mm}
	\begin{equation}\label{E62}
	n|c_1|+n^2a_2'\leq \frac{a}{4}<(n+1)|c_1|+(n+1)^2a_2'.
	\vspace{-2mm}\end{equation}
	
	\noindent Now take $x_2=\pm (n+1), x_1=c_1,x_1x_2=(n+1)|c_1|,~ x_4=x_5=0$,  so that \vspace{-2mm}
	$$F=(n+1)|c_1|+(n+1)^2a_2'-\frac{a}{4}.\vspace{-2mm}$$
	
	\noindent Using (\ref{E62}), we find that $F>0$. And
	\begin{equation*}\begin{array}{lll}
	F&=&(n+1)|c_1|+(n+1)^2a_2'-\frac{a}{4}\vspace{2mm}\\
	&=&\frac{(n+1)^2}{n^2}.(\frac{n^2}{n+1}|c_1|+n^2a_2')-\frac{a}{4}\vspace{2mm}\\
	&\leq &\frac{(n+1)^2}{n^2}.(n|c_1|+n^2a_2')-\frac{a}{4}\vspace{2mm}\\
	&\leq & 4(\frac{a}{4})-\frac{a}{4}
	=3(\frac{a}{4})<d+\frac{3a}{4}.
	\end{array}
	\end{equation*}
	\noindent \textbf{Case 2:} $-\frac{2}{3}+d+a\leq f_1'\leq \frac{a}{4}, ~-1< g_1'\leq -1+\frac{a}{4}$. \vspace{2mm}
	
	\noindent Note that here, $ 0<1+g_1'\leq \frac{a}{4}$. Find an integer $n\geq 1$ such that\vspace{-2mm}
	
	{ \begin{equation}\label{E60}\begin{array}{l}
		(2n-1)n+g_{2n-1}'\leq \frac{a}{4}<(2n+1)(n+1)+g_{2n+1}'.
		\end{array}
		\end{equation}}
	
	\noindent Further $g_{2n-1}'=-(2n-1)|c_1|+(2n-1)^2a_2'\leq -(2n-1)n+\frac{a}{4}$ gives \begin{equation}\label{E61}
	a_2'\leq -\frac{1}{2}+\frac{a}{4}.\frac{1}{(2n-1)^2}.
	\end{equation}
	\noindent   Take $x_2=\pm (2n+1),~ x_1=\pm (n+1)+c_1, ~x_1x_2=(2n+1)(n+1)-(2n+1)|c_1|,~x_4=x_5=0$. So that,
	$$F=(2n+1)(n+1)-(2n+1)|c_1|+(2n+1)^2a_2'-\frac{a}{4}.$$
	
	\noindent Using (\ref{E60}), we find that $F>0$. And from (\ref{E61}), we find that
	\begin{equation*}\begin{array}{lll}
	 F&=&(2n+1)(n+1)+\frac{2n+1}{2n-1}\big(-(2n-1)|c_1|+(2n-1)^2a_2'~\big)\\&&-(2n+1)(2n-1)a_2'+(2n+1)^2a_2'
	-\frac{a}{4}\vspace{2mm}\\
	&\leq& (2n+1)(n+1)+\frac{2n+1}{2n-1}\big(-(2n-1)n+\frac{a}{4}\big)\\&&+2(2n+1)\big(-\frac{1}{2}+\frac{a}{4}
	\frac{1}{(2n-1)^2}\big)-\frac{a}{4}\vspace{2mm}\\
	 &=&\big(\frac{2n+1}{2n-1}+2.\frac{2n+1}{(2n-1)^2}-1\big)\frac{a}{4}=\frac{8n}{(2n-1)^2}\frac{a}{4}\vspace{2mm}\\
	&\leq & 8.\frac{a}{4}<d+\frac{3a}{4}\vspace{2mm}\\
	\end{array}
	\end{equation*}
	\noindent as $\frac{8n}{(2n-1)^2}$ is a decreasing function of $n$ and $n\geq 1$. \vspace{2mm}\hfill $\square$

	\noindent The proof of Theorem \ref{thm6} follows from Lemmas \ref{lem21} - \ref{prob}.\vspace{2mm}
	
	\noindent Thus from the arguments given in the beginning of this section for $(m,K)=(2,1)$ we have proved that $\Gamma_{1,4}<8.486$ instead of $\Gamma_{1,4}<8$.

	\section{Proof of $\Gamma_{1,4}< \frac{32}{3}$ when $(m,K)=(2,2)$}\label{sec9}
	\numberwithin{equation}{section}
	
	If we start  with $d=(8|D|)^{\frac{1}{5}}$ and proceed as in previous sections, we find that Macbeath's Lemma is not applicable for $m=2, K=2$ and $L=4$. This is so because for $d=(8|D|)^{\frac{1}{5}}$, we have
	$a\geq \frac{d}{3}, A\geq \frac{5d}{12}$, $\frac{\delta_{2,2}}{t}>\frac{125}{288d^2}$ and $L<\frac{\delta_{2,2}}{t}\leq L+1$ which gives $d^2\geq \frac{125}{288}\frac{1}{L+1}$.
	It turns out that  $1<\frac{\delta_{2,2}}{t}\leq 5$  which gives $L=1, 2, 3$ and $4$.
	When $L=4$, we get $d>0.2946$, $\delta_2=d+\frac{3}{4}a\geq \frac{5}{4}d> 0.3628$ and
	$\delta_{2,2}= d+\frac{3}{4}a+\frac{3}{4}A\geq \frac{25}{16}d>0.46035$.
	Thus, for $m=2, K=2, L=4$, we can not apply Macbeath Lemma \ref{lem5} either on $G$ or on $F$ or on $Q$. Therefore, we have to work with a constant bigger than $8$. We prove here $\Gamma_{1,4}< \frac{32}{3}$. For $(m,K)=(2,2)$, one may prove the result for a constant slightly smaller than $32/3$ also, but we find that when $(m,K)=(1,1)$, $\delta_m=\delta_{m,K}=d$ turns out to be bigger than $1/2$ only if the constant is $32/3$. So, in the remaining sections, we shall prove that $\Gamma_{1,4}< \frac{32}{3}$.\vspace{2mm}
	
	\noindent  Here, instead of (3.1), (3.5) and (3.6) we have
	\begin{equation}\label{eq7.1} \frac{d}{m+1}\leq a<(2\Delta)^{\frac{1}{3}}=(8|D|)^{\frac{1}{3}}=(3d^5/4)^{\frac{1}{3}}.\end{equation}
	\begin{equation}\label{key2} \frac{(m+3)d}{4(K+1)}\leq A\leq \sqrt{\frac{16|D|}{3a}} = \sqrt{\frac{d^5}{2a}}, ~~ t=\frac{3d^5}{8aA}.\vspace{2mm} \end{equation}
	
	\noindent  Further as $\delta_{m,K}\geq \frac{(m+3)}{4} \frac{(K+3)}{4}d$,  we have,  \begin{equation}\label{Eq55'} \begin{array}{ll}
	\frac{\delta_{m,K}}{t}&=\frac{8aA}{3d^5} \delta_{m,K}
	\geq \frac{8}{3d^5} \frac{d}{m+1}\frac{(m+3)d}{4(K+1)}\frac{(m+3)}{4} \frac{(K+3)}{4}d\vspace{4mm}\\
	&=\frac{1}{24}\frac{(m+3)^2(K+3)}{(m+1)(K+1)}\frac{1}{d^2}.
	\end{array}
	\end{equation}

	\begin{theorem}\label{thm5}
		Let $Q(x_1,x_2, \cdots, x_5)$ be a real indefinite quadratic form of type $(1,4)$ and of determinant $D\neq 0$. Let $d=(\frac{32}{3}|D|)^{\frac{1}{5}}$.  Suppose that $c_2 \equiv 0 \pmod1$, $a<\frac{1}{2}$, $d\leq 1$ and $a+d\leq 1$. Let $(m,K)= (2,2)$.  Then $(\ref{eq101*})$ and hence $(\ref{eq10})$ is soluble with strict inequality.
	\end{theorem}
	
	\noindent \textbf{Proof:} Here, $2 < \frac{d}{a}\leq 3$ and $2<\frac{\delta_{2}}{A}=\frac{d+3a/4}{A} \leq 3$. Since, $\frac{d}{3}+d \leq a+d \leq 1$ we get $d \leq \frac{3}{4}$. Also $d/3 \leq a\leq (3d^5/4)^{1/3}$. Now, $\frac{5d}{12}=\frac{(m+3)d}{4(K+1)} \leq A \leq \sqrt{\frac{d^5}{2a}} \leq \sqrt{\frac{d^4.3}{2}} = \sqrt{\frac{3}{2}}d^2$ gives $d \geq \frac{5\sqrt{6}}{36} > 0.3402$. This together with $d\leq \frac{16}{(m+3)(K+3)}$, gives
	\begin{equation}\label{EqN1}
	0.3402< d \leq 16/25=0.64.
	\vspace{1mm}\end{equation}
	
	\noindent Further we have, from (\ref{Eq55'}), $\frac{\delta_{m,K}}{t}\geq \frac{1}{24}\frac{(m+3)^2(K+3)}{(m+1)(K+1)}\frac{1}{d^2} > 1$,
	for $m=2,K=2$. Also $\frac{\delta_{2,2}}{t}=\frac{d+3a/4+3A/4}{t}\leq \frac{3A+3A/4}{3A/4}=5$  which gives $L=1, 2, 3$ and $4$.  Now as $L<\frac{\delta_{2,2}}{t}\leq L+1$, from (\ref{Eq55'}), we get $d^2\geq \frac{1}{24}\frac{(m+3)^2(K+3)}{(m+1)(K+1)(L+1)}$. Therefore using (\ref{eq26'}), we have
	
	\begin{equation}\label{Eq74}
	\delta_{2,2,L}\leq 1,{\rm ~and~} \sqrt{125/216(L+1)} \leq d \leq 64/25(L+3).\vspace{2mm}
	\end{equation}
	
	\noindent We distinguish the cases $L=1,2,3$ and $4$ in the following subsections.
	
	\subsection{${m=2, K=2, L=4}$}
	\noindent Here $4<\frac{\delta_{m,K}}{t}=\frac{d+3a/4+3A/4}{t}\leq 5$, $\delta_{m,K}=d+\frac{3a}{4}+\frac{3A}{4}\geq d+\frac{d}{4}+\frac{5d}{16}=\frac{25d}{16}$. From (\ref{Eq74}) we have
	\begin{equation}\label{Eq79}
	0.3402< \sqrt{25/216} \leq d \leq 64/175<0.365715.
	\end{equation}

	\noindent Further, since $a\geq d/3$ and $A\geq 5d/12$ we get, $$0.106312< \frac{5d}{16}\leq \frac{d+3a/4+3A/4}{5}\leq t =\frac{3d^5}{8aA}\leq \frac{27d^3}{10}< 0.132067.$$
	
	\noindent Also, $\delta_{2,2,4}=d+\frac{3a}{4}+\frac{3A}{4}+\frac{15t}{4} \leq 1$ gives $t \leq \frac{4}{15}\left(1-\frac{25d}{16}\right)$. Thus, \begin{equation*}
	t \leq \min \left\{ \frac{27d^3}{10},\frac{4}{15}\left(1-\frac{25d}{16}\right) \right\} = \left\{
	\begin{array}{lll}
	&\frac{27d^3}{10} &{\rm~if~} d\leq \alpha\vspace{2mm}\\
	& \frac{4}{15}\left(1-\frac{25d}{16}\right) &{\rm~if~} d> \alpha
	\end{array}
	\right.
	\end{equation*}
	\noindent where $\alpha$ is a root of $\frac{27d^3}{10}-\frac{4}{15}\left(1-\frac{25d}{16}\right)=0$ satisfying $0.353<\alpha<0.354$. This gives $0.106312 < t\leq 0.1194...=\lambda_0(say)$. We divide this range of $t$ into $14$ subintervals $[\lambda_n,\lambda_{n-1}], n=1,2,\cdots 14$ and in each subinterval choose suitable integers $2h_n$ and $k_n$ such that $\vline~ h_n-tk_n^2~\vline +\frac{1}{2} <\frac{25d}{16}\leq d+\frac{3a}{4}+\frac{3A}{4}$. The choice of $(h_n,k_n)$ is done in a manner similar to that in the Case  $m=2,K=3, L=4$.\vspace{2mm}

	{\scriptsize
		\begin{equation*}\begin{array}{lllllllll}
		\hline\\
		n & (h_n,k_n) & \lambda_n&\mbox{Remarks}&~~ &n & (h_n,k_n) & \lambda_n&\mbox{Remarks}\vspace{2mm}\\
		\hline\\
		
		1&(0.5,2)&0.114&\mbox{na} &&2&(1,3)&0.1075&\mbox{tbd} \vspace{2mm}\\
		
		3&(13,11)&0.10717&(15.5,12) &&4&(21,14)&0.10698&\mbox{tbd} \vspace{2mm}\\
		
		5&(78,27)&0.106952&\mbox{na} &&6&(109.5,32)&0.106902&(131,35) \vspace{2mm}\\	
		
		7&(216.5,45)&0.106898&\mbox{na} &&8&(2247.5,145)&0.106896&\mbox{tbd}\vspace{2mm}\\
		
		9&(1178.5,105)&0.106895&\mbox{na} &&10&(171,40)&0.106866&(138.5,36) \vspace{2mm}\\
		
		11&(56.5,23)&0.10681&\mbox{na} &&12&(24,15)&0.106699&\mbox{na}\vspace{2mm}\\
		
		13&(18,13)&0.106322&(492.5,68) &&14&(42.5,20)&0.106314&\mbox{na}\vspace{2mm}\\
		
		\hline
		\end{array}
		\end{equation*}\vspace{1mm}}
	
	\noindent  Thus (\ref{eq101'}) is soluble unless $t=1/9,$ $3/28, 31/290$ and $\beta_{t,G} = \pm a_5''-2tc_5\equiv h_n/k_n \pmod {1/k_n, 2t}$. For each of these values to $t$, we find that (\ref{eq101'}) is soluble by taking $x_2=\pm 1$ and suitably chosing $x_5\equiv c_5 \pmod1$ so that range of $f_1''$ or $p_1''$ is covered.

	\subsection{${m=2, K=2, L=3}$}
	
	\noindent Here $3<\frac{\delta_{m,K}}{t}\leq 4$, $\delta_{m,K}=d+\frac{3a}{4}+\frac{3A}{4}\geq d+\frac{d}{4}+\frac{5d}{16}=\frac{25d}{16}$. From (\ref{Eq74}) we have
	\begin{equation}\label{Eq'79}
	0.38036 =\sqrt{125/864} \leq d \leq 32/75<0.42667.
	\end{equation}

	\noindent Further, since $a\geq d/3$ and $A\geq 5d/12$ we get, $$0.14857< \frac{25d}{64}\leq \frac{d+3a/4+3A/4}{4}\leq t =\frac{3d^5}{8aA}\leq \frac{27d^3}{10}< 0.20973.$$
	
	\noindent Also, $\delta_{2,2,3}=d+\frac{3a}{4}+\frac{3A}{4}+2t \leq 1$ gives $t \leq \frac{1}{2}\left(1-\frac{25d}{16}\right)$. Thus, \begin{equation*}
	t \leq \min \left\{ \frac{27d^3}{10},\frac{1}{2}\left(1-\frac{25d}{16}\right) \right\} = \left\{
	\begin{array}{lll}
	&\frac{27d^3}{10} &{\rm~if~} d\leq \alpha\vspace{2mm}\\
	& \frac{1}{2}\left(1-\frac{25d}{16}\right) &{\rm~if~} d> \alpha
	\end{array}
	\right.
	\end{equation*}
	\noindent where $\alpha$ is a root of $\frac{27d^3}{10}-\frac{1}{2}\left(1-\frac{25d}{16}\right)=0$ satisfying $0.407<\alpha<0.408$. This gives $0.14857< t\leq 0.182032...=\lambda_0(say)$.\vspace{2mm}

	\noindent We divide this range of $t$ into $5$ subintervals $[\lambda_n,\lambda_{n-1}], n=1,2,\cdots 5$ and in each subinterval choose suitable integers $2h_n$ and $k_n$, in the  way already described, satisfying $\vline~ h_n-tk_n^2~\vline +\frac{1}{2} < \frac{25d}{16}\leq d+\frac{3a}{4}+\frac{3A}{4}.$\vspace{4mm}
	
	{\scriptsize
		\begin{equation*}\begin{array}{lllllllll}
		\hline\\
		n & (h_n,k_n) & \lambda_n&\mbox{Remarks}&~~ &n & (h_n,k_n) & \lambda_n&\mbox{Remarks}\vspace{2mm}\\
		\hline\\
		
		1&(3,4)&0.18&\mbox{na}&&2&(1.5,3)&0.156&\mbox{tbd}\vspace{2mm}\\
		
		3&(2.5,4)&0.1503&\mbox{na}&&4&(15,10)&0.1499&\mbox{tbd}\vspace{2mm}\\
		
		5&(9.5,8)&0.14857&\mbox{na}&&&&&\vspace{2mm}\\	
		
		\hline
		\end{array}
		\end{equation*}\vspace{2mm}}
	
	\noindent  Thus (\ref{eq101'}) is soluble unless $t=1/6,$ $3/20$ and $\beta_{t,G} = \pm a_5''-2tc_5\equiv h_n/k_n \pmod {1/k_n, 2t}$.  When $t=3/20$, by Lemma \ref{lem15} we need to consider $(a_5'',c_5)=(0,0), (0,1/3)  {\rm ~and~} (1/20,1/6).$ Here, we find that (\ref{eq101'}) is soluble by taking $x_2=\pm 1$ and suitably choosing $x_5 \equiv c_5 \pmod 1$. The case $t=1/6$ is involved, we discuss it in next lemma.
	
	\begin{lemma} When $t=\frac{1}{6}$ and  $m=2, K=2, L=3$, $(\ref{eq101})$ is soluble. \end{lemma}
	
	\noindent \textbf{Proof:}
	Here, since  $\frac{1}{6}=t=\frac{3d^5}{8aA}\leq \frac{27d^3}{10}$, we get $d^3\geq \frac{5}{81}$ i.e., $d\geq (\frac{5}{81})^{\frac{1}{3}}=0.39521$. Also, $\frac{1}{6}=t=\frac{3d^5}{8aA}$ gives $A=\frac{9d^5}{4a}\leq \frac{27d^4}{4}$. And from (\ref{Eq'79}), we have $d\leq \frac{32}{75}$.\vspace{2mm}
	
	\noindent Therefore,

	\begin{equation}\label{Eq'131}
	\begin{array}{l}
	0.131736< d/3\leq a\leq d/2< 0.213335\vspace{2mm}\\
	
	0.1646708< 5d/12 \leq A \leq 4t/3< 0.222223\vspace{2mm}\\
	0.074101< \frac{a+A}{4} \leq 0.10889.
	\end{array}\vspace{2mm}
	\end{equation}

	\noindent Using Lemma \ref{lem15} we need to consider $(a_5'',c_5)= (0,1/2), {\rm ~and~} (1/6,0).$ If $a_5''=1/6,$ considering the transformation $x_5\rightarrow x_5-x_2$, we may assume that $a_5''=1/2$.\vspace{2mm}

	\noindent When $(a_5'',c_5)=(0,1/2)$, the following table gives a solution to (\ref{eq101'}) i.e., of\\
	 \begin{equation}\label{EQ'55}0<G=(x_1+a_2''x_2)x_2-\frac{1}{6}x_5^2-\frac{a+A}{4}<d+\frac{3a}{4}+\frac{3A}{4}.\end{equation}
	
	{\scriptsize
		\begin{equation*}\begin{array}{lllll}
		\hline\\
		{\rm ~~~~Range} & x_1 & x_1x_2 & x_5& {\rm ~~G}\vspace{2mm}\\
		\hline\\
		
		\frac{3}{8}+\frac{a+A}{4}<f_1''\leq 1&c_1&|c_1|&\frac{3}{2} &f_1''-\frac{3}{8}-\frac{a+A}{4}\vspace{2mm}\\
		
		\frac{1}{24}+\frac{a+A}{4}<f_1''\leq \frac{3}{8}+\frac{a+A}{4}&c_1&|c_1|&\frac{1}{2} &f_1''-\frac{1}{24}-\frac{a+A}{4}\vspace{2mm}\\
		
		-\frac{1}{2}<f_1''<-\frac{5}{8}+d+a+A&\pm 1+c_1&1+|c_1|&\frac{3}{2} &f_1''+\frac{5}{8}-\frac{a+A}{4}\vspace{2mm}\\

		\hline
		\end{array}
		\end{equation*}\\}
	
	{\scriptsize
		\begin{equation*}\begin{array}{lllll}
		\hline\\
		{\rm ~~~~Range} & x_1 & x_1x_2 & x_5& {\rm ~~G}\vspace{2mm}\\
		\hline\\
		
		{\rm When~}  -\frac{5}{8}+d+a+A\leq f_1''\leq \frac{1}{24}+\frac{a+A}{4}: \vspace{2mm}\\
		
		-\frac{5}{8}+d+a+A\leq g_1''\frac{1}{24}+\frac{a+A}{4}&c_1&2|c_1|&\frac{1}{2}&f_2''-\frac{1}{24}-\frac{a+A}{4}\vspace{2mm}\\
		
		-\frac{5}{8}+\frac{a+A}{4}<g_1''<-\frac{5}{8}+d+a+A&\pm 1+c_1&1-|c_1|&\frac{3}{2} &g_1''+\frac{5}{8}-\frac{a+A}{4}\vspace{2mm}\\
		
		-\frac{23}{24}+\frac{a+A}{4}<g_1''\leq -\frac{5}{8}+\frac{a+A}{4}&\pm 1+c_1&1-|c_1|&\frac{1}{2} &g_1''+\frac{23}{24}-\frac{a+A}{4}\vspace{2mm}\\
		
		-1<g_1''\leq -\frac{23}{24}+\frac{a+A}{4}&\pm 1+c_1&2+2|c_1|&\frac{5}{2}&f_2''+\frac{23}{24}-\frac{a+A}{4}\vspace{2mm}\\
		
		\hline
		\end{array}
		\end{equation*}\\}
	\noindent When $(a_5'',c_5)=(1/2,0)$, the following table gives a solution to (\ref{eq101'}) i.e., of
	 \begin{equation}\label{EQ'56}0<G=(x_1+a_2''x_2+\frac{1}{2}x_5)x_2-\frac{1}{6}x_5^2-\frac{a+A}{4}<d+\frac{3a}{4}+\frac{3A}{4}.\end{equation}

	{\scriptsize
		\begin{equation*}\begin{array}{llllll}
		\hline\\
		{\rm ~~~~Range} & x_1 & x_1x_2 & x_5&x_2x_5& {\rm ~~G}\vspace{2mm}\\
		\hline\\

		\frac{a+A}{4}<g_1''\leq \frac{1}{2}&c_1&-|c_1|&0&0 &g_1''-\frac{a+A}{4}\vspace{2mm}\\
		
		-\frac{1}{3}+\frac{a+A}{4}<g_1''\leq \frac{a+A}{4}&c_1&-|c_1|&\pm 1 &1&g_1''+\frac{1}{3}-\frac{a+A}{4}\vspace{2mm}\\
		
		-1+\frac{a+A}{4}<g_1''<-1+d+a+A&\pm 1+c_1&1+|c_1|&0 &0&g_1''+1-\frac{a+A}{4}\vspace{2mm}\\
		
		-1<g_1''\leq -1+\frac{a+A}{4}&\pm 1+c_1&1-|c_1|&\pm 1 &1&g_1''+\frac{4}{3}-\frac{a+A}{4}\vspace{2mm}\\

		{\rm When~}  -1+d+a+A\leq g_1''\leq -\frac{1}{3}+\frac{a+A}{4}: \vspace{2mm}\\
		
		\frac{2}{3}+\frac{a+A}{4}<f_1''\leq 1&\pm 1+c_1&-1+|c_1|&\pm 1&1 &f_1''-\frac{2}{3}-\frac{a+A}{4}\vspace{2mm}\\
		
		\frac{a+A}{4}<f_1''<d+a+A&c_1&|c_1|&0&0 &f_1''-\frac{a+A}{4}\vspace{2mm}\\
		
		-\frac{1}{3}+\frac{a+A}{4}<f_1''\leq \frac{a+A}{4}&c_1&|c_1|&\pm 1 &1&f_1''+\frac{1}{3}-\frac{a+A}{4}\vspace{2mm}\\
		
		-\frac{1}{2}<f_1'' \leq -\frac{1}{3}+\frac{a+A}{4}&c_1&2|c_1|&\pm 2&4&f_2''+\frac{4}{3}-\frac{a+A}{4}\vspace{2mm}\\		
		
		\hline
		\end{array}
		\end{equation*}}
	
	\noindent except when $d+a+A \leq f_1'' \leq \frac{2}{3}+\frac{a+A}{4}$ and $-1+d+a+A \leq g_1'' \leq -\frac{1}{3}+\frac{a+A}{4}$.\vspace{2mm}
	
	\noindent In this exceptional case, we work   in a manner similar to that in Lemma 9 of \cite{RakaRani} by considering the cases $a_2''\leq \frac{1}{4}$ and $a_2''> \frac{1}{4}$ separately and  find that (\ref{eq101'}) is soluble unless $c_1=1/2, a_2''=1/4$ and $d+a+A \leq 3/4, a+A \geq 1/3$ which gives $d \leq 5/12$.\vspace{2mm}
	
	\noindent Here, from (\ref{eq16}) we have
	\begin{equation*}
	\begin{array}{l}
	1/4=a_2''\equiv a_2'+a_4'^2/4A ~({\rm mod~} 1), 1/2=a_5''\equiv a_5'-\lambda a_4'~({\rm mod~} 1) {\rm ~~ and}\vspace{2mm}\\
	
	-1/2<a_2'\leq a_2'+a_4'^2/4A\leq a_2'+1/16A<1, -3/4\leq a_5'-\lambda a_4'\leq 3/4.
	\end{array}
	\end{equation*}
	
	\noindent So we have
	\begin{equation}\label{One}
	a_2'+a_4'^2/4A=1/4 \mbox{~and~} a_5'-\lambda a_4' =1/2 \mbox{~or~} -1/2.
	\end{equation}
	
	\noindent Since here $0.39521 < d$, we have
	$$2<2.223< 45d/8 \leq 9\delta_{2}/2 \leq \delta_{2}/C=(d+3a)/C \leq \delta_{2}/A \leq 3$$
	so that $M=2$, (see Remark \ref{rem1}) and also $A \geq 5d/12$ gives
	$$1/8<A/(3A/2+1)\leq A/(6A\lambda^2+1)=A/6C=At/C\leq 1/6.$$

	\noindent Taking $h^*=3/2, k^*=3$ for $A/6C>(2-\delta_{2,2}^*)/9$ and $h^*=1/2, k^*=2$ for $A/6C\leq (2-\delta_{2,2}^*)/9$ and noting that $\delta_{2,2}^*\geq \delta_{2,2} \geq 25d/16> 0.61751$ we find that (\ref{eq45'}) is satisfied. Therefore, (\ref{e17}) and hence (\ref{eq101}) is soluble unless $A/6C=1/6$ i.e., $A=C$ and $\beta_{G}^*\equiv 1/2 \pmod {1/3}$ i.e., $\pm 3a_4^*-c_4 \equiv 1/2 \pmod 1.$ From Lemma \ref{lem15}, we need to consider $(a_4^*,c_4)=(1/2,0)$ or $(0,1/2).$ Proceeding in a similar way  to that of Lemma 9 of \cite{RakaRani}, now working with $G^*$ instead of $G$, we find that (\ref{e17}) has a solution unless $(a_4^*,c_4)=(1/2,0)$ and $a_2^*=1/4.$ From (\ref{G*}), we get
	\begin{equation}\label{Two}
	a_2'+a_5'^2/4A=1/4 \mbox{~and~} a_4'-\lambda a_5'= 1/2 \mbox{~or~} -1/2.
	\end{equation}
	
	\noindent From (\ref{One}) and (\ref{Two}), we get that $a_4'^2=a_5'^2,$ so $a_4'=\pm a_5'.$ As $A=C=A\lambda^2+1/6$, we get $A=1/6(1-\lambda^2).$ Further $d+a+A  \leq \frac{2}{3}+\frac{a+A}{4}$, gives $A\leq \frac{4}{3}(\frac{2}{3}-\frac{5}{4}d)$. So $0.1646\leq 5d/12 \leq A \leq \mbox{min}\{27d^4/4,\frac{4}{3}(\frac{2}{3}-\frac{5}{4}d)\}\leq 0.1986$ which gives $0\leq \lambda \leq 0.401.$\vspace{4mm}
	
	\noindent When $a_4'=a_5',$ we get from (\ref{Two}) that $a_5'=\pm 1/2(1-\lambda)$ which gives $|a_5'|\geq 1/2$ and hence $a_5'=\pm 1/2$. Therefore, $a_4'=a_5'=1/2$ or $a_4'=a_5'=-1/2$. Taking the unimodular transformation, $x_2 \rightarrow -x_2, x_4 \rightarrow -x_4, x_5 \rightarrow -x_5$, we need to consider $a_4'=a_5'=1/2$ only.\vspace{2mm}
	
	\noindent For $a_4'=a_5'=1/2$, (\ref{Two}) gives $\lambda=0$. Also, $1/4=a_2'+a_5'^2/4A$ gives $a_2'=-1/8$. Then, $x_1=1/2, x_2=1, x_4=x_5=0$ gives a solution of (\ref{eq101}).\vspace{4mm}
	
	\noindent When $a_4'=-a_5',$ we get from (\ref{G*}) that $a_5'=\pm 1/2(1+\lambda).$ So $a_2'=1/4-a_5'^2/4A=(5\lambda -1)/8(1+\lambda)$. When $0 \leq \lambda \leq 0.254$, $x_1=1/2, x_2=1, x_4=0=x_5$ gives a solution of (\ref{eq101}); and when $0.254 \leq \lambda \leq 0.401$, $x_1=1/2, x_2=1, x_4= \pm 1, x_5=\mp 1$ gives a solution of (\ref{eq101}).
	
	\subsection{${m=2, K=2, L=2}$}
	
	\noindent Here $2<\frac{\delta_{m,K}}{t}\leq 3$, $\delta_{m,K}=d+\frac{3a}{4}+\frac{3A}{4}\geq d+\frac{d}{4}+\frac{5d}{16}=\frac{25d}{16}$. From (\ref{Eq74}) we have
	\begin{equation}\label{Eq'80}
	0.4392< \sqrt{125/648} \leq d \leq 64/125=0.512.
	\end{equation}

	\noindent Further, since $a\geq d/3$ and $A\geq 5d/12$ we get, $$0.22875< \frac{25d}{48}\leq \frac{d+3a/4+3A/4}{3}\leq t =\frac{3d^5}{8aA}\leq \frac{27d^3}{10}< 0.36239.$$
	
	\noindent Also, $\delta_{2,2,2}=d+\frac{3a}{4}+\frac{3A}{4}+\frac{2t}{4} \leq 1$ gives $t \leq \frac{4}{3}\left(1-\frac{25d}{16}\right)$. Thus,
	
	\begin{equation*}
	t \leq \min \left\{ \frac{27d^3}{10},\frac{4}{3}\left(1-\frac{25d}{16}\right) \right\} = \left\{
	\begin{array}{lll}
	&\frac{27d^3}{10} &{\rm~if~} d\leq \alpha\vspace{2mm}\\
	& \frac{4}{3}\left(1-\frac{25d}{16}\right) &{\rm~if~} d> \alpha
	\end{array}
	\right.
	\end{equation*}
	\noindent where $\alpha$ is a root of $\frac{27d^3}{10}-\frac{4}{3}\left(1-\frac{25d}{16}\right)=0$ satisfying $0.488<\alpha<0.489$. This gives $0.22875 \leq t\leq 0.31518$.\vspace{3mm}
	
	\noindent Taking $h_t=1$ and $k_t=2$ and using (\ref{Eq'80}) we find that $| h_t-tk_t^2 |+\frac{1}{2}< \delta_{2,2}$ is satisfied. Thus, (\ref{eq101'}) is soluble  unless $t=1/4$ and $\beta_{t,G}=\pm a_5''-c_5 \equiv 1/2 \pmod 1.$ This gives $(a_5'',c_5)=(0,0)$ or $(1/4,1/2).$\vspace{3mm}
	
	\noindent Now, $1/4=t=3d^5/8aA\leq 27d^3/10$ gives $d\geq (5/54)^{1/3}>0.4524$. Also, $1/4=t=3d^5/8aA$ gives $A=3d^5/2a\leq 9d^4/2$. Hence,
	\begin{equation}\label{Eq86}
	\begin{array}{l}
	0.1508< d/3\leq a\leq d/2< 0.256\vspace{2mm}\\
	0.1885< 5d/12 \leq A \leq 9d^4/2< 0.309238\vspace{2mm}\\
	0.08482< \frac{a+A}{4} < 0.14131.
	\end{array}
	\end{equation}
	
	\noindent Thus we need to find the solutions of
	\begin{equation}\label{Eq87}
	0<G=(x_1+a_2''x_2+a_5''x_5)x_2-x_5^2/4-(a+A)/4<d+3a/4+3A/4
	\end{equation}
	for $(a_5'',c_5)=(0,0)$ and $(1/4,1/2)$.\vspace{2mm}
	
	\noindent When  $(a_5'',c_5)=(0,0)$, the following table gives a solution of (\ref{Eq87}).
	
	{\scriptsize
		\begin{equation*}\begin{array}{lllll}
		\hline\\
		{\rm ~~~~Range} & x_1 & x_1x_2 & x_5 & {\rm ~~G}\vspace{2mm}\\
		\hline\\
		
		\frac{1}{4}+\frac{a+A}{4}<f_1''\leq 1&c_1&|c_1|&\pm 1&f_1''-\frac{1}{4}-\frac{a+A}{4}\vspace{2mm}\\
		
		\frac{a+A}{4}<f_1''\leq \frac{1}{4}+\frac{a+A}{4}&c_1&|c_1|&0&f_1''-\frac{a+A}{4}\vspace{2mm}\\
		
		-\frac{1}{2}<f_1''<-\frac{3}{4}+d+a+A&\pm 1+c_1&1+|c_1|&\pm 1&f_1''+\frac{3}{4}-\frac{a+A}{4}\vspace{2mm}\\

		{\rm When~}  -\frac{3}{4}+d+a+A\leq f_1''\leq \frac{a+A}{4}: \vspace{2mm}\\
		
		-\frac{3}{4}+d+a+A\leq g_1''\leq \frac{a+A}{4}&c_1&2|c_1|&0&f_2''-\frac{a+A}{4}\vspace{2mm}\\
		
		-\frac{3}{4}+\frac{a+A}{4}<g_1''<-\frac{3}{4}+d+a+A&\pm 1+c_1&1-|c_1|&\pm 1&g_1''+\frac{3}{4}-\frac{a+A}{4}\vspace{2mm}\\
		
		-1+\frac{a+A}{4}<g_1''\leq -\frac{3}{4}+\frac{a+A}{4}&\pm 1+c_1&1-|c_1|&0&g_1''+1-\frac{a+A}{4}\vspace{2mm}\\
		
		-1<g_1'' \leq -1+ \frac{a+A}{4}&\pm 1+c_1&2+2|c_1|&\pm 2&f_2''+1-\frac{a+A}{4}\vspace{2mm}\\
		
		\hline
		\end{array}
		\end{equation*}\vspace{1mm}}
	
	\noindent When  $(a_5'',c_5)=(1/4,1/2)$, the following table gives a solution of (\ref{Eq87}).
	
	{\scriptsize
		\begin{equation*}\begin{array}{llllll}
		\hline\\
		{\rm ~~~~Range} & x_1 & x_1x_2 & x_5 &x_2x_5 & {\rm ~~G}\vspace{2mm}\\
		\hline\\
		
		\frac{3}{16}+d+a+A\leq f_1''\leq 1&\pm 1+c_1&-2+2|c_1|&\pm \frac{1}{2}&-1&f_2''-\frac{37}{16}-\frac{a+A}{4}\vspace{2mm}\\
		
		\frac{3}{16}+\frac{a+A}{4}<f_1''<\frac{3}{16}+d+a+A&c_1&|c_1|&\pm \frac{1}{2}&-\frac{1}{2}&f_1''-\frac{3}{16}-\frac{a+A}{4}\vspace{2mm}\\
		
		-\frac{1}{16}+\frac{a+A}{4}<f_1''\leq \frac{3}{16}+\frac{a+A}{4}&c_1&|c_1|&\pm \frac{1}{2}&\frac{1}{2}&f_1''+\frac{1}{16}-\frac{a+A}{4}\vspace{2mm}\\
		
		-\frac{1}{2}<f_1''<-\frac{13}{16}+d+a+A&\pm 1+c_1&1+|c_1|&\pm \frac{1}{2}&-\frac{1}{2}&f_1''+\frac{13}{16}-\frac{a+A}{4}\vspace{2mm}\\
		
		{\rm When~}  -\frac{13}{16}+d+a+A\leq f_1''\leq -\frac{1}{16}+\frac{a+A}{4}: \vspace{2mm}\\
		
		-\frac{13}{16}+d+a+A\leq g_1''\leq -\frac{1}{16}+\frac{a+A}{4}&c_1&2|c_1|&\pm \frac{1}{2}&1&f_2''+\frac{3}{16}-\frac{a+A}{4}\vspace{2mm}\\
		
		-\frac{13}{16}+\frac{a+A}{4}<g_1''<-\frac{13}{16}+d+a+A&\pm 1+c_1&1-|c_1|&\pm \frac{1}{2}&-\frac{1}{2}&g_1''+\frac{13}{16}-\frac{a+A}{4}\vspace{2mm}\\
		
		-\frac{17}{16}+\frac{a+A}{4}<g_1''\leq -\frac{13}{16}+\frac{a+A}{4}&\pm 1+c_1&1-|c_1|&\pm \frac{1}{2}&\frac{1}{2}&g_1''+\frac{17}{16}-\frac{a+A}{4}\vspace{2mm}\\
		
		-1<g_1'' \leq -\frac{17}{16}+ \frac{a+A}{4}&\pm 1+c_1&3+3|c_1|&\pm \frac{1}{2}&\frac{3}{2}&f_3''+\frac{53}{16}-\frac{a+A}{4}\vspace{2mm}\\
		
		\hline
		\end{array}
		\end{equation*}}
	
	\subsection{${m=2, K=2, L=1}$}
	
	\noindent Here $1<\frac{\delta_{m,K}}{t}\leq 2$, $\delta_{m,K}\geq \frac{25d}{16}$. From (\ref{Eq74}) we have
	\begin{equation}\label{Eq85}
	0.5379< \sqrt{125/432} \leq d \leq 16/25=0.64.
	\end{equation}
	\noindent Further, since $a\geq d/3$ and $A\geq 5d/12$ we get, $$0.42023< \frac{25d}{32}\leq \frac{d+3a/4+3A/4}{2}\leq t =\frac{3d^5}{8aA}\leq \frac{27d^3}{10}< 0.707789.$$
	
	\noindent Taking $h_t=1/2$ and $k_t=1$ and using (\ref{Eq85})
	we find that $\vline~h_t-tk_t^2~\vline + \frac{1}{2} <\delta_{m,K}$, as $ t < \delta_{m,K}$ and   $t+\delta_{m,K} \geq \frac{25}{32}d+\frac{25}{16}d=\frac{75}{32}d>1$. Therefore,  (\ref{eq101'}) is soluble unless $t=1/2$ and $\beta_{t,G}=\pm a_5''-c_5\equiv 1/2 \pmod 1$. This gives $(a_5'',c_5)=(0,1/2)$ or $(1/2,0)$.\vspace{2mm}
	
	\noindent Now, $1/2=t=3d^5/8aA\leq 27d^3/10$ gives $d\geq (5/27)^{1/3}>0.56999$. Also, $1/2=t=3d^5/8aA$ gives $A=3d^5/4a\leq 9d^4/4$. And, $\delta_{2,2,1}=d+\frac{3a}{4}+\frac{3A}{4}\leq 1$ gives $\frac{a+A}{4}=\frac{1}{3}\left( \frac{3a}{4}+\frac{3A}{4} \right) \leq \frac{1-d}{3}< 0.1433367$. Hence,
	\begin{equation}\label{Eq'86}
	\begin{array}{l}
	0.18999< d/3\leq a\leq d/2< 0.32\vspace{2mm}\\
	0.23749< 5d/12 \leq A \leq 9d^4/4< 0.377488\vspace{2mm}\\
	0.10687< \frac{a+A}{4} < 0.1433367.
	\end{array}
	\end{equation}
	
	\noindent Thus, for $(a_5'',c_5)=(0,1/2)$ and $(1/2,0)$, (\ref{eq101'}) reduces to
	\begin{equation}\label{Eq'87}
	0<G=(x_1+a_2''x_2+a_5''x_5)x_2-x_5^2/2-(a+A)/4<d+3a/4+3A/4.
	\end{equation}
	
	\noindent \textbf{Case (i)}:  $(a_5'',c_5)=(0,1/2)$. The following table gives a solution of (\ref{Eq'87}).\vspace{2mm}
	
	{\scriptsize
		\begin{equation*}\begin{array}{lllll}
		\hline\\
		{\rm ~~~~Range} & x_1 & x_1x_2 & x_5 & {\rm ~~G}\vspace{2mm}\\
		\hline\\
		
		\frac{1}{8}+\frac{a+A}{4}<f_1''\leq 1&c_1&|c_1|&\frac{1}{2}&f_1''-\frac{1}{8}-\frac{a+A}{4}\vspace{2mm}\\
		
		-\frac{1}{2}<f_1''<-\frac{7}{8}+d+a+A&\pm 1+c_1&1+|c_1|&\frac{1}{2}&f_1''+\frac{7}{8}-\frac{a+A}{4}\vspace{2mm}\\
		
		{\rm When~}  -\frac{7}{8}+d+a+A\leq f_1''\leq \frac{1}{8}+\frac{a+A}{4}: \vspace{2mm}\\

		-\frac{7}{8}+d+a+\frac{9A}{4}\leq g_1''\leq \frac{1}{8}+\frac{a+A}{4}&c_1&2|c_1|&\frac{1}{2}&f_2''-\frac{1}{8}-\frac{a+A}{4}\vspace{2mm}\\
		
		-\frac{7}{8}+\frac{a+A}{4}<g_1''<-\frac{7}{8}+d+a+A&\pm 1+c_1&1-|c_1|&\frac{1}{2}&g_1''+\frac{7}{8}-\frac{a+A}{4}\vspace{2mm}\\
		
		-1<g_1'' \leq -\frac{7}{8}+\frac{a+A}{4}&\pm 1+c_1&2+2|c_1|&\frac{3}{2}&f_2''+\frac{7}{8}-\frac{a+A}{4}\vspace{2mm}\\
		
		\hline
		\end{array}
		\end{equation*}\vspace{1mm}}

	\noindent \textbf{Case(ii)}:  $(a_5'',c_5)=(1/2,0)$. Here equation (\ref{Eq'87}) does not have a solution for $c_1=a_2''=\frac{1}{2}.$ This is so because $d+a+A> 0.997$ and $\frac{1}{2}\big((1+x_2+x_5)x_2-x_5^2\big)$ never takes the value $\frac{1}{2}$ for integral values of $x_2,x_5$. Therefore we have to go to $F$. Since $\delta_{2,2,1}\leq 1$, we get $A\leq 4(1-(d+3a/4))/3 \leq 4(1-5d/4)/3$.
	Therefore,
	\begin{equation}\label{Eq'89}
	\begin{array}{ll}
	0.23749< 5d/12 \leq A\leq &\min\{9d^4/4,4(1-5d/4)/3\}\vspace{2mm}\\
	&=\left\{
	\begin{array}{ll}
	3d^4 & {\rm ~if~} d\leq \alpha\\
	4(1-5d/4)/3 &{\rm ~if~} d\geq \alpha
	\end{array}
	\right.\vspace{2mm}\\
	&\leq 0.31437
	\end{array}
	\end{equation}
	
	\noindent where $\alpha$ ia a root of $3d^4-4(1-5d/4)/3=0$ satisfying $0.611<\alpha<0.612$. Taking $h_A=1$ and $k_A=2$ and using $0.56999< d\leq 16/25=0.64$ we find that $\vline~h_A-Ak_A^2~\vline +\frac{1}{2}<\delta_{2}$ is satisfied. Thus, (\ref{eq101}) is soluble unless $A=1/4$ and $\beta_{A,F} = \pm a'_4-(c_4+\lambda x_5)/2 \equiv 1/2 \pmod {1/2}$. Taking $x_5=c_5$ and $1+c_5$, we get $\lambda=0$. By Lemma \ref{lem15}, we need to consider $(a_4',c_4)=(0,0)$ or $(1/4,1/2)$. Also from (\ref{eq16}), we get $(a_5'',c_5)=(a_5',c_5)=(1/2,0)$.\vspace{1mm}
	
	\noindent When $(a_4',c_4)=(0,0)$, on taking $x_5=0$, the following table gives a solution to (\ref{eq101}) i.e., of \vspace{-2mm}
	 \begin{equation}\label{EQN'22}0<F=(x_1+a_2'x_2+\frac{1}{2}x_5)x_2-\frac{1}{4}x_4^2-\frac{1}{2}x_5^2-\frac{a}{4}<d+\frac{3a}{4}.\end{equation}
	
	{\scriptsize
		\begin{equation*}\begin{array}{lllll}
		\hline\\
		{\rm ~~~~Range} & x_1 & x_1x_2 & x_4& {\rm ~~F}\vspace{2mm}\\
		\hline\\
		
		\frac{1}{4}+\frac{a}{4}<f_1'\leq 1&c_1&|c_1|&1&f_1'-\frac{1}{4}-\frac{a}{4}\vspace{2mm}\\
		
		\frac{a}{4}<f_1'\leq \frac{1}{4}+\frac{a}{4}&c_1&|c_1|&0&f_1'-\frac{a}{4}\vspace{2mm}\\
		
		-\frac{1}{2}<f_1'<-\frac{3}{4}+d+a&\pm 1+c_1&1+|c_1|&1&f_1'+\frac{3}{4}-\frac{a}{4}\vspace{2mm}\\
		
		{\rm When~}  -\frac{3}{4}+d+a\leq f_1'\leq \frac{a}{4}: \vspace{2mm}\\
		
		-\frac{3}{4}+d+a\leq g_1'\leq \frac{a}{4}&c_1&3|c_1|&0&f_3'-\frac{a}{4}\vspace{2mm}\\
		
		-\frac{3}{4}+\frac{a}{4}<g_1'<-\frac{3}{4}+d+a&\pm 1+c_1&1-|c_1|&1&g_1'+\frac{3}{4}-\frac{a}{4}\vspace{2mm}\\
		
		-1+\frac{a}{4}<g_1'\leq -\frac{3}{4}+\frac{a}{4}&\pm 1+c_1&1-|c_1|&0&g_1'+1-\frac{a}{4}\vspace{2mm}\\
		
		-1<g_1' \leq -1+\frac{a}{4}&\pm 1+c_1&3+3|c_1|&0&f_3'+3-\frac{a}{4}\vspace{2mm}\\	
		
		\hline
		\end{array}
		\end{equation*}\vspace{1mm}}
	
	\noindent When $(a_4',c_4)=(1/4,1/2)$, on taking $x_5=0$, the following table gives a solution to (\ref{eq101}) i.e., of \vspace{-2mm}
	 \begin{equation}\label{EQN'23}0<F=(x_1+a_2'x_2+\frac{1}{4}x_4+\frac{1}{2}x_5)x_2-\frac{1}{4}x_4^2-\frac{1}{2}x_5^2-\frac{a}{4}<d+\frac{3a}{4}.\end{equation}
	
	{\scriptsize
		\begin{equation*}\begin{array}{llllll}
		\hline\\
		{\rm ~~~~Range} & x_1 & x_1x_2 & x_4&x_2x_4& {\rm ~~F}\vspace{2mm}\\
		\hline\\

		\frac{3}{16}+d+a\leq f_1'\leq 1&\pm 1+c_1&-2+2|c_1|&\pm \frac{1}{2}&-1&f_2'-\frac{37}{16}-\frac{a}{4}\vspace{2mm}\\
		
		\frac{3}{16}+\frac{a}{4}<f_1'<\frac{3}{16}+d+a&c_1&|c_1|&\pm \frac{1}{2}&-\frac{1}{2}&f_1'-\frac{3}{16}-\frac{a}{4}\vspace{2mm}\\
		
		-\frac{1}{16}+\frac{a}{4}<f_1'\leq \frac{3}{16}+\frac{a}{4}&c_1&|c_1|&\pm \frac{1}{2}&\frac{1}{2}&f_1'+\frac{1}{16}-\frac{a}{4}\vspace{2mm}\\
		
		-\frac{1}{2}<f_1'<-\frac{13}{16}+d+a&\pm 1+c_1&1+|c_1|&\pm \frac{1}{2}&-\frac{1}{2}&f_1'+\frac{13}{16}-\frac{a}{4}\vspace{2mm}\\
		
		{\rm When~}  -\frac{13}{16}+d+a\leq f_1'\leq -\frac{1}{16}+\frac{a}{4}: \vspace{2mm}\\
		
		-\frac{13}{16}+d+a\leq g_1'\leq -\frac{1}{16}+\frac{a}{4}&\pm 1+c_1&2+2|c_1|&\pm \frac{3}{2}&-3&f_2'+\frac{11}{16}-\frac{a}{4}\vspace{2mm}\\
		
		-\frac{13}{16}+\frac{a}{4}<g_1'<-\frac{13}{16}+d+a&\pm 1+c_1&1-|c_1|&\pm \frac{1}{2}&-\frac{1}{2}&g_1'+\frac{13}{16}-\frac{a}{4}\vspace{2mm}\\
		
		-1<g_1'\leq -\frac{13}{16}+\frac{a}{4}&\pm 1+c_1&1-|c_1|&\pm \frac{1}{2}&\frac{1}{2}&g_1'+\frac{17}{16}-\frac{a}{4}\vspace{2mm}\\
		
		\hline
		\end{array}
		\end{equation*}}

	\section{Proof of $\Gamma_{1,4}< \frac{32}{3}$ when $(m,K)=(1,2)$}\label{sec10}
	\numberwithin{equation}{section}
	\begin{theorem}\label{thm7}
		Let $Q(x_1,x_2, \cdots, x_5)$ be a real indefinite quadratic form of type $(1,4)$ and of determinant $D\neq 0$. Let $d=(\frac{32}{3}|D|)^{\frac{1}{5}}$.  Suppose that $c_2 \equiv 0 \pmod1$, $a<\frac{1}{2}$, $d\leq 1$ and $a+d\leq 1$. Let $(m,K)= (1,2)$.  Then $(\ref{eq101*})$ and hence $(\ref{eq10})$ is soluble with strict inequality.
	\end{theorem}
	
	\noindent {\bf Proof:} Here, $1 < \frac{d}{a}\leq 2$ and $2<\frac{\delta_{1}}{A}=\frac{d}{A} \leq 3$. Since, $\frac{d}{2}+d \leq a+d \leq 1$ we get $d \leq \frac{2}{3}$. Also $d/2 \leq a\leq (3d^5/4)^{1/3}$ which gives $d \geq \frac{1}{\sqrt{6}} > 0.40824$. Thus,
	\begin{equation}\label{EqN1'}
	0.40824 < d \leq 2/3< 0.6667.
	\end{equation}
	
	\noindent Now, $\frac{3d}{8}\leq \frac{3a}{4} \leq A \leq \sqrt{\frac{d^5}{2a}} \leq \sqrt{\frac{d^4.2}{2}} = d^2$. Thus,
	\begin{equation}\label{EqNn1}
	\delta_{1,2}=d+3A/4 \geq 41d/32.
	\end{equation}

	\noindent Further from (\ref{Eq55'}) we have $\frac{\delta_{m,K}}{t} \geq \frac{1}{24}\frac{(m+3)^2(K+3)}{(m+1)(K+1)}\frac{1}{d^2}$, which is bigger than $1$ for $m=1,K=2$ and $d\leq 2/3.$ Thus, $L\geq 1$ (from the definition of $L$). Also using (\ref{eq22}) we find that $\delta_{m,K,L}>1$ for $L\geq 5$ and $d\geq {1}/{\sqrt{6}}$. Therefore, we get that $L\leq 4$. Further, $\delta_{1,2,L}\leq 1$ and (\ref{EqNn1}) gives $d \leq \frac{128}{41(L+3)}$; and $\frac{1}{L+1} . \frac{41d}{32} \leq \frac{\delta_{1,2}}{L+1} \leq t=\frac{3d^5}{8aA}\leq 2d^3$ gives $d^2\geq \frac{41}{64(L+1)}$. Thus,
	
	\begin{equation}\label{Eq74'}
	\sqrt{41/64(L+1)} \leq d \leq 128/41(L+3)
	\end{equation}
	\noindent and for $L \neq 1$, as $\delta_{1,2,L}=\delta_{1,2}+ \frac{L^2-1}{4}t \leq 1$, we have
	
	\begin{equation}\label{EqNn2}
	t\leq \frac{4}{L^2-1}(1-\delta_{1,2}) \leq \frac{4}{L^2-1}\left(1-\frac{41d}{32}\right).
	\end{equation}

	\noindent We distinguish the cases $L=1,2,3$ and $4$ in the following subsections:
	
	\subsection{${m=1, K=2, L=4}$}	
	
	\noindent Here $4<\frac{\delta_{m,K}}{t}=\frac{d+3A/4}{t}\leq 5$, $\delta_{m,K}=d+\frac{3A}{4}\geq d+\frac{9d}{32}=\frac{41d}{32}$. From (\ref{EqN1'}) and (\ref{Eq74'}), we have
	
	\begin{equation}\label{Eq'85''}
	0.40824< 1/\sqrt{6} \leq d \leq 128/287< 0.445994.
	\end{equation}
	
	\noindent Further since $A\geq 3d/8$,
	from (\ref{EqNn2}), we get
	\begin{equation*}\label{Eq''89''}
	\begin{array}{ll}
	0.114817< 9d/32 \leq 3A/4 \leq t\leq
	4(1-41d/32)/15  \leq  0.127185 =\lambda_0 {\rm ~(say)}.
	\end{array}
	\end{equation*}
	
	\noindent Thus, we have
	\begin{equation}\label{Eq90}
	\begin{array}{l}
	0.20412< d/2\leq a\leq 4A/3\leq 16t/9 < 0.226107\vspace{2mm}\\
	0.15309< 3d/8 \leq A \leq 4t/3 < 0.16958\vspace{2mm}\\
	0.118006875 \leq (a+A+t)/4 \leq 0.130718.
	\end{array}
	\end{equation}

	\noindent When, $0.1148175 \leq t \leq 0.1159$, we have $d\leq 4t/3<0.4121,$ $a\leq 16t/9<0.2061$ and $A\leq 4t/3<0.1546$. For this range of $t$, the following table gives a solution to (\ref{eq21}) with strict inequality, i.e., of
	\begin{equation}\label{Eq91}
	0<H=(x_1+a_2'''x_2)x_2-\frac{a+A+t}{4}< d+\frac{3A}{4}+\frac{15t}{4}.
	\end{equation}
	
	{\scriptsize
		\begin{equation*}\begin{array}{llll}
		\hline\\
		{\rm ~~~~Range} & x_1 & x_1x_2 & {\rm ~~H}\vspace{2mm}\\
		\hline\\
		
		\frac{a+A+t}{4}<f_1'''\leq 1&c_1&|c_1|&f_1'''-\frac{a+A+t}{4} \vspace{2mm}\\
		
		-\frac{1}{2}<f_1'''<-1+d+\frac{a}{4}+A+4t &\pm 1+c_1&1+|c_1|&f_1'''+1-\frac{a+A+t}{4}\vspace{2mm}\\
		
		{\rm When~}  -1+d+\frac{a}{4}+A+4t\leq f_1'''\leq \frac{a+A+t}{4}: \vspace{2mm}\\
		
		-1+d+\frac{a}{4}+A+4t\leq g_1'''\leq \frac{a+A+t}{4}&c_1&2|c_1|&f_2'''-\frac{a+A+t}{4}\vspace{2mm}\\
		
		-1+\frac{a+A+t}{4}<g_1'''<-1+d+\frac{a}{4}+A+4t &\pm 1+c_1&1-|c_1|&g_1'''+1-\frac{a+A+t}{4}\vspace{2mm}\\
		
		-1<g_1'''\leq -1+\frac{a+A+t}{4}&\pm 1+c_1&3+3|c_1|&f_3'''+3-\frac{a+A+t}{4}\vspace{2mm}\\
		
		\hline
		\end{array}
		\end{equation*}\vspace{1mm}}
	
	\noindent When $0.1159<t\leq 0.127185=\lambda_0(\mbox{say})$, we divide the range of $t$ into $62$ subintervals $[\lambda_n,\lambda_{n-1}], n=1,2,\cdots 62$ and in each subinterval choose suitable integers $2h_n$ and $k_n$ such that $\vline~ h_n-tk_n^2~\vline +\frac{1}{2} <d+\frac{3A}{4}.$
	
	{\scriptsize
		\begin{equation*}\begin{array}{lllllllllll}
		\hline\\
		n & h_n & k_n & \lambda_n&\mbox{Remarks}&~~~~ &n & h_n & k_n & \lambda_n&\mbox{Remarks}\vspace{2mm}\\
		\hline\\
		
		1&0.5&2&0.11924&\mbox{tbd}&&2&3&5&0.11908&(0.5,2)\vspace{2mm}\\
		
		3&43&19&0.11905&(3,5)&&4&52.5&21&0.118996&\mbox{tbd}\vspace{2mm}\\
		
		5&181&39&0.118986&(52.5,21)&&6&200&41&0.118964&(583,70)\vspace{2mm}\\
		
		7&137.5&34&0.11893&(68.5,24)&&8&68.5&24&0.118884&(100,29)\vspace{2mm}\\
		
		9&38.5&18&0.11876&(57.5,22)&&10&47.5&20&0.118693&\mbox{tbd}\vspace{2mm}\\
		
		11&162.5&37&0.118683&(47.5,20)&&12&86.5&27&0.11863&(121.5,32)\vspace{2mm}\\
		
		13&93&28&0.118594&(114,31)&&14&685&76&0.118591&(296.5,50)\vspace{2mm}\\

		15&721.5&78&0.118586&\mbox{tbd}&&16&982&91&0.118582&(501,65)\vspace{2mm}\\
		
		17&501&65&0.118575&(667,75)&&18&516.5&66&0.118567&(564.5,69)\vspace{2mm}\\
		
		19&564.5&69&0.118563&(516.5,66)&&20&1003.5&92&0.11856&(1070,95)\vspace{2mm}\\
		
		21&1070&95&0.118558&(1003.5,92)&&22&1115.5&97&0.118555&\mbox{tbd}\vspace{2mm}\\
		
		23&333&53&0.11854&(229.5,44)&&24&229.5&44&0.118532&(470.5,63)\vspace{2mm}\\
		
		25&137&34&0.118493&(240,45)&&26&209&42&0.11847&(273,48)\vspace{2mm}\\
		
		27&20&13&0.11821&(358,55)&&28&17&12&0.1179&\mbox{tbd}\vspace{2mm}\\
		
		29&52&21&0.11787&(249.5,46)&&30&26.5&15&0.1177&(1249.5,103)\vspace{2mm}\\
		
		31&34&17&0.117568&\mbox{tbd}&&32&9.5&9&0.117&(7.5,8)\vspace{2mm}\\
		
		33&7.5&8&0.11683&(9.5,9)&&34&51.5&21&0.11673&(56.5,22)\vspace{2mm}\\
		
		35&56.5&22&0.1167&(91.5,28)&&36&91.5&28&0.11668&(105,30)\vspace{2mm}\\
		
		37&105&30&0.1166411&\mbox{tbd}&&38&434&61&0.1166292&(406,59)\vspace{2mm}\\
		
		39&85&27&0.116567&(196,41)&&40&112&31&0.116522&(215.5,43)\vspace{2mm}\\
		
		41&159.5&37&0.116492&(151,36)&&42&205.5&42&0.116484&(246.5,46)\vspace{2mm}\\
		
		43&225.5&44&0.116466&\mbox{tbd}&&44&554.5&69&0.1164623&(225.5,44)\vspace{2mm}\\
		
		45&538.5&68&0.116453&(477,64)&&46&477&64&0.11645&(492,65)\vspace{2mm}\\
		
		47&492&65&0.116445&(477,64)&&48&587&71&0.116441&(492,65)\vspace{2mm}\\
		
		49&620.5&73&0.1164341&\mbox{tbd}&&50&339.5&54&0.116419&(570.5,70)\vspace{2mm}\\
		
		51&279.5&49&0.116401&(327,53)&&52&42&19&0.1163&(67,24)\vspace{2mm}\\
		
		53&46.5&20&0.1162&(61.5,23)&&54&119&32&0.11619&(205,42)\vspace{2mm}\\
		
		\hline
		\end{array}
		\end{equation*}}
	
	{\scriptsize
		\begin{equation*}\begin{array}{lllllllllll}
		\hline\\
		n & h_n & k_n & \lambda_n&\mbox{Remarks}&~~~~ &n & h_n & k_n & \lambda_n&\mbox{Remarks}\vspace{2mm}\\
		\hline\\
		
		55&377.5&57&0.116183&(1281,105)&&56&1431.5&111&0.116182&(377.5,57)\vspace{2mm}\\
		
		57&126.5&33&0.11615&\mbox{tbd}&&58&78.5&26&0.1161&(150.5,36)\vspace{2mm}\\
		
		59&91&28&0.116043&\mbox{tbd}&&60&111.5&31&0.11601&(245.5,46)\vspace{2mm}\\
		
		61&72.5&25&0.11597&(195,41)&&62&33.5&17&0.1159&(84.5,27)\vspace{2mm}\\	
		
		\hline
		\end{array}
		\end{equation*}\vspace{2mm}}
	
	\noindent The choice of $(h_n,k_n)$ is done in a manner similar to that in the Case  $m=2,K=3, L=4$. Thus (\ref{eq101*}) is soluble unless $t=1/8,$ $525/4410,$ $19/160,$ $37/312,$ $23/194,$ $17/144,$ $2/17,$ $7/60,$ $41/352,$ $17/146,$ $23/198,$ $13/112$  and $\beta_{t,G} = \pm a_5''-2tc_5\equiv h_n/k_n \pmod {1/k_n, 2t}$.
	\noindent  For all these values of $t$ and the corresponding $(a_5'',c_5)$, using Lemma \ref{lem15}, and the  bounds $1/\sqrt{6}\leq d \leq 128/287,$ $d/2\leq a\leq 4A/3\leq 16t/9,$ $3d/8\leq A\leq 4t/3,$ we find that (\ref{eq101'}) is soluble on taking $x_2=\pm 1$ and suitably choosing $x_5 \equiv c_5 \pmod 1$.
	
	\subsection{${m=1, K=2, L=3}$}
	
	\noindent Here $3<\frac{\delta_{m,K}}{t}=\frac{d+3A/4}{t}\leq 4$, $\delta_{m,K}=d+\frac{3A}{4}\geq d+\frac{9d}{32}=\frac{41d}{32}$. From (\ref{EqN1'}) and (\ref{Eq74'}), we have
	\begin{equation}\label{Eq''85''}
	0.40824< 1/\sqrt{6} \leq d \leq 128/246<0.520326.
	\end{equation}
	\noindent Further, since $a\geq d/2$ and $A\geq 3d/8$ we get, $$0.130764< \frac{41d}{128}\leq \frac{d+3A/4}{4}=\frac{\delta_{1,2}}{4}\leq t =\frac{3d^5}{8aA}\leq 2d^3< 0.28175.$$
	
	\noindent Also from (\ref{EqNn2}), we get $t \leq [1-41d/32]/2$. Therefore,
	\begin{equation*}\label{Eq''89'''}
	\begin{array}{ll}
	0.130764< 41d/128 \leq t\leq &\min\{2d^3,(1-41d/32)/2\}\vspace{2mm}\\
	&=\left\{
	\begin{array}{ll}
	2d^3 & {\rm ~if~} d\leq \alpha\\
	(1-41d/32)/2 &{\rm ~if~} d\geq \alpha
	\end{array}
	\right.\vspace{2mm}\\
	&\leq 0.20174=\lambda_0(say)
	\end{array}
	\end{equation*}

	\noindent where $\alpha$ is a root of $2d^3-(1-41d/32)/2=0$ satisfying $0.465 < \alpha <0.466$. We divide this range of $t$ into $272$ subintervals $[\lambda_n,\lambda_{n-1}], n=1,2,\cdots 272$ and in each subinterval choose suitable integers $2h_n$ and $k_n$ such that $\vline~ h_n-tk_n^2~\vline +\frac{1}{2} <d+\frac{3A}{4}.$\vspace{2mm}
		
	{\scriptsize
		\begin{equation*}\begin{array}{lllllllll}
		
		\hline\\
		n & (h_n,k_n)& \lambda_n & \mbox{Remarks} &&n & (h_n,k_n)& \lambda_n & \mbox{Remarks} \vspace{2mm}\\
		\hline\\
		
		1 & (5,5) & 0.1964 & \mbox{tbd} & ~ & 2 & (7,6) & 0.19297 & (163.5,29) \\
		3 & (3,4) & 0.1827 & \mbox{tbd} & ~ & 4 & (4.5,5) & 0.17716 & (6.5,6) \\
		5 & (25.5,12) & 0.1766 & \mbox{tbd} & ~ & 6 & (51,17) & 0.17637 & \mbox{tbd} \\
		7 & (34.5,14) & 0.17569 & (119,26) & ~ & 8 & (17.5,10) & 0.17487 & \mbox{tbd} \\
		9 & (8.5,7) & 0.1742 & \mbox{na} & ~ & 10 & (1.5,3) & 0.1608 & \mbox{tbd} \\
		11 & (4,5) & 0.15946 & (219,37) & ~ & 12 & (2.5,4) & 0.1535 & \mbox{tbd} \\
		13 & (5.5,6) & 0.15161 & (128.5,29) & ~ & 14 & (102.5,26) & 0.15157 & \mbox{na} \\
		15 & (110.5,27) & 0.15152 & \mbox{na} & ~ & 16 & (165,33) & 0.151477 & \mbox{tbd} \\
		
		\hline
		\end{array}
		\end{equation*}}
	
	{\scriptsize
		\begin{equation*}\begin{array}{lllllllll}
		
		\hline\\
		n & (h_n,k_n)& \lambda_n & \mbox{Remarks} &&n & (h_n,k_n)& \lambda_n & \mbox{Remarks} \vspace{2mm}\\
		\hline\\
		
		17 & (349,48) & 0.151458 & (425.5,53) & ~ & 18 & (185.5,35) & 0.151448 & \mbox{na} \\
		19 & (145.5,31) & 0.151362 & (441.5,54) & ~ & 20 & (49,18) & 0.15129 & \mbox{na} \\
		
		21 & (34,15) & 0.15113 & \mbox{na} & ~ & 22 & (25.5,13) & 0.15065 & (73,22) \\
		23 & (29.5,14) & 0.15055 & \mbox{na} & ~ & 24 & (38.5,16) & 0.1504 & \mbox{na} \\
		25 & (15,10) & 0.14961 & \mbox{tbd} & ~ & 26 & (48.5,18) & 0.14958 & \mbox{na} \\
		27 & (21.5,12) & 0.14904 & (79,23) & ~ & 28 & (9.5,8) & 0.14785 & (78.5,23) \\
		29 & (12,9) & 0.14769 & \mbox{na} & ~ & 30 & (71.5,22) & 0.14765 & \mbox{na} \\
		31 & (299,45) & 0.147636 & \mbox{na} & ~ & 32 & (273,43) & 0.14763 & \mbox{na} \\
		33 & (85,24) & 0.14759 & \mbox{na} & ~ & 34 & (59,20) & 0.147476 & (170.5,34) \\
		35 & (65,21) & 0.14731 & (743,71) & ~ & 36 & (115.5,28) & 0.147275 & \mbox{na} \\
		37 & (141.5,31) & 0.147259 & \mbox{na} & ~ & 38 & (92,25) & 0.14718 & (413.5,53) \\
		39 & (42.5,17) & 0.14694 & \mbox{tbd} & ~ & 40 & (160,33) & 0.146915 & (180,35) \\
		41 & (53,19) & 0.14682 & \mbox{na} & ~ & 42 & (33,15) & 0.14651 & (47.5,18) \\
		
		43 & (37.5,16) & 0.14635 & (381,51) & ~ & 44 & (91.5,25) & 0.146343 & \mbox{na} \\
		
		45 & (246,41) & 0.146339 & \mbox{tbd} & ~ & 46 & (58.5,20) & 0.146161 & (64.5,21) \\
		47 & (617.5,65) & 0.14615 & \mbox{tbd} & ~ & 48 & (131.5,30) & 0.146139 & \mbox{na} \\
		49 & (106.5,27) & 0.14608 & (1055.5,85) & ~ & 50 & (21,12) & 0.145591 & \mbox{tbd} \\
		
		51 & (77,23) & 0.145586 & \mbox{na} & ~ & 52 & (28.5,14) & 0.14534 & (106,27) \\
		53 & (14.5,10) & 0.14467 & (24.5,13) & ~ & 54 & (17.5,11) & 0.14436 & (574,63) \\
		55 & (32.5,15) & 0.144297 & \mbox{na} & ~ & 56 & (187,36) & 0.14428 & (855.5,77) \\
		57 & (97.5,26) & 0.144182 & \mbox{tbd} & ~ & 58 & (157,33) & 0.144175 & \mbox{na} \\
		59 & (113,28) & 0.144136 & \mbox{na} & ~ & 60 & (52,19) & 0.143954 & (83,24) \\	
		61 & (63.5,21) & 0.143917 & \mbox{na} & ~ & 62 & (186.5,36) & 0.143915 & \mbox{na} \\
		63 & (121,29) & 0.143838 & (685,69) & ~ & 64 & (766.5,73) & 0.143835 & \mbox{tbd} \\
		65 & (404,53) & 0.14383 & \mbox{na} & ~ & 66 & (57.5,20) & 0.14371 & \mbox{tbd} \\
		67 & (41.5,17) & 0.14351 & (46.5,18) & ~ & 68 & (7,7) & 0.14223 & \mbox{tbd} \\
		69 & (11.5,9) & 0.14161 & (205,38) & ~ & 70 & (88.5,25) & 0.14159 & (145,32) \\
		71 & (68.5,22) & 0.14147 & (136,31) & ~ & 72 & (81.5,24) & 0.141442 & \mbox{na} \\
		73 & (163.5,34) & 0.141411 & (154,33) & ~ & 74 & (154,33) & 0.141387 & \mbox{na} \\
		75 & (509,60) & 0.141381 & \mbox{na} & ~ & 76 & (1773.5,112) & 0.14138 & \mbox{na} \\
		77 & (2972.5,145) & 0.1413785 & \mbox{tbd} & ~ & 78 & (1588.5,106) & 0.1413733 & (1649,108) \\
		79 & (1968.5,118) & 0.141373 & \mbox{na} & ~ & 80 & (997.5,84) & 0.14137 & \mbox{na} \\
		81 & (927.5,81) & 0.141365 & (193.5,37) & ~ & 82 & (193.5,37) & 0.14135 & \mbox{na} \\
		83 & (51,19) & 0.1412 & (56.5,20) & ~ & 84 & (56.5,20) & 0.14118 & \mbox{na} \\
		85 & (3.5,5) & 0.13894 & \mbox{tbd} & ~ & 86 & (5,6) & 0.13818 & \mbox{tbd} \\
		87 & (141.5,32) & 0.138159 & \mbox{na} & ~ & 88 & (199.5,38) & 0.138141 & \mbox{tbd} \\
		89 & (221,40) & 0.138137 & \mbox{na} & ~ & 90 & (179,36) & 0.138098 & (221,40) \\
			
		91 & (657.5,69) & 0.138096 & \mbox{na} & ~ & 92 & (1522.5,105) & 0.138093 & \mbox{tbd} \\
		93 & (2298,129) & 0.1380927 & \mbox{na} & ~ & 94 & (1465,103) & 0.138091 & \mbox{na} \\
		95 & (601.5,66) & 0.138083 & (464.5,58) & ~ & 96 & (210,39) & 0.138075 & \mbox{na} \\
		97 & (189,37) & 0.138065 & \mbox{na} & ~ & 98 & (79.5,24) & 0.138044 & \mbox{na} \\
		99 & (73,23) & 0.137949 & (997,85) & ~ & 100 & (116,29) & 0.137901 & \mbox{tbd} \\
		101 & (132.5,31) & 0.13789 & \mbox{na} & ~ & 102 & (31,15) & 0.13788 & \mbox{na} \\
		103 & (27,14) & 0.13763 & (31,15) & ~ & 104 & (86,25) & 0.1375601 & (386.5,53) \\
		105 & (55,20) & 0.13744 & \mbox{tbd} & ~ & 106 & (66.5,22) & 0.13742 & \mbox{na} \\
		107 & (44.5,18) & 0.13727 & (115.5,29) & ~ & 108 & (149.5,33) & 0.13726 & \mbox{na} \\
		109 & (357,51) & 0.137249 & \mbox{tbd} & ~ & 110 & (123.5,30) & 0.137243 & \mbox{na} \\
		
		111 & (60.5,21) & 0.137186 & (79,24) & ~ & 112 & (49.5,19) & 0.137052 & (79,24) \\
		113 & (72.5,23) & 0.137006 & (290,46) & ~ & 114 & (342.5,50) & 0.1369904 & (833.5,78) \\
		115 & (730,73) & 0.1369818 & \mbox{tbd} & ~ & 116 & (770.5,75) & 0.1369815 & \mbox{na} \\
		117 & (690.5,71) & 0.136978 & \mbox{na} & ~ & 118 & (177.5,36) & 0.136942 & (187.5,37) \\
		119 & (302.5,47) & 0.136929 & (543.5,63) & ~ & 120 & (315.5,48) & 0.1369253 & \mbox{na} \\
		121 & (509.5,61) & 0.1369191 & \mbox{na} & ~ & 122 & (578.5,65) & 0.136918 & \mbox{na} \\
		123 & (241.5,42) & 0.1368912 & \mbox{tbd} & ~ & 124 & (265,44) & 0.13689 & \mbox{na} \\
		125 & (219,40) & 0.136869 & (265,44) & ~ & 126 & (92.5,26) & 0.136838 & \mbox{na} \\
		127 & (85.5,25) & 0.136812 & \mbox{na} & ~ & 128 & (35,16) & 0.136626 & (302,47) \\
		129 & (39.5,17) & 0.136597 & \mbox{na} & ~ & 130 & (187,37) & 0.136592 & (459.5,58) \\ 	
		131 & (177,36) & 0.136559 & (218.5,40) & ~ & 132 & (16.5,11) & 0.136173 & \mbox{tbd} \\
		133 & (23,13) & 0.13608 & (60,21) & ~ & 134 & (11,9) & 0.13553 & (157,34) \\
		135 & (19.5,12) & 0.135264 & \mbox{tbd} & ~ & 136 & (26.5,14) & 0.13521 & \mbox{na} \\
		137 & (13.5,10) & 0.13479 & (39,17) & ~ & 138 & (34.5,16) & 0.134683 & (129.5,31) \\
		139 & (215.5,40) & 0.134675 & \mbox{na} & ~ & 140 & (249,43) & 0.134661 & (297.5,47) \\
		
		141 & (174.5,36) & 0.134646 & \mbox{na} & ~ & 142 & (91,26) & 0.134585 & \mbox{tbd} \\
		143 & (77.5,24) & 0.134513 & (105.5,28) & ~ & 144 & (155.5,34) & 0.134498 & \mbox{na} \\
		
		\hline
		\end{array}
		\end{equation*}}

	{\scriptsize
		\begin{equation*}\begin{array}{lllllllll}
		
		\hline\\
		n & (h_n,k_n)& \lambda_n & \mbox{Remarks} &&n & (h_n,k_n)& \lambda_n & \mbox{Remarks} \vspace{2mm}\\
		\hline\\
		
		145 & (517,62) & 0.13449 & (678,71) & ~ & 146 & (659,70) & 0.1344856 & (756.5,75) \\
		
		147 & (1657,111) & 0.1344842 & \mbox{na} & ~ & 148 & (1969,121) & 0.134484 & \mbox{na} \\
		
		149 & (4261,178) & 0.1344837 & \mbox{na} & ~ & 150 & (2827.5,145) & 0.1344818 & \mbox{tbd} \\
		151 & (2750,143) & 0.1344817 & \mbox{na} & ~ & 152 & (1778.5,115) & 0.1344803 & \mbox{na} \\
		153 & (1511,106) & 0.1344786 & \mbox{na} & ~ & 154 & (1454.5,104) & 0.134477 & \mbox{na} \\	
		155 & (1291.5,98) & 0.1344731 & (1318,99) & ~ & 156 & (1568.5,108) & 0.1344718 & \mbox{na} \\
		157 & (2711.5,142) & 0.1344713 & \mbox{na} & ~ & 158 & (3272.5,156) & 0.1344709 & \mbox{na} \\
		159 & (2343,132) & 0.1344695 & \mbox{tbd} & ~ & 160 & (1399,102) & 0.1344688 & \mbox{na} \\
		161 & (1163,93) & 0.134467 & \mbox{na} & ~ & 162 & (121,30) & 0.134459 & \mbox{na} \\	
		163 & (98,27) & 0.134433 & \mbox{na} & ~ & 164 & (84,25) & 0.134406 & \mbox{na} \\
		165 & (48.5,19) & 0.13432 & (4209,177) & ~ & 166 & (43.5,18) & 0.1342 & (8257.5,248) \\
		167 & (71,23) & 0.13418 & \mbox{na} & ~ & 168 & (349,51) & 0.134172 & (3435,160) \\
		169 & (483,60) & 0.134162 & (4742,188) & ~ & 170 & (467,59) & 0.134158 & \mbox{na} \\
		171 & (225.5,41) & 0.134135 & \mbox{tbd} & ~ & 172 & (204,39) & 0.134109 & (3306,157) \\
		173 & (322,49) & 0.134103 & \mbox{na} & ~ & 174 & (499,61) & 0.1341 & \mbox{na} \\
		175 & (155,34) & 0.134086 & \mbox{na} & ~ & 176 & (146,33) & 0.13405 & (6548,221) \\
		177 & (183.5,37) & 0.134025 & (9342,264) & ~ & 178 & (482.5,60) & 0.1340222 & \mbox{na} \\
		179 & (1261,97) & 0.134019 & \mbox{tbd} & ~ & 180 & (466.5,59) & 0.134016 & \mbox{na} \\181 & (193.5,38) & 0.1339887 & (5200.5,197) & ~ & 182 & (283.5,46) & 0.13397 & (515,62) \\183 & (498.5,61) & 0.133964 & (2778,144) & ~ & 184 & (566,65) & 0.13396 & \mbox{na} \\
		185 & (815,78) & 0.133959 & \mbox{na} & ~ & 186 & (753.5,75) & 0.133957 & \mbox{na} \\
		187 & (583.5,66) & 0.133954 & \mbox{na} & ~ & 188 & (105,28) & 0.133905 & \mbox{tbd} \\
		189 & (90.5,26) & 0.133846 & (1476,105) & ~ & 190 & (225,41) & 0.133837 & \mbox{na} \\
		191 & (271,45) & 0.13383 & \mbox{na} & ~ & 192 & (59,21) & 0.1338 & \mbox{na} \\
		193 & (53.5,20) & 0.133701 & (112.5,29) & ~ & 194 & (77,24) & 0.133646 & (183,37) \\
		195 & (154.5,34) & 0.133634 & \mbox{na} & ~ & 196 & (834,79) & 0.13363244 & \mbox{na} \\
		197 & (813,78) & 0.133631 & \mbox{na} & ~ & 198 & (83.5,25) & 0.133569 & (145.5,33) \\
		199 & (497,61) & 0.133563 & (5668,206) & ~ & 200 & (224.5,41) & 0.13354 & (5130.5,196) \\201 & (295,47) & 0.1335357 & \mbox{na} & ~ & 202 & (530,63) & 0.1335329 & (3504.5,162) \\
		203 & (258.5,44) & 0.133514 & \mbox{tbd} & ~ & 204 & (235.5,42) & 0.133502 & (4668.5,187) \\
		205 & (173,36) & 0.133485 & (4718,188) & ~ & 206 & (163.5,35) & 0.133454 & (4568,185) \\
		207 & (418.5,56) & 0.133449 & (4919.5,192) & ~ & 208 & (213.5,40) & 0.133426 & (2504.5,137) \\
		209 & (433.5,57) & 0.13342 & (213.5,40) & ~ & 210 & (30,15) & 0.13325 & \mbox{tbd} \\
		211 & (22.5,13) & 0.13304 & (213,40) & ~ & 212 & (6.5,7) & 0.13229 & (58.5,21) \\
		
		213 & (16,11) & 0.132084 & (244.5,43) & ~ & 214 & (233,42) & 0.132077 & \mbox{na} \\
		215 & (371,53) & 0.132074 & \mbox{tbd} & ~ & 216 & (222,41) & 0.132068 & \mbox{na} \\
		217 & (19,12) & 0.131823 & \mbox{tbd} & ~ & 218 & (135,32) & 0.131819 & \mbox{na} \\
		219 & (200.5,39) & 0.13181 & \mbox{na} & ~ & 220 & (232.5,42) & 0.131793 & (4318,181) \\
		
		221 & (1115.5,92) & 0.131792 & \mbox{na} & ~ & 222 & (822.5,79) & 0.131788 & (4270,180) \\
		223 & (143.5,33) & 0.131757 & (221.5,41) & ~ & 224 & (721.5,74) & 0.1317549 & (143.5,33) \\
		225 & (683,72) & 0.1317482 & (5376,202) & ~ & 226 & (1091,91) & 0.1317453 & (2545.5,139) \\
		227 & (1452.5,105) & 0.131745 & \mbox{na} & ~ & 228 & (1189,95) & 0.1317433 & \mbox{na} \\
		229 & (2058.5,125) & 0.131743 & \mbox{na} & ~ & 230 & (2401,135) & 0.1317417 & (2694,143) \\
		
		231 & (291,47) & 0.131726 & (1239.5,97) & ~ & 232 & (303.5,48) & 0.1317238 & \mbox{na} \\
		233 & (255,44) & 0.13171 & (2158,128) & ~ & 234 & (96,27) & 0.131682 & (1772,116) \\
		235 & (89,26) & 0.131632 & (591,67) & ~ & 236 & (126.5,31) & 0.131626 & \mbox{na} \\
		237 & (47.5,19) & 0.131547 & \mbox{tbd} & ~ & 238 & (38,17) & 0.131429 & (180,37) \\
		
		239 & (161,35) & 0.131415 & (243,43) & ~ & 240 & (489,61) & 0.131412 & \mbox{na} \\
		241 & (69.5,23) & 0.131348 & (5467.5,204) & ~ & 242 & (134.5,32) & 0.1313312 & (383,54) \\
		243 & (2085,126) & 0.1313308 & \mbox{na} & ~ & 244 & (799,78) & 0.1313299 & \mbox{na} \\
		245 & (643.5,70) & 0.1313247 & (1017,88) & ~ & 246 & (589.5,67) & 0.1313173 & (143,33) \\
		247 & (1767,116) & 0.131317 & \mbox{na} & ~ & 248 & (572,66) & 0.13131 & \mbox{tbd} \\
		249 & (819.5,79) & 0.1313066 & (1706.5,114) & ~ & 250 & (778.5,77) & 0.131302 & (1736.5,115) \\
		251 & (341.5,51) & 0.1313 & \mbox{na} & ~ & 252 & (302.5,48) & 0.131292 & (341.5,51) \\
		253 & (52.5,20) & 0.131224 & \mbox{tbd} & ~ & 254 & (42.5,18) & 0.13118 & \mbox{na} \\
		255 & (29.5,15) & 0.131038 & (4634,188) & ~ & 256 & (1258.5,98) & 0.1310377 & \mbox{na} \\
		257 & (1614.5,111) & 0.1310352 & (3481.5,163) & ~ & 258 & (2755,145) & 0.1310338 & \mbox{tbd} \\
		259 & (160.5,35) & 0.131024 & \mbox{na} & ~ & 260 & (95.5,27) & 0.13098 & (1732.5,115) \\
		261 & (536.5,64) & 0.130978 & \mbox{na} & ~ & 262 & (679,72) & 0.13097678 & \mbox{na} \\
		263 & (471.5,60) & 0.130968 & (756.5,76) & ~ & 264 & (425.5,57) & 0.130961 & (4290.5,181) \\
		265 & (231,42) & 0.130948 & \mbox{tbd} & ~ & 266 & (253.5,44) & 0.130941 & \mbox{na} \\
		267 & (88.5,26) & 0.130923 & \mbox{na} & ~ & 268 & (33.5,16) & 0.1307966 & (142.5,33) \\
		269 & (110,29) & 0.130778 & \mbox{na} & ~ & 270 & (314,49) & 0.130773 & \mbox{na} \\
		271 & (552.5,65) & 0.130767 & \mbox{tbd} & ~ & 272 & (519,63) & 0.130764 & \mbox{na} \\
		\hline
		\end{array}
		\end{equation*}\vspace{5mm}}
	
	\noindent The choice of $(h_n,k_n)$ is done in a manner similar to that in Section \ref{subsec7.2} above $(m=2,~K=3,~L=4)$. Thus by Macbeath Lemma, (\ref{eq101'}) is soluble unless $t=\frac{h_n}{k_n^2}$  and $\beta_{t,G} = \pm a_5''-2tc_5\equiv h_n/k_n \pmod {1/k_n, 2t}$, where $(h_n,k_n)$ to be discussed are listed in the above table.
	For a fixed $t$ we get
	\begin{equation}\label{Eqn1}
	\begin{array}{llll}
	{\rm max~}\{(\frac{t}{2})^{1/3}, \frac{1}{\sqrt{6}}\} &\leq &d &\leq {\rm ~min~}\{\frac{128t}{41}, \frac{32}{41}[1-2t] \}\vspace{2mm}\\
	
	{\rm ~~~~~~~~~~~}d/2 &\leq &a&\leq {\rm ~min~} \{\frac{16t}{9}, \frac{4A}{3}, (\frac{3d^5}{4})^{1/3} \} \vspace{2mm}\\
	
	{\rm ~~~~~~~~~~~}3d/8 &\leq & A&\leq {\rm ~min~} \{\frac{4t}{3}, d^2, \frac{3d^4}{4t} \}.
	
	\end{array}
	\end{equation}
	
	\noindent The cases  $t=1/5$ and $1/6$, i.e., when  $n=1, 10$, are discussed  separately in Lemmas \ref{lem81}, \ref{lem82} and \ref{lem82'}. For all other values of $t$,   Lemma \ref{lem15} fixes the values of $(a_5'',c_5)$ to be considered and  for each value of $(a_5'',c_5)$, using (\ref{Eqn1}) we give suitable values to $x_2$ and $x_5$, which covers the range of $f_1''$ or $p_1''$ and hence (\ref{eq101'}) is soluble.
	
	\begin{lemma}\label{lem81}
		If  $t=\frac{1}{5}$ and $\beta_{t,G} = \pm a_5''-2tc_5\equiv 1 \pmod {1/5, 2/5}$, then inequality  $(\ref{eq101'})$ is soluble. \end{lemma}
	
	\noindent {\bf Proof:} Here, since $\frac{d}{a}\leq 2$, and $\frac{1}{5}=t=\frac{3d^5}{8aA}\leq 2d^3$, we get $d^3\geq \frac{1}{10}$. And $t\leq (1-41d/32)/2$ gives $d\leq 32(1-2t)/41$. Also,  $A=15d^5/8a\leq 15d^4/4$. Therefore,
	
	\begin{equation}\label{Eqn2}
	\begin{array}{l}
	0.46415 \leq (1/10)^{1/3} \leq d \leq 32[1-2t]/41 \leq 0.4683\vspace{2mm}\\
	0.232\leq d/2\leq a\leq 4A/3\leq 0.2406\vspace{2mm}\\
	0.174 \leq 3d/8 \leq A \leq 15d^4/4\leq 0.1804\vspace{2mm}\\
	0.1015 \leq \frac{a+A}{4} \leq 0.10525.
	\end{array}
	\end{equation}
	
	\noindent Using Lemma \ref{lem15} we need to consider $(a_5'',c_5)=(0,0),$ $(0,1/2),$ $(1/10,1/4).$\vspace{3mm}
	
	\noindent When $(a_5'',c_5)=(0,0)$, the following table gives a solution to (\ref{eq101'}), i.e., of
	\begin{equation}\label{EQ'''58'}0<G=(x_1+a_2''x_2)x_2-\frac{1}{5}x_5^2-\frac{a+A}{4}
	<d+\frac{3A}{4}.\end{equation}

	{\scriptsize
		\begin{equation*}\begin{array}{lllll}
		\hline\\
		{\rm ~~~~Range} & x_1 & x_1x_2 & x_5& {\rm ~~G}\vspace{2mm}\\
		\hline\\
		
		\frac{1}{5}+\frac{a+A}{4}<g_1''\leq \frac{1}{2}&c_1&-|c_1|&\pm 1&g_1''-\frac{1}{5}-\frac{a+A}{4}\vspace{2mm}\\
		
		-\frac{1}{5}+\frac{a+A}{4}<g_1''\leq \frac{1}{5}+\frac{a+A}{4}&\pm 1+c_1&1-|c_1|&\pm 2&g_1''+\frac{1}{5}-\frac{a+A}{4}\vspace{2mm}\\
		
		-\frac{4}{5}+\frac{a+A}{4}<g_1''<-\frac{4}{5}+d+\frac{a}{4}+A&\pm 1+c_1&1-|c_1|&\pm 1&g_1''+\frac{4}{5}-\frac{a+A}{4}\vspace{2mm}\\
		
		-1<g_1''\leq -\frac{4}{5}+\frac{a+A}{4}&\pm 2+c_1&2-|c_1|&\pm 2&g_1''+\frac{6}{5}-\frac{a+A}{4}\vspace{2mm}\\
		
		{\rm When~}  -\frac{4}{5}+d+\frac{a}{4}+A\leq g_1''\leq -\frac{1}{5}+\frac{a+A}{4}: \vspace{2mm}\\
		
		\frac{1}{5}+d+\frac{a}{4}+A\leq f_1''\leq 1&\pm 1+c_1&-2+2|c_1|&\pm 1&f_2''-\frac{11}{5}-\frac{a+A}{4}\vspace{2mm}\\
		
		\frac{1}{5}+\frac{a+A}{4}<f_1''<\frac{1}{5}+d+\frac{a}{4}+A&c_1&|c_1|&\pm 1&f_1''-\frac{1}{5}-\frac{a+A}{4}\vspace{2mm}\\
		
		-\frac{1}{5}+\frac{a+A}{4}<f_1''\leq \frac{1}{5}+\frac{a+A}{4}&\pm 1+c_1&1+|c_1|&\pm 2&f_1''+\frac{1}{5}-\frac{a+A}{4}\vspace{2mm}\\
		
		-\frac{1}{2}<f_1''\leq -\frac{1}{5}+\frac{a+A}{4}&\pm 2+c_1&4+2|c_1|&\pm 4&f_2''+\frac{4}{5}-\frac{a+A}{4}\vspace{2mm}\\
		
		\hline
		\end{array}
		\end{equation*}\vspace{2mm}}

	\noindent When $(a_5'',c_5)=(0,1/2)$, the following table gives a solution to (\ref{EQ'''58'}).\vspace{2mm}
	
	{\scriptsize
		\begin{equation*}\begin{array}{lllll}
		\hline\\
		{\rm ~~~~Range} & x_1 & x_1x_2 & x_5& {\rm ~~G}\vspace{2mm}\\
		\hline\\
		
		\frac{9}{20}+\frac{a+A}{4}<f_1''\leq 1&c_1&|c_1|&\pm \frac{3}{2}&f_1''-\frac{9}{20}-\frac{a+A}{4}\vspace{2mm}\\
		
		\frac{1}{20}+\frac{a+A}{4}<f_1''\leq \frac{9}{20}+\frac{a+A}{4}&c_1&|c_1|&\pm \frac{1}{2}&f_1''-\frac{1}{20}-\frac{a+A}{4}\vspace{2mm}\\
		
		-\frac{11}{20}+\frac{a+A}{4}<f_1''<-\frac{11}{20}+d+\frac{a}{4}+A&\pm 1+c_1&1+|c_1|&\pm \frac{3}{2}&f_1''+\frac{11}{20}-\frac{a+A}{4}\vspace{2mm}\\
		
		-\frac{1}{2}<f_1''\leq -\frac{11}{20}+\frac{a+A}{4}&\pm 1+c_1&1+|c_1|&\pm \frac{1}{2}&f_1''+\frac{19}{20}-\frac{a+A}{4}\vspace{2mm}\\
		
		{\rm When~}  -\frac{11}{20}+d+\frac{a}{4}+A\leq f_1''\leq \frac{1}{20}+\frac{a+A}{4}: \vspace{2mm}\\
		
		-\frac{11}{20}+d+\frac{a}{4}+A\leq g_1''\leq \frac{1}{20}+\frac{a+A}{4}&c_1&2|c_1|&\pm \frac{3}{2}&f_2''-\frac{9}{20}-\frac{a+A}{4}\vspace{2mm}\\
		
		-\frac{11}{20}+\frac{a+A}{4}<g_1''<-\frac{11}{20}+d+\frac{a}{4}+A&\pm 1+c_1&1-|c_1|&\pm \frac{3}{2}&g_1''+\frac{11}{20}-\frac{a+A}{4}\vspace{2mm}\\
		
		-\frac{19}{20}+\frac{a+A}{4}<g_1''\leq -\frac{11}{20}+\frac{a+A}{4}&\pm 1+c_1&1-|c_1|&\pm \frac{1}{2}&g_1''+\frac{19}{20}-\frac{a+A}{4}\vspace{2mm}\\
		
		-1<g_1''\leq -\frac{19}{20}+\frac{a+A}{4}&\pm 1+c_1&2+2|c_1|&\pm \frac{5}{2}&f_2''+\frac{3}{4}-\frac{a+A}{4}\vspace{2mm}\\
		
		\hline
		\end{array}
		\end{equation*}\vspace{1mm}}
	
	\noindent When $(a_5'',c_5)=(1/10,1/4)$, the following table gives a solution to (\ref{eq101'}), i.e., of \vspace{-2mm}
	\begin{equation}\label{EQ'''58'''}0<G=(x_1+a_2''x_2+\frac{1}{10}x_5)x_2-\frac{1}{5}x_5^2-
	\frac{a+A}{4}<d+\frac{3A}{4}.\end{equation}
	{\scriptsize
		\begin{equation*}\begin{array}{lllll}
		\hline\\
		{\rm ~~~~Range} & x_1 & x_2 & x_5& {\rm ~~G}\vspace{2mm}\\
		\hline\\
		
		\frac{63}{80}+\frac{a+A}{4}<p_1''\leq 1&c_1&1&\frac{9}{4}&p_1''-\frac{63}{80}-\frac{a+A}{4}\vspace{2mm}\\
		
		\frac{3}{16}+d+\frac{a}{4}+A\leq p_1''\leq \frac{63}{80}+\frac{a+A}{4}&-1+c_1&-1&\frac{5}{4}&q_1''+\frac{9}{16}-\frac{a+A}{4}\vspace{2mm}\\
		
		 \frac{3}{16}+\frac{a+A}{4}<p_1''<\frac{3}{16}+d+\frac{a}{4}+A&c_1&1&\frac{5}{4}&p_1''-\frac{3}{16}-\frac{a+A}{4}\vspace{2mm}\\
		
		-\frac{17}{80}+\frac{a+A}{4}<p_1''\leq \frac{3}{16}+\frac{a+A}{4}&1+c_1&1&\frac{9}{4}&p_1''+\frac{17}{80}-\frac{a+A}{4}\vspace{2mm}\\
		
		 -\frac{13}{16}+\frac{a+A}{4}<p_1''<-\frac{13}{16}+d+\frac{a}{4}+A&1+c_1&1&\frac{5}{4}&p_1''+\frac{13}{16}-\frac{a+A}{4}\vspace{2mm}\\
		
		-1<p_1''\leq -\frac{13}{16}+\frac{a+A}{4}&2+c_1&1&\frac{9}{4}&p_1''+\frac{97}{80}-\frac{a+A}{4}\vspace{2mm}\\
				
		{\rm When~}  -\frac{13}{16}+d+\frac{a}{4}+A\leq p_1''\leq -\frac{17}{80}+\frac{a+A}{4}: \vspace{2mm}\\
		
		\frac{7}{16}+\frac{a+A}{4}<q_1''\leq 1&c_1&-1&\frac{5}{4}&q_1''-\frac{7}{16}-\frac{a+A}{4}\vspace{2mm}\\
		
		\frac{3}{80}+\frac{a+A}{4}<q_1''\leq \frac{7}{16}+\frac{a+A}{4}&c_1&-1&\frac{1}{4}&q_1''-\frac{3}{80}-\frac{a+A}{4}\vspace{2mm}\\
		
		-\frac{9}{16}+d+\frac{a}{4}+A\leq q_1''\leq \frac{3}{80}+\frac{a+A}{4}&1+c_1&2&\frac{13}{4}&p_2''+\frac{43}{80}-\frac{a+A}{4}\vspace{2mm}\\
		
		 -\frac{9}{16}+\frac{a+A}{4}<q_1''<-\frac{9}{16}+d+\frac{a}{4}+A&-1+c_1&-1&\frac{5}{4}&q_1''+\frac{9}{16}-\frac{a+A}{4}\vspace{2mm}\\
		
		-\frac{77}{80}+\frac{a+A}{4}<q_1''\leq -\frac{9}{16}+\frac{a+A}{4}&-1+c_1&-1&\frac{1}{4}&q_1''+\frac{77}{80}-\frac{a+A}{4}\vspace{2mm}\\
		
		-1<q_1''\leq -\frac{77}{80}+\frac{a+A}{4}&1+c_1&2&\frac{9}{4}&p_2''+\frac{23}{16}-\frac{a+A}{4}\vspace{2mm}\\
		
		\hline
		\end{array}
		\end{equation*}}

	\begin{lemma}\label{lem82}
		If  $t=\frac{1}{6}$ and $ \pm 3a_5''-c_5\equiv 1/2 \pmod 1$, then inequality $(\ref{eq101'})$ is soluble  unless $(a_5'',c_5)=(\frac{1}{2},0)$, $c_1=\frac{1}{2}$ and $a_2''=\frac{1}{4}$. \end{lemma}
	
	\noindent {\bf Proof:}  Here, since $\frac{d}{a}\leq 2$, and $\frac{1}{6}=t=\frac{3d^5}{8aA}\leq 2d^3$, we get $d^3\geq \frac{1}{12}$ i.e., $d\geq (\frac{1}{12})^{\frac{1}{3}}>0.43679$. And from (\ref{Eq''85''}) we have $d\leq 128/246 \leq 0.52032$. Therefore,\vspace{-2mm}
	
	\begin{equation}\label{Eqn3}
	\begin{array}{l}
	0.21835\leq d/2\leq a\leq 4A/3\leq 16t/9\leq 0.2963\vspace{2mm}\\
	0.1637 \leq 3d/8 \leq A \leq 4t/3\leq 0.2223\vspace{2mm}\\
	0.0955125 \leq \frac{a+A}{4} \leq 0.12965.
	\end{array}
	\end{equation}
	\noindent Using Lemma \ref{lem15} we need to consider $(a_5'',c_5)=(0,1/2),$ and $(1/6,0).$ If $a_5''=1/6$, considering the transformation $x_5 \rightarrow x_5-x_2$, we may assume that $a_5''=1/2$.\vspace{2mm}
	
	\noindent When $(a_5'',c_5)=(0,1/2)$, the following table gives a solution to (\ref{eq101'}), i.e., of
	
	\begin{equation}\label{Eqn4}
	0<G=(x_1+a_2''x_2)x_2-\frac{1}{6}x_5^2-\frac{a+A}{4}<d+\frac{3A}{4}.
	\end{equation}

	{\scriptsize
		\begin{equation*}\begin{array}{lllll}
		\hline\\
		{\rm ~~~~Range} & x_1 & x_1x_2 & x_5& {\rm ~~G}\vspace{1mm}\\
		\hline\\
		
		\frac{3}{8}+\frac{a+A}{4}<f_1''\leq 1&c_1&|c_1|&\frac{3}{2}&f_1''-\frac{3}{8}-\frac{a+A}{4}\vspace{2mm}\\
		
		\frac{1}{24}+\frac{a+A}{4}<f_1''\leq \frac{3}{8}+\frac{a+A}{4}&c_1&|c_1|&\frac{1}{2}&f_1''-\frac{1}{24}-\frac{a+A}{4}\vspace{2mm}\\
		
		-\frac{5}{8}+\frac{a+A}{4}<f_1''<-\frac{5}{8}+d+\frac{a}{4}+A&\pm 1+c_1&1+|c_1|&\frac{3}{2}&f_1''+\frac{5}{8}-\frac{a+A}{4}\vspace{2mm}\\
		
		-\frac{1}{2}<f_1''\leq -\frac{5}{8}+\frac{a+A}{4}&\pm 1+c_1&1+|c_1|&\frac{1}{2}&f_1''+\frac{23}{24}-\frac{a+A}{4}\vspace{2mm}\\
		
		{\rm When~}  -\frac{5}{8}+d+\frac{a}{4}+A\leq f_1''\leq \frac{1}{24}+\frac{a+A}{4}: \vspace{2mm}\\
		
		\left.\begin{array}{l}-\frac{5}{8}+d+\frac{a}{4}+A\leq g_1''\leq \frac{1}{24}+\frac{a+A}{4}\vspace{1mm}\\\frac{1}{24}+\frac{a+A}{4}<f_2''\end{array}\right\}& c_1&2|c_1|&\frac{1}{2}&f_2''-\frac{1}{24}-\frac{a+A}{4}\vspace{2mm}\\
		
		\left.\begin{array}{l}-\frac{5}{8}+d+\frac{a}{4}+A\leq g_1''\leq \frac{1}{24}+\frac{a+A}{4}\vspace{1mm}\\f_2''\leq \frac{1}{24}+\frac{a+A}{4}\end{array}\right\}& c_1&3|c_1|&\frac{1}{2}&f_3''-\frac{1}{24}-\frac{a+A}{4}\vspace{2mm}\\

		-\frac{5}{8}+\frac{a+A}{4}<g_1''<-\frac{5}{8}+d+\frac{a}{4}+A&\pm 1+c_1&1-|c_1|&\frac{3}{2}&g_1''+\frac{5}{8}-\frac{a+A}{4}\vspace{2mm}\\
		
		-\frac{23}{24}+\frac{a+A}{4}<g_1''<-\frac{5}{8}+\frac{a+A}{4}&\pm 1+c_1&1-|c_1|&\frac{1}{2}&g_1''+\frac{23}{24}-\frac{a+A}{4}\vspace{2mm}\\
		
		\left.\begin{array}{l}-1<g_1'' \leq -\frac{23}{24}+\frac{a+A}{4}\vspace{1mm}\\-\frac{23}{24}+\frac{a+A}{4}<f_2''\end{array}\right\}& \pm 1+c_1&2+2|c_1|&\frac{5}{2}&f_2''+\frac{23}{24}-\frac{a+A}{4}\vspace{2mm}\\
		
		\left.\begin{array}{l}-1<g_1'' \leq -\frac{23}{24}+\frac{a+A}{4}\vspace{1mm}\\f_2''\leq -\frac{23}{24}+\frac{a+A}{4}\end{array}\right\}& \pm 1+c_1&3+3|c_1|&\frac{1}{2}&f_3''+\frac{71}{24}-\frac{a+A}{4}\vspace{2mm}\\
		
		\hline
		\end{array}
		\end{equation*}\vspace{1mm}}

	\noindent When $(a_5'',c_5)=(1/2,0)$, the following table gives a solution to (\ref{eq101'}) i.e., of
	
	\begin{equation}\label{Eqn5}
	0<G=(x_1+a_2''x_2+\frac{1}{2}x_5)x_2-\frac{1}{6}x_5^2-\frac{a+A}{4}<d+\frac{3A}{4}.
	\end{equation}
	
	{\scriptsize
		\begin{equation*}\begin{array}{llllll}
		\hline\\
		{\rm ~~~~Range} & x_1 & x_1x_2 & x_5&x_2x_5& {\rm ~~G}\vspace{1mm}\\
		\hline\\
		
		\frac{a+A}{4}<g_1''\leq \frac{1}{2}&c_1&-|c_1|&0&0&g_1''-\frac{a+A}{4}\vspace{2mm}\\
		
		-\frac{1}{3}+\frac{a+A}{4}<g_1''\leq \frac{a+A}{4}&c_1&-|c_1|&\pm 1&1&g_1''+\frac{1}{3}-\frac{a+A}{4}\vspace{2mm}\\
		
		-1+\frac{a+A}{4}<g_1''<-1+d+\frac{a}{4}+A&\pm 1+c_1&1-|c_1|&0&0&g_1''+1-\frac{a+A}{4}\vspace{2mm}\\
		
		-1<g_1''\leq -1+\frac{a+A}{4}&\pm 1+c_1&1-|c_1|&\pm 1&1&g_1''+\frac{4}{3}-\frac{a+A}{4}\vspace{2mm}\\
		
		{\rm When~}  -1+d+\frac{a}{4}+A\leq g_1''\leq -\frac{1}{3}+\frac{a+A}{4}: \vspace{2mm}\\
		
		\frac{a+A}{4}<f_1''<d+\frac{a}{4}+A&c_1&|c_1|&0&0&f_1''-\frac{a+A}{4}\vspace{2mm}\\
		
		-\frac{1}{3}+\frac{a+A}{4}<f_1''\leq \frac{a+A}{4}&c_1&|c_1|&\pm 1&1&f_1''+\frac{1}{3}-\frac{a+A}{4}\vspace{2mm}\\
		
		-\frac{1}{2}<f_1''\leq -\frac{1}{3}+\frac{a+A}{4}&c_1&2|c_1|&\pm 3&6&f_2''+\frac{3}{2}-\frac{a+A}{4}\vspace{2mm}\\
		
		\hline
		\end{array}
		\end{equation*}\vspace{1mm}}
	
	\noindent  except  when $ d+\frac{a}{4}+A\leq f_1''\leq 1$ and $-1+d+\frac{a}{4}+A\leq  g_1''\leq -\frac{1}{3}+\frac{a+A}{4}$.\vspace{10mm}
	
	\noindent So, let now
	\begin{equation}\label{Eqn6}
	\begin{array}{l}
	g_1''\leq -\frac{1}{3}+\frac{a+A}{4}, ~~f_1'' \geq d+\frac{a}{4}+A.
	\end{array}
	\end{equation}
	
	\noindent We distinguish the cases $a_2''>\frac{1}{4}$ and $a_2''\leq \frac{1}{4}$ and work as in Lemma 9 of Raka and Rani \cite{RakaRani}.\vspace{2mm}
	
	\noindent \textbf{Case (i):} $a_2''>\frac{1}{4}.$\vspace{2mm}
	
	\noindent Let $|c_1|=\frac{1}{2}-\epsilon_1$ and $a_2''=\frac{1}{4}+\epsilon_2$, where $\epsilon_1 \geq 0,\epsilon_2>0$. Then $g_1''\leq -\frac{1}{3}+\frac{a+A}{4}$ gives $\epsilon_1+\epsilon_2\leq -\frac{1}{12}+\frac{a+A}{4}$. Chose an integer $\ell(\geq 1)$ such that\vspace{-2mm}
	\begin{equation}\label{Eqn7}
	\epsilon_1.3^{\ell-1}+\epsilon_2.(3^{\ell-1})^2 \leq \frac{a+A}{4}-\frac{1}{12}<\epsilon_1.3^\ell + \epsilon_2.(3^\ell)^2.
	\end{equation}
	
	\noindent Let $p$ be an integer such that (such an integer exists)\vspace{-2mm}
	\begin{equation}\label{Eqn8}
	2p^2+1 \equiv 0 \pmod {3^{\ell+1}}.
	\end{equation}
	
	\noindent Clearly $f(p,\ell)=1+2p^2+6.3^\ell-3.3^{2\ell}+6.3^\ell p\equiv 0 \pmod {3^{\ell+1}}$. Also, $f(p,\ell)=2p(p+3^{\ell+1})-3(3^\ell-1)^2+4 \equiv 0 \pmod 4$. Therefore, there exists an integer $q$ such that $f(p,\ell)=3^{\ell+1}4q$, so that\vspace{-2mm} $$q3^\ell-\frac{3^\ell}{2}+\frac{3^{2\ell}}{4}-\frac{3^\ell p}{2}-\frac{p^2}{6}=\frac{1}{12}.$$
	
	\noindent Take $x_2=\pm 3^\ell, x_5=\pm p, x_1=\pm q+c_1$ such that $x_1x_2=q.3^\ell-3^\ell|c_1|, x_2x_5=-p3^\ell$, so that
	\begin{equation*}\begin{array}{ll}
	 G&=q.3^\ell-3^\ell|c_1|+3^{2\ell}a_2''-p.\frac{3^\ell}{2}-\frac{p^2}{6}-\frac{a+A}{4}\vspace{2mm}\\
	
	&=\frac{1}{12}+3^\ell\epsilon_1+3^{2\ell}\epsilon_2-\frac{a+A}{4}\vspace{2mm}\\
	
	&>0.
	\end{array}
	\end{equation*}
	
	\noindent Also,
	\begin{equation*}\begin{array}{ll}
	 G&=\big(\frac{3^\ell}{3^{\ell-1}}\big)^2\big(\epsilon_1\frac{(3^{\ell-1})^2}{3^\ell}+\epsilon_2(3^{\ell-1})^2\big)
	-(\frac{a+A}{4}-\frac{1}{12})\vspace{2mm}\\
	
	 &<9\big(\epsilon_13^{\ell-1}+\epsilon_2(3^{\ell-1})^2\big)-(\frac{a+A}{4}-\frac{1}{12})\vspace{2mm}\\
	
	&=2(a+A)-2/3 < d+\frac{3A}{4}.
	\end{array}
	\end{equation*}

	\noindent \textbf{Case (ii):} $a_2''\leq \frac{1}{4}.$\vspace{2mm}
	
	\noindent Here we modify the technique used in Lemma 9 of Raka and Rani \cite{RakaRani}.\vspace{2mm}
	
	\noindent Let $|c_1|=\frac{1}{2}-\epsilon_1$ and $a_2''=\frac{1}{4}-\epsilon_2$, where   $\epsilon_1\geq 0,\epsilon_2\geq 0$. Then $f_1''\geq d+\frac{a}{4}+A$ gives $\epsilon_1+\epsilon_2\leq \frac{3}{4}-(d+\frac{a}{4}+A)< \frac{5}{6}-(d+\frac{a}{4}+A)$. Suppose $(\epsilon_1,\epsilon_2)\neq (0,0)$. Choose an integer $\ell(\geq 1)$ such that \vspace{-2mm}
	\begin{equation}\label{Eqn9}
	\epsilon_1.3^{\ell-1}+\epsilon_2.(3^{\ell-1})^2 \leq \frac{5}{6}-(d+\frac{a}{4}+A)<\epsilon_1.3^\ell + \epsilon_2.(3^\ell)^2.
	\end{equation}
	
	\noindent We  divide the above range into two parts and will take  different values of $x_1, x_2, x_5$ for a solution of (\ref{eq101'}).
	\begin{equation}\label{1}\begin{array}{l}
	{\rm \textbf{I}~~}\epsilon_1.3^{\ell-1}+\epsilon_2.(3^{\ell-1})^2\leq \frac{5}{6}-(d+\frac{a}{4}+A)<2.\epsilon_1.3^{\ell-1}+4.\epsilon_2.(3^{\ell-1})^2\end{array}
	\end{equation}
	\begin{equation}\label{2}\begin{array}{l}
	{\rm \textbf{II}~~~~~~~}2.\epsilon_1.3^{\ell-1}+4.\epsilon_2.(3^{\ell-1})^2 \leq \frac{5}{6}-(d+\frac{a}{4}+A)< \epsilon_1.3^\ell+\epsilon_2.(3^{\ell})^2.
	\end{array}
	\end{equation}

	\noindent \noindent \textbf{Part I:} ~~If (\ref{1}) holds, take $x_2=2.3^{\ell-1}, x_1=\pm q+c_1,  x_5=\pm p$ such that $x_1x_2= q.2.3^{\ell-1}+|c_1|.2.3^{\ell-1}, x_2x_5=-2.3^{\ell-1}.p$, where $p$ is an \textbf{odd} integer  satisfying $p^2+5 \equiv 0 \pmod {3^\ell}$ (such an integer exists). Now choose an  integer $q$ such that
	$$ q.2.3^{\ell-1}+ 3^{\ell-1}+(3^{\ell-1})^2-3^{\ell-1}p-\frac{p^2}{6}=\frac{5}{6}.$$
	This is possible as $f(p,\ell)=5+p^2-6.3^{\ell-1}-6(3^{\ell-1})^2-6.3^{\ell-1}p\equiv 0 \pmod {3^\ell}$ and $f(p,\ell) \equiv 0 \pmod 4$. Then
	\begin{equation*}\begin{array}{ll}
	G&= q.2.3^{\ell-1}+|c_1|.2.3^{\ell-1}+a_2''.4.(3^{\ell-1})^2-3^{\ell-1}p-\frac{p^2}{6}-\frac{a+A}{4}\vspace{2mm}\\&=
	q.2.3^{\ell-1}+ 3^{\ell-1}+(3^{\ell-1})^2-3^{\ell-1}p-\frac{p^2}{6}-\big(2.\epsilon_1.3^{\ell-1}+4.\epsilon_2.(3^{\ell-1})^2\big)-\frac{a+A}{4}
	\vspace{2mm}\\
	
	 &=\frac{5}{6}-\big(2.\epsilon_1.3^{\ell-1}+4.\epsilon_2.(3^{\ell-1})^2)-\frac{a+A}{4}\vspace{2mm}\\
	
	&<d+\frac{a}{4}+A-\frac{a+A}{4}=d+\frac{3A}{4}.
	\end{array}
	\end{equation*}
	\noindent Also,
	\begin{equation*}\begin{array}{ll}
	 G&=\frac{5}{6}-(2.\epsilon_1.3^{\ell-1}+4.\epsilon_2.(3^{\ell-1})^2)-\frac{a+A}{4}\vspace{2mm}\\
	
	 &=\frac{5}{6}-\frac{a+A}{4}-4\Big(\frac{\epsilon_1}{2}.3^{\ell-1}+\epsilon_2.(3^{\ell-1})^2\Big)\vspace{2mm}\\
	
	 &>\frac{5}{6}-\frac{a+A}{4}-4\Big(\epsilon_1.3^{\ell-1}+\epsilon_2.(3^{\ell-1})^2\Big)\vspace{2mm}\\
	
	&\geq \frac{5}{6}-\frac{a+A}{4}-4(\frac{5}{6}-(d+\frac{a}{4}+A)) > 0.
	\end{array}
	\end{equation*}
	
	\noindent \textbf{Part II:}~~ If (\ref{2}) holds, take  $x_2=\pm 3^\ell,~ x_1=\pm q+c_1, ~ x_5=\pm p$ such that $x_1x_2=q.3^\ell+3^\ell.|c_1|,~ x_2x_5=-p.3^\ell$,  where $p$ is an  integer  satisfying $2p^2+9 \equiv 0 \pmod {3^{\ell+1}}$ (such an integer exists). Now choose an  integer $q$ such that $$ q.3^\ell+\frac{1}{2}3^\ell+\frac{1}{4}3^{2\ell}-\frac{1}{2}p.3^\ell-\frac{p^2}{6}=\frac{3}{4}.$$ This is possible as $f(p,\ell)=9+2p^2-12(\frac{1}{2}3^\ell+\frac{1}{4}3^{2\ell}-\frac{1}{2}p.3^\ell) \equiv 0 \pmod {3^{\ell+1}}$ and $f(p,\ell) \equiv 0 \pmod 4$. Then
	\begin{equation*}\begin{array}{ll}
	G&= q.3^\ell+3^\ell.|c_1|+3^{2\ell}a_2''-\frac{1}{2}p.3^\ell-\frac{p^2}{6}-\frac{a+A}{4}\vspace{2mm}\\
	&=q.3^\ell+\frac{1}{2}3^\ell+\frac{1}{4}3^{2\ell}-\frac{1}{2}p.3^\ell-\frac{p^2}{6}-
	(3^\ell\epsilon_1+3^{2\ell}\epsilon_2)-\frac{a+A}{4}\vspace{2mm}\\
	&=\frac{3}{4}-3^\ell\epsilon_1-3^{2\ell}\epsilon_2-\frac{a+A}{4}\vspace{2mm}\\
	
	&<\frac{5}{6}-(3^\ell\epsilon_1+3^{2\ell}\epsilon_2)-\frac{a+A}{4}\vspace{2mm}\\
	
	&<d+\frac{3A}{4}.
	\end{array}
	\end{equation*}
	\noindent Also,
	\begin{equation*}\begin{array}{ll}
	G&=\frac{3}{4}-(3^\ell\epsilon_1+3^{2\ell}\epsilon_2)-\frac{a+A}{4}\vspace{2mm}\\
	
	 &=\frac{3}{4}-\frac{a+A}{4}-\frac{9}{4}\left(\frac{4}{3}(3^{\ell-1}\epsilon_1)+4(3^{\ell-1})^2\epsilon_2\right)\vspace{2mm}\\
	
	&>\frac{3}{4}-\frac{a+A}{4}-\frac{9}{4}\left( 2\epsilon_1.3^{\ell-1} +4\epsilon_2.(3^{\ell-1})^2 \right)\vspace{2mm}\\
	
	&> \frac{3}{4}-\frac{a+A}{4}-\frac{9}{4}\left( \frac{5}{6}-(d+\frac{a}{4}+A) \right) \vspace{2mm}\\&>\frac{101}{32}d-\frac{9}{8}>0.
	
	\end{array}
	\end{equation*}
	
	\noindent Thus, (\ref{Eqn5}) is soluble unless $(\epsilon_1,\epsilon_2)=(0,0)$, i.e., unless $c_1=\frac{1}{2}$ and $a_2''=\frac{1}{4}$.\hfill $\square$
	
	\begin{lemma}\label{lem82'}  If  $t=\frac{1}{6}$, $(a_5'',c_5)=(\frac{1}{2},0)$,  $c_1=\frac{1}{2}$ and $a_2''=\frac{1}{4}$, then inequality $(\ref{eq101})$ is soluble.\end{lemma}
	
	\noindent {\bf Proof:}~ From (\ref{eq16}), we have
	$$\frac{1}{4}=a_2''\equiv a_2'+\frac{a_4'^2}{4A}\hspace{-2mm}\pmod 1,~~~ \frac{1}{2}=a_5''\equiv a_5'-\lambda a_4'\hspace{-2mm}\pmod 1,$$
	
	\noindent and since $0.43679 \leq d\leq 0.52032$, we have $0.16379\leq 3d/8\leq A\leq 4t/3=2/9<0.2223$ and $0.218395\leq d/2\leq a\leq (3d^5/4)^{1/3} \leq 0.3059$. So that,
	$$-1/2<a_2'\leq a_2'+a_4'^2/4A \leq a_2'+/16A <1,~ -3/4\leq a_5'-\lambda a_4'\leq 3/4$$
	\noindent which gives
	\begin{equation}\label{Eqn12}
	a_2'+\frac{a_4'^2}{4A}=\frac{1}{4} {\rm ~
		~and~} a_5'-\lambda a_4'=\frac{1}{2} {\rm ~ or} -\frac{1}{2}.
	\end{equation}
	
	\noindent Note that,
	\begin{equation*}\label{Eqn13}
	\begin{array}{l}
	A\leq \left\{
	\begin{array}{ll}
	\frac{4t}{3}=\frac{2}{9} & {\rm ~if~} d>\frac{4}{9}\vspace{2mm}\\
	\frac{9d^4}{2}<{\tiny 0.1756} &{\rm ~if~} d\leq \frac{4}{9}.
	\end{array}
	\right.\vspace{2mm}
	\end{array}
	\end{equation*}
	
	\noindent This gives
	\begin{equation*}\label{Eqn14}
	\begin{array}{l}
	C=t+A\lambda^2 \leq \left\{
	\begin{array}{ll}
	\frac{2}{9} & {\rm ~if~} d>\frac{4}{9}\vspace{2mm}\\
	{\tiny 0.211} &{\rm ~if~} d\leq \frac{4}{9}
	\end{array}
	\right.\vspace{2mm}
	\end{array}
	\end{equation*}
	
	\noindent so that, $2<\frac{d}{C}=\frac{\delta_m}{C}\leq \frac{d}{A}<3$, i.e., $M=2$. And also, using $A\leq C$ $$\frac{1}{8}<\frac{A}{\frac{3A}{2}+1}\leq \frac{A}{6A\lambda^2 +1}=\frac{A}{6C}=\frac{At}{C}\leq t=\frac{1}{6}=\lambda_0~( {\rm say}).$$
	
	\noindent We divide this range of $\frac{At}{C}$ into $17$ subintervals $[\lambda_n,\lambda_{n-1}], n=1,2,\cdots 17$ and in each subinterval choose suitable integers $2h_n$ and $k_n$ such that $\vline~ h_n-\frac{At}{C}k_n^2~\vline +\frac{1}{2} <\delta_{m,M}^*=\delta_{1,2}^*=d+\frac{3C}{4}$.
	The choice of $(h_n,k_n)$ is done in a manner described earlier.\vspace{1mm}

	{\scriptsize
		\begin{equation*}\begin{array}{lllllllll}
		\hline\\
		n & (h_n, k_n) & \lambda_n&\mbox{Remarks}&&n & (h_n, k_n) & \lambda_n&\mbox{Remarks}\vspace{2mm}\\
		\hline\\
		
		1&(1.5,3)&0.16&\mbox{tbd}&&2&(2.5,4)&0.153&\mbox{tbd}\vspace{2mm}\\
		
		3&(5.5,6)&0.1511&\mbox{na}&&4&(25.5,13)&0.1506&(34,15)\vspace{2mm}\\
		
		5&(15,10)&0.1493&\mbox{tbd}&&6&(9.5,8)&0.1474&(12,9)\vspace{2mm}\\
		
		7&(65,21)&0.14726&(12,9)&&8&(42.5,17)&0.1468&\mbox{tbd}\vspace{2mm}\\
		
		9&(33,15)&0.1463&(37.5,16)&&10&(58.5,20)&0.14609&(37.5,16)\vspace{2mm}\\
		
		11&(21,12)&0.1454&\mbox{tbd}&&12&(14.5,10)&0.1443&(17.5,11)\vspace{2mm}\\
		
		13&(17.5,11)&0.14412&\mbox{na}&&14&(52,19)&0.1438&\mbox{na}\vspace{2mm}\\
		
		15&(7,7)&0.1415&(0.5,2)&&16&(3.5,5)&0.1375&\mbox{na}\vspace{2mm}\\
		
		17&(0.5,2)&0.125&\mbox{na}\vspace{2mm}\\

		\hline
		\end{array}
		\end{equation*}\vspace{1mm}}
	
	\noindent Thus, by Macbeath's Lemma, (\ref{e17}) and hence (\ref{eq101}) is soluble unless $\frac{At}{C}=1/6,$ $5/32,$ $3/20,$ $5/34,$ $7/48$ and $\beta_{G}^* = \pm a_4^*-\frac{2At}{C}c_4\equiv h_n/k_n \pmod {1/k_n, \frac{2At}{C}}$. Further, for a particular pair $(h,k)$ we have $\frac{A}{6C}=\frac{h}{k^2}$ i.e., $C=\frac{A}{6}.\frac{k^2}{h}\geq \frac{d}{16}.\frac{k^2}{h}$. Hence, $\delta_{1,2}^*\geq (1+\frac{3\cdot k^2}{64\cdot h})d$. For all these values of $\frac{At}{C}$, we have $\frac{d}{16}.\frac{k^2}{h}\leq C\leq \frac{2}{9}$ and Lemma \ref{lem15} fixes the values of $(a_4^*,c_4)$ to be considered.
	For each of these values of $\frac{At}{C}$, except when $\frac{At}{C}=1/6$, we find that (\ref{e17}) is soluble, by giving suitable values to $x_2$ and $x_4$. \vspace{2mm}
	
	\noindent When $A/6C=1/6$, i.e., $A=C$ and $\beta_G^*\equiv 1/2 \pmod {1/3}$,  we need to consider $(a_4^*,c_4)=(1/2,0)$ or $(0,1/2)$. Proceeding similar to that in Lemma \ref{lem82} (now working for $G^*$ instead of $G$), we find that (\ref{e17}) is soluble unless $(a_4^*,c_4)=(1/2,0)$ and $a_2^*=1/4$. Here also, we have, by (\ref{G*})
	
	\begin{equation}\label{Eqn16}
	a_2'+\frac{a_5'^2}{4A}=\frac{1}{4} {\rm ~
		~and~} a_4'-\lambda a_5'=\frac{1}{2} {\rm ~ or} -\frac{1}{2}.
	\end{equation}
	
	\noindent From (\ref{Eqn12}) and (\ref{Eqn16}), we get $ a_4'^2=a_5'^2$, i.e., $a_4'=\pm a_5'$. As, $A=C=A\lambda^2+1/6$, we get $A=C=1/6(1-\lambda^2)$. Also, we have $0\leq \lambda \leq 1/2$.\vspace{2mm}
	
	\noindent \textbf{(i)} When $a_4'=a_5'$, we have from (\ref{G*}) that $a_5'=\pm 1/2(1-\lambda)$ giving thereby $|a_5'|\geq 1/2$, i.e., $a_5'=\pm 1/2=a_4'$. Taking the unimodular transformation $x_2\rightarrow -x_2,~ x_4\rightarrow -x_4,~ x_5\rightarrow -x_5$, we need to consider only $a_4'=1/2=a_5'$.  In this case, from (\ref{Eqn16}), we get $\lambda=0$ and $a_2'=-1/8$. Then, $x_1=1/2, x_2=1, x_4=0=x_5$ gives a solution of (\ref{eq101}).\vspace{2mm}
	
	\noindent \textbf{(ii)} When $a_4'=-a_5'=\pm 1/2(1+\lambda)$, consider without loss of generality, $a_4'=1/2(1+\lambda), a_5'=-1/2(1+\lambda)$. Then, $a_2'=1/4-a_5'^2/4A=(5\lambda-1)/8(1+\lambda)$. In different subintervals of  $\lambda$, the following  gives a solution of (\ref{eq101}).
	\begin{itemize}
		\item When $0\leq \lambda \leq 0.182$, take $x_1=1/2, x_2=1, x_4=0, x_5=0$.
		\item When $0.182\leq \lambda \leq 0.275$, take $x_1=1/2, x_2=1, x_4=1, x_5=1$.
		\item When $0.275\leq \lambda \leq 0.425$, take $x_1=1/2, x_2=2, x_4=0, x_5=1$.
		\item When $0.425\leq \lambda \leq 0.5$, take $x_1=-1/2, x_2=3, x_4=1, x_5=0$. \hfill $\square$
	\end{itemize}
	
	\subsection{${m=1, K=2, L=2}$}
	
	\noindent Here $2<\frac{\delta_{1,2}}{t}=\frac{d+3A/4}{t}\leq 3$, $\delta_{1,2}=d+\frac{3A}{4}\geq d+\frac{9d}{32}=\frac{41d}{32}$. From (\ref{Eq74'}) we have
	\begin{equation}\label{Eq85''}
	0.4621< \sqrt{41/192} \leq d \leq 128/205<0.624391.
	\end{equation}
	\noindent Further, since $a\geq d/2$ and $A\geq 3d/8$ we get, $\frac{41d}{96}\leq \frac{\delta_{1,2}}{3}\leq t =\frac{3d^5}{8aA}\leq 2d^3$.\vspace{2mm}
	
	\noindent Also from (\ref{EqNn2}), we get $t\leq 4(1-41d/32)/3$. Therefore,
	\begin{equation}\label{Eq'89''}
	\begin{array}{ll}
	0.19735< 41d/96 \leq t\leq & \min\{2d^3,4(1-41d/32)/3\}\vspace{2mm}\\
	&=\left\{
	\begin{array}{ll}
	2d^3 & {\rm ~if~} d\leq \alpha\\
	4(1-41d/32)/3 &{\rm ~if~} d\geq \alpha
	\end{array}
	\right.\vspace{2mm}\\
	&\leq 0.3646=\lambda_0~({\rm say})
	\end{array}
	\end{equation}
	
	\noindent where $\alpha$ is a root of $2d^3-4(1-41d/32)/3=0$ satisfying $0.567<\alpha<0.568$. We divide this range of $t$ into $9$ subintervals $[\lambda_n,\lambda_{n-1}], n=1,2,\cdots 9$ and in each subinterval choose suitable integers $2h_n$ and $k_n$ satisfying $\vline~ h_n-tk_n^2~\vline +\frac{1}{2} <d+\frac{3A}{4}.$\vspace{1mm}
		
	{\scriptsize
		\begin{equation*}\begin{array}{lllllllll}
		\hline\\
		n & (h_n, k_n) & \lambda_n&\mbox{Remarks}&&n & (h_n, k_n) & \lambda_n&\mbox{Remarks}\vspace{2mm}\\
		\hline\\
		
		1&(0.5,1)&0.3111&\mbox{na}&&2&(5,4)&0.3012&\mbox{na}\vspace{2mm}\\
		
		3&(7.5,5)&0.2976&\mbox{tbd}&&4&(2.5,3)&0.294&\mbox{na}\vspace{2mm}\\
		
		5&(1,2)&0.222&\mbox{tbd}&&6&(2,3)&0.211&\mbox{na}\vspace{2mm}\\	
		
		7&(7.5,6)&0.2056&\mbox{tbd}&&8&(10,7)&0.2039&(13,8)\vspace{2mm}\\
		
		9&(5,5)&0.1973&\mbox{tbd}&&&&&\vspace{2mm}\\
		
		\hline
		\end{array}
		\end{equation*}\vspace{1mm}}
	
	\noindent The choice of $(h_n,k_n)$ is done in a manner already described. Thus (\ref{eq101'}) is soluble unless $t=3/10,$ $1/4,$ $5/24,$ $1/5$ and $\beta_{t,G} = \pm a_5''-2tc_5\equiv h_n/k_n \pmod {1/k_n, 2t}$. For a fixed $t$ we get \begin{equation}\label{Eqn_1'}
	\begin{array}{l}
	(\frac{t}{2})^{1/3} \leq d \leq \frac{128t}{41},~~
	d/2 \leq a\leq  \frac{4A}{3},~~
	3d/8 \leq  A\leq \frac{4t}{3}.
	
	\end{array}
	\end{equation}
	
	\noindent For $t=3/10$ and $5/24$, we find that (\ref{eq101'}) is soluble by taking $x_2=\pm 1$ and choosing $x_5 \equiv c_5 \pmod 1$ suitably. For $t=\frac{1}{5}$, we need to consider $(a_5'',c_5)=(0,0),$ $(0,1/2),$ $(1/10,1/4),$ then tables similar to those in Lemma \ref{lem81} give a solution to (\ref{eq101'}).
	
	\begin{lemma}\label{lem28} When $t=\frac{1}{4}$ and $\beta_{t,G} = \pm a_5''-2tc_5\equiv 1/2 \pmod {1/2, 2t}$, $(\ref{eq101})$ is soluble.
		
	\end{lemma}
	
	\noindent {\bf Proof:} Here, since $\frac{d}{a}\leq 2$, and $\frac{1}{4}=t=\frac{3d^5}{8aA}\leq 2d^3$, we get $d^3\geq \frac{1}{8}$ i.e., $d\geq ({1}/{8})^{{1}/{3}}=0.5$. And from (\ref{Eq85''}), we have $d\leq \frac{128}{205}$. Also, $2<\frac{d}{A}=\frac{\delta_1}{A}\leq 3$ gives $A< d/2< 0.3121955$. Therefore,
	\begin{equation}\label{Eq''134''}
	\begin{array}{l}
	0.25< d/2\leq a\leq 4A/3< 0.416261\vspace{2mm}\\
	0.1875< 3d/8 \leq A < d/2 < 0.3121955\vspace{2mm}\\
	0.109375< \frac{a+A}{4} < 0.18212.
	\end{array}
	\end{equation}
	
	\noindent Using Lemma \ref{lem15} we need to consider $(a_5'',c_5)=(0,0),$ $(1/4,1/2).$\vspace{2mm}

	\noindent When $(a_5'',c_5)=(1/4,1/2)$, the following table gives a solution to (\ref{eq101'}) i.e., of\\
	 \begin{equation}\label{EQ'''58''}0<G=(x_1+a_2''x_2+\frac{1}{4}x_5)x_2-\frac{1}{4}x_5^2-\frac{a+A}{4}<d+\frac{3A}{4}.\end{equation}

	{\scriptsize
		\begin{equation*}\begin{array}{llllll}
		\hline\\
		{\rm ~~~~Range} & x_1 & x_1x_2 & x_5&x_2x_5& {\rm ~~G}\vspace{2mm}\\
		\hline\\

		\frac{3}{16}+d+\frac{a}{4}+A\leq f_1''\leq 1& \pm 1+c_1&-2+2|c_1|&\pm \frac{5}{2}&5&f_2''-\frac{37}{16}-\frac{a+A}{4}\vspace{2mm}\\
		
		\frac{3}{16}+\frac{a+A}{4}<f_1''<\frac{3}{16}+d+\frac{a}{4}+A&c_1&|c_1|&\pm \frac{3}{2}&\frac{3}{2}&f_1''-\frac{3}{16}-\frac{a+A}{4}\vspace{2mm}\\
		
		-\frac{1}{16}+\frac{a+A}{4}<f_1''\leq \frac{3}{16}+\frac{a+A}{4}&c_1&|c_1|&\pm \frac{1}{2}&\frac{1}{2}&f_1''+\frac{1}{16}-\frac{a+A}{4}\vspace{2mm}\\	
		
		-\frac{1}{2}<f_1''<-\frac{13}{16}+d+\frac{a}{4}+A&\pm 1+c_1&1+|c_1|&\pm \frac{3}{2}&\frac{3}{2}&f_1''+\frac{13}{16}-\frac{a+A}{4}\vspace{2mm}\\

		\hline
		\end{array}
		\end{equation*}\vspace{1mm}}
	
	{\scriptsize
		\begin{equation*}\begin{array}{llllll}
		\hline\\
		{\rm ~~~~Range} & x_1 & x_1x_2 & x_5&x_2x_5& {\rm ~~G}\vspace{2mm}\\
		\hline\\
		
		{\rm When~}  -\frac{13}{16}+d+\frac{a}{4}+A\leq f_1''\leq -\frac{1}{16}+\frac{a+A}{4}: \vspace{2mm}\\

		\left.\begin{array}{l}-\frac{13}{16}+d+\frac{a}{4}+A\leq g_1''\leq -\frac{1}{16}+\frac{a+A}{4}\vspace{1mm}\\-\frac{3}{16}+\frac{a+A}{4}<f_2''\end{array}\right\}& c_1&2|c_1|&\pm \frac{1}{2}&1&f_2''+\frac{3}{16}-\frac{a+A}{4}\vspace{2mm}\\
		
		\left.\begin{array}{l}-\frac{13}{16}+d+\frac{a}{4}+A\leq g_1''\leq -\frac{1}{16}+\frac{a+A}{4}\vspace{1mm}\\f_2''\leq -\frac{3}{16}+\frac{a+A}{4}\end{array}\right\}& \pm 1+c_1&2+2|c_1|&\pm \frac{3}{2}&-3&f_2''+\frac{11}{16}-\frac{a+A}{4}\vspace{2mm}\\	
		
		-\frac{13}{16}+\frac{a+A}{4}<g_1''<-\frac{13}{16}+d+\frac{a}{4}+A&\pm 1+c_1&1-|c_1|&\pm \frac{3}{2}&\frac{3}{2}&g_1''+\frac{13}{16}-\frac{a+A}{4}\vspace{2mm}\\
		
		-\frac{17}{16}+\frac{a+A}{4}<g_1''\leq -\frac{13}{16}+\frac{a+A}{4}&\pm 1+c_1&1-|c_1|&\pm \frac{1}{2}&\frac{1}{2}&g_1''+\frac{17}{16}-\frac{a+A}{4}\vspace{2mm}\\
		
		\left.\begin{array}{l}-1<g_1''\leq -\frac{17}{16}+\frac{a+A}{4}\vspace{1mm}\\-\frac{19}{16}+\frac{a+A}{4}<f_2''\end{array}\right\}&\pm 2+ c_1&4+2|c_1|&\pm \frac{5}{2}&-5&f_2''+\frac{19}{16}-\frac{a+A}{4}\vspace{2mm}\\
		
		\left.\begin{array}{l}-1<g_1''\leq -\frac{17}{16}+\frac{a+A}{4}\vspace{1mm}\\f_2''\leq -\frac{19}{16}+\frac{a+A}{4}\end{array}\right\}& \pm 1+c_1&2+2|c_1|&\pm \frac{1}{2}&-1&f_2''+\frac{27}{16}-\frac{a+A}{4}\vspace{2mm}\\	
		
		\hline
		\end{array}
		\end{equation*}\vspace{1mm}}
	
	\noindent When $(a_5'',c_5)=(0,0)$, the following table gives a solution to (\ref{eq101'}) i.e., of\vspace{-2mm}
	 \begin{equation}\label{EQ59'}0<G=(x_1+a_2''x_2)x_2-\frac{1}{4}x_5^2-\frac{a+A}{4}<d+\frac{3A}{4},\end{equation}
	{\scriptsize
		\begin{equation*}\begin{array}{lllll}
		\hline\\
		{\rm ~~~~Range} & x_1 & x_1x_2 & x_5& {\rm ~~G}\vspace{2mm}\\
		\hline\\
		
		\frac{1}{4}+\frac{a+A}{4}<f_1''<1&c_1&|c_1|&\pm 1&f_1''-\frac{1}{4}-\frac{a+A}{4}\vspace{2mm}\\
		
		 \frac{a+A}{4}<f_1''\leq\frac{1}{4}+\frac{a+A}{4}&c_1&|c_1|&0&f_1''-\frac{a+A}{4}\vspace{2mm}\\
		
		-\frac{1}{2}<f_1''<-\frac{3}{4}+d+\frac{a}{4}+A&\pm 1+c_1&1+|c_1|&\pm 1&f_1''+\frac{3}{4}-\frac{a+A}{4}\vspace{2mm}\\
		
		\frac{1}{4}+\frac{a+A}{4}<g_1''\leq \frac{1}{2}&c_1&-|c_1|&\pm 1&g_1''-\frac{1}{4}-\frac{a+A}{4}\vspace{2mm}\\
		
		\frac{a+A}{4}<g_1''\leq \frac{1}{4}+\frac{a+A}{4}&c_1&-|c_1|&0&g_1''-\frac{a+A}{4}\vspace{2mm}\\
		
		-\frac{3}{4}+\frac{a+A}{4}<g_1''<-\frac{3}{4}+d+\frac{a}{4}+A&\pm 1+c_1&1-|c_1|&\pm 1&g_1''+\frac{3}{4}-\frac{a+A}{4}\vspace{2mm}\\
		
		-1+\frac{a+A}{4}<g_1''\leq-\frac{3}{4}+\frac{a+A}{4}&\pm 1+c_1&1-|c_1|&0&g_1''+1-\frac{a+A}{4}\vspace{2mm}\\
		
		\hline
		\end{array}
		\end{equation*}\vspace{1mm}}
	
	\noindent except in the following three cases:
	\begin{enumerate}
		\item $-\frac{3}{4}+d+\frac{a}{4}+A\leq f_1''\leq \frac{a+A}{4}, ~ g_1''\leq -1+\frac{a+A}{4}$
		\item $-\frac{3}{4}+d+\frac{a}{4}+A\leq f_1''\leq \frac{a+A}{4}$,~ $-\frac{3}{4}+d+\frac{a}{4}+A\leq g_1''\leq \frac{a+A}{4}$
		\item $f_1''=\frac{1}{4}+d+\frac{a}{4}+A=1$.
	\end{enumerate}
	
	\noindent \textbf{Case 1:} $-\frac{3}{4}+d+\frac{a}{4}+A\leq f_1''\leq \frac{a+A}{4}, ~ g_1''\leq -1+\frac{a+A}{4}$.\vspace{2mm}
	
	\noindent Note that here, $g_1''\leq -1+\frac{a+A}{4} \Rightarrow a_2''\leq -\frac{1}{2}+\frac{a+A}{4}<0$.  Take $x_2=\pm 2, x_1=\pm 2+c_1, x_1x_2=4-2|c_1|,~ x_5=2$. Then, $$\begin{array}{ll}G&=3+g_2''-\frac{a+A}{4}\\&
	=3+2g_1''+2a_2''-\frac{a+A}{4}\\ &\leq 3+2(-1+\frac{a+A}{4})+2(-\frac{1}{2}+\frac{a+A}{4})-\frac{a+A}{4}\\& =3(\frac{a+A}{4})<d+3A/4.\end{array}$$.
	
	\noindent And, $G>0$ for $g_2''>-3+\frac{a+A}{4}$. So, let now $g_2''\leq -3+\frac{a+A}{4}$ which gives $a_2''\leq -\frac{1}{2}+\frac{a+A}{4}.\frac{1}{4}$.\vspace{2mm}
	
	\noindent Now take $x_2=\pm 3, x_1=\pm 2+c_1, x_1x_2=6-3|c_1|, x_5=0$. Then, $$\begin{array}{ll}G&=6+g_3''-\frac{a+A}{4}\\&=6+\frac{3}{2}(-2|c_1|+4a_2'')+3a_2''-\frac{a+A}{4}\\&
	=6+\frac{3}{2}g_2''+3a_2''-\frac{a+A}{4}\\ &\leq 6+\frac{3}{2}(-3+\frac{a+A}{4})+3(-\frac{1}{2}+\frac{a+A}{16})-\frac{a+A}{4}\\&=5(\frac{a+A}{16})<d+3A/4.
	\end{array}$$
	\noindent And, $G>0$ for $g_3''>-6+\frac{a+A}{4}$. So, let now \begin{equation*}g_3''=-3|c_1|+9a_2''\leq -6+\frac{a+A}{4}, {\rm ~which ~ gives~} a_2''\leq -\frac{1}{2}+\frac{a+A}{4}.\frac{1}{9}.\end{equation*}
	
	\noindent Note that $a_2''<0$. Find an integer $n\geq 1$ such that\vspace{-3mm}
	
	{\footnotesize \begin{equation}\label{Eq60}\begin{array}{l}
		(2n+1)\{(n+1)-|c_1|+(2n+1)a_2''\}\leq \frac{a+A}{4}
		<(2n+3)\{(n+2)-|c_1|+(2n+3)a_2''\}
		\end{array}
		\end{equation}}
	i.e., \begin{equation*}\begin{array}{l}(2n+1)(n+1)+g_{2n+1}''\leq \frac{a+A}{4}<(2n+3)(n+2)+g_{2n+3}''.\end{array}
	\end{equation*}
	\noindent Further $g_{2n+1}''=-(2n+1)|c_1|+(2n+1)^2a_2''\leq -(2n+1)(n+1)+\frac{a+A}{4}$ gives \begin{equation}\label{Eq61}
	a_2''\leq -\frac{1}{2}+\frac{a+A}{4}.\frac{1}{(2n+1)^2}.
	\end{equation}
	\noindent   Take $x_2=\pm (2n+3), x_1=\pm (n+2)+c_1, ~x_1x_2=(2n+3)(n+2)-(2n+3)|c_1|,~x_5=0$. So that, $$G=(2n+3)(n+2)-(2n+3)|c_1|+(2n+3)^2a_2''-\frac{a+A}{4}.$$
	
	\noindent Using (\ref{Eq60}), we find that $G>0$. And from (\ref{Eq61}), we find that
	\begin{equation*}\begin{array}{lll}
	 G&=&(2n+3)(n+2)+\frac{2n+3}{2n+1}\big(-(2n+1)|c_1|+(2n+1)^2a_2''\big)\\&&+\big(-(2n+3)(2n+1)+(2n+3)^2\big)a_2''
	-\frac{a+A}{4}\vspace{2mm}\\
	&\leq& (2n+3)(n+2)+\frac{2n+3}{2n+1}\big(-(2n+1)(n+1)+\frac{a+A}{4}\big)\\&&+2(2n+3)\big(-\frac{1}{2}+\frac{a+A}{4}
	\frac{1}{(2n+1)^2}\big)-\frac{a+A}{4}\vspace{2mm}\\
	&=&\big(\frac{2n+3}{2n+1}+2.\frac{2n+3}{(2n+1)^2}-1\big)\frac{a+A}{4}\vspace{2mm}\\
	&\leq & \frac{16}{9}.\frac{a+A}{4}<d+\frac{3A}{4}\vspace{2mm}\\
	\end{array}
	\end{equation*}
	\noindent as $\frac{2n+3}{2n+1}+2.\frac{2n+3}{(2n+1)^2}-1$ is a decreasing function of $n$ and $n\geq 1$. \vspace{4mm}
	
	\noindent \textbf{Case 2:} $-\frac{3}{4}+d+\frac{a}{4}+A\leq f_1''\leq \frac{a+A}{4}$,~ $-\frac{3}{4}+d+\frac{a}{4}+A\leq g_1''\leq \frac{a+A}{4}$.\vspace{2mm}
	
	\noindent Here as $g_1''\geq -\frac{3}{4}+d+\frac{a}{4}+A \geq -\frac{3}{4}+d+\frac{d}{8}+\frac{3d}{8}\geq -\frac{3}{4}+\frac{3d}{2}\geq 0$, we have
	$ g_1''=-|c_1|+a_2'' \geq 0$. Which  gives $a_2'' \geq 0$.\vspace{2mm}
		
	\noindent  First suppose $a_2''\neq 0$ or $c_1\neq 0$. Then find an integer $n\geq 1$ such that
	\begin{equation}\label{Eq62}
	n|c_1|+n^2a_2''\leq \frac{a+A}{4}<(n+1)|c_1|+(n+1)^2a_2''.
	\end{equation}
	
	\noindent Now take $x_2=\pm (n+1), x_1=c_1,x_1x_2=(n+1)|c_1|, x_5=0$ so that
	$$G=(n+1)|c_1|+(n+1)^2a_2''-\frac{a+A}{4}.$$
	
	\noindent Using (\ref{Eq62}), we find that $G>0$. And
	\begin{equation*}\begin{array}{lll}
	G&=&(n+1)|c_1|+(n+1)^2a_2''-\frac{a+A}{4}\vspace{2mm}\\
	&=&\frac{(n+1)^2}{n^2}.(\frac{n^2}{n+1}|c_1|+n^2a_2'')-\frac{a+A}{4}\vspace{2mm}\\
	&\leq &\frac{(n+1)^2}{n^2}.(n|c_1|+n^2a_2'')-\frac{a+A}{4}\vspace{2mm}\\
	&\leq & 4(\frac{a+A}{4})-\frac{a+A}{4}\vspace{2mm}\\
	&=&3(\frac{a+A}{4})<d+\frac{3A}{4}.
	\end{array}
	\end{equation*}
	\noindent If $a_2''= 0$ and $c_1= 0$, we get $0=g_1''\geq -\frac{3}{4}+d+\frac{a}{4}+A \geq -\frac{3}{4}+d+\frac{d}{8}+\frac{3d}{8}\geq -\frac{3}{4}+\frac{3d}{2}\geq 0$, which gives $d=\frac{1}{2}$; $a=\frac{d}{2}=\frac{1}{4}$ and $A=\frac{3d}{8}=\frac{3}{16}$. \vspace{2mm}
	
	\noindent Further since, $-\frac{3}{4}<a_5'-\lambda a_4'\leq \frac{3}{4}$ and $a_5''=0$, and $-\frac{1}{2}< a_2'+ a_4'^2/4A = a_2'+ \frac{4}{3}a_4'^2<1$ we get from (\ref{eq16}) that
	\begin{equation}\label{Eq63}
	a_5'-\lambda a_4'=a_5''=0, {\rm ~ and ~}\end{equation}
	\begin{equation}\label{Eq64}a_2'+ \frac{4}{3}a_4'^2=a_2''=0.
	\end{equation}
	\noindent Taking $h=3, k=4$, we find by Macbeath Lemma that (\ref{eq101}) is soluble unless $A=3/16$ and $\pm a_4'-\frac{3}{8}(c_4+\lambda x_5) \equiv \frac{3}{4} \pmod {\frac{1}{4},\frac{3}{8}}$ i.e., $ \pm 8a_4'-3c_4-3\lambda x_5 \equiv 0\pmod 1$.  Taking $x_5=c_5$ and $1+c_5$ simultaneously, we get $3\lambda \equiv 0 \pmod 1$, i.e., $\lambda=0$ or, $\frac{1}{3}$ and $ \pm 8a_4'-3c_4 \equiv 0\pmod 1$ (as $c_5=0$ here).
	Therefore using Lemma \ref{lem15}, we need to consider $(a_4',c_4)=(0,0),$ $(0,1/3),$ $(1/16,1/6)$ and $(1/16,1/2)$. \vspace{2mm}\\
	\noindent (i) $\lambda=0$:\\
	\noindent Here, since $(a_5'',c_5)=(0,0)$, we get from (\ref{Eq63}) that $a_5'=0$.
	\begin{itemize}
		\item When $(a_4',c_4)=(0,0)$, (\ref{Eq64}) gives $a_2'=0$. Then $x_1=1, x_2=2, x_4=3, x_5=0$ gives a solution to (\ref{eq101}).
		\item When $(a_4',c_4)=(0,1/3)$, (\ref{Eq64}) gives $a_2'=0$. Then $x_1=1, x_2=1, x_4=\frac{4}{3}, x_5=1$ gives a solution to (\ref{eq101}).
		\item When $(a_4',c_4)=(1/16,1/6)$, (\ref{Eq64}) gives $a_2'=-\frac{1}{192}$. Then $x_1=-1, x_2=-1, x_4=-\frac{11}{6}, x_5=1$ gives a solution to (\ref{eq101}).
		\item When $(a_4',c_4)=(1/16,1/2)$, (\ref{Eq64}) gives $a_2'=-\frac{1}{192}$. Then $x_1=-1, x_2=-1, x_4=-\frac{3}{2}, x_5=1$ gives a solution to (\ref{eq101}).
	\end{itemize}
	
	\noindent (ii) $\lambda=\frac{1}{3}$:\\
	\noindent Here, since $(a_5'',c_5)=(0,0)$, we get from (\ref{Eq63}) that
	$a_5'=\frac{1}{3}a_4'$ and $a_2'=-\frac{4}{3} a_4'^2$.
	
	\begin{itemize}
		\item When $(a_4',c_4)=(0,0)$, (\ref{Eq63}) gives $a_2'=0, a_5'=0$. Then $x_1=1, x_2=1, x_4=1, x_5=1$ gives a solution to (\ref{eq101}).
		\item When $(a_4',c_4)=(0,1/3)$, (\ref{Eq63}) gives $a_2'=0,a_5'=0$. Then $x_1=1, x_2=1, x_4=\frac{4}{3}, x_5=1$ gives a solution to (\ref{eq101}).
		\item When $(a_4',c_4)=(1/16,1/6)$, (\ref{Eq63}) gives $a_2'=-\frac{1}{192}, a_5'=\frac{1}{48}$. Then $x_1=1, x_2=1, x_4=\frac{7}{6}, x_5=1$ gives a solution to (\ref{eq101}).
		\item When $(a_4',c_4)=(1/16,1/2)$, (\ref{Eq63}) gives $a_2'=-\frac{1}{192},a_5'=\frac{1}{48}$. Then $x_1=1, x_2=1, x_4=\frac{3}{2}, x_5=1$ gives a solution to (\ref{eq101}).
	\end{itemize}
	
	\noindent \textbf{Case 3:} $f_1''=\frac{1}{4}+d+\frac{a}{4}+A=1$.\\
	
	\noindent This gives $c_1=\frac{1}{2}=a_2''$; $d=\frac{1}{2}$; $a=\frac{d}{2}=\frac{1}{4}$ and $A=\frac{3d}{8}=\frac{3}{16}$. Here in place of (\ref{Eq64}) we have
	\begin{equation}\label{Eq65}
	a_2'+ \frac{4}{3}a_4'^2=a_2''=\frac{1}{2}.
	\end{equation}
	Working as in previous case we get $\lambda=0$ or, $\frac{1}{3}$.
	and we need to consider $(a_4',c_4)=(0,0),$ $(0,1/3),$ $(1/16,1/6)$ and $(1/16,1/2)$.
	\begin{itemize}
		\item When $\lambda=0$, $(a_4',c_4)=(0,0)$, (\ref{Eq65}) gives $a_2'=\frac{1}{2}$. (\ref{Eq63}) gives $a_5'=0$. Then $x_1=\frac{-3}{2}, x_2=-1, x_4=-2, x_5=-2$ gives a solution to (\ref{eq101}).
		\item When  $\lambda=0$, $(a_4',c_4)=(0,1/3)$, (\ref{Eq65}) gives $a_2'=\frac{1}{2}$. (\ref{Eq63}) gives $a_5'=0$. Then $x_1=\frac{1}{2}, x_2=1, x_4=\frac{4}{3}, x_5=1$ gives a solution to (\ref{eq101}).
		\item When  $\lambda=0$, $(a_4',c_4)=(1/16,1/6)$, (\ref{Eq65}) gives $a_2'=\frac{95}{192}$. (\ref{Eq63}) gives  $a_5'=0$. Then $x_1=\frac{-1}{2}, x_2=-1, x_4=\frac{-11}{6}, x_5=-1$ gives a solution to (\ref{eq101}).
		
		\item When  $\lambda=0$, $(a_4',c_4)=(1/16,1/2)$, (\ref{Eq65}) gives $a_2'=\frac{95}{192}$. (\ref{Eq63}) gives $a_5'=0$. Then $x_1=\frac{1}{2}, x_2=1, x_4=\frac{3}{2}, x_5=1$ gives a solution to (\ref{eq101}).
		
		\item When $\lambda=\frac{1}{3}$, $(a_4',c_4)=(0,0)$, (\ref{Eq65}) gives $a_2'=\frac{1}{2}$. (\ref{Eq63}) gives $a_5'=0$. Then $x_1=\frac{1}{2}, x_2=1, x_4=0, x_5=0$ gives a solution to (\ref{eq101}).
		
		\item When  $\lambda=\frac{1}{3}$, $(a_4',c_4)=(0,1/3)$, (\ref{Eq65}) gives $a_2'=\frac{1}{2}$. (\ref{Eq63}) gives $a_5'=0$. Then $x_1=\frac{1}{2}, x_2=1, x_4=\frac{4}{3}, x_5=1$ gives a solution to (\ref{eq101}).
		
		\item When  $\lambda=\frac{1}{3}$, $(a_4',c_4)=(1/16,1/6)$, (\ref{Eq65}) gives $a_2'=\frac{95}{192}$. (\ref{Eq63}) gives  $a_5'=\frac{1}{48}$. Then $x_1=\frac{1}{2}, x_2=1, x_4=\frac{7}{6}, x_5=1$ gives a solution to (\ref{eq101}).
		
		\item When  $\lambda=\frac{1}{3}$, $(a_4',c_4)=(1/16,1/2)$, (\ref{Eq65}) gives $a_2'=\frac{95}{192}$. (\ref{Eq63}) gives $a_5'=\frac{1}{48}$. Then $x_1=\frac{1}{2}, x_2=1, x_4=\frac{3}{2}, x_5=1$ gives a solution to (\ref{eq101}).
	\end{itemize}
	
	\subsection{${m=1, K=2, L=1}$}
	
	\noindent Here $1<\frac{\delta_{m,K}}{t}=\frac{d+3A/4}{t}\leq 2$, $\delta_{m,K}=d+\frac{3A}{4}\geq d+\frac{9d}{32}=\frac{41d}{32}$. From (\ref{Eq74'}) and (\ref{EqN1'}) we have
	\begin{equation}\label{Eq85'}
	\begin{array}{l}~~~~~~~~~~0.5659< \sqrt{41/128} \leq d \leq 2/3<0.6667,\vspace{2mm}\\0.3625< \frac{41d}{64}\leq \frac{d+3A/4}{2}=\frac{\delta_{1,2}}{2}\leq t =\frac{3d^5}{8aA}\leq 2d^3< 0.5927.\end{array}
	\end{equation}

	\noindent Taking $h_t=1/2$ and $k_t=1$ and using (\ref{Eq85'})
	we find that $\vline~h_t-tk_t^2~\vline + \frac{1}{2} <\delta_{1,2}$, as $ t < \delta_{1,2}$ and   $t+\delta_{1,2} \geq \frac{41}{64}d+\frac{41}{32}d=\frac{123}{64}d>1$. Therefore,  (\ref{eq101'}) is soluble unless $t=1/2$ and $\beta_{t,G}=\pm a_5''-c_5\equiv 1/2 \pmod 1$. This gives $(a_5'',c_5)=(0,1/2)$ or $(1/2,0)$.\vspace{2mm}
	
	\noindent  We shall find solutions to (\ref{eq101}). For $t=1/2$ we have $1/2= 3d^5/8aA\leq 2d^3$ and $a\geq d/2$, which gives $3d/8\leq A\leq 3d^4/2$. Thus we have
	\begin{equation}\label{Eq'88'}
	d>0.62996,~a\leq 1-d<0.371 ~~{\rm and }~~ 0.2362< A < 0.297.
	\end{equation}

	\noindent Take $(h_A,k_A) = (7.5,5), (2.5,3)$ and $ (1,2)$ respectively for $A\in [0.296,0.297),~ [0.29,0.296]$ and $[0.236,0.29]$ to get that $\vline~ h_A-Ak_A^2~\vline +\frac{1}{2} <d$. As $h_A/k_A^2=7.5/25$ and $2.5/9$ lies outside the range of $A$, (\ref{eq101}) is soluble unless $A=1/4$ and $\beta_{A,F} = \pm a'_4-(c_4+\lambda x_5)/2 \equiv 1/2 \pmod {1/2}$. Taking $x_5=c_5$ and $1+c_5$, we get $\lambda=0$. By Lemma \ref{lem15}, we need to consider $(a_4',c_4)=(0,0)$ or $(1/4,1/2)$. Also from (\ref{eq16}), we get $(a_5'',c_5)=(a_5',c_5)=(0,1/2) {\rm ~or~} (1/2,0)$.\\
	
	\noindent \textbf{Case (i):} $(a_5',c_5)=(0,1/2)$, $(a_4',c_4)=(0,0)$. Here, the following table gives a solution to (\ref{eq101}) i.e., of
	 \begin{equation}\label{EQN'22'}0<F=(x_1+a_2'x_2)x_2-\frac{1}{4}x_4^2-\frac{1}{2}x_5^2-\frac{a}{4}<d.\end{equation}
	
	{\scriptsize
		\begin{equation*}\begin{array}{llllll}
		\hline\\
		{\rm ~~~~Range} & x_1 & x_1x_2 & x_4&x_5& {\rm ~~F}\vspace{2mm}\\
		\hline\\
		
		\frac{3}{8}+\frac{a}{4}<f_1'\leq 1&c_1&|c_1|&\pm 1&\pm \frac{1}{2}&f_1'-\frac{3}{8}-\frac{a}{4}\vspace{2mm}\\
		
		\frac{1}{8}+\frac{a}{4}<f_1'\leq \frac{3}{8}+\frac{a}{4}&c_1&|c_1|&0&\pm \frac{1}{2}&f_1'-\frac{1}{8}-\frac{a}{4}\vspace{2mm}\\
		
		-\frac{5}{8}+\frac{a}{4}<f_1'<-\frac{5}{8}+d+\frac{a}{4}&\pm 1+c_1&1+|c_1|&\pm 1&\pm \frac{1}{2}&f_1'+\frac{5}{8}-\frac{a}{4}\vspace{2mm}\\
		
		-\frac{1}{2}<f_1'\leq -\frac{5}{8}+\frac{a}{4}&\pm 1+c_1&1+|c_1|&0&\pm \frac{1}{2}&f_1'+\frac{7}{8}-\frac{a}{4}\vspace{2mm}\\
		
		{\rm When~}  -\frac{5}{8}+d+\frac{a}{4}\leq f_1'\leq \frac{1}{8}+\frac{a}{4}: \vspace{2mm}\\
		
		\left.\begin{array}{l}-\frac{5}{8}+d+\frac{a}{4}\leq g_1'\leq \frac{1}{8}+\frac{a}{4}\vspace{1mm}\\\frac{3}{8}+\frac{a}{4}< f_2'\end{array}\right\}& c_1&2|c_1|&\pm 1&\pm \frac{1}{2}&f_2'-\frac{3}{8}-\frac{a}{4}\vspace{2mm}\\
		
		\left.\begin{array}{l}-\frac{5}{8}+d+\frac{a}{4}\leq g_1'\leq \frac{1}{8}+\frac{a}{4}\vspace{1mm}\\f_2'\leq \frac{3}{8}+\frac{a}{4}\end{array}\right\}& c_1&2|c_1|&0&\pm \frac{1}{2}&f_2'-\frac{1}{8}-\frac{a}{4}\vspace{2mm}\\	
		
		-\frac{5}{8}+\frac{a}{4}<g_1'<-\frac{5}{8}+d+\frac{a}{4}&\pm 1+c_1&1-|c_1|&\pm 1&\pm \frac{1}{2}&g_1'+\frac{5}{8}-\frac{a}{4}\vspace{2mm}\\
		
		-\frac{7}{8}+\frac{a}{4}<g_1'\leq -\frac{5}{8}+\frac{a}{4}&\pm 1+c_1&1-|c_1|&0&\pm \frac{1}{2}&g_1'+\frac{7}{8}-\frac{a}{4}\vspace{2mm}\\
		
		\left.\begin{array}{l}-1<g_1' \leq -\frac{7}{8}+\frac{a}{4}\vspace{1mm}\\f_2'< -\frac{7}{8}+d+\frac{a}{4}\end{array}\right\}& \pm 1+ c_1&2+2|c_1|&\pm 2&\pm \frac{1}{2}&f_2'+\frac{7}{8}-\frac{a}{4}\vspace{2mm}\\
		
		\left.\begin{array}{l}-1<g_1' \leq -\frac{7}{8}+\frac{a}{4}\vspace{1mm}\\-\frac{7}{8}+d+\frac{a}{4}\leq f_2'\end{array}\right\}& \pm 1+c_1&2+2|c_1|&\pm 1&\pm \frac{3}{2}&f_2'+\frac{5}{8}-\frac{a}{4}\vspace{2mm}\\
		
		\hline
		\end{array}
		\end{equation*}\vspace{2mm}}

	\noindent \textbf{Case (ii):} $(a_5',c_5)=(0,1/2)$, $(a_4',c_4)=(1/4,1/2)$. Here, the following table gives a solution to (\ref{eq101}) i.e., of\\
	 \begin{equation}\label{EQN'22''}0<F=(x_1+a_2'x_2+\frac{1}{4}x_4)x_2-\frac{1}{4}x_4^2-\frac{1}{2}x_5^2-\frac{a}{4}<d.\end{equation}
	
	{\scriptsize
		\begin{equation*}\begin{array}{lllllll}
		\hline\\
		{\rm ~~~~Range} & x_1 & x_1x_2 & x_4&x_2x_4&x_5& {\rm ~~F}\vspace{2mm}\\
		\hline\\
		
		\frac{5}{16}+\frac{a}{4}<f_1'\leq 1&c_1&|c_1|&\pm \frac{1}{2}&-\frac{1}{2}&\pm \frac{1}{2}&f_1'-\frac{5}{16}-\frac{a}{4}\vspace{2mm}\\
		
		\frac{1}{16}+\frac{a}{4}<f_1'\leq \frac{5}{16}+\frac{a}{4}&c_1&|c_1|&\pm \frac{1}{2}&\frac{1}{2}&\pm \frac{1}{2}&f_1'-\frac{1}{16}-\frac{a}{4}\vspace{2mm}\\
		
		-\frac{1}{2}<f_1'<-\frac{11}{16}+d+\frac{a}{4}&\pm 1+c_1&1+|c_1|&\pm \frac{1}{2}&-\frac{1}{2}&\pm \frac{1}{2}&f_1'+\frac{11}{16}-\frac{a}{4}\vspace{2mm}\\
		
		{\rm When~}  -\frac{11}{16}+d+\frac{a}{4}\leq f_1'\leq \frac{1}{16}+\frac{a}{4}: \vspace{2mm}\\
		
		\left.\begin{array}{l}-\frac{11}{16}+d+\frac{a}{4}\leq g_1'\leq \frac{1}{16}+\frac{a}{4}\vspace{1mm}\\\frac{7}{16}+\frac{a}{4}<f_2'\end{array}\right\}& c_1&2|c_1|&\pm \frac{1}{2}&-1&\pm \frac{1}{2}&f_2'-\frac{7}{16}-\frac{a}{4}\vspace{2mm}\\
		
		\left.\begin{array}{l}-\frac{11}{16}+d+\frac{a}{4}\leq g_1'\leq \frac{1}{16}+\frac{a}{4}\vspace{1mm}\\f_2'\leq \frac{7}{16}+\frac{a}{4}\end{array}\right\}& c_1&2|c_1|&\pm \frac{1}{2}&1&\pm \frac{1}{2}&f_2'+\frac{1}{16}-\frac{a}{4}\vspace{2mm}\\	
		
		-\frac{11}{16}+\frac{a}{4}<g_1'<-\frac{11}{16}+d+\frac{a}{4}&\pm 1+c_1&1-|c_1|&\pm \frac{1}{2}&-\frac{1}{2}&\pm \frac{1}{2}&g_1'+\frac{11}{16}-\frac{a}{4}\vspace{2mm}\\
		
		-\frac{15}{16}+\frac{a}{4}<g_1'\leq -\frac{11}{16}+\frac{a}{4}&\pm 1+c_1&1-|c_1|&\pm \frac{1}{2}&\frac{1}{2}&\pm \frac{1}{2}&g_1'+\frac{15}{16}-\frac{a}{4}\vspace{2mm}\\
		
		\left.\begin{array}{l}-1<g_1' \leq -\frac{15}{16}+\frac{a}{4}\vspace{1mm}\\f_2'< -\frac{17}{16}+d+\frac{a}{4}\end{array}\right\}& \pm 1+c_1&2+2|c_1|&\pm \frac{1}{2}&1&\pm \frac{3}{2}&f_2'+\frac{17}{16}-\frac{a}{4}\vspace{2mm}\\
		
		\left.\begin{array}{l}-1<g_1' \leq -\frac{15}{16}+\frac{a}{4}\vspace{1mm}\\-\frac{17}{16}+d+\frac{a}{4} \leq f_2'\end{array}\right\}& \pm 1+c_1&2+2|c_1|&\pm \frac{1}{2}&-1&\pm \frac{3}{2}&f_2'+\frac{9}{16}-\frac{a}{4}\vspace{2mm}\\	 
		
		\hline
		\end{array}
		\end{equation*}\vspace{2mm}}

	\noindent \textbf{Case (iii):} $(a_5',c_5)=(1/2,0)$, $(a_4',c_4)=(0,0)$. Here, the following table gives a solution to (\ref{eq101}) i.e., of \begin{equation}\label{EQN'22'''}0<F=(x_1+a_2'x_2+\frac{1}{2}x_5)x_2-\frac{1}{4}x_4^2-\frac{1}{2}x_5^2-\frac{a}{4}<d.\end{equation}
	
	{\scriptsize
		\begin{equation*}\begin{array}{lllllll}
		\hline\\
		{\rm ~~~~Range} & x_1 & x_1x_2 & x_4&x_5&x_2x_5& {\rm ~~F}\vspace{2mm}\\
		\hline\\
		
		\frac{1}{4}+\frac{a}{4}<f_1'\leq 1&c_1&|c_1|&\pm 1&0&0&f_1'-\frac{1}{4}-\frac{a}{4}\vspace{2mm}\\
		
		\frac{a}{4}<f_1'\leq \frac{1}{4}+\frac{a}{4}&c_1&|c_1|&0&0&0&f_1'-\frac{a}{4}\vspace{2mm}\\
		
		-\frac{1}{2}<f_1'<-\frac{3}{4}+d+\frac{a}{4}&\pm 1+c_1&1+|c_1|&\pm 1&0&0&f_1'+\frac{3}{4}-\frac{a}{4}\vspace{2mm}\\
		
		{\rm When~}  -\frac{3}{4}+d+\frac{a}{4}\leq f_1'\leq \frac{a}{4}: \vspace{2mm}\\
		
		\left.\begin{array}{l}-\frac{3}{4}+d+\frac{a}{4}\leq g_1'\leq \frac{a}{4}\vspace{1mm}\\\frac{a}{4}< f_2'\end{array}\right\}& c_1&2|c_1|&0&0&0&f_2'-\frac{a}{4}\vspace{2mm}\\
		
		\left.\begin{array}{l}-\frac{3}{4}+d+\frac{a}{4}\leq g_1'\leq \frac{a}{4}\vspace{1mm}\\f_2'\leq \frac{a}{4}\end{array}\right\}& c_1&2|c_1|&0&\pm 1&2&f_2'+\frac{1}{2}-\frac{a}{4}\vspace{2mm}\\	
		
		-\frac{3}{4}+\frac{a}{4}<g_1'<-\frac{3}{4}+d+\frac{a}{4}&\pm 1+c_1&1-|c_1|&\pm 1&0&0&g_1'+\frac{3}{4}-\frac{a}{4}\vspace{2mm}\\
		
		-1+\frac{a}{4}<g_1'\leq -\frac{3}{4}+\frac{a}{4}&\pm 1+c_1&1-|c_1|&0&0&0&g_1'+1-\frac{a}{4}\vspace{2mm}\\
		
		\left.\begin{array}{l}-1<g_1' \leq -1+\frac{a}{4}\vspace{1mm}\\-1+\frac{a}{4}<f_2'\end{array}\right\}& \pm 1+c_1&2+2|c_1|&\pm 2&0&0&f_2'+1-\frac{a}{4}\vspace{2mm}\\
		
		\left.\begin{array}{l}-1<g_1' \leq -1+\frac{a}{4}\vspace{1mm}\\f_2'\leq -1+\frac{a}{4}\end{array}\right\}& \pm 1+c_1&2+2|c_1|&\pm 2&\pm 1&2&f_2'+\frac{3}{2}-\frac{a}{4}\vspace{2mm}\\	
		
		\hline
		\end{array}
		\end{equation*}\vspace{1mm}}

	\noindent \textbf{Case (iv):} $(a_5',c_5)=(1/2,0)$, $(a_4',c_4)=(1/4,1/2)$. On taking $x_5=0$, the following table gives a solution to (\ref{eq101}) i.e., of\\
	 \begin{equation}\label{EQN'22''''}0<F=(x_1+a_2'x_2+\frac{1}{4}x_4+\frac{1}{2}x_5)x_2-\frac{1}{4}x_4^2-\frac{1}{2}x_5^2-\frac{a}{4}<d.\end{equation}
	
	{\scriptsize
		\begin{equation*}\begin{array}{llllll}
		\hline\\
		{\rm ~~~~Range} & x_1 & x_1x_2 & x_4&x_2x_4& {\rm ~~F}\vspace{2mm}\\
		\hline\\
		
		\frac{3}{16}+d+\frac{a}{4}\leq f_1'\leq 1&\pm 1+c_1&-2+2|c_1|&\pm \frac{1}{2}&-1&f_2'-\frac{37}{16}-\frac{a}{4}\vspace{2mm}\\
		
		\frac{3}{16}+\frac{a}{4}<f_1'<\frac{3}{16}+d+\frac{a}{4}&c_1&|c_1|&\pm \frac{1}{2}&-\frac{1}{2}&f_1'-\frac{3}{16}-\frac{a}{4}\vspace{2mm}\\
		
		-\frac{1}{16}+\frac{a}{4}<f_1'\leq \frac{3}{16}+\frac{a}{4}&c_1&|c_1|&\pm \frac{1}{2}&\frac{1}{2}&f_1'+\frac{1}{16}-\frac{a}{4}\vspace{2mm}\\
		
		-\frac{1}{2}<f_1'<-\frac{13}{16}+d+\frac{a}{4}&\pm 1+c_1&1+|c_1|&\pm \frac{1}{2}&-\frac{1}{2}&f_1'+\frac{13}{16}-\frac{a}{4}\vspace{2mm}\\
		
		{\rm When~}  -\frac{13}{16}+d+\frac{a}{4}\leq f_1'\leq -\frac{1}{16}+\frac{a}{4}: \vspace{2mm}\\
		
		\left.\begin{array}{l}-\frac{13}{16}+d+\frac{a}{4}\leq g_1'\leq -\frac{1}{16}+\frac{a}{4}\vspace{1mm}\\f_2'<-\frac{11}{16}+d+\frac{a}{4}\end{array}\right\}& \pm 1+c_1&2+2|c_1|&\pm \frac{3}{2}&-3&f_2'+\frac{11}{16}-\frac{a}{4}\vspace{2mm}\\
		
		\left.\begin{array}{l}-\frac{13}{16}+d+\frac{a}{4}\leq g_1'\leq -\frac{1}{16}+\frac{a}{4}\vspace{1mm}\\-\frac{11}{16}+d+\frac{a}{4}\leq f_2'\end{array}\right\}& c_1&2|c_1|&\pm \frac{1}{2}&1&f_2'+\frac{3}{16}-\frac{a}{4}\vspace{2mm}\\	
		
		-\frac{13}{16}+\frac{a}{4}<g_1'<-\frac{13}{16}+d+\frac{a}{4}&\pm 1+c_1&1-|c_1|&\pm \frac{1}{2}&-\frac{1}{2}&g_1'+\frac{13}{16}-\frac{a}{4}\vspace{2mm}\\
		
		-\frac{17}{16}+\frac{a}{4}<g_1'\leq -\frac{13}{16}+\frac{a}{4}&\pm 1+c_1&1-|c_1|&\pm \frac{1}{2}&\frac{1}{2}&g_1'+\frac{17}{16}-\frac{a}{4}\vspace{2mm}\\
		
		\left.\begin{array}{l}-1<g_1'\leq -\frac{17}{16}+\frac{a}{4}\vspace{1mm}\\-\frac{19}{16}+\frac{a}{4}<f_2'\end{array}\right\}& \pm 2+c_1&4+2|c_1|&\pm \frac{5}{2}&-5&f_2'+\frac{19}{16}-\frac{a}{4}\vspace{2mm}\\
		
		\left.\begin{array}{l}-1<g_1'\leq -\frac{17}{16}+\frac{a}{4}\vspace{1mm}\\f_2'\leq -\frac{19}{16}+\frac{a}{4}\end{array}\right\}& \pm 1+c_1&2+2|c_1|&\pm \frac{1}{2}&-1&f_2'+\frac{27}{16}-\frac{a}{4}\vspace{2mm}\\	
		
		\hline
		\end{array}
		\end{equation*}}

	\section{Proof of $\Gamma_{1,4}< \frac{32}{3}$ when $(m,K)=(1,1)$}\label{sec11}
	\numberwithin{equation}{section}
	\begin{theorem}\label{thm8}
		Let $Q(x_1,x_2, \cdots, x_5)$ be a real indefinite quadratic form of type $(1,4)$ and of determinant $D\neq 0$. Let $d=(\frac{32}{3}|D|)^{\frac{1}{5}}$.  Suppose that $c_2 \equiv 0 \pmod1$, $a<\frac{1}{2}$, $d\leq 1$ and $a+d\leq 1$. Let $(m,K)= (1,1)$.  Then $(\ref{eq101*})$ and hence $(\ref{eq10})$ is soluble with strict inequality.
	\end{theorem}
	\noindent \textbf{Proof:} Here, $1<\frac{d}{a}\leq 2, {\rm ~and~} 1<\frac{d}{A}\leq 2$, $ \frac{d}{2}\leq A \leq \sqrt{\frac{d^5}{2a}}\leq d^2$.  Also $\frac{d}{2}+d\leq a+d \leq 1$ gives $d\leq \frac{2}{3}$. From  (\ref{key2}) and (\ref{eq09}) we get $\frac{3d}{8}\leq \frac{3A}{4}\leq t=\frac{\delta}{A}=\frac{3d^5}{8aA}\leq \frac{3d^3}{2}$. Therefore,
	
	\begin{equation}\label{one}
	\frac{1}{2}\leq d\leq \frac{2}{3}, ~~0.1875=\frac{3}{16}\leq\frac{3d}{8}\leq t \leq \frac{3d^3}{2}\leq \frac{4}{9}<0.44445=\lambda_0.
	\end{equation}
	
	\noindent Here along with Macbeath's Lemma \ref{lem5}, we take help of Second Lemma of Macbeath \cite{Macbeath}  to deal with the whole range of $t$. We divide the range of $t$ into two parts.  For $0.195965 \leq t < 0.44445=\lambda_0 $ we use Lemma \ref{lem5}, and for $0.1875 \leq t \leq 0.195965$, we  use the following Second Lemma of Macbeath \cite{Macbeath}.
	
	\begin{lemma}\label{lem6}
		Let $t, \beta, d$ be real numbers with $t > 0, d > 0.$ Let $h,k$ be integers such that $$|h-k^2t|\leq (d/2)^3.$$ Suppose that either $t \neq h/k^2$ or $\beta \neq p/q$ with $q\leq 2/d$. Then for any real number $\nu$, there exist integers $x,y$ satisfying $$0< \pm x + \beta y\pm t y^2 +\nu \leq d.$$
	\end{lemma}
	
	\begin{remark}\label{rem2} \normalfont Macbeath's first Lemma \ref{lem5} is applicable only when $d \geq \frac{1}{2}$. This restriction is not required in Macbeath's second Lemma \ref{lem6}. We find that   $(\frac{d}{2})^3>d-\frac{1}{2}$, for $d <0.517304$, so that for the integral values of $h$ and $k$, the length of interval of $t$ covered by Lemma \ref{lem6} is larger than that covered by Lemma \ref{lem5}. For $d=(\frac{32}{3}|D|)^{\frac{1}{5}}$, we do have $d \geq \frac{1}{2}$. But when $d$ is very very close to half, Lemma \ref{lem5} is not practical to apply. We shall use Lemma \ref{lem6} for $\frac{1}{2}\leq d< 0.5075$. Once $t$ is fixed by Lemma \ref{lem6}, we go back to Lemma \ref{lem5}, as $d> \frac{1}{2}$, and then fix $(a_5'',c_5)$.\vspace{2mm}
		
		\noindent Note that in Lemma \ref{lem6}, $h$ has to be an integer; whereas in Lemma \ref{lem5}, $h$ can be half an integer also. This restricts the choices of pairs $(h,k)$ in Lemma \ref{lem6} to half compared to its choices in Lemma \ref{lem5}.  Moreover, Lemma \ref{lem6}  is very weak in fixing the values of $(a_5'',c_5)$ or $(a_4',c_4)$. When $t = h/k^2$ we get  $\beta = p/q$, where $p,q$ are any integers with $0<q\leq 2/d$. Here integers $p$ and $q$ are not related to $(h,k)$.
	\end{remark}
	
	\noindent Here when $(m,K)=(1,1)$ we will show that (\ref{eq101'}) is soluble. We divide the range of $t$ into $\textbf{3121}$ subintervals $[\lambda_n,\lambda_{n-1}]$. For   $n=1,2,\cdots 1775$, we use Macbeath first Lemma (Lemma \ref{lem5}) and for $n=1776, \cdots, 3121$ use Macbeath's second Lemma (Lemma \ref{lem6}).\vspace{4mm}
	
	\noindent For   $n=1,2,\cdots 1775$, choose suitable integers $2h_n$ and $k_n$  such that
	$\vline~ h_n-tk_n^2~\vline +\frac{1}{2} \leq d$.
	(see Table I). The choice of $(h_n,k_n)$ is done in a manner similar to that in the Case  $m=2,K=3, L=4$.
	{\scriptsize
		\begin{equation*}
\begin{array}{lllllllll}
		&&&&\mbox{\normalsize{Table I.}}&&&&\vspace{2mm}\\
		\hline\\
		n & (h_n, k_n) & \lambda_n&\mbox{Remarks}&&n & (h_n, k_n) & \lambda_n&\mbox{Remarks}\vspace{2mm}\\
		\hline\\
		1 & (4,3) & 0.42692 & (44.5,10) & ~ & 2 & (34.5,9) & 0.42629 & \mbox{na} \\
		3 & (10.5,5) & 0.41396 & (15,6) & ~ & 4 & (15,6) & 0.412581 & \mbox{na} \\
		5 & (1.5,2) & 0.346595 & (8325.5,149) & ~ & 6 & (5.5,4) & 0.34591 & \mbox{na} \\
		7 & (3,3) & 0.32234 & \mbox{tbd} & ~ & 8 & (8,5) & 0.31853 & (26,9) \\
		9 & (5,4) & 0.30693 & \mbox{tbd} & ~ & 10 & (11,6) & 0.30348 & (15,7) \\
		11 & (7.5,5) & 0.2967 & \mbox{tbd} & ~ & 12 & (14.5,7) & 0.29427 & (20004,260) \\
		13 & (29.5,10) & 0.29419 & \mbox{na} & ~ & 14 & (85,17) & 0.29405 & \mbox{tbd} \\
		15 & (35.5,11) & 0.2939 & \mbox{na} & ~ & 16 & (10.5,6) & 0.28951 & \mbox{tbd} \\
		17 & (18.5,8) & 0.28787 & (35,11) & ~ & 18 & (41.5,12) & 0.28767 & \mbox{na} \\
		19 & (115,20) & 0.28743 & \mbox{tbd} & ~ & 20 & (48.5,13) & 0.28727 & \mbox{na} \\
		21 & (14,7) & 0.28615 & \mbox{na} & ~ & 22 & (2.5,3) & 0.270556 & (54.5,14) \\
		23 & (27,10) & 0.26936 & (97.5,19) & ~ & 24 & (45.5,13) & 0.26912 & \mbox{tbd} \\
		25 & (32.5,11) & 0.26808 & (118.5,21) & ~ & 26 & (52.5,14) & 0.26779 & \mbox{tbd} \\
		27 & (38.5,12) & 0.26693 & (10271,196) & ~ & 28 & (60,15) & 0.26657 & \mbox{tbd} \\
		29 & (13,7) & 0.26533 & \mbox{na} & ~ & 30 & (1,2) & 0.239395 & \mbox{tbd} \\
		31 & (6,5) & 0.23834 & \mbox{na} & ~ & 32 & (61,16) & 0.23812 & (229,31) \\
		33 & (105,21) & 0.238002 & \mbox{tbd} & ~ & 34 & (173.5,27) & 0.23796 & (4799,142) \\
		35 & (53.5,15) & 0.2376 & (325.5,37) & ~ & 36 & (95,20) & 0.23745 & \mbox{tbd} \\
		\hline
		\end{array}
		\end{equation*}}
	{\scriptsize
		\begin{equation*}\begin{array}{lllllllll}
		\hline\\
		n & (h_n, k_n) & \lambda_n&\mbox{Remarks}&&n & (h_n, k_n) & \lambda_n&\mbox{Remarks}\vspace{2mm}\\
		\hline\\	
		37 & (46.5,14) & 0.2372 & (12769.5,232) & ~ & 38 & (8.5,6) & 0.2351 & \mbox{tbd} \\
		39 & (11.5,7) & 0.23392 & (357,39) & ~ & 40 & (39.5,13) & 0.23351 & (10209.5,209) \\
		41 & (67.5,17) & 0.2335 & \mbox{na} & ~ & 42 & (52.5,15) & 0.23317 & \mbox{tbd} \\
		43 & (170,27) & 0.233145 & \mbox{na} & ~ & 44 & (269.5,34) & 0.23314 & \mbox{na} \\
		45 & (75.5,18) & 0.23291 & (4122,133) & ~ & 46 & (238.5,32) & 0.23289 & \mbox{na} \\
\hline
		\end{array}
		\end{equation*}}
	{\scriptsize
		\begin{equation*}\begin{array}{lllllllll}
		\hline\\
		n & (h_n, k_n) & \lambda_n&\mbox{Remarks}&&n & (h_n, k_n) & \lambda_n&\mbox{Remarks}\vspace{2mm}\\
		\hline\\
		169 & (3004.5,119) & 0.2121669 & (6647,177) & ~ & 170 & (1643,88) & 0.2121633 & (5499.5,161) \\
		171 & (869,64) & 0.2121577 & (5229.5,157) & ~ & 172 & (738.5,59) & 0.2121457 & (9090.5,207) \\
		173 & (815.5,62) & 0.2121434 & \mbox{na} & ~ & 174 & (842,63) & 0.21214 & \mbox{na} \\
		175 & (231,33) & 0.212102 & \mbox{tbd} & ~ & 176 & (429.5,45) & 0.2120885 & (4965,153) \\
		177 & (356.5,41) & 0.212065 & (17347,286) & ~ & 178 & (93.5,21) & 0.211971 & (13786.5,255) \\
\hline
		\end{array}
		\end{equation*}}
	{\scriptsize
		\begin{equation*}\begin{array}{lllllllll}
		\hline\\
		n & (h_n, k_n) & \lambda_n&\mbox{Remarks}&&n & (h_n, k_n) & \lambda_n&\mbox{Remarks}\vspace{2mm}\\
		\hline\\
		301 & (253.5,35) & 0.2069252 & (9925,219) & ~ & 302 & (365,42) & 0.206916 & (7787.5,194) \\
		303 & (174,29) & 0.206877 & \mbox{tbd} & ~ & 304 & (331,40) & 0.206866 & (400.5,44) \\
		305 & (382.5,43) & 0.20686 & \mbox{na} & ~ & 306 & (1291,79) & 0.206859 & \mbox{na} \\
		307 & (956.5,68) & 0.206857 & \mbox{na} & ~ & 308 & (821,63) & 0.206849 & (4171,142) \\
		309 & (538,51) & 0.20684 & (1459.5,84) & ~ & 310 & (67,18) & 0.206748 & (1986,98) \\
\hline
		\end{array}
		\end{equation*}}
	{\scriptsize
		\begin{equation*}\begin{array}{lllllllll}
		\hline\\
		n & (h_n, k_n) & \lambda_n&\mbox{Remarks}&&n & (h_n, k_n) & \lambda_n&\mbox{Remarks}\vspace{2mm}\\
		\hline\\
		433 & (688.5,58) & 0.204663 & (13835.5,260) & ~ & 434 & (891.5,66) & 0.2046569 & (6631,180) \\
		435 & (361,42) & 0.20465 & \mbox{na} & ~ & 436 & (295.5,38) & 0.20463 & (4666,151) \\
		437 & (433,46) & 0.2046244 & \mbox{na} & ~ & 438 & (619,55) & 0.2046241 & \mbox{na} \\
		439 & (452,47) & 0.204611 & (12483.5,247) & ~ & 440 & (918.5,67) & 0.204608 & \mbox{na} \\
		441 & (786.5,62) & 0.204604 & (8598.5,205) & ~ & 442 & (209.5,32) & 0.2045755 & (2170.5,103) \\
\hline
		\end{array}
		\end{equation*}}
	{\scriptsize
		\begin{equation*}\begin{array}{lllllllll}
		\hline\\
		n & (h_n, k_n) & \lambda_n&\mbox{Remarks}&&n & (h_n, k_n) & \lambda_n&\mbox{Remarks}\vspace{2mm}\\
		\hline\\
		565 & (8326.5,203) & 0.20205538 & \mbox{na} & ~ & 566 & (4307,146) & 0.2020544 & \mbox{tbd} \\
		567 & (2489.5,111) & 0.2020529 & (45973,477) & ~ & 568 & (1261,79) & 0.202051 & \mbox{na} \\
		569 & (247.5,35) & 0.202031 & (465.5,48) & ~ & 570 & (220,33) & 0.202009 & \mbox{tbd} \\
		571 & (633.5,56) & 0.2020078 & (19538.5,311) & ~ & 572 & (373.5,43) & 0.202 & (505,50) \\
		573 & (233.5,34) & 0.201979 & (2060.5,101) & ~ & 574 & (276.5,37) & 0.2019631 & (8992,211) \\
\hline
		\end{array}
		\end{equation*}}
	{\scriptsize
		\begin{equation*}\begin{array}{lllllllll}
		\hline\\
		n & (h_n, k_n) & \lambda_n&\mbox{Remarks}&&n & (h_n, k_n) & \lambda_n&\mbox{Remarks}\vspace{2mm}\\
		\hline\\
		697 & (407,45) & 0.200982 & (875.5,66) & ~ & 698 & (653,57) & 0.2009815 & \mbox{na} \\
		699 & (2091,102) & 0.2009798 & (86753,657) & ~ & 700 & (1130.5,75) & 0.200976 & (4582.5,151) \\
		701 & (146.5,27) & 0.200945 & (3140,125) & ~ & 702 & (169,29) & 0.200938 & \mbox{na} \\
		703 & (321.5,40) & 0.200931 & (5208.5,161) & ~ & 704 & (389,44) & 0.2009257 & (1254,79) \\
		705 & (371.5,43) & 0.20091312 & (4765,154) & ~ & 706 & (7329.5,191) & 0.2009128 & (4581,151) \\
\hline
		\end{array}
		\end{equation*}}
	{\scriptsize
		\begin{equation*}\begin{array}{lllllllll}
		\hline\\
		n & (h_n, k_n) & \lambda_n&\mbox{Remarks}&&n & (h_n, k_n) & \lambda_n&\mbox{Remarks}\vspace{2mm}\\
		\hline\\
		829 & (538.5,52) & 0.1991457 & (9903.5,223) & ~ & 830 & (6596.5,182) & 0.19914533 & (9205.5,215) \\
		831 & (7572.5,195) & 0.1991451 & (2280,107) & ~ & 832 & (5891.5,172) & 0.199145 & \mbox{na} \\
		833 & (3212,127) & 0.199144 & (2280,107) & ~ & 834 & (765.5,62) & 0.199142 & \mbox{na} \\
		835 & (624.5,56) & 0.199136 & (6524,181) & ~ & 836 & (1613,90) & 0.1991346 & (16061.5,284) \\
		837 & (11091,236) & 0.19913441 & (647,57) & ~ & 838 & (16402.5,287) & 0.19913438 & (647,57) \\
\hline
		\end{array}
		\end{equation*}}
	{\scriptsize
		\begin{equation*}\begin{array}{lllllllll}
		\hline\\
		n & (h_n, k_n) & \lambda_n&\mbox{Remarks}&&n & (h_n, k_n) & \lambda_n&\mbox{Remarks}\vspace{2mm}\\
		\hline\\
		961 & (763.5,62) & 0.19861878 & (4292,147) & ~ & 962 & (6945.5,187) & 0.1986185 & (3408.5,131) \\
		963 & (16360,287) & 0.1986183 & (1755,94) & ~ & 964 & (1335.5,82) & 0.1986155 & (4350.5,148) \\
		965 & (1435,85) & 0.1986146 & \mbox{na} & ~ & 966 & (1792.5,95) & 0.1986139 & \mbox{na} \\
		967 & (2765.5,118) & 0.1986138 & \mbox{na} & ~ & 968 & (715,60) & 0.1986103 & \mbox{tbd} \\
		969 & (384.5,44) & 0.1986004 & (2107,103) & ~ & 970 & (496.5,50) & 0.1985962 & (10597.5,231) \\
\hline
		\end{array}
		\end{equation*}}
	{\scriptsize
		\begin{equation*}\begin{array}{lllllllll}
		\hline\\
		n & (h_n, k_n) & \lambda_n&\mbox{Remarks}&&n & (h_n, k_n) & \lambda_n&\mbox{Remarks}\vspace{2mm}\\
		\hline\\
		1093 & (6486.5,181) & 0.1979943 & \mbox{na} & ~ & 1094 & (1430.5,85) & 0.1979919 & (8731.5,210) \\
		1095 & (6133,176) & 0.1979917 & \mbox{na} & ~ & 1096 & (6063.5,175) & 0.19799154 & \mbox{na} \\
		1097 & (7147.5,190) & 0.19799144 & \mbox{na} & ~ & 1098 & (18660.5,307) & 0.1979914 & \mbox{na} \\
		1099 & (96186,697) & 0.19799139 & \mbox{tbd} & ~ & 1100 & (23429.5,344) & 0.19799137 & \mbox{na} \\
		1101 & (1901.5,98) & 0.1979909 & \mbox{na} & ~ & 1102 & (1940.5,99) & 0.1979903 & \mbox{na} \\
\hline
		\end{array}
		\end{equation*}}
	{\scriptsize
		\begin{equation*}\begin{array}{lllllllll}
		\hline\\
		n & (h_n, k_n) & \lambda_n&\mbox{Remarks}&&n & (h_n, k_n) & \lambda_n&\mbox{Remarks}\vspace{2mm}\\
		\hline\\
		1225 & (50.5,16) & 0.1972325 & \mbox{na} & ~ & 1226 & (57,17) & 0.197203 & (4738.5,155) \\
		1227 & (331.5,41) & 0.1971992 & \mbox{na} & ~ & 1228 & (189.5,31) & 0.1971817 & (2386,110) \\
		1229 & (994,71) & 0.19718143 & \mbox{na} & ~ & 1230 & (17510.5,298) & 0.19718134 & \mbox{na} \\
		1231 & (47730.5,492) & 0.19718133 & \mbox{na} & ~ & 1232 & (16355,288) & 0.19718121 & (994,71) \\
		1233 & (9457,219) & 0.19718116 & \mbox{na} & ~ & 1234 & (7043.5,189) & 0.1971807 & (8449,207) \\
\hline
		\end{array}
		\end{equation*}}
	{\scriptsize
		\begin{equation*}\begin{array}{lllllllll}
		\hline\\
		n & (h_n, k_n) & \lambda_n&\mbox{Remarks}&&n & (h_n, k_n) & \lambda_n&\mbox{Remarks}\vspace{2mm}\\
		\hline\\
		1357 & (1229.5,79) & 0.1970033 & (2255.5,107) & ~ & 1358 & (492.5,50) & 0.1970001 & \mbox{na} \\
		1359 & (347.5,42) & 0.196991 & (4373.5,149) & ~ & 1360 & (733,61) & 0.1969879 & (3861,140) \\
		1361 & (1491,87) & 0.19698662 & (11346.5,240) & ~ & 1362 & (2471,112) & 0.1969862 & (3484.5,133) \\
		1363 & (640,57) & 0.1969812 & (2383.5,110) & ~ & 1364 & (1357,83) & 0.1969795 & (19298,313) \\
		1365 & (2836.5,120) & 0.1969792 & \mbox{na} & ~ & 1366 & (1108,75) & 0.1969773 & (4028,143) \\
\hline
		\end{array}
		\end{equation*}}
	{\scriptsize
		\begin{equation*}\begin{array}{lllllllll}
		\hline\\
		n & (h_n, k_n) & \lambda_n&\mbox{Remarks}&&n & (h_n, k_n) & \lambda_n&\mbox{Remarks}\vspace{2mm}\\
		\hline\\
		1489 & (71,19) & 0.1966538 & (3745.5,138) & ~ & 1490 & (2127,104) & 0.1966524 & (3272.5,129) \\
		1491 & (6091.5,176) & 0.196652 & (9175,216) & ~ & 1492 & (9955.5,225) & 0.1966517 & (15417.5,280) \\
		1493 & (3799.5,139) & 0.196651 & (11804,245) & ~ & 1494 & (780.5,63) & 0.1966471 & (4724.5,155) \\
		1495 & (4134.5,145) & 0.1966465 & (6585.5,183) & ~ & 1496 & (5289,164) & 0.19664605 & (36191,429) \\
		1497 & (9691.5,222) & 0.19664582 & (9779,223) & ~ & 1498 & (6876.5,187) & 0.1966454 & (26486,367) \\
\hline
		\end{array}
		\end{equation*}}
	{\scriptsize
		\begin{equation*}\begin{array}{lllllllll}
		\hline\\
		n & (h_n, k_n) & \lambda_n&\mbox{Remarks}&&n & (h_n, k_n) & \lambda_n&\mbox{Remarks}\vspace{2mm}\\
		\hline\\
		1621 & (3793.5,139) & 0.19634037 & \mbox{na} & ~ & 1622 & (3959,142) & 0.1963401 & \mbox{na} \\
		1623 & (3631.5,136) & 0.1963394 & (5541.5,168) & ~ & 1624 & (5808.5,172) & 0.19633898 & (10206.5,228) \\
		1625 & (9851.5,224) & 0.19633882 & (12667,254) & ~ & 1626 & (3421,132) & 0.19633794 & (36472,431) \\
		1627 & (4128,145) & 0.19633733 & (11215,239) & ~ & 1628 & (4717,155) & 0.19633714 & (13170.5,259) \\
		1629 & (3473,133) & 0.19633627 & (4839.5,157) & ~ & 1630 & (50867,509) & 0.196336252 & \mbox{na} \\
\hline
		\end{array}
		\end{equation*}}
	{\scriptsize
		\begin{equation*}\begin{array}{lllllllll}
		\hline\\
		n & (h_n, k_n) & \lambda_n&\mbox{Remarks}&&n & (h_n, k_n) & \lambda_n&\mbox{Remarks}\vspace{2mm}\\
		\hline\\
		1753 & (11104.5,238) & 0.19604008 & (23605,347) & ~ & 1754 & (9836.5,224) & 0.1960398 & (11198,239) \\
		1755 & (10370.5,230) & 0.19603965 & (41482,460) & ~ & 1756 & (9924.5,225) & 0.19603936 & (15369.5,280) \\
		1757 & (10280.5,229) & 0.19603922 & (19452,315) & ~ & 1758 & (13150.5,259) & 0.19603918 & \mbox{na} \\
		1759 & (11863.5,246) & 0.19603909 & \mbox{na} & ~ & 1760 & (5533,168) & 0.1960386 & (77561,629) \\
		1761 & (1073.5,74) & 0.1960364 & (12154.5,249) & ~ & 1762 & (880,67) & 0.1960355 & \mbox{na} \\
		1763 & (593,55) & 0.1960306 & (2683.5,117) & ~ & 1764 & (2039.5,102) & 0.19602966 & (5867,173) \\
		1765 & (2459,112) & 0.19602938 & \mbox{na} & ~ & 1766 & (4236,147) & 0.19602909 & \mbox{na} \\
		1767 & (7454,195) & 0.19602897 & \mbox{na} & ~ & 1768 & (2729.5,118) & 0.19602791 & (21347.5,330) \\
		1769 & (4649,154) & 0.19602789 & \mbox{na} & ~ & 1770 & (2329,109) & 0.19602665 & (4352,149) \\
		1771 & (4178.5,146) & 0.1960265 & \mbox{na} & ~ & 1772 & (379.5,44) & 0.1960189 & (9748,223) \\
		
		1773 & (433,47) & 0.196016 & (6854.5,187) & ~ & 1774 & (132.5,26) & 0.19601 & \mbox{na} \\
		1775 & (63.5,18) & 0.195965 & (283,38) & ~ & ~ & ~ & ~ & ~ \\	
		\hline
		\end{array}
		\end{equation*}\vspace{1mm}}

	\noindent For $n=1776,1777,\cdots 3121$, we use Macbeath second Lemma (Lemma \ref{lem6}) and  choose suitable integers $h_n$ and $k_n$ (see Table II), such that
	\begin{equation}\label{Eqn21}\vline~ h_n-tk_n^2~\vline <d^3/8.\end{equation}
	
	\noindent The choice of $(h_n,k_n)$ is done exactly as  for Lemma \ref{lem5}. Here $a_n=a_n(d)= (h_n-\frac{d^3}{8})/k^2_n$, $ b_n=b_n(d)=(h_n+\frac{d^3}{8})/k^2_n$. For a fixed value of the pair $(h_n,k_n)$, $h_n,k_n \in \mathbb{Z}$, (\ref{Eqn21}) is satisfied if $t$ lies in the interval $I_n=(a_n,b_n)\supseteq[\lambda_n,\lambda_{n-1}]$ and $h_n/k^2_n $ is the center of the interval $I_n$.
	{\scriptsize
		\begin{equation*}
\begin{array}{lllllllll}
		&&&&\mbox{\normalsize{Table II.}}&&&&\vspace{2mm}\\
		\hline\\
		n & (h_n, k_n) & \lambda_n&\mbox{Remarks}&&n & (h_n, k_n) & \lambda_n&\mbox{Remarks}\vspace{2mm}\\
		\hline\\
		1776 & (933,69) & 0.19596381 & \mbox{na} & ~ & 1777 & (1223,79) & 0.1959596 & (4708,155) \\
		1778 & (1882,98) & 0.19595832 & \mbox{na} & ~ & 1779 & (6070,176) & 0.19595764 & (5730,171) \\
		1780 & (4954,159) & 0.1959576 & \mbox{na} & ~ & 1781 & (2869,121) & 0.1959555 & (4177,146) \\
		1782 & (3013,124) & 0.19595368 & (11287,240) & ~ & 1783 & (10276,229) & 0.19595361 & \mbox{na} \\
\hline
		\end{array}
		\end{equation*}}
	{\scriptsize
		\begin{equation*}\begin{array}{lllllllll}
		\hline\\
		n & (h_n, k_n) & \lambda_n&\mbox{Remarks}&&n & (h_n, k_n) & \lambda_n&\mbox{Remarks}\vspace{2mm}\\
		\hline\\
		1836 & (853,66) & 0.1958182 & (2118,104) & ~ & 1837 & (571,54) & 0.1958106 & (11279,240) \\
		1838 & (1131,76) & 0.1958075 & (3729,138) & ~ & 1839 & (1551,89) & 0.1958066 & \mbox{na} \\
		1840 & (8069,203) & 0.1958064 & \mbox{na} & ~ & 1841 & (10001,226) & 0.1958063 & \mbox{na} \\
		1842 & (8553,209) & 0.1958056 & (9477,220) & ~ & 1843 & (9477,220) & 0.1958055 & \mbox{na} \\
		1844 & (7599,197) & 0.1958047 & (9391,219) & ~ & 1845 & (550,53) & 0.1957935 & (802,64) \\
\hline
		\end{array}
		\end{equation*}}
	{\scriptsize
		\begin{equation*}\begin{array}{lllllllll}
		\hline\\
		n & (h_n, k_n) & \lambda_n&\mbox{Remarks}&&n & (h_n, k_n) & \lambda_n&\mbox{Remarks}\vspace{2mm}\\
		\hline\\
		1968 & (10313,230) & 0.1949529 & \mbox{na} & ~ & 1969 & (2232,107) & 0.1949508 & (29047,386) \\
		1970 & (1343,83) & 0.1949462 & (5634,170) & ~ & 1971 & (1126,76) & 0.1949418 & (7489,196) \\
		1972 & (1650,92) & 0.1949414 & \mbox{na} & ~ & 1973 & (4270,148) & 0.1949411 & \mbox{na} \\
		1974 & (3244,129) & 0.19494021 & \mbox{na} & ~ & 1975 & (1186,78) & 0.1949349 & (6246,179) \\
		1976 & (1686,93) & 0.194934 & \mbox{na} & ~ & 1977 & (2316,109) & 0.1949318 & (15612,283) \\
\hline
		\end{array}
		\end{equation*}}
	{\scriptsize
		\begin{equation*}\begin{array}{lllllllll}
		\hline\\
		n & (h_n, k_n) & \lambda_n&\mbox{Remarks}&&n & (h_n, k_n) & \lambda_n&\mbox{Remarks}\vspace{2mm}\\
		\hline\\
		2100 & (4825,158) & 0.1932784 & \mbox{na} & ~ & 2101 & (2737,119) & 0.1932762 & \mbox{tbd} \\
		2102 & (14723,276) & 0.19327611 & \mbox{na} & ~ & 2103 & (6904,189) & 0.19327579 & \mbox{na} \\
		2104 & (5455,168) & 0.1932747 & (17279,299) & ~ & 2105 & (21950,337) & 0.1932745 & (21950,337) \\
		2106 & (21432,333) & 0.19327436 & \mbox{na} & ~ & 2107 & (14937,278) & 0.1932742 & \mbox{na} \\
		2108 & (4764,157) & 0.19327291 & (5919,175) & ~ & 2109 & (27179,375) & 0.1932729 & \mbox{na} \\
\hline
		\end{array}
		\end{equation*}}
	{\scriptsize
		\begin{equation*}\begin{array}{lllllllll}
		\hline\\
		n & (h_n, k_n) & \lambda_n&\mbox{Remarks}&&n & (h_n, k_n) & \lambda_n&\mbox{Remarks}\vspace{2mm}\\
		\hline\\
		2232 & (1085,75) & 0.1928861 & (3412,133) & ~ & 2233 & (1361,84) & 0.1928833 & (5509,169) \\
		2234 & (2824,121) & 0.19288191 & (4168,147) & ~ & 2235 & (5907,175) & 0.19288111 & (9939,227) \\
		2236 & (6389,182) & 0.19288081 & \mbox{na} & ~ & 2237 & (6963,190) & 0.19288045 & \mbox{na} \\
		2238 & (13443,264) & 0.19288028 & \mbox{na} & ~ & 2239 & (13957,269) & 0.19288025 & \mbox{na} \\
		2240 & (13545,265) & 0.1928801 & \mbox{na} & ~ & 2241 & (3310,131) & 0.1928783 & (38539,447) \\
\hline
		\end{array}
		\end{equation*}}
	{\scriptsize
		\begin{equation*}\begin{array}{lllllllll}
		\hline\\
		n & (h_n, k_n) & \lambda_n&\mbox{Remarks}&&n & (h_n, k_n) & \lambda_n&\mbox{Remarks}\vspace{2mm}\\
		\hline\\
		2364 & (1809,97) & 0.1922612 & (9820,226) & ~ & 2365 & (889,68) & 0.1922544 & (8398,209) \\
		2366 & (5108,163) & 0.19225353 & (12306,253) & ~ & 2367 & (6160,179) & 0.19225318 & \mbox{na} \\
		2368 & (6229,180) & 0.1922526 & (10527,234) & ~ & 2369 & (36379,435) & 0.19225253 & \mbox{na} \\
		2370 & (93398,697) & 0.19225249 & \mbox{tbd} & ~ & 2371 & (35713,431) & 0.19225247 & \mbox{na} \\
		2372 & (13603,266) & 0.19225239 & \mbox{na} & ~ & 2373 & (9305,220) & 0.1922518 & (66470,588) \\
\hline
		\end{array}
		\end{equation*}}
	{\scriptsize
		\begin{equation*}\begin{array}{lllllllll}
		\hline\\
		n & (h_n, k_n) & \lambda_n&\mbox{Remarks}&&n & (h_n, k_n) & \lambda_n&\mbox{Remarks}\vspace{2mm}\\
		\hline\\
		2496 & (1385,85) & 0.1916933 & \mbox{na} & ~ & 2497 & (9965,228) & 0.191693291 & \mbox{na} \\
		2498 & (18780,313) & 0.191693253 & (193615,1005) & ~ & 2499 & (91795,692) & 0.191693251 & \mbox{na} \\
		2500 & (60115,560) & 0.191693228 & (11695,247) & ~ & 2501 & (99650,721) & 0.191693223 & \mbox{na} \\
		2502 & (53441,528) & 0.191693205 & (11695,247) & ~ & 2503 & (27535,379) & 0.191693154 & (11695,247) \\
		2504 & (32381,411) & 0.191693137 & \mbox{na} & ~ & 2505 & (119333,789) & 0.19169313 & (11695,247) \\
\hline
		\end{array}
		\end{equation*}}
	{\scriptsize
		\begin{equation*}\begin{array}{lllllllll}
		\hline\\
		n & (h_n, k_n) & \lambda_n&\mbox{Remarks}&&n & (h_n, k_n) & \lambda_n&\mbox{Remarks}\vspace{2mm}\\
		\hline\\
		2628 & (3802,141) & 0.1912371 & (5592,171) & ~ & 2629 & (1654,93) & 0.191236 & \mbox{na} \\
		2630 & (759,63) & 0.1912281 & (6059,178) & ~ & 2631 & (2529,115) & 0.1912277 & \mbox{na} \\
		2632 & (937,70) & 0.1912213 & (7958,204) & ~ & 2633 & (3485,135) & 0.1912208 & (5206,165) \\
		2634 & (1019,73) & 0.1912149 & (2356,111) & ~ & 2635 & (3182,129) & 0.1912139 & \mbox{tbd} \\
		2636 & (2485,114) & 0.1912115 & (2846,122) & ~ & 2637 & (2846,122) & 0.1912112 & \mbox{na} \\
\hline
		\end{array}
		\end{equation*}}
	{\scriptsize
		\begin{equation*}\begin{array}{lllllllll}
		\hline\\
		n & (h_n, k_n) & \lambda_n&\mbox{Remarks}&&n & (h_n, k_n) & \lambda_n&\mbox{Remarks}\vspace{2mm}\\
		\hline\\
		2760 & (1613,92) & 0.19057 & (9140,219) & ~ & 2761 & (1944,101) & 0.190568 & (17961,307) \\
		2762 & (2101,105) & 0.1905655 & (12198,253) & ~ & 2763 & (7025,192) & 0.19056538 & \mbox{na} \\
		2764 & (6952,191) & 0.19056452 & (8891,216) & ~ & 2765 & (7699,201) & 0.1905642 & \mbox{na} \\
		2766 & (19392,319) & 0.19056401 & (23746,353) & ~ & 2767 & (23746,353) & 0.19056396 & \mbox{na} \\
		2768 & (26371,372) & 0.1905639 & (40851,463) & ~ & 2769 & (1793,97) & 0.1905606 & (12882,260) \\
\hline
		\end{array}
		\end{equation*}}
	{\scriptsize
		\begin{equation*}\begin{array}{lllllllll}
		\hline\\
		n & (h_n, k_n) & \lambda_n&\mbox{Remarks}&&n & (h_n, k_n) & \lambda_n&\mbox{Remarks}\vspace{2mm}\\
		\hline\\
		2892 & (493,51) & 0.1895365 & \mbox{tbd} & ~ & 2893 & (2170,107) & 0.18953548 & (11846,250) \\
		2894 & (2050,104) & 0.1895337 & (5871,176) & ~ & 2895 & (1213,80) & 0.1895319 & \mbox{na} \\
		2896 & (1010,73) & 0.1895261 & (10556,236) & ~ & 2897 & (3203,130) & 0.1895257 & \mbox{na} \\
		2898 & (3930,144) & 0.18952471 & (9007,218) & ~ & 2899 & (5413,169) & 0.18952363 & (9766,227) \\
		2900 & (6914,191) & 0.18952357 & \mbox{na} & ~ & 2901 & (5738,174) & 0.1895227 & (13209,264) \\
\hline
		\end{array}
		\end{equation*}}
	{\scriptsize
		\begin{equation*}\begin{array}{lllllllll}
		\hline\\
		n & (h_n, k_n) & \lambda_n&\mbox{Remarks}&&n & (h_n, k_n) & \lambda_n&\mbox{Remarks}\vspace{2mm}\\
		\hline\\
		3024 & (4308,151) & 0.1889384 & (35424,433) & ~ & 3025 & (823,66) & 0.1889319 & (7482,199) \\
		3026 & (703,61) & 0.1889235 & (2952,125) & ~ & 3027 & (1634,93) & 0.1889224 & \mbox{na} \\
		3028 & (1463,88) & 0.1889194 & (4897,161) & ~ & 3029 & (1333,84) & 0.1889151 & (1599,92) \\
		3030 & (3047,127) & 0.1889145 & \mbox{na} & ~ & 3031 & (1179,79) & 0.1889097 & (12284,255) \\
		3032 & (848,67) & 0.188906 & (22746,347) & ~ & 3033 & (170,30) & 0.1888715 & \mbox{tbd} \\
\hline
		\end{array}
		\end{equation*}\vspace{2mm}}
	
	\noindent For example take $(h_{1776},k_{1776})=(933,69)$ in (\ref{Eqn21}) to get $|933-t.69^2|<d^3/8$. This means for $t\in I_{1776}= (a_{1776}, b_{1776})$, $a_{1776}=(933-\frac{d^3}{8})/69^2$, $ b_{1776}=(933+\frac{d^3}{8})/69^2$,   (\ref{eq101'}) is soluble by Lemma \ref{lem6}. Now
	
	\begin{equation*}
	%
	\right.\vspace{2mm}\\
	&<0.19596381=\lambda_{1776}\end{array}
	\end{equation*}
	
	\noindent  where $\alpha $ is a root of $\frac{3}{2}d^3- (933-\frac{d^3}{8})/69^2 =0$ satisfying $0.5074<\alpha<0.5075$. This means (\ref{eq101'}) is soluble for $t\in[0.19596381,0.195965]$, $\lambda_{1775}=0.195965.$ We repeat this process 1346 times.
	Take $(h_{3121},k_{3121})=(3,4)$ in (\ref{Eqn21}) to get $|3-16t|<d^3/8$. This means for $t\in I_{3121}= (a_{3121}, b_{3121})$, $a_{3121}=(3-\frac{d^3}{8})/16$, $ b_{3121}=(3+\frac{d^3}{8})/16$,   (\ref{eq101'}) is soluble by Lemma \ref{lem6}. As $\frac{3}{16}<b_{3121}$, and $t\geq \frac{3}{16}$, we are done.\vspace{2mm}

	\noindent For $t=\frac{h_n}{k_n^2}$ obtained by Lemma \ref{lem6}, we go back to Lemma \ref{lem5}, which is possible as $d\geq \frac{1}{2}$ here. When $t=\frac{h_n}{k_n^2}$, in most of the  cases, alternate pair of $(h'_n,k'_n)$ are obtained to satisfy $|h_n-tk_n^2|+\frac{1}{2}<d$, and in many cases the ratio $h_n/k_n^2$ is outside the region $[\lambda_n,\lambda_{n-1}]$. In total, $203=(114 \mbox{~in table I and~} 89 \mbox{~in table II})$ values of $t$ are left to be discussed.  For these  values of $t=\frac{h_n}{k_n^2}$, $\beta_{t,G}=\pm a_5''-2tc_5\equiv h/k\pmod {\frac{1}{k},2t}$, Lemma \ref{lem15} fixes the values of $(a_5'',c_5)$ to be considered and we have,
	\begin{equation}\label{Eqn22}
	\begin{array}{llll}
	(\frac{2t}{3})^{1/3} &\leq &d &\leq {\rm ~min~}\{\frac{8t}{3}, \frac{2}{3} \}\vspace{1mm}\\
	
	d/2 &\leq &a&\leq {\rm ~min~} \{(\frac{d^5}{2t})^{1/2}, (\frac{3d^5}{4})^{1/3} \} \vspace{1mm}\\
	
	3d/8 &\leq & A&\leq \frac{4t}{3}.
	
	\end{array}
	\end{equation}
	
	\noindent For each value of $t$ and $(a_5'',c_5)$, except for $t=1/3,$ $1/4,$ $1/5$,  and $3/16$, using (\ref{Eqn22}) we take $x_2=\pm 1$ and give suitable value to  $x_5$, which covers the range of $f_1''$ or, $p_1''$ and hence (\ref{eq101'}) is soluble.  We are not giving the details here, the proof being routine. For $t=\frac{1}{5}$, we need to consider $(a_5'',c_5)=(0,0),$ $(0,1/2),$ $(1/10,1/4),$ then tables similar to those in Lemma \ref{lem81} give a solution to (\ref{eq101'}).	We  discuss here $t=3/16$, $t=1/3$ and $t=1/4$ in following lemmas.
	
	\begin{lemma} For $t=\frac{3}{16}$, $\pm 8a_5''-3c_5 \equiv 0 \pmod 1$,  $m=1, K=1$, $(\ref{eq101'})$ is soluble.\end{lemma}
	
	\noindent {\bf Proof:} ~ Here, since $\frac{d}{a}\leq 2$, and $\frac{3}{16}=t=\frac{3d^5}{8aA}\leq 3d^3/2$, we get $d^3\geq \frac{1}{8}$ i.e., $d\geq 1/2$. And $d\leq 8t/3=1/2$. Together, we get $d=1/2$. Now, $1/4 =d/2 \leq A \leq 4t/3=1/4$ and, $t=3d^5/8aA$ gives $A=1/4, a=1/4$. Also, from $A\leq C=t+A\lambda^2$, we get $\lambda^2 \geq 1/4$. But, $0\leq \lambda \leq 1/2$, so that $\lambda=1/2$.   Therefore, we have \vspace{-1mm}
	\begin{equation}\label{Eqn30}
	\begin{array}{l}
	d=1/2; ~ a=1/4; ~A=1/4;~ t=3/16; {\rm ~and~} \lambda=1/2.
	\end{array}
	\end{equation}
	
	\noindent Using Lemma \ref{lem15} we need to consider $(a_5'',c_5)=(0,0),$ $(0,1/3),$ $(1/16,1/6),$ and $(1/16,1/2).$\vspace{2mm}

	\noindent When $(a_5'',c_5)=(0,0)$, the following table gives a solution to (\ref{eq101'}) i.e., of\vspace{-1mm}
	\begin{equation}\label{Eqn31}0<G=(x_1+a_2''x_2)x_2-\frac{3}{16}x_5^2-\frac{a+A}{4}<d. \vspace{-1mm}\end{equation}

	{\scriptsize
		\begin{equation*}\begin{array}{lllll}
		\hline\\
		{\rm ~~~~Range} & x_1 & x_1x_2 & x_5& {\rm ~~G}\vspace{2mm}\\
		\hline\\
		
		\frac{13}{16}<f_1''\leq 1&\pm 1+c_1&1+|c_1|&3 &f_1''-\frac{13}{16}\vspace{2mm}\\
		
		\frac{5}{16}<f_1''<\frac{13}{16}&c_1&|c_1|&1&f_1''-\frac{5}{16}\vspace{2mm}\\
		
		\frac{1}{8}<f_1''\leq \frac{5}{16}&c_1&|c_1|&0&f_1''-\frac{1}{8}\vspace{2mm}\\
		
		-\frac{3}{16}<f_1''\leq \frac{1}{8}&\pm 2+c_1&2+|c_1|&3&f_1''+\frac{3}{16}\vspace{2mm}\\
		
		-\frac{1}{2}<f_1''<-\frac{3}{16}&\pm 1+c_1&1+|c_1|&1&f_1''+\frac{11}{16}\vspace{2mm}\\
		\hline
		\end{array}
		\end{equation*}}

	{\scriptsize
		\begin{equation*}\begin{array}{lllll}
		\hline\\
		{\rm ~~~~Range} & x_1 & x_1x_2 & x_5& {\rm ~~G}\vspace{2mm}\\
		\hline\\

		\frac{1}{8}<g_1''\leq \frac{1}{2}&c_1&-|c_1|&0&g_1''-\frac{1}{8}\vspace{2mm}\\
		
		-\frac{3}{16}<g_1''\leq \frac{1}{8}&\pm 2+c_1&2+|c_1|&3&g_1''+\frac{3}{16}\vspace{2mm}\\
		
		-\frac{11}{16}<g_1''<-\frac{3}{16}&\pm 1+c_1&1-|c_1|&1&g_1''+\frac{11}{16}\vspace{2mm}\\
		
		-1<g_1''<-\frac{11}{16}&\pm 3+c_1&3+|c_1|&3&g_1''+\frac{19}{16}\vspace{2mm}\\
		
		\hline
		\end{array}
		\end{equation*}}
	
	\begin{itemize}
		\item When $f_1''=\frac{13}{16}, g_1''= -\frac{3}{16}$, then $c_1=\frac{1}{2}, a_2''=\frac{5}{16}$. Here $x_1=-\frac{1}{2}, x_2=3, x_5=2$ gives a solution of (\ref{Eqn31}).\vspace{-2mm}
		\item When $f_1''=-\frac{3}{16}=g_1''$, then $c_1=0, a_2''=-\frac{3}{16}$. Here $x_1=8, x_2=14, x_5=20$ gives a solution of (\ref{Eqn31}).\vspace{-2mm}
		\item When $f_1''=-\frac{3}{16}, g_1''=-\frac{11}{16}$, then $c_1=\pm \frac{1}{4}$, $a_2''=-\frac{7}{16}$.
		If $c_1= \frac{1}{4}$, then $x_1=-\frac{7}{4}, x_2=-1, x_5=2$ gives a solution of (\ref{Eqn31}); and if $c_1=-\frac{1}{4}$, then $x_1=\frac{3}{4}, x_2=1, x_5=0$ gives a solution of (\ref{Eqn31}).\vspace{-1mm}
		
	\end{itemize}
	
	\noindent That $f_1''=\frac{13}{16}$ and $ g_1''= -\frac{11}{16}$ can not hold. The other cases  $(a_5'',c_5)=(0,1/3)$, $(1/16,1/6)$ and  $(1/16,1/2)$ are dealt with similarly.

	\begin{lemma} For $t=\frac{1}{3}$,  $m=1, K=1$, $(\ref{eq101})$ is soluble.\end{lemma}
	
	\noindent {\bf Proof: } Here, since $\frac{d}{a}\leq 2$, and $\frac{1}{3}=t=\frac{3d^5}{8aA}\leq 3d^3/2$, we get $d^3\geq \frac{2}{9}$ i.e., $d\geq (\frac{2}{9})^{\frac{1}{3}}>0.6057$. Therefore,\vspace{-1mm}
	\begin{equation}\label{Eqn23}
	\begin{array}{l}
	0.6057 < (2/9)^{1/3} \leq d \leq 2/3<0.6667\vspace{1mm}\\
	0.30285< d/2\leq a\leq 4A/3\leq 16t/9< 0.5926\vspace{1mm}\\
	0.30285 < d/2 \leq A \leq 4t/3< 0.4445=\lambda_0{\rm (say)}.
	\end{array}
	\end{equation}

	\noindent We divide this range of $A$ into $6$ subintervals $[\lambda_n,\lambda_{n-1}], n=1,2,\cdots 6$ and in each subinterval choose suitable integers $2h_n$ and $k_n$ such that
	$\vline~ h_n-Ak_n^2~\vline +\frac{1}{2} <d.$\vspace{2mm}
	
	{\scriptsize
		\begin{equation*}\begin{array}{lllllllll}
		\hline\\
		n & (h_n,k_n) & \lambda_n&\mbox{Remarks}&~~ &n & (h_n,k_n) & \lambda_n&\mbox{Remarks}\vspace{2mm}\\
		\hline\\
		
		1&(0.5,1)&0.366&\mbox{na}&&2&(1.5,2)&0.3438&(0.5,1)\vspace{2mm}\\
		
		3&(3,3)&0.3207&\mbox{tbd}&&4&(8,5)&0.3196&(11.5,6)\vspace{2mm}\\
		
		5&(5,4)&0.3059&\mbox{tbd}&&6&(11,6)&0.30285&(15,7)\vspace{2mm}\\
		
		\hline
		\end{array}
		\end{equation*}\vspace{1mm}}
	
	\noindent The choice of $(h_n,k_n)$ is done in the usual manner. Thus (\ref{eq101}) is soluble unless $A=1/3,$ $5/16$ and $\beta_{A,F} = \pm a'_4-2A(c_4+\lambda x_5)\equiv h_n/k_n \pmod {1/k_n, 2A}$. For each of these values to $A$ we will give solutions of (\ref{eq101}) depending upon values of $a_2'$ and $c_1$.\vspace{2mm}
	
	\noindent \textbf{{Case (i): $A=\frac{1}{3}$, $ \pm 3a'_4-2(c_4+\lambda x_5)\equiv 0 \pmod 1$.}}\vspace{2mm}
	
	\noindent Taking $x_5=c_5$ and $1+c_5$ simultaneously and from (\ref{eq8}), we get $\lambda=0$, or $1/2$.
	
	\noindent For $\lambda=1/2$, we have $C=A\lambda^2+t=5/12$ and $|0.5-C.1^2|+1/2<d$, so that (\ref{eq45}) is satisfied.\vspace{1mm}
	
	\noindent For $\lambda=0$, we have $C=A\lambda^2+t=t=1/3=A$. By Lemma (\ref{lem15}), we need to consider $(a_4',c_4)=(0,0)$, $(0,1/2)$ or $(1/6,1/4)$. Because of symmetry in $x_4$ and $x_5$ we need to consider $(a_4',c_4,a_5',c_5)=(0,0,0,0)$, $(0,0,0,1/2)$, $(0,0,1/6,1/4)$, $(0,1/2,0,1/2)$, $(0,1/2,1/6,1/4)$, and $(1/6,1/4,1/6,1/4)$.
	In each of these, we find that (\ref{eq101}) is soluble by taking $x_2=\pm 1$ and choosing $x_5\equiv c_5 \pmod 1$. \vspace{4mm}
	
	\noindent \textbf{{Case (ii)}:} $A=\frac{5}{16}$, $ \pm 8a'_4-5(c_4+\lambda x_5)\equiv 0 \pmod 1$. \vspace{2mm}
	
	\noindent Taking $x_5=c_5$ and $1+c_5$ simultaneously and using (\ref{eq8}) and (\ref{eq09}), we get $\lambda=0$, or $1/5$, or $2/5$.
	
	\noindent	 For $\lambda=0$, we have $C=A\lambda^2+t=1/3$ and $|12-C.6^2|+1/2<d$, so that (\ref{eq45}) is satisfied.\vspace{1mm}
	
	\noindent For $\lambda=1/5$, we have $C=A\lambda^2+t=83/240$ and $|5.5-C.4^2|+1/2<d$ so that (\ref{eq45}) is satisfied.\vspace{1mm}
	
	\noindent For $\lambda=2/5$, we have $C=A\lambda^2+t=23/60$ and $|1.5-C.2^2|+1/2<d$ so that (\ref{eq45}) is satisfied.
	
	\begin{lemma} For $t=\frac{1}{4}$, $\pm 2a_5''-c_5 \equiv 0 \pmod 1$, $m=1, K=1$, $(\ref{eq101*})$ is soluble.\end{lemma}
	
	\noindent {\bf Proof: } Here, using Lemma \ref{lem15} we need to consider $(a_5'',c_5)=(0,0),$ $(1/4,1/2)$. Since $\frac{d}{a}\leq 2$, and $\frac{1}{4}=t=\frac{3d^5}{8aA}\leq 3d^3/2$, we get $d^3\geq \frac{1}{6}$. Therefore,
	\begin{equation}\label{Eqn35}
	\begin{array}{l}
	0.5503 \leq (1/6)^{1/3} \leq d \leq 2/3 \leq 0.6667\vspace{2mm}\\
	0.2751\leq d/2\leq a\leq 4A/3\leq 16t/9= 4/9\vspace{2mm}\\
	0.2751 \leq d/2 \leq A \leq 4t/3=1/3\vspace{2mm}\\
	0.13758 \leq \frac{a+A}{4} \leq 7/36<0.1945.
	\end{array}
	\end{equation}
	
	\noindent When $(a_5'',c_5)=(1/4,1/2)$, a table similar to that in Lemma \ref{lem28} gives a solution to (\ref{eq101'}).
	
	\noindent When $(a_5'',c_5)=(0,0)$, (\ref{eq101'}) has no solution for $(a_2'',c_1)=(0,0)$. This is so because $0<G=x_1x_2-\frac{1}{4}x_5^2-\frac{a+A}{4}\leq d$ with $x_1,x_2,x_5 \in \mathbb{Z}$ is not soluble for $d<0.556$, as $0.5502<a+A< 4x_1x_2-x_5^2\leq 4(d+\frac{a+A}{4})<3$ has no solution. We note that $4x_1x_2-x_5^2$ never takes the value $1$ or $2$ for integers $x_1,x_2,x_5$, as $-x_5^2 \not \equiv 1 {\rm ~ or ~} 2 \pmod 4$. Therefore we have to go to $F$. For this, we have from (\ref{Eqn35}) \begin{equation}\label{Eqn37}
	0.2751\leq d/2 \leq A\leq 1/3=\lambda_0{\rm (say)}.
	\end{equation}

	\noindent We divide the range of $A$ in $33$ subintervals $[\lambda_n,\lambda_{n-1}], n=1,2,\cdots 33$ and in each subinterval choose suitable integers $2h_n$ and $k_n$ such that
	$\vline~h_n-Ak_n^2~\vline +\frac{1}{2} <d.$
	
	{\scriptsize
		\begin{equation*}\begin{array}{lllllllll}
		\hline\\
		n & (h_n, k_n) & \lambda_n&\mbox{Remarks}&&n & (h_n, k_n) & \lambda_n&\mbox{Remarks}\vspace{2mm}\\
		\hline\\
		
		1&(3,3)&0.32514&\mbox{tbd}&&2&(32.5,10)&0.32427&\mbox{tbd}\vspace{2mm}\\
		
		3&(63.5,14)&0.32361&(83,16)&&4&(93.5,17)&0.32342&\mbox{tbd}\vspace{2mm}\\
		
		5&(46.5,12)&0.32291&\mbox{tbd}&&6&(8,5)&0.31719&(11.5,6)\vspace{2mm}\\
		
		7&(15.5,7)&0.31688&\mbox{na}&&8&(5,4)&0.30837&\mbox{tbd}\vspace{2mm}\\
		
		9&(25,9)&0.30783&\mbox{na}&&10&(52,13)&0.30738&\mbox{tbd}\vspace{2mm}\\
		11&(11,6)&0.3057&\mbox{na}&&12&(19.5,8)&0.30369&\mbox{tbd}\vspace{2mm}\\
		13&(59.5,14)&0.30326&\mbox{tbd}&&14&(24.5,9)&0.30254&\mbox{na}\vspace{2mm}\\
		
		15&(7.5,5)&0.29756&\mbox{tbd}&&16&(36,11)&0.29716&(1673.5,75)\vspace{2mm}\\
		
		17&(14.5,7)&0.2947&(24,9)&&18&(29.5,10)&0.2944&\mbox{na}\vspace{2mm}\\
		19&(142.5,22)&0.29432&\mbox{na}&&20&(85,17)&0.29391&\mbox{tbd}\vspace{2mm}\\
		
		21&(155.5,23)&0.29388&\mbox{na}&&22&(35.5,11)&0.29331&(6778.5,152)\vspace{2mm}\\
		
		23&(10.5,6)&0.29007&\mbox{tbd}&&24&(23.5,9)&0.28996&\mbox{na}\vspace{2mm}\\
		
		25&(18.5,8)&0.28818&(35,11)&&26&(41.5,12)&0.2878&\mbox{na}\vspace{2mm}\\
		
		27&(194.5,26)&0.28764&(13300,215)&&28&(115,20)&0.28736&\mbox{tbd}\vspace{2mm}\\
		
		29&(152,23)&0.28732&(4931,131)&&30&(48.5,13)&0.28686&(73.5,16)\vspace{2mm}\\
		
		31&(14,7)&0.28469&\mbox{tbd}&&32&(4.5,4)&0.28385&\mbox{na}\vspace{2mm}\\
		
		33&(2.5,3)&0.27515&\mbox{tbd}\vspace{2mm}\\
		
		\hline
		\end{array}
		\end{equation*}\vspace{1mm}}
	
	\noindent The choice of $(h_n,k_n)$ is done in a manner already described. Thus, (\ref{eq101}) is soluble unless $A=1/3,$ $13/40,$ $11/34,$ $93/288,$ $5/16,$ $4/13,$ $39/128,$ $17/56,$ $3/10,$ $5/17,$ $7/24,$ $23/80,$ $2/7,$ $5/18$ and $\beta_{A,F} = \pm a'_4-2A(c_4+\lambda x_5)\equiv h_n/k_n \pmod {1/k_n, 2A}$. On taking $x_5=c_5$ and $1+c_5$ simultaneously and using (\ref{eq8}) and (\ref{eq09}), we find that (\ref{eq101}) is soluble unless
	$A=\frac{1}{3}$, $\lambda=\frac{1}{2}$;~~
	$A=\frac{93}{288}$, $\lambda=\frac{15}{31}$;~~
	$A=\frac{4}{13}$, $\lambda=\frac{1}{2}$;~~ $A=\frac{39}{128}$, $\lambda=\frac{17}{39}, \frac{18}{39}, \frac{19}{39}$;~~
	$A=\frac{17}{56}$, $\lambda=\frac{8}{17}$;~~
	$A=\frac{5}{17}$, $\lambda=\frac{2}{5}, \frac{1}{2}$;~~ $A=\frac{7}{24}$, $\lambda=\frac{3}{7}$;~~
	$A=\frac{23}{80}$, $\lambda=\frac{9}{23}, \frac{10}{23}, \frac{11}{23}$;~~
	$A=\frac{2}{7}$, $\lambda=\frac{1}{2}$;~~
	$A=\frac{5}{18}$, $\lambda=\frac{2}{5}$.\vspace{2mm}
	
	\noindent For $A=\frac{93}{288}$, $\lambda=\frac{15}{31}$; $A=\frac{4}{13}$, $\lambda=\frac{1}{2}$; $A=\frac{39}{128}$, $\lambda=\frac{17}{39}, \frac{18}{39}, \frac{19}{39}$; $A=\frac{17}{56}$, $\lambda=\frac{8}{17}$; $A=\frac{5}{17}$, $\lambda=\frac{2}{5}, \frac{1}{2}$; $A=\frac{23}{80}$, $\lambda=\frac{9}{23}, \frac{10}{23}, \frac{11}{23}$; $A=\frac{2}{7}$, $\lambda=\frac{1}{2}$; and $A=\frac{5}{18}$, $\lambda=\frac{2}{5}$, we find that (\ref{eq45}) is satisfied on taking $(h,k)=(55,31),$ $(16,7),$ $(52,13),$ $(25.5,9),$ $(39,11),$ $(114.5,19),$ $(20.5,9),$ $(93.5,17),$ $(85,17),$ $(19.5,8),$ $(45.5,12),$ $(116,19),$ $(26.5,9)$ respectively and noting that $d\geq (A/3)^{1/4}$ as $A=3d^5/2a \leq 3d^4$. When $A=7/24$ and $\lambda = 3/7$, we have $At/C=49/204$ and $|19.5 - (49/204)9|+1/2<d$; proving thereby (\ref{eq45'}). \vspace{2mm}
	
	\noindent When $A=1/3$ and $\lambda=1/2$ we have $C=A$. Here,
	$F=(x_1+a_2'x_2+a_4'x_4+a_5'x_5)x_2-(x_4^2+x_4x_5+x_5^2)/3-a/4$ is symmetrical in $x_4$ and $x_5$.
	\noindent We already had $c_5=0$. Because of symmetry in $x_4$ and $x_5$, we must have $c_4=0$. Also $\beta_{A,F} \equiv 3/3 \pmod {1/3,2A}$. Taking $c_4=0, x_5=0$ in $\beta_F$, we must have $3a_4'\equiv 0 \pmod 1$ i.e., $a_4'=0$ or $1/3$. By symmetry we have $a_5'=0$ or $1/3$. But from (\ref{eq16}) we have $0=a_5''\equiv a_5'-a_4'/2 \pmod 1$ which gives $a_5'-a_4'/2=0$, as $|a_5'-a_4'/2| \leq 3/4$. So we must have $a_4'=a_5'=0$.\vspace{2mm}

	\noindent Now, we are left with $A=1/3, \lambda=1/2, C=A=1/3$; $(a_4',c_4,a_5',c_5)=(0,0,0,0)$, $\left(1/9\right)^{1/4}=\left(A/3\right)^{1/4}\leq d\leq 2/3 {\rm ~and, ~} 0.2886< d/2\leq a= 9d^5/2\leq 4A/3=4/9=\lambda_0({\rm say}).$ We find that (\ref{eq13}) has no solution for $(a_2',c_1)=(0,0)$. This is so because $0<F=x_1x_2-\frac{1}{3}(x_4^2+x_4x_5+x_5^2)-\frac{a}{4}\leq d$ has no solution for $d<0.556$ with $x_1,x_2,x_4,x_5 \in \mathbb{Z}$, as $0.216<\frac{3a}{4}<3x_1x_2-(x_4^2+x_4x_5+x_5^2)\leq 3(d+a/4)<2$  is not possible because $3x_1x_2-(x_4^2+x_4x_5+x_5^2)$ never takes the value $1$ for integers $x_1,x_2,x_4,x_5$.  So we will find the solution of (\ref{eq101*}).\vspace{2mm}
	
	\noindent Take $x_1=x+c_1$, $x_2=\pm 1$, $x_3=y+c_3$, and choose $(x_4,x_5)\equiv (c_4,c_5)\pmod 1$ arbitrary to get $Q=\pm x-ay^2+\beta_{a,Q}y+\nu$, where $\beta_{a,Q}=\pm a_3-2a(c_3+h_4x_4+h_5x_5)$. For a solution to $0<Q<d$, by Macbeath's Lemma (Lemma \ref{lem5}), we divide the range of $a$ into $17$ subintervals $[\lambda_n, \lambda_{n-1}], n=1, 2, \cdots, 17$ and in each subinterval choose suitable integers $2h_n$ and $k_n$ such that
	$\vline~h_n-ak_n^2~\vline + \frac{1}{2} <d.$\vspace{2mm}

	{\scriptsize
		\begin{equation*}\begin{array}{lllllllll}
		\hline\\
		n & (h_n, k_n) & \lambda_n&\mbox{Remarks}&&n & (h_n, k_n) & \lambda_n&\mbox{Remarks}\vspace{2mm}\\
		\hline\\
		
		1&(7,4)&0.42968&(4,3)&&2&(15.5,6)&0.42711&\mbox{na}\vspace{2mm}\\
		
		3&(34.5,9)&0.4249&\mbox{tbd}&&4&(10.5,5)&0.4152&(15,6)\vspace{2mm}\\
		
		5&(15,6)&0.4137&\mbox{na}&&6&(6.5,4)&0.3991&\mbox{tbd}\vspace{2mm}\\
		
		7&(1.5,2)&0.35&\mbox{tbd}&&8&(5.5,4)&0.3378&\mbox{na}\vspace{2mm}\\
		
		9&(3,3)&0.3233&\mbox{tbd}&&10&(8,5)&0.318&(11.5,6)\vspace{2mm}\\
		
		11&(5,4)&0.3073&\mbox{tbd}&&12&(11,6)&0.3033&(19.5,8)\vspace{2mm}\\
		
		13&(7.5,5)&0.2968&\mbox{tbd}&&14&(14.5,7)&0.2943&(24,9)\vspace{2mm}\\
		
		15&(85,17)&0.29387&\mbox{tbd}&&16&(10.5,6)&0.28951&(3790.5,114)\vspace{2mm}\\
		
		17&(18.5,8)&0.2886&(35,11)\vspace{2mm}\\
		
		\hline
		\end{array}
		\end{equation*}\vspace{1mm}}
	
	\noindent The choice of $(h_n,k_n)$ is done in the prescribed manner. Thus, (\ref{eq101*}) is soluble unless $a=23/54,$ $13/32,$ $3/8,$ $1/3,$ $5/16,$ $3/10,$ $5/17,$  and $\beta_{a,Q} \equiv h/k \pmod {1/k, 2a}$.
	
	\noindent On taking $x_4=c_4$ and $1+c_4$ simultaneously and using $a\leq A+h_4^2 a$ and $-1/2<h_4\leq 1/2$, we find that (\ref{eq101*}) is soluble unless
	$a=\frac{23}{54}$, $h_4=\frac{11}{23}$;~~
	$a=\frac{13}{32}$, $h_4=\frac{6}{13}$;~~
	$a=\frac{3}{8}$, $h_4=\frac{1}{3}$;~~
	$a=\frac{1}{3}$, $h_4=0,\frac{1}{2}$;~~
	$a=\frac{5}{16}$, $h_4=0,\frac{1}{5}, \frac{2}{5}$;~~
	$a=\frac{3}{10}$, $h_4=0,\frac{1}{3}$;~~
	$a=\frac{5}{17}$, $h_4=0,\frac{1}{10},\frac{2}{10},\frac{3}{10},\frac{4}{10},\frac{5}{10}$.\vspace{2mm}
	
	\noindent Now take $x_1=x+c_1$, $x_2=\pm 1$, $x_4=y+c_4=y$, $x_5=0$ and choose $x_3\equiv c_3\pmod 1$ arbitrary to get $Q=\pm x+\beta_a y-(A+ah_4^2)y^2+\nu$ where, $\beta_{a}=\pm a_4-2ah_4(x_3+h_4c_4)$ and $\nu$ is some real number. We find suitable integers $2h, k$ satisfying \begin{equation}\label{x}\vline~h-(A+ah_4^2).k^2~\vline+\frac{1}{2}<d,\end{equation} and apply Lemma \ref{lem5}. For $a=\frac{23}{54}$, $h_4=\frac{11}{23}$, we find that (\ref{x}) is satisfied for $(h,k)=(0.5,1)$. Similarly we find that for  $a=\frac{13}{32}$, $h_4=\frac{6}{13}$; $a=\frac{1}{3}$, $h_4=\frac{1}{2}$; $a=\frac{5}{16}$, $h_4=\frac{1}{5}, \frac{2}{5}$; $a=\frac{3}{10}$, $h_4=\frac{1}{3}$; and $a=\frac{5}{17}$, $h_4=\frac{1}{10},\frac{2}{10},\frac{3}{10},\frac{4}{10},\frac{5}{10}$; (\ref{x}) is satisfied for the values $(h,k)=$ $(0.5,1),$ $(0.5,1),$ $(28,9),$ $(169,21),$ $(62,13),$ $(16.5,7),$ $(28,9),$ $(23,8),$ $(46,11)$ and $(6.5,4)$ respectively.  Therefore we are left with \vspace{2mm}
	
	\noindent $\begin{array}{l}
	\mbox{(i)}~ a=\frac{3}{8}, h_4=\frac{1}{3};~ \mbox{(ii)}~ ~a=\frac{1}{3}, h_4=0; ~\mbox{(iii)}~ ~a=\frac{5}{16}, h_4=0;~\mbox{(iv)}~ ~a=\frac{3}{10}, h_4=0;\vspace{2mm}\\ \mbox{(v)}~ a=\frac{5}{17}, h_4=0.
	\end{array}$\vspace{2mm}
	
	\noindent  Interchanging the role of $x_4$ and $x_5$ and taking $x_1=x+c_1$, $x_2=\pm 1$, $x_5=y+c_5=y$, $x_4=0$, choosing $x_3\equiv c_3\pmod 1$ arbitrary we get $Q=\pm x+\beta_a y-(A+ah_5^2)y^2+\nu$.\vspace{2mm}
	
	\noindent Now we find that (\ref{eq101*}) is soluble unless
	\begin{itemize}
		\item $a=\frac{3}{8}$, with $(h_4,h_5)= (\frac{1}{3},\frac{1}{3})$
		\item $a=\frac{1}{3}, \frac{5}{16}, \frac{3}{10}, \frac{5}{17}$, with $(h_4,h_5)=(0,0)$.
	\end{itemize}
	
	\noindent Finally, Lemma \ref{lem15} fixes the values of $(a_3,c_3)$ to be considered. Using (\ref{eq12}), values of $(a_4,a_5)$ are also fixed. For these particular values, by choosing suitable values of $x_3$, $x_4$ and $x_5$, we find that the range of $f_1$ or, $p_1$ is covered and hence (\ref{eq101*})
	is soluble.\vspace{4mm}\hfill $\Box$
	
	\noindent This completes the proof of Theorem \ref{thm8}, i.e. $\Gamma_{1,4}< \frac{32}{3}$ when $(m,K)=(1,1)$.


\begin{thebibliography}{99}
		
		\bibitem{AggarwalGupta1988} S.K. Aggarwal, D.P. Gupta, Positive values of inhomogeneous quadratic forms of signature 2, \emph{J. Number Theory} \textbf{29} (1988), 138-165.
		
		\bibitem{AggarwalGupta1991} S.K. Aggarwal, D.P. Gupta, Least positive value of non-homogeneous indefinite quadratic forms of signature 3, \emph{J. Number Theory} \textbf{37} (1991), 260-278.
		
		\bibitem{Bambah1981} R.P. Bambah, V.C. Dumir and R.J. Hans-Gill, Positive values of non homogeneous indefinite quadratic forms, \emph{Proceedings Colloqu, Classical Number Theory, Budapest(Hungary)}, (1981), 111-170 .
		
		\bibitem{Bambah1983} R.P. Bambah, V.C. Dumir, R.J. Hans-Gill, On a conjecture of Jackson on non-homogeneous quadratic forms, \emph{J. Number Theory} \textbf{16} (1983), 403-419.
		
		\bibitem{Bambah1984} R.P. Bambah, V.C. Dumir, R.J. Hans-Gill, Positive values of non-homogeneous indefinite quadratic forms II, \emph{J. Number Theory} \textbf{18} (1984), 313-341.
		
		\bibitem{Barnes} E.S. Barnes, The positive values of inhomogeneous ternary quadratic forms, \emph{J. Austral. Math. Soc.} \textbf{2} (1961), 127-132.
		
		\bibitem{Birch} B.J. Birch, The inhomogeneous minimum of quadratic forms of signature zero, \emph{Acta} \textbf{3} (1958), 85-98.
		
		\bibitem{Blaney1948} H. Blaney, Indefinite quadratic forms in n variables, \emph{J. London Math. Soc.} \textbf{23} (1948), 153-160.
		
		\bibitem{Blaney1950} H. Blaney, Some asymmetric inequalities, \emph{Proc. Cambridge Philos. Soc.} \textbf{46} (1950), 359-376.
		
		\bibitem{Davenport1947} H. Davenport, H. Heilbronn, Asymmetric inequalities for non-homogeneous linear forms, \emph{J. London Math. Soc.} \textbf{22} (1947), 53-61.
		
		\bibitem{Dumir1967} V.C. Dumir, Asymmetric inequality for non-homogeneous ternary quadratic forms, \emph{Math. Proc. Cambridge Philos. Soc.} \textbf{63} (1967), 291-303.
		
		\bibitem{Dumir1968} V.C. Dumir, Positive values of inhomogeneous quadratic forms I and II, \emph{J. Aust. Math. Soc.} \textbf{8} (1968), 87-101, 287-303.
		
		\bibitem{DumirGill} V.C. Dumir, R.J. Hans-Gill, On positive values of non-homogeneous quaternary quadratic forms of type $(1,3)$, \emph{Indian J. Pure Appl. Math.} \textbf{12} (1981), 814-825.
		
		\bibitem{DGS1995} V.C. Dumir, R.J. Hans-Gill, R. Sehmi, Positive values of non-homogeneous indefinite quadratic forms of type  (2,4), \emph{J. Number Theory} \textbf{55}(2) (1995), 261-284.
		
		\bibitem{DumirSehmi} V.C. Dumir, R. Sehmi, Positive Values of non homogeneous indefinite quadratic forms of type (1,4), \emph{Proc. Indian Acad. Sci.} \textbf{104} (1994), 557-579.
		
		\bibitem{DumirWoods} V.C. Dumir, R.J. Hans-Gill, A.C. Woods, Values of non-homogeneous indefinite quadratic forms, \emph{J. Number Theory} \textbf{47} (1994), 190-197.
		
		\bibitem{Flahive1988} M.Flahive, Indefinite quadratic forms in many variables, \emph{Indian J. Pure Appl. Math} \textbf{19}(10) (1988), 931-959.
		
		\bibitem{Gruber} P.M. Gruber, C.G. Lekkerkerker, Geometry of numbers, Second Edition, North Holland Elsevier Science Publishers, (1987).
		
		\bibitem{GillRaka1980} R.J. Hans-Gill, M. Raka, Positive values of inhomogeneous 5-ary quadratic forms of type $(3,2)$, \emph{J. Austral. Math. Soc. Ser. A} \textbf{29} (1980), 439-453.
		
		\bibitem{GillRaka1981} R.J. Hans-Gill, M. Raka, Positive values of inhomogeneous quinary quadratic forms of type $(4,1)$, \emph{J. Austral. Math. Soc. Ser. A} \textbf{31} (1981), 175-188.
		
		\bibitem{Jackson} T.H. Jackson, Gaps between values of quadratic forms, \emph{J. London Math.Soc.} \textbf{3} (1971), 47-58.
		
		\bibitem{Macbeath} A.M. Macbeath, A new sequence of minima in the geometry of numbers, \emph{Math. Proc. Cambridge Philos. Soc.} \textbf{47}(2) (1951), 266-273.
		
		\bibitem{Marguils} G.A. Marguils, Indefinite quadratic forms and unipotent flows on homogeneous spaces, \emph{Comp. Rend. Acad. Sci.} \textbf{304} (1987), 249-253.
		
		\bibitem{RakaRani} M. Raka, U. Rani, Positive Values of non-homogeneous quadratic forms of type (1,4), \emph{Proc. Indian Acad. Sci.} \textbf{107}(4) (1997), 329-361.
		
		\bibitem{DumirSehmi1994} R. Sehmi, V.C. Dumir, Positive Values of non-homogeneous indefinite quadratic forms of type (2,5), \emph{Journal of Number Theory} \textbf{48}(1) (1994), 1-35. 
		
		\bibitem{Watson} G.L. Watson, Indefinite quadratic polynomials, \emph{Mathematika} \textbf{7} (1960), 141-144.
		
		
	\end{thebibliography}
\end{document}